\documentclass[12pt,a4paper]{article}
\usepackage{amsmath,amssymb,bm,ascmac,bbm,mathtools}
\usepackage[dvipdfmx, usenames]{color}
\usepackage[dvipdfmx]{graphicx}
\usepackage{xcolor}
\usepackage{here}
\usepackage{authblk}
\usepackage[hang,small,bf]{caption}
\usepackage[subrefformat=parens]{subcaption}
\usepackage{dcolumn}

\newcommand{\tr}{\text{tr}}

\newcommand{\mcal}{\mathcal}
\newcommand{\mbb}{\mathbb}
\newcommand{\mfrak}{\mathfrak}

\newcommand{\rank}{\text{rank}}

\newcommand{\mc}{\mathcal}

\newcommand{\wf}[1]{\widehat{\mfrak{#1}}}
\newcommand{\vecG}{\text{Vec}}
\newcommand{\vect}{\text{Vect}_{\mbb C}}
\newcommand{\fp}{\text{FPdim}}
\newcommand{\fib}{\text{Fib}}
\newcommand{\ising}{\text{Ising}}
\newcommand{\tc}{\text{ToricCode}}

\newcommand{\mods}{\text{mod }}

\begin{document}
\title{Classification of connected \'etale algebras in multiplicity-free modular fusion categories up to rank nine}
\author{Ken KIKUCHI}
\affil{Department of Physics, National Taiwan University, Taipei 10617, Taiwan}
\date{}
\maketitle

\begin{abstract}
We classify connected \'etale algebras $A$'s in multiplicity-free modular fusion categories $\mathcal B$'s with $\text{rank}(\mathcal B)\le9$. We also identify categories $\mathcal B_A$'s of right $A$-modules. The results have physical applications in constraining renormalization group flows. As demonstration, we study massive renormalization group flows from non-unitary minimal models to predict ground state degeneracies and prove spontaneous $\mathcal B$-symmetry breaking.
\end{abstract}


\makeatletter
\renewcommand{\theequation}
{\arabic{section}.\arabic{equation}}
\@addtoreset{equation}{section}
\makeatother

\section{Introduction}
Throughout the paper, $\mc C$ denotes a fusion category over the field $\mbb C$ of complex numbers (see \cite{EGNO15} for definitions). Its rank and simple objects are denoted $\rank(\mc C)$ and $c_i$'s with $i=1,2,\dots,\rank(\mc C)$, respectively. When a fusion category admits a braiding $c$, the braided fusion category (BFC) and its simple objects are denoted $\mc B$ and $b_i$'s, respectively. A BFC with non-degenerate braiding $c$ is called modular or modular fusion category (MFC). In a BFC, we can define an étale algebra $A\in\mc B$. (The definitions are collected in section \ref{def}.) The goal of this paper is to classify connected étale algebras in multiplicity-free MFCs up to rank nine. Our main results are summarized as the
\newpage

\textbf{Theorem.} \textit{Connected étale algebras in multiplicity-free modular fusion categories up to rank nine are given by}
\begin{table}[H]
\begin{center}
\makebox[1 \textwidth][c]{       
\resizebox{1.2 \textwidth}{!}{\begin{tabular}{c|c|c|c}
    Rank&$\mcal B$&Results&Completely anisotropic?\\\hline\hline
    6&$su(6)_1\simeq\vecG_{\mbb Z/6\mbb Z}^\alpha$&\cite{G23,KKH24} &Yes\\
    &$\vecG_{\mbb Z/2\mbb Z}^{-1}\boxtimes\ising$&e.g. \cite{KKH24} &Yes\\
    &$su(3)_2\simeq\fib\boxtimes\vecG_{\mbb Z/3\mbb Z}^\alpha$&\cite{EP09,KKH24}&Yes\\
    &$\text{TriCritIsing}$&\cite{KL02,KKH24}&Yes\\
    &$su(2)_5\simeq\vecG_{\mbb Z/2\mbb Z}^{-1}\boxtimes psu(2)_5$&\cite{KO01,KKH24}&Yes\\
    &$so(5)_2$&\cite{SW86,BB87,DMNO10,G23,KKH24}&No\\
    &$\fib\boxtimes psu(2)_5$&e.g. \cite{KKH24} &Yes\\
    &$psu(2)_{11}$&e.g. \cite{KKH24}&Yes\\\hline
    7&$su(7)_1\simeq\vecG_{\mbb Z/7\mbb Z}^1$&\cite{G23}, Table \ref{rank7Z7results}&Yes\\
    &$su(2)_6$&\cite{KO01}, Table \ref{rank7su26results}&Yes\\
    &$so(7)_2$&\cite{SW86,BB87,DMNO10,G23}, Table \ref{rank7so72results}&No\\
    &$psu(2)_{13}$&Table \ref{rank7psu213results}&Yes\\\hline
    8&$\vecG_{\mbb Z/2\mbb Z\times\mbb Z/2\mbb Z\times\mbb Z/2\mbb Z}^\alpha\simeq\vecG_{\mbb Z/2\mbb Z}^{-1}\boxtimes\tc$&Table \ref{rank8Z2Z2Z2results}&No (16 with $d,h$ in (\ref{Z2Z2Z2etale}))/ Yes (the other 24)\\
    &$\vecG_{\mbb Z/2\mbb\times\mbb Z/4\mbb Z}^\alpha$&Table \ref{rank8Z2Z4results}&Yes\\
    &$su(8)_1\simeq\vecG_{\mbb Z/8\mbb Z}^\alpha$&\cite{SW86,BB87,DMNO10}, Table \ref{rank8Z8results}&No\\
    &$\fib\boxtimes\vecG_{\mbb Z/2\mbb Z\times\mbb Z/2\mbb Z}^\alpha\simeq\begin{cases}\fib\boxtimes\vecG_{\mbb Z/2\mbb Z}^{-1}\boxtimes\vecG_{\mbb Z/2\mbb Z}^{-1}\\\fib\boxtimes\tc\end{cases}$&$\begin{cases}\text{Table }\ref{rank8fibZ2Z2results}\\\text{Table }\ref{rank8fibtoriccoderesults}\end{cases}$&$\begin{cases}\text{No }(\text{16 with }d,h\text{ in }(\ref{fibZ2Z21Zetale}))/\text{ Yes (the other 64)}\\\text{No }(\text{16 with }d,h\text{ in }(\ref{fibtcetale}))/\text{ Yes (the other 24)}\end{cases}$\\
    &$\fib\boxtimes\vecG_{\mbb Z/4\mbb Z}^\alpha$&Table \ref{rank8fibZ4results}&Yes\\
    &$\vecG_{\mbb Z/2\mbb Z}^{-1}\boxtimes\fib\boxtimes\fib$&Table \ref{rank8Z2fibfibresults}&No (16 with $d,h$ in (\ref{Z2fibfibetale})/ Yes (the other 64)\\
    &$so(9)_2$&\cite{SW86,BB87,DMNO10,G23}, Table \ref{rank8so92results}&No\\
    &$\text{Rep}(D(D_3))$&Table \ref{rank8RepDD3results}&No\\
    &$su(2)_7$&\cite{KO01}, Table \ref{rank8su27results}&Yes\\
    &$\fib\boxtimes\fib\boxtimes\fib$&\cite{BD11}, Table \ref{rank8fibfibfibresults}&No (16 with $d,h$ in (\ref{FibFibFib1Tetale},\ref{FibFibFib1UVetale}))/ Yes (the other 24)\\
    &$\fib\boxtimes psu(2)_7$&Table \ref{rank8fibpsu27results}&Yes\\
    &$psu(2)_{15}$&Table \ref{rank8psu215results}&Yes\\\hline
    9&$su(9)_1\simeq\vecG_{\mbb Z/9\mbb Z}^1$&\cite{SW86,BB87,DMNO10}, Table \ref{rank9Z9results}&No\\
    &$\vecG_{\mbb Z/3\mbb Z}^1\boxtimes\vecG_{\mbb Z/3\mbb Z}^1$&Table \ref{rank9Z3Z3results}&No (two with the 1st $h$)/ Yes (the other two)\\
    &$\vecG_{\mbb Z/3\mbb Z}^1\boxtimes\ising$&Table \ref{rank9Z3isingresults}&Yes\\
    &$\ising\boxtimes\ising$&Table \ref{rank9isingisingresults}&No\\
    &$\vecG_{\mbb Z/3\mbb Z}^1\boxtimes psu(2)_5$&Table \ref{rank9Z3psu25results}&Yes\\
    &$\ising\boxtimes psu(2)_5$&Table \ref{rank9isingpsu25results}&Yes\\
    &$so(11)_2$&\cite{SW86,BB87,DMNO10,G23}, Table \ref{rank9so112results}&No\\
    &$su(2)_8$&\cite{SW86,BB87,KO01,DMNO10}, Table \ref{rank9su28results}&No\\
    &$psu(2)_5\boxtimes psu(2)_5$&Table \ref{rank9psu25psu25results}&No (six with $d,h$ in (\ref{psu25psu25etale}))/ Yes (the other 36)\\
    &$psu(2)_{17}$&Table \ref{rank9psu217results}&Yes
\end{tabular}.}}
\end{center}
\caption{Connected étale algebras in multiplicity-free MFC $\mcal B$ up to rank nine}\label{results}
\end{table}

\textbf{Remark.} For a reader's convenience, we also included known results \cite{KKH24} at rank six.\newline

\textbf{Remark.} Some MFCs are realized by rational conformal field theories (RCFTs) such as Wess-Zumino-Witten (WZW) models or minimal models. In those cases, we collectively denote the MFCs sharing the same fusion ring by the realization, e.g., $su(2)_6$ or $\ising$. Other MFCs are realized by subcategories of simple objects invariant under centers. We denote the MFCs by the realizations with $p$ in their head, e.g., $psu(2)_{13}$. Note that this is \textit{not} a `gauging' of their centers. For instance, the $\mbb Z/2\mbb Z$ center symmetry of $su(2)_{13}$ is anomalous, and cannot be gauged. (Mathematically, this means there is no $\mbb Z/2\mbb Z$ algebra.)\newline

\textbf{Remark.} MFCs up to rank five have been classified in \cite{GK94,RSW07,BNRW15}. Classification of larger MFCs is not complete yet. Results we know are limited to multiplicity-free ones up to rank nine \cite{LPR20,VS22}, summarized in AnyonWiki \cite{anyonwiki}. We study these MFCs; there are four rank seven, 12 rank eight, and 10 rank nine multiplicity-free fusion rings.\newline

\textbf{Remark.} The classification problem was initiated in \cite{DMNO10,DMO11}, but some results had been known before these works. Especially when MFCs admit realizations by WZW models, some classification results are known. (In this context, an MFC describing $\wf g_k$ WZW model is denoted $\mc C(\mfrak g,k)$, and connected étale algebras in it are called quantum subgroups \cite{O00}. In another context, say, \cite{CZW18}, connected étale algebras are also called condensable or normal algebras.) For example, connected étale algebras were classified in \cite{KO01} (for $\wf{su}(2)_k$), in \cite{EP09} (for $\wf{su}(3)_k$), in \cite{CEM23} (for $\wf{su}(4)_k$), and many more in \cite{G23}. Connected étale algebras in (pre-)MFCs up to rank five have been classified in \cite{KK23GSD,KK23preMFC,KK23MFC}. When available, our results are consistent with them; $\mc C(A_6,1),\mc C(A_1,6),\mc C(A_1,7)$ are known to be completely anisotropic. Also, $\fib\boxtimes\fib\boxtimes\fib$ is known to be completely anisotropic when braiding structures of all factors are the same, while it can admit nontrivial connected étale algebras when two factors have the opposite braidings. Examples with nontrivial quantum subgroups are also consistent. The MFCs $\mc C(B_3,2),\mc C(B_5,2),\mc C(A_1,8)$ are known to have two quantum subgroups $1,1\oplus X$ where $X$ is a $\mbb Z/2\mbb Z$ simple object in each MFC. The nontrivial connected étale algebras give conformal embeddings $\wf{so}(7)_2\subset\wf{su}(7)_1,\wf{so}(11)_2\subset\wf{su}(11)_1,(\wf{g_2})_1\times\wf{su}(2)_8\subset(\wf{f_4})_1$, respectively. The MFC $\mc C(B_4,2)$ has $1,1\oplus X$ and exotic\footnote{A quantum subgroup which is not a direct sum of invertible simple objects was called exotic in \cite{G23}.} quantum subgroups $1\oplus Z,1\oplus X\oplus2Z,1\oplus Z\oplus V$.\footnote{We thank Terry Gannon for teaching us this fact.} The nontrivial connected étale algebras give $\wf{so}(9)_2\subset\wf{su}(9)_1,\wf{so}(9)_2\subset\wf{so}(16)_1,\wf{so}(9)_2\subset(\wf{e_8})_1$, and another $\wf{so}(9)_2\subset(\wf{e_8})_1$, respectively.\newline

\textbf{Remark.} Connected étale algebras have also been studied in physics. For example, in the context of anyon condensation \cite{BSS02,BSS02',BS08}, physically natural conditions demand condensing object $A$ be connected étale \cite{K13}.\newline

\textbf{Remark.} Algebras in symmetry categories have another aspect; they give anomaly-free subsymmetries \cite{FRS02,CR12,BT17}. Therefore, these results also classify certain symmetries which can be gauged.

\section{Classification}
In this section, we classify connected étale algebras in multiplicity-free MFCs up to rank nine. Before we study the problem, we first collect relevant definitions in section \ref{def} and explain our classification method in section \ref{method}.
\subsection{Definitions}\label{def}
Let $\mc C$ be a fusion category over $\mbb C$, i.e., a finite $\mbb C$-linear semisimple abelian rigid monoidal category with bilinear monoidal product $\otimes$ and simple unit object $1$. The monoidal product of $\mc C$ is specified by fusion matrices $(N_i)_{jk}:={N_{ij}}^k$ with $\mbb N$-coefficients
\[ c_i\otimes c_j\cong\bigoplus_{k=1}^{\rank(\mc C)}{N_{ij}}^kc_k. \]
We denote the fusion ring of $\mc C$ as $K(\mc C)$. On the other hand, we denote a fusion category with a fusion ring $K$ as $\mc C(K)$. Since the entries of fusion matrices are non-negative, we can apply the Perron-Frobenius theorem to obtain the largest eigenvalue $\fp_{\mc C}(c_i)$ called the Frobenius-Perron dimension of a simple object $c_i$.\footnote{We add ambient categories in the subscript because the Frobenius-Perron dimension of a given object depend on the ambient categories.} The Frobenius-Perron dimension of the category $\mc C$ is defined as
\begin{equation}
    \fp(\mc C):=\sum_{i=1}^{\rank(\mc C)}\left(\fp_{\mc C}(c_i)\right)^2.\label{FPdimC}
\end{equation}
In a spherical fusion category (including MFCs), one can also define quantum dimension $d_i$ of $c_i$ by the quantum (or categorical) trace
\[ d_i:=\tr(a_{c_i}), \]
where $a:id_{\mc C}\cong(-)^{**}$ is a pivotal structure. Its multiplication rules are the same as the fusion rules
\begin{equation}
    d_id_j=\sum_{k=1}^{\rank(\mc C)}{N_{ij}}^kd_k.\label{quantumdimmultiplication}
\end{equation}
When various fusion categories are involved, in order to avoid confusion, we denote a quantum dimension of $c_i\in\mc C$ as $d_{\mc C}(c_i)$. The squared sum of quantum dimensions define the categorical (or global) dimension
\begin{equation}
    D^2(\mc C):=\sum_{i=1}^{\rank(\mc C)}d_i^2.\label{categoricaldim}
\end{equation}
Note that there are two $D(\mc C)$'s, one positive and one negative, for each categorical dimension. A fusion category is called pseudo-unitary if $D^2(\mc C)=\fp(\mc C)$ and unitary if $\forall c_i\in\mc C,\ d_i=\fp_{\mc C}(c_i)$.

An MFC $\mc B$ has additional structure, braiding $c$. It is a natural isomorphism between two bifunctors
\begin{equation}
    c_{-,-}:-\otimes-\stackrel\sim\Rightarrow-\otimes-\label{braiding}
\end{equation}
subject to hexagon axioms. More concretely, for $b,b'\in\mc B$ (not necessarily simple), it is a family of natural isomorphisms
\[ c_{b,b'}:b\otimes b'\stackrel\sim\Rightarrow b'\otimes b. \]
A fusion category with a braiding is called braided fusion category (BFC). A BFC is a pair $(\mc B,c)$ of a fusion category $\mc B$ and a braiding $c$, but we often write $\mc B$. For a BFC $\mc B=(\mc B,c)$, there exists another BFC $\tilde{\mc B}$ with the reverse (or opposite) braiding $\tilde c_{b,b'}=c_{b',b}^{-1}$. The BFC $\tilde{\mc B}=(\mc B,\tilde c)$ is called reverse BFC. The structure is specified by conformal dimensions $h_i$'s of $b_i$'s. For instance, the double braiding of two simple objects $b_i,b_j\in\mc B$ is given by
\begin{equation}
    c_{b_j,b_i}\cdot c_{b_i,b_j}\cong\sum_{k=1}^{\rank(\mc B)}{N_{ij}}^k\frac{e^{2\pi ih_k}}{e^{2\pi i(h_i+h_j)}}id_k,\label{doublebraiding}
\end{equation}
where $id_k$ is the identity morphism at $b_k\in\mc B$. (A BFC $\mc B$ is called symmetric if $\forall b,b'\in\mc B$, $c_{b',b}\cdot c_{b,b'}\cong id_{b\otimes b'}$.) Its quantum trace defines (unnormalized) $S$-matrix
\begin{equation}
    \widetilde S_{i,j}:=\tr(c_{b_j,b_i}\cdot c_{b_i,b_j})=\sum_{k=1}^{\rank(\mc B)}{N_{ij}}^k\frac{e^{2\pi ih_k}}{e^{2\pi i(h_i+h_j)}}d_k.\label{tildeS}
\end{equation}
A normalized $S$-matrix is defined by
\begin{equation}
    S_{i,j}:=\frac{\widetilde S_{i,j}}{D(\mc B)}.\label{normalizedS}
\end{equation}
An MFC is defined as a spherical BFC (called pre-modular fusion category, pre-MFC) with non-degenerate $S$-matrix. It squares to charge conjugation matrix
\[ S^2=C. \]
The charge conjugation matrix is defined by
\[ C_{i,j}=\delta_{i,j}\quad(b_i^*\cong b_j), \]
where $b_i^*\in\mc B$ is the dual of $b_i\in\mc B$. It obeys\
\[ \widetilde S_{i,j^*}=\left(\widetilde S_{i,j}\right)^*, \]
where the RHS is the complex conjugate of $\widetilde S_{i,j}$. Another modular matrix $T$ is also defined with conformal dimensions
\begin{equation}
    T_{i,j}:=e^{2\pi ih_i}\delta_{i,j}.\label{Tmatrix}
\end{equation}
The two modular matrices define an additive central charge $c(\mc B)$ mod 8 by
\begin{equation}
    (ST)^3=e^{2\pi ic(\mc B)/8}S^2.\label{additivecentralcharge}
\end{equation}

Given a BFC, we can define commutative algebras. An algebra in a fusion category $\mc C$ is a triple $(A,\mu,\eta)$ of an object $A\in\mc C$, multiplication morphism $\mu:A\otimes A\to A$, and unit morphism $\eta:1\to A$ subject to associativity and unit axioms. (We abuse the notation and an algebra is often denoted by $A$.) A category $\mc C_A$ of right $A$-modules consists of pairs $(m,p)$ of an object $m\in\mc C$ and morphism $p:m\otimes A\to m$ subject to consistency conditions. A category $_A\mc C$ of left $A$-modules are defined analogously. (We also abuse the notation and an $A$-module is often denoted by $m$.) An algebra $A\in\mc C$ is called separable if $\mc C_A$ is semisimple. An algebra $A$ in a BFC with braiding $c$ is called commutative if
\begin{equation}
    \mu\cdot c_{A,A}=\mu.\label{commutativealg}
\end{equation}
A commutative separable algebra is called étale. Any étale algebra decomposes to a direct sum of connected ones \cite{DMNO10}. Here, an algebra $A$ in a fusion category $\mc C$ is called connected if $\dim_{\mbb C}\mc C(1,A)=1$. A connected étale algebra $A\in\mc B$ is called Lagrangian if $(\fp_{\mc B}(A))^2=\fp(\mc B)$. An example of a connected étale algebra is the unit object $A\cong1\in\mc B$ giving $\mc B_A\simeq\mc B$. The connected étale algebra always exists, and is thus called trivial. A BFC is called completely anisotropic if it has no nontrivial connected étale algebra.

The category of right $A$-modules has an important subcategory $\mc B_A^0\subset\mc B_A$. It consists of dyslectic (or local) modules \cite{P95} $(m,p)\in\mc B_A$ obeying
\begin{equation}
    p\cdot c_{A,m}\cdot c_{m,A}=p.\label{dyslectic}
\end{equation}
(In the context of anyon condensation \cite{BSS02,BSS02',BS08}, $\mc B_A,\mc B_A^0$ are called broken and deconfined phase, respectively.) The latter, the category $\mc B_A$ of right $A$-modules, is a left $\mc B$-module category. Here, for a fusion category $\mc C$, a left $\mc C$-module category (or module category over $\mc C$) \cite{O01} is a quadruple $(\mc M,\triangleright,m,l)$ of a category $\mc M$, an action (or module product) bifunctor $\triangleright:\mc C\times\mc M\to\mc M$, and natural isomorphisms $m_{-,-,-}:(-\otimes-)\triangleright-\stackrel\sim\Rightarrow-\triangleright(-\triangleright-)$ and $l:1\triangleright\mc M\simeq\mc M$ called module associativity constraint and unit constraint, respectively. They are subject to associativity and unit axioms. Let $\mc M_{1,2}$ be $\mc C$-module categories. The category $\mc M\simeq\mc M_1\oplus\mc M_2$ is called a direct sum of the module categories $\mc M_{1,2}$. A $\mc C$-module category is called indecomposable if it is not equivalent to a direct sum of nontrivial module categories. In our setup, the categories $\mc M\simeq\mc B_A$ (and also $\mc B_A^0$) are fusion categories with monoidal products $\otimes_A$.\footnote{For an explanation, see the footnote 2 of \cite{KK23preMFC}.} Then, the action of (left) $\mc B$-module categories forms non-negative integer matrix representations (NIM-reps). Here is the reason. For any $b\in\mc B,m\in\mc M$, $b\triangleright m$ is an object of $\mc M$. Hence, it can be decomposed to a direct sum of simple objects in $\mc M$ with $\mbb N$-coefficients. The natural numbers assemble to $r\times r$-matrices where $r=\rank(\mc M)$. The NIM-rep is called $r$-dimensional.

\subsection{Method}\label{method}
In this section, we explain our classification method. It is based on \cite{KK23GSD,KK23preMFC}, but developed further.

Let $\mc B$ be an MFC and $A\in\mc B$ a connected étale algebra. The category $\mc B_A^0$ of dyslectic right $A$-modules is modular \cite{P95,KO01} obeying \cite{KO01,DMNO10}
\begin{equation}
\begin{split}
    \fp(\mc B_A^0)&=\frac{\fp(\mc B)}{(\fp_{\mc B}(A))^2},\\
    e^{2\pi ic(\mc B_A^0)/8}&=e^{2\pi ic(\mc B)/8}.\label{propBA0}
\end{split}
\end{equation}
Since we have \cite{ENO02,EGNO15}
\begin{equation}
    \forall c\in\mc C,\quad\fp_{\mc C}(c)\ge1,\label{FPdimge1}
\end{equation}
we obtain an inequality
\begin{equation}
    1\le\left(\fp_{\mc B}(A)\right)^2\le\fp(\mc B).\label{FPdimA2bound}
\end{equation}
In addition, since $A$ consists of simple objects of $\mc B$, its general form is given by
\begin{equation}
    A\cong\bigoplus_{i=1}^{\rank(\mc B)}n_ib_i,\label{generalA}
\end{equation}
where $n_i\in\mbb N$ counts the number of $b_i$ in $A$. As the direct sum is defined as (co)limit, the object is equipped with product projections $p_i:A\to b_i$ and coproduct injections $\iota_i:b_i\to A$. Its Frobenius-Perron dimension is given by the linear sum of those of simple objects:
\begin{equation}
    \fp_{\mc B}(A)=\sum_{i=1}^{\rank(\mc B)}n_i\fp_{\mc B}(b_i).\label{FPdimAlinearsum}
\end{equation}
In classifying connected algebras, we can set $n_i=1$ for $b_i\cong1$ at the outset. Thus, for a given MFC $\mc B$, we solve (\ref{FPdimA2bound}) with an ansatz
\[ A\cong1\oplus\bigoplus_{b_i\not\cong1}n_ib_i\quad(n_i\in\mbb N), \]
or
\[ \fp_{\mc B}(A)=1+\sum_{b_i\not\cong1}n_i\fp_{\mc B}(b_i). \]
As a result, we get a set of natural numbers $n_i$'s. Each element gives a candidate for $A$. Then, for each candidate, we check whether it satisfies the axioms of connected étale algebra. Our strategy here is to check various necessary conditions to reduce the number of candidates. In particular, we check the commutativity (\ref{commutativealg}). In checking this axiom, it turns out useful to study a necessary condition
\begin{equation}
    \mu\cdot c_{A,A}\cdot c_{A,A}=\mu.\label{commutativealgnecessary}
\end{equation}
If $A$ satisfies (\ref{commutativealg}), then by substituting the RHS in the LHS, we obtain the necessary condition. Since the condition is written in terms of double braiding, we can use the formula (\ref{doublebraiding}). Here, the double braiding of the direct sum (\ref{generalA}) is computed by chasing commuting diagrams \cite{KK23GSD}:
\begin{equation}
    c_{A,A}\cdot c_{A,A}\cong\sum_{i,j=1}^{\rank(\mc B)}n_in_j(\iota_i\otimes\iota_j)\cdot c_{b_j,b_i}\cdot c_{b_i,b_j}\cdot(p_i\otimes p_j).\label{doublecAA}
\end{equation}
(The morphisms $\iota,p$ are coproduct injection and product projection, respectively, introduced below (\ref{generalA}).) If there exists a term with nontrivial phase in (\ref{doublecAA}), it means the candidate is non-commutative, and it is ruled out. In particular, if a simple object $b_i\in\mc B$ has $1\in b_i\otimes b_i$ and $e^{-4\pi ih_i}\neq1$, $b_i$ cannot enter a commutative algebra. Thus, we can set $n_i=0$ for such simple objects at the outset. We will see the necessary condition (\ref{commutativealgnecessary}) is strong enough to rule out most of the simple objects. For the remaining candidates which pass the necessary condition, we check other necessary conditions such as an existence of MFC with the Frobenius-Perron dimension (\ref{propBA0}) which matches the additive central charges. For a handful final candidates passing all necessary conditions, we check axioms of an algebra and directly compute the braiding \cite{KK23GSD}
\begin{equation}
    c_{A,A}\cong\sum_{i,j=1}^{\rank(\mc B)}n_in_j(\iota_j\otimes\iota_i)\cdot c_{b_i,b_j}\cdot(p_i\otimes p_j)\label{cAA}
\end{equation}
to check (\ref{commutativealg}). Here, some facts save us from tedious computations.\newline

\textbf{Lemma 1.} \textit{Let $(\mc C,\otimes,\alpha,1,\iota)$ be a monoidal category, $(\mc C',\otimes,\alpha,1,\iota)$ with $\mc C'\subset\mc C$ a monoidal subcategory, and $(A,\mu,\eta)$ an algebra in $(\mc C',\otimes,\alpha,1,\iota)$. The algebra in the monoidal subcategory $\mc C'$ is also an algebra in the larger monoidal category $\mc C$.}\newline

\textit{Proof.} By definition, the algebra obeys associativity
\[ \mu\cdot(id_A'\otimes'\mu)\cdot\alpha_{A,A,A}'=\mu\cdot(\mu\otimes'id_A'), \]
and unit
\[ id_A'\cdot l_A'=\mu\cdot(\eta\otimes'id_A'),\quad id_A'\cdot r_A'=\mu\cdot(id_A'\otimes'\eta) \]
axioms where quantities with primes belong to $\mc C'$. Here, note that the monoidal subcategory has the same structure as the larger monoidal category by definition \cite{EGNO15}. (We already used this fact in the statement.) Furthermore, by uniqueness of identity morphisms (up to isomorphism), we have
\[ id_A'\cong id_A. \]
Since left and right unit constraints $l_A,r_A$ are defined with $\alpha_{1,1,A},\alpha_{A,1,1},\iota,id_A$, they are also isomorphic on $A$:
\[ l_A'\cong l_A,\quad r_A'\cong r_A. \]
The axioms of the algebra $(A,\mu,\eta)$ now reduce to
\[ \mu\cdot(id_A\otimes\mu)\cdot\alpha_{A,A,A}=\mu\cdot(\mu\otimes id_A),\quad id_A\cdot l_A=\mu\cdot(\eta\otimes id_A),\quad id_A\cdot r_A=\mu\cdot(id_A\otimes\eta). \]
This is nothing but the definition of $(A,\mu,\eta)$ being an algebra in the larger monoidal category $(\mc C,\otimes,\alpha,1,\iota)$. $\square$\newline

\textbf{Remark.} Physically, this means anomaly-free symmetries in $\mc C'$ remain anomaly-free in a larger symmetry $\mc C\supset\mc C'$, as physicists know.\newline

Furthermore, commutative algebras in braided fusion subcategory remain commutative in a larger BFC. This follows from the definition of braided fusion subcategory. Let $(\mc B',c')$ be a braided fusion subcategory of a larger BFC $(\mc B,c)$. Then, the braidings are equivalent on $\mc B'$. Namely, $\forall b,b'\in\mc B'$, we have
\[ c'_{b,b'}\cong c_{b,b'}. \]
Therefore, a commutative algebra $(A,\mu,\eta)$ in a BFC $(\mc B',c')$ obeying $\mu\cdot c'_{A,A}=\mu$ is a commutative algebra in a larger BFC $(\mc B,c)$:
\[ \mu\cdot c_{A,A}=\mu\cdot c'_{A,A}=\mu. \]
Sometimes, a direct check of commutativity axiom is difficult especially when braidings $c_{b_i,b_j}$ are not known. In those cases, we indirectly find connected étale algebras manipulating the following fact:\newline

\textbf{Lemma 2.} \cite{CR12,BT17} \textit{The operations (not all) $\mc B\to\mc B_A^0$ are composable, associative, and invertible.}\newline

\textbf{Remark.} Note that not all operations are composable. For example, an operation $\mc B\to\vect$ given by a Lagrangian algebra (if exists) cannot be composed with other operations. Therefore, the operations and their compositions fail to form a group.\newline

How do we use this fact? Let $A\in\mc B$ and $A'\in\mc B_A^0$ be connected étale algebras. The second algebra gives another MFC $(\mc B_A^0)_{A'}^0$. We get a sequence $\mc B\to\mc B_A^0\to(\mc B_A^0)_{A'}^0$. Since the two operations can be composed, we learn there should exist a connected étale algebra $\tilde A\in\mc B$ giving $\mc B_{\tilde A}^0\simeq(\mc B_A^0)_{A'}^0$. Since $\mc B_A^0$ has smaller Frobenius-Perron dimension (hence usually smaller rank) than $\mc B$ thanks to (\ref{propBA0}), this can reduce the classification problem in larger MFC $\mc B$ to that in smaller MFC $\mc B_A^0$. For an example of this lemma in action, see section \ref{so92}.

Another fact we use is the\newline

\textbf{Lemma 3.} \cite{DMNO10,M12} \textit{Let $\mc B$ be a modular fusion category. The modular fusion category\footnote{This special MFC is equivalent to the Drinfeld center (or quantum double) $Z(\mc B)$ \cite{M01,DGNO09}.} $\mc B\boxtimes\tilde{\mc B}$ has a Lagrangian algebra $A\in\mc B\boxtimes\tilde{\mc B}$.}\newline

Even when we cannot check axioms directly, this fact helps us to ensure an existence of connected étale algebra. See section \ref{psu25psu25} for our use of this lemma. We manipulate all these facts at our disposal.

In order to discuss (physical) applications of classification results, we also identify categories $\mc B_A^0,\mc B_A$ of right $A$-modules. We do this in two steps. First, we identify the category $\mc B_A^0$ of dyslectic right $A$-modules, and next we identify $\mc B_A$ which contains $\mc B_A^0$ as a subcategory. In identifying $\mc B_A^0$, the following fact turns out useful; the MFCs $\mc B$ and $\mc B_A^0$ have the same topological twists \cite{FFRS03}
\begin{equation}
    e^{2\pi ih_b^{\mc B}}=e^{2\pi ih_b^{\mc B_A^0}},\label{invtopologicaltwist}
\end{equation}
where $h_b^{\mc B}$ is the conformal dimension of $b\in\mc B$ and $h_b^{\mc B_A^0}$ is that of $b\in\mc B_A^0$. Physically, this means anyon condensation preserves conformal dimensions (mod 1) of deconfined particles. We see this usually fixes $\mc B_A^0$ uniquely.

When the category of dyslectic right $A$-modules is specified, the category $\mc B_A$ of right $A$-modules is highly constrained. The remaining contributions to the Frobenius-Perron dimensions come from confined particles. This sets an upper bound on $\rank(\mc B_A)$. Furthermore, we can get all candidate simple objects as follows. The free module functor
\begin{equation}
    F_A:=-\otimes A:\mc B\to\mc B_A\label{FA}
\end{equation}
gives objects of $\mc B_A$. The functor is surjective (or dominant) \cite{DMNO10}. Namely, $\forall m\in\mc B_A$, $\exists b\in\mc B$ such that $m$ is a subobject of $b\otimes A$. Since $\mc B_A$ is a fusion category in our setup, $\forall b\in\mc B,\ b\otimes A\in\mc B_A$ can be decomposed as a direct sum of simple objects of $\mc B_A$ (definition of semisimplicity). All possible sets of simple objects are given by consistent collection of some of them. Some consistencies are correct Frobenius-Perron dimensions and closedness under actions $b_j\otimes-$ or monoidal products $\otimes_A$.

Given simple objects of $\mc B_A$, we identify the category by computing monoidal products $\otimes_A$ as follows. The free module functor satisfies two nice properties: 1) $F_A$ is a tensor functor, and 2) it admits the forgetful functor $U:\mc B_A\to\mc B$ as a right adjoint \cite{O01}\footnote{This second property is sometimes called the Frobenius reciprocity (see, say, \cite{CZW18}).}
\begin{equation}
    \forall b\in\mc B,\ \forall m\in\mc B_A,\quad\mc B_A(F_A(b),m)\cong\mc B(b,U(m)).\label{Frobeniusreciprocity}
\end{equation}
The first property in particular implies preservation of Frobenius-Perron dimensions
\begin{equation}
    \forall b\in\mc B,\quad\fp_{\mc B}(b)=\fp_{\mc B_A}(F_A(b)),\label{FPpreserve}
\end{equation}
and
\[ \forall b,b'\in\mc B,\quad F_A(b)\otimes_AF_A(b')\cong F_A(b\otimes b'). \]
When $b\otimes b'$ is a (finite) direct sum, we can distribute using the biexactness of monoidal product
\[ (b_1\oplus b_2\oplus\cdots\oplus b_n)\otimes A\cong(b_1\otimes A)\oplus(b_2\otimes A)\oplus\cdots\oplus(b_n\otimes A). \]
With these facts, we can compute monoidal products $\otimes_A$ to identify $\mc B_A$. In addition, we can compute quantum dimensions. A right $A$-module $m\in\mc B_A$ has \cite{KO01}
\begin{equation}
    d_{\mc B_A}(m)=\frac{d_{\mc B}(m)}{d_{\mc B}(A)}.\label{dBAm}
\end{equation}
\newline

Before we start classifying connected étale algebras, let us see how our method works in a simple example. We take a symmetric pre-MFC $\text{Rep}(S_3)$ as our ambient category $\mc B$. (Note that this is not modular, but most of our classification method except (\ref{propBA0}) still works here.) The pre-MFCs have three simple objects $\{1,X,Y\}$ obeying monoidal products
\begin{table}[H]
\begin{center}
\begin{tabular}{c|c|c|c}
    $\otimes$&1&$X$&$Y$\\\hline
    $1$&1&$X$&$Y$\\\hline
    $X$&&$1$&$Y$\\\hline
    $Y$&&&$1\oplus X\oplus Y$
\end{tabular}.
\end{center}
\end{table}
\hspace{-17pt}They have
\[ \fp_{\mcal B}(1)=1=\fp_{\mcal B}(X),\quad\fp_{\mcal B}(Y)=2, \]
and
\[ \fp(\mcal B)=6. \]
It is known \cite{CZW18,KK23preMFC} that the pre-MFCs with $(h_X,h_Y)=(0,0)$ mod 1 have a connected étale algebra $A\cong1\oplus Y$ giving $\mc B_A\simeq\vecG_{\mbb Z/2\mbb Z}^\alpha$. Let us see how our classification method reproduces this result.

Calculating $b_j\otimes A$, we get
\[ 1\otimes A\cong A\cong1\oplus Y,\quad X\otimes A\cong X\oplus Y,\quad Y\otimes A\cong(1\oplus Y)\oplus(X\oplus Y). \]
The result tells us $\mc B_A$ has two invertible simple objects\footnote{Logically, $F_A(Y)=Y\otimes A$ can be simple. Let us rule out this possibility. Since $F_A$ preserves Frobenius-Perron dimensions, we have $\fp_{\mc B_A}(F_A(Y))=2$. If this is simple, it contributes $2^2=4$ to $\fp(\mc B_A)$, which is a contradiction because it exceeds $\fp(\mc B_A)=2$ we get below. Hence, $F_A(Y)=Y\otimes A$ should be a direct sum of two simple objects with Frobenius-Perron dimensions one.}
\[ 1\oplus Y,\quad X\oplus Y. \]
In $\mc B_A$, they both have Frobenius-Perron dimensions one because $\fp_{\mc B}(1)=1=\fp_{\mc B}(X)$ and the free module functor $F_A:=-\otimes A$ preserves Frobenius-Perron dimensions. On the other hand, for a generic BFC $\mc B$ and a connected étale algebra $A\in\mc B$, we have the formula \cite{KO01,ENO02,DMNO10}
\begin{equation}
    \fp_{\mc B}(A)=\frac{\fp(\mc B)}{\fp(\mc B_A)}.\label{FPdimBA}
\end{equation}
In our current example, this imposes
\[ \fp(\mc B_A)=\frac{\fp(\mc B)}{\fp_{\mc B}(A)}=2. \]
Therefore, no other simple objects can appear. We can identify $\mc B_A\simeq\{1\oplus Y,X\oplus Y\}$. While we know rank two fusion category with Frobenius-Perron dimension two should be $\vecG_{\mbb Z/2\mbb Z}^\alpha$, let us check this fact by computing monoidal products $\otimes_A$. First, we know \cite{EGNO15} $A$ is the monoidal unit of $\otimes_A$. Namely, $\forall m\in\mc B_A$, $m\otimes_AA\cong m\cong A\otimes_Am$. Thus, we only have to compute $(X\oplus Y)\otimes_A(X\oplus Y)$. We have
\begin{align*}
    (X\oplus Y)\otimes_A(X\oplus Y)&=F_A(X)\otimes_AF_A(X)\\
    &\cong F_A(X\otimes X)\\
    &\cong F_A(1)\cong 1\oplus Y.
\end{align*}
Therefore, the two right $A$-modules
\[ m_1\cong1\oplus Y,\quad m_2\cong X\oplus Y \]
obey monoidal products
\begin{table}[H]
\begin{center}
\begin{tabular}{c|c|c}
    $\otimes_A$&$m_1$&$m_2$\\\hline
    $m_1$&$m_1$&$m_2$\\\hline
    $m_2$&&$m_1$
\end{tabular},
\end{center}
\end{table}
\hspace{-17pt}showing $\mc B_A\simeq\vecG_{\mbb Z/2\mbb Z}^\alpha$. Below, we apply this method to classify connected étale algebras in MFCs.

\subsection{Rank seven}
\subsubsection{$\mc B\simeq su(7)_1\simeq\vecG_{\mbb Z/7\mbb Z}^1$}
The MFCs have seven simple objects $\{1,X,Y,Z,U,V,W\}$ obeying monoidal products
\begin{table}[H]
\begin{center}
\begin{tabular}{c|c|c|c|c|c|c|c}
    $\otimes$&$1$&$X$&$Y$&$Z$&$U$&$V$&$W$\\\hline
    $1$&$1$&$X$&$Y$&$Z$&$U$&$V$&$W$\\\hline
    $X$&&$W$&$1$&$V$&$Z$&$Y$&$U$\\\hline
    $Y$&&&$V$&$U$&$W$&$Z$&$X$\\\hline
    $Z$&&&&$X$&$1$&$W$&$Y$\\\hline
    $U$&&&&&$Y$&$X$&$V$\\\hline
    $V$&&&&&&$U$&$1$\\\hline
    $W$&&&&&&&$Z$
\end{tabular}.
\end{center}
\end{table}
\hspace{-17pt}Thus, they have
\[ \hspace{-30pt}\fp_{\mc B}(1)=\fp_{\mc B}(X)=\fp_{\mc B}(Y)=\fp_{\mc B}(Z)=\fp_{\mc B}(U)=\fp_{\mc B}(V)=\fp_{\mc B}(W)=1, \]
and
\[ \fp(\mc B)=7. \]
Their quantum dimensions $d_j$'s are solutions of the same multiplication rules $d_id_j=\sum_{k=1}^7{N_{ij}}^kd_k$. There is only one solution
\[ (d_X,d_Y,d_Z,d_U,d_V,d_W)=(1,1,1,1,1,1). \]
Thus, all MFCs are unitary. The only categorical dimension is
\[ D^2(\mc B)=7. \]
They have conformal dimensions\footnote{Naively, one also finds
\[ (h_X,h_Y,h_Z,h_U,h_V,h_W)=(\frac17,\frac17,\frac27,\frac27,\frac47,\frac47),(\frac27,\frac27,\frac47,\frac47,\frac17,\frac17),(\frac57,\frac57,\frac37,\frac37,\frac67,\frac67),(\frac67,\frac67,\frac57,\frac57,\frac37,\frac37)\quad(\mods1). \]
However, they do not give new conformal dimensions because these are related to those in the main text by permutations $(XUWYZV)$ of simple objects.}
\[ (h_X,h_Y,h_Z,h_U,h_V,h_W)=(\frac37,\frac37,\frac67,\frac67,\frac57,\frac57),(\frac47,\frac47,\frac17,\frac17,\frac27,\frac27)\quad(\mods1). \]
The $S$-matrices are given by
\[ \widetilde S=\begin{pmatrix}1&1&1&1&1&1&1\\1&e^{\mp2\pi i/7}&e^{\pm2\pi i/7}&e^{\pm6\pi i/7}&e^{\mp6\pi i/7}&e^{\pm4\pi i/7}&e^{\mp4\pi i/7}\\1&e^{\pm2\pi i/7}&e^{\mp2\pi i/7}&e^{\mp6\pi i/7}&e^{\pm6\pi i/7}&e^{\mp4\pi i/7}&e^{\pm4\pi i/7}\\1&e^{\pm6\pi i/7}&e^{\mp6\pi i/7}&e^{\mp4\pi i/7}&e^{\pm4\pi i/7}&e^{\pm2\pi i/7}&e^{\mp2\pi i/7}\\1&e^{\mp6\pi i/7}&e^{\pm6\pi i/7}&e^{\pm4\pi i/7}&e^{\mp4\pi i/7}&e^{\mp2\pi i/7}&e^{\pm2\pi i/7}\\1&e^{\pm4\pi i/7}&e^{\mp4\pi i/7}&e^{\pm2\pi i/7}&e^{\mp2\pi i/7}&e^{\pm6\pi i/7}&e^{\mp6\pi i/7}\\1&e^{\mp4\pi i/7}&e^{\pm4\pi i/7}&e^{\mp2\pi i/7}&e^{\pm2\pi i/7}&e^{\mp6\pi i/7}&e^{\pm6\pi i/7}\end{pmatrix}. \]
(All signs are correlated. In other words, one $S$-matrix is given by choosing upper signs in all elements, and the other $S$-matrix is its complex conjugate.) They have additive central charges
\[ c(\mc B)=\begin{cases}-2&(\text{1st }h),\\+2&(\text{2nd }h).\end{cases}\quad(\mods8) \]
There are
\[ 1(\text{quantum dimension})\times2(\text{conformal dimensions})\times2(\text{categorical dimensions})=4 \]
MFCs, and all of them are unitary. We study connected étale algebras in all four MFCs simultaneously.

The most general form of a connected algebra is given by
\[ A\cong1\oplus n_XX\oplus n_YY\oplus n_ZZ\oplus n_UU\oplus n_VV\oplus n_WW \]
with $n_j\in\mbb N$. It has
\[ \fp_{\mc B}(A)=1+n_X+n_Y+n_Z+n_U+n_V+n_W. \]
For this to obey (\ref{FPdimA2bound}), the natural numbers $n_j$'s can only take seven values
\begin{align*}
    (n_X,n_Y,n_Z,n_U,n_V,n_W)=&(0,0,0,0,0,0),(1,0,0,0,0,0),(0,1,0,0,0,0),(0,0,1,0,0,0),\\
    &(0,0,0,1,0,0),(0,0,0,0,1,0),(0,0,0,0,0,1).
\end{align*}
The first solution is nothing but the trivial connected étale algebra $A\cong1$ giving $\mc B_A^0\simeq\mc B_A\simeq\mc B$. The other six solutions do not give connected étale algebras because they fail to satisfy the necessary condition (\ref{commutativealgnecessary}).

We conclude
\begin{table}[H]
\begin{center}
\begin{tabular}{c|c|c|c}
    Connected étale algebra $A$&$\mc B_A$&$\rank(\mc B_A)$&Lagrangian?\\\hline
    $1$&$\mc B$&$7$&No
\end{tabular}.
\end{center}
\caption{Connected étale algebras in rank seven MFC $\mcal B\simeq\vecG_{\mbb Z/7\mbb Z}^1$}\label{rank7Z7results}
\end{table}
\hspace{-17pt}That is, all the four MFCs $\mc B\simeq\vecG_{\mbb Z/7\mbb Z}^1$'s are completely anisotropic.

\subsubsection{$\mc B\simeq su(2)_6$}
The MFCs have seven simple objects $\{1,X,Y,Z,U,V,W\}$ obeying monoidal products
\begin{table}[H]
\begin{center}
\begin{tabular}{c|c|c|c|c|c|c|c}
    $\otimes$&$1$&$X$&$Y$&$Z$&$U$&$V$&$W$\\\hline
    $1$&$1$&$X$&$Y$&$Z$&$U$&$V$&$W$\\\hline
    $X$&&$1$&$Z$&$Y$&$V$&$U$&$W$\\\hline
    $Y$&&&$1\oplus V$&$X\oplus U$&$Z\oplus W$&$Y\oplus W$&$U\oplus V$\\\hline
    $Z$&&&&$1\oplus V$&$Y\oplus W$&$Z\oplus W$&$U\oplus V$\\\hline
    $U$&&&&&$1\oplus U\oplus V$&$X\oplus U\oplus V$&$Y\oplus Z\oplus W$\\\hline
    $V$&&&&&&$1\oplus U\oplus V$&$Y\oplus Z\oplus W$\\\hline
    $W$&&&&&&&$1\oplus X\oplus U\oplus V$
\end{tabular}.
\end{center}
\end{table}
\hspace{-17pt}Thus, they have
\begin{align*}
    \fp_{\mc B}(1)=1=\fp_{\mc B}(X),\quad\fp_{\mc B}(Y)=\sqrt{2+\sqrt2}&=\fp_{\mc B}(Z),\\
    \fp_{\mc B}(U)=1+\sqrt2=\fp_{\mc B}(V),\quad\fp_{\mc B}(W)=&\sqrt{4+2\sqrt2},
\end{align*}
and
\[ \fp(\mc B)=8(2+\sqrt2)\approx27.3. \]
Their quantum dimensions $d_j$'s are solutions of the same multiplication rules $d_id_j=\sum_{k=1}^7{N_{ij}}^kd_k$. There are four (nonzero) solutions
\begin{align*}
    (d_X,d_Y,d_Z,d_U,d_V,d_W)=&(1,\sqrt{2-\sqrt2},\sqrt{2-\sqrt2},1-\sqrt2,1-\sqrt2,-\sqrt{4-2\sqrt2}),\\
    &(1,-\sqrt{2-\sqrt2},-\sqrt{2-\sqrt2},1-\sqrt2,1-\sqrt2,\sqrt{4-2\sqrt2}),\\
    &(1,-\sqrt{2+\sqrt2},-\sqrt{2+\sqrt2},1+\sqrt2,1+\sqrt2,-\sqrt{4+2\sqrt2}),\\
    &(1,\sqrt{2+\sqrt2},\sqrt{2+\sqrt2},1+\sqrt2,1+\sqrt2,\sqrt{4+2\sqrt2})
\end{align*}
with categorical dimensions
\[ D^2(\mc B)=8(2-\sqrt2)(\approx4.7),\quad8(2+\sqrt2), \]
respectively for each pair. Only the last quantum dimensions give unitary MFCs. Each pair has eight conformal dimensions:
\begin{align*}
    (h_X,h_Y,h_Z,h_U,h_V,h_W)=&(\frac12,\frac1{32},\frac1{32},\frac14,\frac34,\frac{5}{32}),(\frac12,\frac7{32},\frac7{32},\frac34,\frac14,\frac{3}{32}),\\
    &(\frac12,\frac{9}{32},\frac{9}{32},\frac14,\frac34,\frac{13}{32}),(\frac12,\frac{15}{32},\frac{15}{32},\frac34,\frac14,\frac{11}{32}),\\
    &(\frac12,\frac{17}{32},\frac{17}{32},\frac14,\frac34,\frac{21}{32}),(\frac12,\frac{23}{32},\frac{23}{32},\frac34,\frac14,\frac{19}{32}),\\
    &(\frac12,\frac{25}{32},\frac{25}{32},\frac14,\frac34,\frac{29}{32}),(\frac12,\frac{31}{32},\frac{31}{32},\frac34,\frac14,\frac{27}{32})\quad(\mods1)
\end{align*}
for the first and second quantum dimensions, and
\begin{align*}
    (h_X,h_Y,h_Z,h_U,h_V,h_W)=&(\frac12,\frac3{32},\frac3{32},\frac34,\frac14,\frac{15}{32}),(\frac12,\frac5{32},\frac5{32},\frac14,\frac34,\frac{25}{32}),\\
    &(\frac12,\frac{11}{32},\frac{11}{32},\frac34,\frac14,\frac{23}{32}),(\frac12,\frac{13}{32},\frac{13}{32},\frac14,\frac34,\frac1{32}),\\&(\frac12,\frac{19}{32},\frac{19}{32},\frac34,\frac14,\frac{31}{32}),(\frac12,\frac{21}{32},\frac{21}{32},\frac14,\frac34,\frac9{32}),\\
    &(\frac12,\frac{27}{32},\frac{27}{32},\frac34,\frac14,\frac7{32}),(\frac12,\frac{29}{32},\frac{29}{32},\frac14,\frac34,\frac{17}{32})\quad(\mods1)
\end{align*}
for the third and fourth quantum dimensions. The $S$-matrices are given by
\[ \widetilde S=\begin{pmatrix}1&d_X&d_Y&d_Z&d_U&d_V&d_W\\d_X&1&-d_Y&-d_Z&d_U&d_V&-d_W\\d_Y&-d_Y&d_W&-d_W&\pm d_Y&\mp d_Y&0\\d_Z&-d_Z&-d_W&d_W&\pm d_Z&\mp d_Z&0\\d_U&d_U&\pm d_Y&\pm d_Z&-1&-1&\mp d_W\\d_V&d_V&\mp d_Y&\mp d_Z&-1&-1&\pm d_W\\d_W&-d_W&0&0&\mp d_W&\pm d_W&0\end{pmatrix}. \]
There are
\[ 4(\text{quantum dimensions})\times8(\text{conformal dimensions})\times2(\text{categorical dimensions})=64 \]
MFCs, among which those 16 with the last quantum dimensions give unitary MFCs. We classify connected étale algebras in all 64 MFCs simultaneously.

An ansatz
\[ A\cong1\oplus n_XX\oplus n_YY\oplus n_ZZ\oplus n_UU\oplus n_VV\oplus n_WW \]
with $n_j\in\mbb N$ has
\[ \fp_{\mc B}(A)=1+n_X+\sqrt{2+\sqrt2}(n_Y+n_Z)+(1+\sqrt2)(n_U+n_V)+\sqrt{4+2\sqrt2}n_W. \]
For this to obey (\ref{FPdimA2bound}), the natural numbers can take only 20 values
\begin{align*}
    (n_X,n_Y,n_Z,n_U,n_V,n_W)=&(0,0,0,0,0,0),(1,0,0,0,0,0),(2,0,0,0,0,0),(3,0,0,0,0,0),\\
    &(4,0,0,0,0,0),(2,1,0,0,0,0),(2,0,1,0,0,0),(1,1,0,0,0,0),\\
    &(1,0,1,0,0,0),(1,0,0,1,0,0),(1,0,0,0,1,0),(1,0,0,0,0,1),\\
    &(0,1,0,0,0,0),(0,2,0,0,0,0),(0,1,1,0,0,0),(0,0,1,0,0,0),\\
    &(0,0,2,0,0,0),(0,0,0,1,0,0),(0,0,0,0,1,0),(0,0,0,0,0,1).
\end{align*}
The first is nothing but the trivial connected étale algebra $A\cong1$ giving $\mc B_A^0\simeq\mc B_A\simeq\mc B$. Those with $X$ do not give commutative algebra because $X$ has $(d_X,h_X)=(1,\frac12)$ (mod 1 for $h$) and $c_{X,X}\cong-id_1$ \cite{KK23preMFC}. The others also fail to be commutative because they contain simple object(s) $b_j\not\cong1$ with nontrivial conformal dimensions.

We conclude
\begin{table}[H]
\begin{center}
\begin{tabular}{c|c|c|c}
    Connected étale algebra $A$&$\mc B_A$&$\rank(\mc B_A)$&Lagrangian?\\\hline
    $1$&$\mc B$&$7$&No
\end{tabular}.
\end{center}
\caption{Connected étale algebras in rank seven MFC $\mcal B\simeq su(2)_6$}\label{rank7su26results}
\end{table}
\hspace{-17pt}That is, all the 64 MFCs $\mc B\simeq su(2)_6$'s are completely anisotropic.

\subsubsection{$\mc B\simeq so(7)_2$}
The MFCs have seven simple objects $\{1,X,Y,Z,U,V,W\}$ obeying monoidal products
\begin{table}[H]
\begin{center}
\begin{tabular}{c|c|c|c|c|c|c|c}
    $\otimes$&$1$&$X$&$Y$&$Z$&$U$&$V$&$W$\\\hline
    $1$&$1$&$X$&$Y$&$Z$&$U$&$V$&$W$\\\hline
    $X$&&$1$&$Y$&$Z$&$U$&$W$&$V$\\\hline
    $Y$&&&$1\oplus X\oplus U$&$Z\oplus U$&$Y\oplus Z$&$V\oplus W$&$V\oplus W$\\\hline
    $Z$&&&&$1\oplus X\oplus Y$&$Y\oplus U$&$V\oplus W$&$V\oplus W$\\\hline
    $U$&&&&&$1\oplus X\oplus Z$&$V\oplus W$&$V\oplus W$\\\hline
    $V$&&&&&&$1\oplus Y\oplus Z\oplus U$&$X\oplus Y\oplus Z\oplus U$\\\hline
    $W$&&&&&&&$1\oplus Y\oplus Z\oplus U$
\end{tabular}.
\end{center}
\end{table}
\hspace{-17pt}Thus, they have
\begin{align*}
    \fp_{\mc B}(1)=1=\fp_{\mc B}(X),\quad&\fp_{\mc B}(Y)=\fp_{\mc B}(Z)=\fp_{\mc B}(U)=2,\\
    \fp_{\mc B}(V)=&\sqrt7=\fp_{\mc B}(W),
\end{align*}
and
\[ \fp(\mc B)=28. \]
Their quantum dimensions $d_j$'s are solutions of the same multiplication rules $d_id_j=\sum_{k=1}^7{N_{ij}}^kd_k$. There are two (nonzero) solutions
\[ (d_X,d_Y,d_Z,d_U,d_V,d_W)=(1,2,2,2,-\sqrt7,-\sqrt7),(1,2,2,2,\sqrt7,\sqrt7). \]
Only the second quantum dimensions give unitary MFCs. They both have the categorical dimension
\[ D^2(\mc B)=28. \]
They have four conformal dimensions\footnote{Naively, one finds 24 conformal dimensions, but the others are equivalent to the four in the main text by permutation $(VW)$ or $(YZU)$ of simple objects. For example, one also finds $(h_X,h_Y,h_Z,h_U,h_V,h_W)=(0,\frac17,\frac27,\frac47,\frac58,\frac18)$ (mod 1) is consistent, but it is the same as our first conformal dimension under a permutation $(VW)$.}
\[ \hspace{-30pt}(h_X,h_Y,h_Z,h_U,h_V,h_W)=(0,\frac17,\frac27,\frac47,\frac18,\frac58),(0,\frac17,\frac27,\frac47,\frac38,\frac78),(0,\frac37,\frac67,\frac57,\frac18,\frac58),(0,\frac37,\frac67,\frac57,\frac38,\frac78).\quad(\mods1) \]
The $S$-matrices are given by
\[ \widetilde S=\begin{pmatrix}1&1&2&2&2&d_V&d_W\\1&1&2&2&2&-d_V&-d_W\\2&2&s&s'&s''&0&0\\2&2&s'&s''&s&0&0\\2&2&s''&s&s'&0&0\\d_V&-d_V&0&0&0&\pm d_V&\mp d_W\\d_W&-d_W&0&0&0&\mp d_W&\pm d_W\end{pmatrix} \]
with
\[ s=-4\sin\frac\pi{14},\quad s'=4\sin\frac{3\pi}{14},\quad s''=-4\cos\frac\pi7, \]
or their permutations $(ss's'')$. They have additive central charges
\[ c(\mc B)=\begin{cases}2&(\text{1st\&2nd }h),\\-2&(\text{3rd\&4th }h).\end{cases}\quad(\mods8) \]
There are
\[ 2(\text{quantum dimensions})\times4(\text{conformal dimensions})\times2(\text{categorical dimensions})=16 \]
MFCs, among which those eight with the second quantum dimensions give unitary MFCs. We study connected étale algebras in all 16 MFCs simultaneously.

We work with an ansatz
\[ A\cong1\oplus n_XX\oplus n_YY\oplus n_ZZ\oplus n_UU\oplus n_VV\oplus n_WW \]
with $n_j\in\mbb N$. It has
\[ \fp_{\mc B}(A)=1+n_X+2(n_Y+n_Z+n_U)+\sqrt7(n_V+n_W). \]
For this to obey (\ref{FPdimA2bound}), the natural numbers can take only 24 values 
\begin{align*}
    (n_X,n_Y,n_Z,n_U,n_V,n_W)=&(0,0,0,0,0,0),(1,0,0,0,0,0),(2,0,0,0,0,0),(3,0,0,0,0,0),\\
    &(4,0,0,0,0,0),(2,1,0,0,0,0),(2,0,1,0,0,0),(2,0,0,1,0,0),\\
    &(1,1,0,0,0,0),(1,0,1,0,0,0),(1,0,0,1,0,0),(1,0,0,0,1,0),\\
    &(1,0,0,0,0,1),(0,1,0,0,0,0),(0,2,0,0,0,0),(0,1,1,0,0,0),\\
    &(0,1,0,1,0,0),(0,0,1,0,0,0),(0,0,2,0,0,0),(0,0,1,1,0,0),\\
    &(0,0,0,1,0,0),(0,0,0,2,0,0),(0,0,0,0,1,0),(0,0,0,0,0,1).
\end{align*}
The first is nothing but the trivial connected étale algebra $A\cong1$ giving $\mc B_A^0\simeq\mc B_A\simeq\mc B$. Next, the 19 candidates with $Y,Z,U,V,W$ have nontrivial conformal dimensions, and they fail to be commutative. Thus, we are left with those four with only $X$'s:
\[ n_X=1,2,3,4. \]
All but the first solution is ruled out by studying Frobenius-Perron dimension. The latter three have $\fp_{\mc B}(A)=3,4,5$, and demand $\fp(\mc B_A^0)=\frac{28}9,\frac74,\frac{28}{25}$. However, there is no MFC with such Frobenius-Perron dimensions. Thus, the candidates are ruled out. We are left with just $A\cong1\oplus X$.

The candidate is indeed a commutative algebra by the lemma 1 and $c_{X,X}\cong id_1$ because $X$ has $(d_X,h_X)=(1,0)$ \cite{KK23preMFC}. Furthermore, it turns out to be separable, hence connected étale. Let us check this fact by identifying $\mc B_A$.

It has $\fp_{\mc B}(A)=2$, and demands
\[ \fp(\mc B_A^0)=7,\quad\fp(\mc B_A)=14. \]
We can identify the MFC as
\[ \mc B_A^0\simeq\vecG_{\mbb Z/7\mbb Z}^1 \]
from its Frobenius-Perron dimension. Taking the invariance of topological twists (\ref{invtopologicaltwist}) into account, we see central charges are also matched. What is a fusion category $\mc B_A$ containing the rank seven MFC $\vecG_{\mbb Z/7\mbb Z}^1$? It turns out that the category of right $A$-modules is
\[ \mc B_A\simeq\text{TY}(\mbb Z/7\mbb Z), \]
a $\mbb Z/7\mbb Z$ Tambara-Yamagami category \cite{TY98}. One of the easiest ways to see this fact is to perform anyon condensation. It `identifies' $1$ and $X$, and hence $V$ and $W$. The other invariant simple objects $Y,Z,U$ `split' into two each. As a result, we get seven invertible objects. The deconfined seven particles have the same conformal dimensions (mod 1) as simple objects in $\vecG_{\mbb Z/7\mbb Z}^1$ MFCs. Thus, they form $\mc B_A^0\simeq\vecG_{\mbb Z/7\mbb Z}^1$. The category $\mc B_A$ of right $A$-modules have one additional non-invertible simple object with quantum dimension $\pm\sqrt7$. This is nothing but a $\mbb Z/7\mbb Z$ Tambara-Yamagami category.\footnote{More rigorously, we have to find NIM-reps. Indeed, we find an eight-dimensional NIM-rep
\begin{align*}
    n_1=1_8=n_X,\quad n_Y=\begin{pmatrix}0&1&1&0&0&0&0&0\\1&0&0&0&0&1&0&0\\1&0&0&0&0&0&1&0\\0&0&0&0&1&1&0&0\\0&0&0&1&0&0&1&0\\0&1&0&1&0&0&0&0\\0&0&1&0&1&0&0&0\\0&0&0&0&0&0&0&2\end{pmatrix},\quad n_Z=\begin{pmatrix}0&0&0&1&1&0&0&0\\0&0&0&0&1&0&1&0\\0&0&0&1&0&1&0&0\\1&0&1&0&0&0&0&0\\1&1&0&0&0&0&0&0\\0&0&1&0&0&0&1&0\\0&1&0&0&0&1&0&0\\0&0&0&0&0&0&0&2\end{pmatrix},\\
    n_U=\begin{pmatrix}0&0&0&0&0&1&1&0\\0&0&1&1&0&0&0&0\\0&1&0&0&1&0&0&0\\0&1&0&0&0&0&1&0\\0&0&1&0&0&1&0&0\\1&0&0&0&1&0&0&0\\1&0&0&1&0&0&0&0\\0&0&0&0&0&0&0&2\end{pmatrix},\quad n_V=\begin{pmatrix}0&0&0&0&0&0&0&1\\0&0&0&0&0&0&0&1\\0&0&0&0&0&0&0&1\\0&0&0&0&0&0&0&1\\0&0&0&0&0&0&0&1\\0&0&0&0&0&0&0&1\\0&0&0&0&0&0&0&1\\1&1&1&1&1&1&1&0\end{pmatrix}=n_W.
\end{align*}
Denoting a basis of $\mc B_A$ by $\{m_1,m_2,m_3,m_4,m_5,m_6,m_7,m_8\}$, we obtain a multiplication table
\begin{table}[H]
\begin{center}
\makebox[1 \textwidth][c]{       
\resizebox{1.0 \textwidth}{!}{\begin{tabular}{c|c|c|c|c|c|c|c|c}
    $b_j\otimes\backslash$&$m_1$&$m_2$&$m_3$&$m_4$&$m_5$&$m_6$&$m_7$&$m_8$\\\hline
    $1,X$&$m_1$&$m_2$&$m_3$&$m_4$&$m_5$&$m_6$&$m_7$&$m_8$\\
    $Y$&$m_2\oplus m_3$&$m_1\oplus m_6$&$m_1\oplus m_7$&$m_5\oplus m_6$&$m_4\oplus m_7$&$m_2\oplus m_4$&$m_3\oplus m_5$&$2m_8$\\
    $Z$&$m_4\oplus m_5$&$m_5\oplus m_7$&$m_4\oplus m_6$&$m_1\oplus m_3$&$m_1\oplus m_2$&$m_3\oplus m_7$&$m_2\oplus m_6$&$2m_8$\\
    $U$&$m_6\oplus m_7$&$m_3\oplus m_4$&$m_2\oplus m_5$&$m_2\oplus m_7$&$m_3\oplus m_6$&$m_1\oplus m_5$&$m_1\oplus m_4$&$2m_8$\\
    $V,W$&$m_8$&$m_8$&$m_8$&$m_8$&$m_8$&$m_8$&$m_8$&$\bigoplus_{j=1}^7m_j$
\end{tabular}.}}
\end{center}
\end{table}
\hspace{-11pt}From the multiplication rules, in this basis, we can identify
\[ m_1\cong1\oplus X,\quad m_2\cong Y\cong m_3,\quad m_4\cong Z\cong m_5,\quad m_6\cong U\cong m_7,\quad m_8\cong V\oplus W. \]
In the category of right $A$-modules (or broken phase) $\mc B_A$, they have quantum dimensions (\ref{dBAm})
\[ \hspace{-50pt}d_{\mc B_A}(m_1)=d_{\mc B_A}(m_2)=d_{\mc B_A}(m_3)=d_{\mc B_A}(m_4)=d_{\mc B_A}(m_5)=d_{\mc B_A}(m_6)=d_{\mc B_A}(m_7)=1,\quad d_{\mc B_A}(m_8)=\pm\sqrt7. \]
Furthermore, employing the free module functor $F_A$ (\ref{FA}), one can compute monoidal products $\otimes_A$:
\[ m_j\otimes_Am_8\cong m_8\cong m_8\otimes_Am_j\ (j=1,2,\dots,7),\quad m_8\otimes_Am_8\cong\bigoplus_{j=1}^7m_j. \]
This shows $\mc B_A\simeq\text{TY}(\mbb Z/7\mbb Z)$.} Note that this example has
\[ \rank(\mc B_A)>\rank(\mc B). \]
While $\mc B_A$ consists of objects of $\mc B$ and the former has Frobenius-Perron dimension no larger than the latter, in general, the former can have larger rank as in this example.

We conclude
\begin{table}[H]
\begin{center}
\begin{tabular}{c|c|c|c}
    Connected étale algebra $A$&$\mc B_A$&$\rank(\mc B_A)$&Lagrangian?\\\hline
    $1$&$\mc B$&$7$&No\\
    $1\oplus X$&$\text{TY}(\mbb Z/7\mbb Z)$&8&No
\end{tabular}.
\end{center}
\caption{Connected étale algebras in rank seven MFC $\mcal B\simeq so(7)_2$}\label{rank7so72results}
\end{table}
\hspace{-17pt}All the 16 MFCs $\mc B\simeq so(7)_2$'s fail to be completely anisotropic.

\subsubsection{$\mc B\simeq psu(2)_{13}$}
The MFCs have seven simple objects $\{1,X,Y,Z,U,V,W\}$ obeying monoidal products
\begin{table}[H]
\begin{center}
\makebox[1 \textwidth][c]{       
\resizebox{1.2 \textwidth}{!}{\begin{tabular}{c|c|c|c|c|c|c|c}
    $\otimes$&$1$&$X$&$Y$&$Z$&$U$&$V$&$W$\\\hline
    $1$&$1$&$X$&$Y$&$Z$&$U$&$V$&$W$\\\hline
    $X$&&$1\oplus Y$&$X\oplus Z$&$Y\oplus U$&$Z\oplus V$&$U\oplus W$&$V\oplus W$\\\hline
    $Y$&&&$1\oplus Y\oplus U$&$X\oplus Z\oplus V$&$Y\oplus U\oplus W$&$Z\oplus V\oplus W$&$U\oplus V\oplus W$\\\hline
    $Z$&&&&$1\oplus Y\oplus U\oplus W$&$X\oplus Z\oplus V\oplus W$&$Y\oplus U\oplus V\oplus W$&$Z\oplus U\oplus V\oplus W$\\\hline
    $U$&&&&&$1\oplus Y\oplus U\oplus V\oplus W$&$X\oplus Z\oplus U\oplus V\oplus W$&$Y\oplus Z\oplus U\oplus V\oplus W$\\\hline
    $V$&&&&&&$1\oplus Y\oplus Z\oplus U\oplus V\oplus W$&$X\oplus Y\oplus Z\oplus U\oplus V\oplus W$\\\hline
    $W$&&&&&&&$1\oplus X\oplus Y\oplus Z\oplus U\oplus V\oplus W$
\end{tabular}.}}
\end{center}
\end{table}
\hspace{-17pt}Thus, they have
\begin{align*}
    \fp_{\mc B}(1)=1,\quad\fp_{\mc B}(X)=&\frac{\sin\frac{2\pi}{15}}{\sin\frac\pi{15}},\quad\fp_{\mc B}(Y)=\frac{\sin\frac{3\pi}{15}}{\sin\frac\pi{15}},\quad\fp_{\mc B}(Z)=\frac{\sin\frac{4\pi}{15}}{\sin\frac\pi{15}},\\
    \fp_{\mc B}(U)=\frac{\sin\frac{5\pi}{15}}{\sin\frac\pi{15}},&\quad\fp_{\mc B}(V)=\frac{\sin\frac{6\pi}{15}}{\sin\frac\pi{15}},\quad\fp_{\mc B}(W)=\frac{\sin\frac{7\pi}{15}}{\sin\frac\pi{15}},
\end{align*}
and
\[ \fp(\mc B)=\frac{15}{4\sin^2\frac\pi{15}}\approx86.8. \]
Their quantum dimensions are solutions of the same multiplication rules $d_id_j=\sum_{k=1}^7{N_{ij}}^kd_k$. There are four (nonzero) solutions
\begin{align*}
    (d_X,d_Y,d_Z,d_U,d_V,d_W)=&(\frac{\sin\frac{\pi}{15}}{\cos\frac\pi{30}},-\frac{\sin\frac{6\pi}{15}}{\cos\frac\pi{30}},-\frac{\sin\frac{2\pi}{15}}{\cos\frac\pi{30}},\frac{\sin\frac{5\pi}{15}}{\cos\frac\pi{30}},\frac{\sin\frac{3\pi}{15}}{\cos\frac\pi{30}},-\frac{\sin\frac{4\pi}{15}}{\cos\frac\pi{30}}),\\
    &(-\frac{\sin\frac{7\pi}{15}}{\cos\frac{7\pi}{30}},\frac{\sin\frac{3\pi}{15}}{\cos\frac{7\pi}{30}},\frac{\sin\frac{\pi}{15}}{\cos\frac{7\pi}{30}},-\frac{\sin\frac{5\pi}{15}}{\cos\frac{7\pi}{30}},\frac{\sin\frac{6\pi}{15}}{\cos\frac{7\pi}{30}},-\frac{\sin\frac{2\pi}{15}}{\cos\frac{7\pi}{30}}),\\
    &(-\frac{\sin\frac{4\pi}{15}}{\cos\frac{11\pi}{30}},\frac{\sin\frac{6\pi}{15}}{\cos\frac{11\pi}{30}},-\frac{\sin\frac{7\pi}{15}}{\cos\frac{11\pi}{30}},\frac{\sin\frac{5\pi}{15}}{\cos\frac{11\pi}{30}},-\frac{\sin\frac{3\pi}{15}}{\cos\frac{11\pi}{30}},\frac{\sin\frac{\pi}{15}}{\cos\frac{11\pi}{30}}),\\
    &(\frac{\sin\frac{2\pi}{15}}{\sin\frac\pi{15}},\frac{\sin\frac{3\pi}{15}}{\sin\frac\pi{15}},\frac{\sin\frac{4\pi}{15}}{\sin\frac\pi{15}},\frac{\sin\frac{5\pi}{15}}{\sin\frac\pi{15}},\frac{\sin\frac{6\pi}{15}}{\sin\frac\pi{15}},\frac{\sin\frac{7\pi}{15}}{\sin\frac\pi{15}})
\end{align*}
with categorical dimensions
\[ D^2(\mc B)=\frac{15}{4\cos^2\frac\pi{30}}(\approx3.8),\quad\frac{15}{4\cos^2\frac{7\pi}{30}}(\approx6.8),\quad\frac{15}{4\cos^2\frac{11\pi}{30}}(\approx22.7),\quad\frac{15}{4\sin^2\frac\pi{15}}, \]
respectively. Each quantum dimension has two conformal dimensions
\[ (h_X,h_Y,h_Z,h_U,h_V,h_W)=\begin{cases}(\frac25,\frac1{15},0,\frac15,\frac23,\frac25),(\frac35,\frac{14}{15},0,\frac45,\frac13,\frac35)&(\text{1st quantum dimension}),\\(\frac15,\frac8{15},0,\frac35,\frac13,\frac15),(\frac45,\frac7{15},0,\frac25,\frac23,\frac45)&(\text{2nd quantum dimension}),\\
(\frac25,\frac{11}{15},0,\frac15,\frac13,\frac25),(\frac35,\frac4{15},0,\frac45,\frac23,\frac35)&(\text{3rd quantum dimension}),\\(\frac15,\frac{13}{15},0,\frac35,\frac23,\frac15),(\frac45,\frac2{15},0,\frac25,\frac13,\frac45)&(\text{4th quantum dimension}).\end{cases}\quad(\mods1) \]
The $S$-matrices are given by
\[ \widetilde S=\begin{pmatrix}1&d_X&d_Y&d_Z&d_U&d_V&d_W\\d_X&-d_Z&d_V&-d_W&d_U&-d_Y&1\\d_Y&d_V&d_V&d_Y&0&-d_Y&-d_V\\d_Z&-d_W&d_Y&1&-d_U&d_V&-d_X\\d_U&d_U&0&-d_U&-d_U&0&d_U\\d_V&-d_Y&-d_Y&d_V&0&-d_V&d_Y\\d_W&1&-d_V&-d_X&d_U&d_Y&-d_Z\end{pmatrix}.\]
There are
\[ 4(\text{quantum dimensions})\times2(\text{conformal dimensions})\times2(\text{categorical dimensions})=16 \]
MFCs, among which those four with the last quantum dimensions give unitary MFCs. We study connected étale algebras in all 16 MFCs simultaneously.

An ansatz
\[ A\cong1\oplus n_XX\oplus n_YY\oplus n_ZZ\oplus n_UU\oplus n_VV\oplus n_WW \]
with $n_j\in\mbb N$ has
\[ \fp_{\mc B}(A)=1+\frac{\sin\frac{2\pi}{15}}{\sin\frac\pi{15}}n_X+\frac{\sin\frac{3\pi}{15}}{\sin\frac\pi{15}}n_Y+\frac{\sin\frac{4\pi}{15}}{\sin\frac\pi{15}}n_Z+\frac{\sin\frac{5\pi}{15}}{\sin\frac\pi{15}}n_U+\frac{\sin\frac{6\pi}{15}}{\sin\frac\pi{15}}n_V+\frac{\sin\frac{7\pi}{15}}{\sin\frac\pi{15}}n_W. \]
For this to obey (\ref{FPdimA2bound}), the natural number can take only 27 values
\begin{align*}
    (n_X,n_Y,n_Z,n_U,n_V,n_W)=&(0,0,0,0,0,0),(1,0,0,0,0,0),(2,0,0,0,0,0),\\
    &(3,0,0,0,0,0),(4,0,0,0,0,0),(2,1,0,0,0,0),\\
    &(2,0,1,0,0,0),(2,0,0,1,0,0),(1,1,0,0,0,0),\\
    &(1,2,0,0,0,0),(1,0,1,0,0,0),(1,0,0,1,0,0),\\
    &(1,0,0,0,1,0),(1,0,0,0,0,1),(0,1,0,0,0,0),\\
    &(0,2,0,0,0,0),(0,1,1,0,0,0),(0,1,0,1,0,0),\\
    &(0,1,0,0,1,0),(0,1,0,0,0,1),(0,0,1,0,0,0),\\
    &(0,0,2,0,0,0),(0,0,1,1,0,0),(0,0,1,0,1,0),\\
    &(0,0,0,1,0,0),(0,0,0,0,1,0),(0,0,0,0,0,1).
\end{align*}
The first is nothing but the trivial connected étale algebra $A\cong1$ giving $\mc B_A^0\simeq\mc B_A\simeq\mc B$. Others with nontrivial simple objects $b_j\not\cong Z$ have nontrivial conformal dimensions, and do not give commutative algebras. Since $Z$ has trivial conformal dimension $h_Z=0$ mod 1, it can give commutative algebra. However, those with $Z$'s are also ruled out as follows. Apart from the trivial one, the natural number $n_Z$ can take one or two. They have $\fp_{\mc B}=1+\frac{\sin\frac{4\pi}{15}}{\sin\frac\pi{15}},1+2\frac{\sin\frac{4\pi}{15}}{\sin\frac\pi{15}}$, and demand $\fp(\mc B_A^0)\approx4.1,1.3$, but there are no MFC with these Frobenius-Perron dimensions. Thus, the two candidates are ruled out.

We conclude
\begin{table}[H]
\begin{center}
\begin{tabular}{c|c|c|c}
    Connected étale algebra $A$&$\mc B_A$&$\rank(\mc B_A)$&Lagrangian?\\\hline
    $1$&$\mc B$&$7$&No
\end{tabular}.
\end{center}
\caption{Connected étale algebras in rank seven MFC $\mcal B\simeq psu(2)_{13}$}\label{rank7psu213results}
\end{table}
\hspace{-17pt}Namely, all the 16 MFCs $\mc B\simeq psu(2)_{13}$'s are completely anisotropic.

\subsection{Rank eight}
\subsubsection{$\mc B\simeq\vecG_{\mbb Z/2\mbb Z\times\mbb Z/2\mbb Z\times\mbb Z/2\mbb Z}^\alpha\simeq\vecG_{\mbb Z/2\mbb Z}^{-1}\boxtimes\tc$}\label{Z2Z2Z2}
The MFCs have eight simple objects $\{1,X,Y,Z,T,U,V,W\}$ obeying monoidal products
\begin{table}[H]
\begin{center}
\begin{tabular}{c|c|c|c|c|c|c|c|c}
    $\otimes$&$1$&$X$&$Y$&$Z$&$T$&$U$&$V$&$W$\\\hline
    $1$&$1$&$X$&$Y$&$Z$&$T$&$U$&$V$&$W$\\\hline
    $X$&&$1$&$W$&$U$&$V$&$Z$&$T$&$Y$\\\hline
    $Y$&&&$1$&$V$&$U$&$T$&$Z$&$X$\\\hline
    $Z$&&&&$1$&$W$&$X$&$Y$&$T$\\\hline
    $T$&&&&&$1$&$Y$&$X$&$Z$\\\hline
    $U$&&&&&&$1$&$W$&$V$\\\hline
    $V$&&&&&&&$1$&$U$\\\hline
    $W$&&&&&&&&$1$
\end{tabular}.
\end{center}
\end{table}
\hspace{-17pt}Thus, they have
\begin{align*}
    &\fp_{\mc B}(1)=\fp_{\mc B}(X)=\fp_{\mc B}(Y)=\fp_{\mc B}(Z)\\
    &=\fp_{\mc B}(T)=\fp_{\mc B}(U)=\fp_{\mc B}(V)=\fp_{\mc B}(W)=1,
\end{align*}
and
\[ \fp(\mc B)=8. \]
Their quantum dimensions $d_j$'s are solutions of the same multiplication rules $d_id_j=\sum_{k=1}^8{N_{ij}}^kd_k$. There are eight solutions
\begin{align*}
    (d_X,d_Y,d_Z,d_T,d_U,d_V,d_W)=&(-1,-1,-1,-1,1,1,1),(-1,-1,1,1,-1,-1,1),\\
    &(-1,1,-1,1,1,-1,-1),(-1,1,1,-1,-1,1,-1),\\
    &(1,-1,-1,1,-1,1,-1),(1,-1,1,-1,1,-1,-1),\\
    &(1,1,-1,-1,-1,-1,1),(1,1,1,1,1,1,1).
\end{align*}
Only the last solution gives unitary MFCs. All quantum dimensions have the same categorical dimension
\[ D^2(\mc B)=8. \]

In order to count the number of MFCs, let us list up conformal dimensions. Naively, the MFCs can have two structures $\vecG_{\mbb Z/2\mbb Z}^{-1}\boxtimes\vecG_{\mbb Z/2\mbb Z}^{-1}\boxtimes\vecG_{\mbb Z/2\mbb Z}^{-1}$ or $\vecG_{\mbb Z/2\mbb Z}^{-1}\boxtimes\tc$. However, by computing conformal dimensions, one finds they all have the structure $\vecG_{\mbb Z/2\mbb Z}^{-1}\boxtimes\tc$. (We have already specified this fact in the name.) In order to avoid double-counting, we choose $\vecG_{\mbb Z/2\mbb Z}^{-1}=\{1,X\}$ and $\tc=\{1,U,V,W\}$. Then the other simple objects are given by
\[ Y\cong X\otimes W,\quad Z\cong X\otimes U,\quad T\cong X\otimes V. \]
We have $h_X=\frac14,\frac34$ (mod 1) and $(h_U,h_V,h_W)=(0,0,\frac12),(0,\frac12,0),(\frac12,0,0),(\frac12,\frac12,\frac12)$ (mod 1). Depending on quantum dimensions, we have different symmetries.

\paragraph{$(d_X,d_U,d_V,d_W)=(1,1,1,1)$.} This quantum dimension gives unitary MFCs. Even after our choice, we still have permutations $(UV),(UW),(VW)$ of simple objects. Different MFCs are given by
\[ (h_X,h_U,h_V,h_W)=(\frac14,0,0,\frac12),(\frac14,\frac12,\frac12,\frac12),(\frac34,0,0,\frac12),(\frac34,\frac12,\frac12,\frac12)\quad(\mods1). \]
Including the two signs of categorical dimensions, we have eight unitary MFCs.

\paragraph{$(d_X,d_U,d_V,d_W)=(1,1,-1,-1)$.} We have smaller symmetry $(VW)$. Thus, different MFCs are given by
\begin{align*}
    (h_X,h_U,h_V,h_W)&=(\frac14,0,0,\frac12),(\frac14,\frac12,0,0),(\frac14,\frac12,\frac12,\frac12),\\
    &~~~~(\frac34,0,0,\frac12),(\frac34,\frac12,0,0),(\frac34,\frac12,\frac12,\frac12)\quad(\mods1).
\end{align*}
With two signs of categorical dimensions, we have 12 MFCs.

\paragraph{$(d_X,d_U,d_V,d_W)=(-1,1,1,1)$.} Different MFCs are given by
\[ (h_X,h_U,h_V,h_W)=(\frac14,0,0,\frac12),(\frac14,\frac12,\frac12,\frac12),(\frac34,0,0,\frac12),(\frac34,\frac12,\frac12,\frac12)\quad(\mods1). \]
There are eight MFCs.

\paragraph{$(d_X,d_U,d_V,d_W)=(-1,1,-1,-1)$.} Different MFCs are given by
\begin{align*}
    (h_X,h_U,h_V,h_W)&=(\frac14,0,0,\frac12),(\frac14,\frac12,0,0),(\frac14,\frac12,\frac12,\frac12),\\
    &~~~~(\frac34,0,0,\frac12),(\frac34,\frac12,0,0),(\frac34,\frac12,\frac12,\frac12)\quad(\mods1).
\end{align*}
There are 12 MFCs.\newline

In total, there are
\[ 8+12+8+12=40 \]
MFCs, among which those eight in the first case are unitary. The $S$-matrices are given by
\[ \widetilde S=\begin{pmatrix}1&d_X&d_Xd_W&d_Xd_U&d_Xd_V&d_U&d_V&d_W\\d_X&-1&-d_W&-d_U&-d_V&d_Xd_U&d_Xd_V&d_Xd_W\\d_Xd_W&-d_W&-1&d_W&1&-d_Xd_W&-d_X&d_X\\d_Xd_U&-d_U&d_W&-1&d_V&d_X&-d_Xd_V&-d_Xd_W\\d_Xd_V&-d_V&1&d_V&-1&-d_Xd_V&d_X&-d_X\\d_U&d_Xd_U&-d_Xd_W&d_X&-d_Xd_V&1&-d_V&-d_W\\d_V&d_Xd_V&-d_X&-d_Xd_V&d_X&-d_V&1&-1\\d_W&d_Xd_W&d_X&-d_Xd_W&-d_X&-d_W&-1&1\end{pmatrix}. \]
They have additive central charges
\[ c(\mc B)=c(\vecG_{\mbb Z/2\mbb Z}^{-1})+c(\tc)\quad(\mods8), \]
where
\[ c(\vecG_{\mbb Z/2\mbb Z}^{-1})=\begin{cases}1&(h_X=\frac14),\\-1&(h_X=\frac34),\end{cases}\quad c(\tc)=\begin{cases}0&(\text{one }h=\frac12\&\text{the other }h=0),\\4&(\text{all }h=\frac12).\end{cases}\quad(\mods8) \]
We classify connected étale algebras in all 40 MFCs simultaneously.

An ansatz
\[ A\cong1\oplus n_XX\oplus n_YY\oplus n_ZZ\oplus n_TT\oplus n_UU\oplus n_VV\oplus n_WW \]
with $n_j\in\mbb N$ has
\[ \fp_{\mc B}(A)=1+n_X+n_Y+n_Z+n_T+n_U+n_V+n_W. \]
For this to obey (\ref{FPdimA2bound}), the natural numbers can take only eight values
\begin{align*}
    (n_X,n_Y,n_Z,n_T,n_U,n_V,n_W)=&(0,0,0,0,0,0,0),(1,0,0,0,0,0,0),(0,1,0,0,0,0,0),(0,0,1,0,0,0,0),\\
    &(0,0,0,1,0,0,0),(0,0,0,0,1,0,0),(0,0,0,0,0,1,0),(0,0,0,0,0,0,1).
\end{align*}
The first is nothing but the trivial connected étale algebra $A\cong1$ giving $\mc B_A^0\simeq\mc B_A\simeq\mc B$. Among the other seven candidates, the second with $X$ cannot be commutative due to our choice $\vecG_{\mbb Z/2\mbb Z}^{-1}=\{1,X\}$. Accordingly, those three with $Y,Z,T$ cannot be commutative either because they have $h=\frac14$ (mod $\frac12$). Thus, we are left with three nontrivial candidates
\[ 1\oplus U,\quad1\oplus V,\quad1\oplus W. \]
These are $\mbb Z/2\mbb Z$ algebras by the lemma 1, and they can also be commutative depending on quantum and conformal dimensions. A $\mbb Z/2\mbb Z$ algebra is commutative iff the nontrivial simple object has $(d,h)=(1,0),(-1,\frac12)$ \cite{KK23preMFC}. Thus, the $\mbb Z/2\mbb Z$ algebras are commutative when
\begin{equation}
\begin{split}
    \hspace{-40pt}c_{U,U}\cong id_1\iff(d_X,d_U,d_V,d_W,h_X,h_U,h_V,h_W)&=(1,1,1,1,\frac14,0,0,\frac12),(1,1,1,1,\frac34,0,0,\frac12),\\
    &~~~~(1,1,-1,-1,\frac14,0,0,\frac12),(1,1,-1,-1,\frac34,0,0,\frac12),\\
    &~~~~(-1,1,1,1,\frac14,0,0,\frac12),(-1,1,1,1,\frac34,0,0,\frac12),\\
    &~~~~(-1,1,-1,-1,\frac14,0,0,\frac12),(-1,1,-1,-1,\frac34,0,0,\frac12),\\
    \hspace{-40pt}c_{V,V}\cong id_1\iff(d_X,d_U,d_V,d_W,h_X,h_U,h_V,h_W)&=(1,1,1,1,\frac14,0,0,\frac12),(1,1,1,1,\frac34,0,0,\frac12),\\
    &~~~~(1,1,-1,-1,\frac14,\frac12,\frac12,\frac12),(1,1,-1,-1,\frac34,\frac12,\frac12,\frac12),\\
    &~~~~(-1,1,1,1,\frac14,0,0,\frac12),(-1,1,1,1,\frac34,0,0,\frac12),\\
    &~~~~(-1,1,-1,-1,\frac14,\frac12,\frac12,\frac12),(-1,1,-1,-1,\frac34,\frac12,\frac12,\frac12),\\
    \hspace{-40pt}c_{W,W}\cong id_1\iff(d_X,d_U,d_V,d_W,h_X,h_U,h_V,h_W)&=(1,1,-1,-1,\frac14,0,0,\frac12),(1,1,-1,-1,\frac14,\frac12,\frac12,\frac12),\\
    &~~~~(1,1,-1,-1,\frac34,0,0,\frac12),(1,1,-1,-1,\frac34,\frac12,\frac12,\frac12),\\
    &~~~~(-1,1,-1,-1,\frac14,0,0,\frac12),(-1,1,-1,-1,\frac14,\frac12,\frac12,\frac12),\\
    &~~~~(-1,1,-1,-1,\frac34,0,0,\frac12),(-1,1,-1,-1,\frac34,\frac12,\frac12,\frac12).\\
    &\hspace{200pt}(\text{mod }1\text{ for }h)
\end{split}\label{commutativeAZ2Z2Z2}
\end{equation}
Note that not all of these are separable. The obstruction is in quantum dimensions. For example, let us look at the third commutative algebra of $1\oplus V$. If we naively use the formula (\ref{dBAm}), a right $A$-module $m\cong1\oplus V$ has $d_{\mc B_A}(m)=0$, a contradiction.\footnote{In view of anyon condensation, this fact leads to the following observation. In this case, $1$ and $V$ which are naively `identified' under anyon condensation have different quantum dimensions. While anyon condensation can `identify' two (or more) simple objects with different conformal dimensions -- in this case, the resulting object is confined -- it may not be allowed to `identify' simple objects with quantum dimensions of different signs. Since most literature on anyon condensation studies only unitary MFCs, it seems this problem was unknown.}  This line of reasoning rules out the second and fourth lines of $1\oplus V$, and all of $1\oplus W$. Therefore, we find three connected étale algebras
\begin{equation}
    A\cong\begin{cases}1&(\text{all MFCs}),\\1\oplus U&(\text{those in }(\ref{commutativeAZ2Z2Z2})),\\1\oplus V&(d_U,d_V,d_W,h_U,h_V,h_W)=(1,1,1,0,0,\frac12).\end{cases}\quad(\mods1\text{ for }h)\label{Z2Z2Z2etale}
\end{equation}

Let us determine the categories of (dyslectic) right $A$-modules $\mc B_A^0,\mc B_A$. The nontrivial connected étale algebras have $\fp_{\mc B}(A)=2$ and demand
\[ \fp(\mc B_A^0)=2,\quad\fp(\mc B_A)=4. \]
Since the only MFC with $\fp=2$ is $\vecG_{\mbb Z/2\mbb Z}^{-1}$, we get
\[ \mc B_A^0\simeq\vecG_{\mbb Z/2\mbb Z}^{-1}. \]
This identification also matches central charges because when the algebras are separable, we have $c(\tc)=0$ mod 8 and $c(\mc B_A^0)$ is determined by $h_X$ as we will see. The category $\mc B_A$ of right $A$-modules contain this MFC as a subcategory and has $\fp=4$. There are three candidates, $\ising$, $\vecG_{\mbb Z/2\mbb Z\times\mbb Z/2\mbb Z}^\alpha$, and $\vecG_{\mbb Z/4\mbb Z}^\alpha$. Here, note that all simple objects of $\mc B$ has $\fp_{\mc B}=1$ because the free module functor preserves Frobenius-Perron dimensions (\ref{FPpreserve}). This observation rules out $\ising$, and we are left with rank four candidates $\vecG_{\mbb Z/2\mbb Z\times\mbb Z/2\mbb Z}^\alpha,\vecG_{\mbb Z/4\mbb Z}^\alpha$. In order to find out the correct category, we search for four-dimensional NIM-reps. We start from $A\cong1\oplus U$. We find a four-dimensional NIM-rep
\[ \hspace{-20pt}n_1=1_4=n_U,\quad n_X=\begin{pmatrix}0&0&0&1\\0&0&1&0\\0&1&0&0\\1&0&0&0\end{pmatrix}=n_Z,\quad n_Y=\begin{pmatrix}0&0&1&0\\0&0&0&1\\1&0&0&0\\0&1&0&0\end{pmatrix}=n_T,\quad n_V=\begin{pmatrix}0&1&0&0\\1&0&0&0\\0&0&0&1\\0&0&1&0\end{pmatrix}=n_W. \]
Denoting a basis of $\mc B_A$ by $\{m_1,m_2,m_3,m_4\}$, we obtain a multiplication table
\begin{table}[H]
\begin{center}
\begin{tabular}{c|c|c|c|c}
    $b_j\otimes\backslash$&$m_1$&$m_2$&$m_3$&$m_4$\\\hline
    $1,U$&$m_1$&$m_2$&$m_3$&$m_4$\\
    $X,Z$&$m_4$&$m_3$&$m_2$&$m_1$\\
    $Y,T$&$m_3$&$m_4$&$m_1$&$m_2$\\
    $V,W$&$m_2$&$m_1$&$m_4$&$m_3$
\end{tabular}.
\end{center}
\end{table}
\hspace{-17pt}We can identify
\[ m_1\cong1\oplus U,\quad m_2\cong V\oplus W,\quad m_3\cong Y\oplus T,\quad m_4\cong X\oplus Z. \]
They have
\[ d_{\mc B_A}(m_1)=1,\quad d_{\mc B_A}(m_2)=d_{\mc B}(V),\quad d_{\mc B_A}(m_3)=d_{\mc B}(X)d_{\mc B}(V),\quad d_{\mc B_A}(m_4)=d_{\mc B}(X). \]
One can also obtain these identifications via the free module functor $F_A$:
\[ \hspace{-40pt}F_A(1)\cong1\oplus U\cong F_A(U),\quad F_A(X)\cong X\oplus Z\cong F_A(Z),\quad F_A(Y)\cong Y\oplus T\cong F_A(T),\quad F_A(V)\cong V\oplus W\cong F_A(W). \]
They obey the $\mbb Z/2\mbb Z\times\mbb Z/2\mbb Z$ monoidal products
\begin{table}[H]
\begin{center}
\begin{tabular}{c|c|c|c|c}
    $\otimes_A$&$m_1$&$m_2$&$m_3$&$m_4$\\\hline
    $m_1$&$m_1$&$m_2$&$m_3$&$m_4$\\\hline
    $m_2$&&$m_1$&$m_4$&$m_3$\\\hline
    $m_3$&&&$m_1$&$m_2$\\\hline
    $m_4$&&&&$m_1$
\end{tabular}.
\end{center}
\end{table}
\hspace{-17pt}This shows $\mc B_A\simeq\vecG_{\mbb Z/2\mbb Z\times\mbb Z/2\mbb Z}^\alpha$.

For the other nontrivial connected étale algebra, just the names of matrices change. For $A\cong1\oplus V$, we have
\[ \hspace{-20pt}n_1=1_4=n_V,\quad n_X=\begin{pmatrix}0&0&0&1\\0&0&1&0\\0&1&0&0\\1&0&0&0\end{pmatrix}=n_T,\quad n_Y=\begin{pmatrix}0&0&1&0\\0&0&0&1\\1&0&0&0\\0&1&0&0\end{pmatrix}=n_Z,\quad n_U=\begin{pmatrix}0&1&0&0\\1&0&0&0\\0&0&0&1\\0&0&1&0\end{pmatrix}=n_W. \]
It has identifications
\[ m_1\cong1\oplus V,\quad m_2\cong X\oplus T,\quad m_3\cong Y\oplus Z,\quad m_4\cong U\oplus W. \]
They obey the same $\mbb Z/2\mbb\times\mbb Z/2\mbb Z$ monoidal products
\begin{table}[H]
\begin{center}
\begin{tabular}{c|c|c|c|c}
    $\otimes_A$&$m_1$&$m_2$&$m_3$&$m_4$\\\hline
    $m_1$&$m_1$&$m_2$&$m_3$&$m_4$\\\hline
    $m_2$&&$m_1$&$m_4$&$m_3$\\\hline
    $m_3$&&&$m_1$&$m_2$\\\hline
    $m_4$&&&&$m_1$
\end{tabular}.
\end{center}
\end{table}
\hspace{-17pt}This shows $\mc B_A\simeq\vecG_{\mbb Z/2\mbb Z\times\mbb Z/2\mbb Z}^\alpha$.

We conclude
\begin{table}[H]
\begin{center}
\begin{tabular}{c|c|c|c}
    Connected étale algebra $A$&$\mc B_A$&$\rank(\mc B_A)$&Lagrangian?\\\hline
    $1$&$\mc B$&$8$&No\\
    $1\oplus U$ for $(d,h)$ in (\ref{Z2Z2Z2etale})&$\vecG_{\mbb Z/2\mbb Z\times\mbb Z/2\mbb Z}^\alpha$&4&No\\
    $1\oplus V$ for $(d,h)$ in (\ref{Z2Z2Z2etale})&$\vecG_{\mbb Z/2\mbb Z\times\mbb Z/2\mbb Z}^\alpha$&4&No
\end{tabular}.
\end{center}
\caption{Connected étale algebras in rank eight MFC $\mcal B\simeq\vecG_{\mbb Z/2\mbb Z\times\mbb Z/2\mbb Z\times\mbb Z/2\mbb Z}^\alpha\simeq\vecG_{\mbb Z/2\mbb Z}^{-1}\boxtimes\tc$}\label{rank8Z2Z2Z2results}
\end{table}
\hspace{-17pt}Namely, 16 MFCs $\mc B\simeq\vecG_{\mbb Z/2\mbb Z\times\mbb Z/2\mbb Z\times\mbb Z/2\mbb Z}^\alpha\simeq\vecG_{\mbb Z/2\mbb Z}^{-1}\boxtimes\tc$'s in (\ref{Z2Z2Z2etale}) fail to be completely anisotropic, while the other 24 with
\begin{align*}
    (d_X,d_U,d_V,d_W,h_X,h_U,h_V,h_W)&=(1,1,1,1,\frac14,\frac12,\frac12,\frac12),(1,1,1,1,\frac34,\frac12,\frac12,\frac12),\\
    &~~~~(1,1,-1,-1,\frac14,\frac12,0,0),(1,1,-1,-1,\frac14,\frac12,\frac12,\frac12),\\
    &~~~~(1,1,-1,-1,\frac34,\frac12,0,0),(1,1,-1,-1,\frac34,\frac12,\frac12,\frac12),\\
    &~~~~(-1,1,1,1,\frac14,\frac12,\frac12,\frac12),(-1,1,1,1,\frac34,\frac12,\frac12,\frac12),\\
    &~~~~(-1,1,-1,-1,\frac14,\frac12,0,0),(-1,1,-1,-1,\frac14,\frac12,\frac12,\frac12),\\
    &~~~~(-1,1,-1,-1,\frac34,\frac12,0,0),(-1,1,-1,-1,\frac34,\frac12,\frac12,\frac12)\quad(\mods1\text{ for }h)
\end{align*}
are completely anisotropic.

\subsubsection{$\mc B\simeq\vecG_{\mbb Z/2\mbb Z\times\mbb Z/4\mbb Z}^\alpha$}
The MFCs have eight simple objects $\{1,X,Y,Z,T,U,V,W\}$ obeying monoidal products
\begin{table}[H]
\begin{center}
\begin{tabular}{c|c|c|c|c|c|c|c|c}
    $\otimes$&$1$&$X$&$Y$&$Z$&$T$&$U$&$V$&$W$\\\hline
    $1$&$1$&$X$&$Y$&$Z$&$T$&$U$&$V$&$W$\\\hline
    $X$&&$1$&$Z$&$Y$&$W$&$V$&$U$&$T$\\\hline
    $Y$&&&$1$&$X$&$V$&$W$&$T$&$U$\\\hline
    $Z$&&&&$1$&$U$&$T$&$W$&$V$\\\hline
    $T$&&&&&$Z$&$1$&$X$&$Y$\\\hline
    $U$&&&&&&$Z$&$Y$&$X$\\\hline
    $V$&&&&&&&$Z$&$1$\\\hline
    $W$&&&&&&&&$Z$
\end{tabular}.
\end{center}
\end{table}
\hspace{-17pt}(One can identify $\vecG_{\mbb Z/2\mbb Z}^{-1}=\{1,X\},\vecG_{\mbb Z/4\mbb Z}^\alpha=\{1,Z,V,W\}$, and $Y\cong X\otimes Z,T\cong X\otimes W,U\cong X\otimes V$.) Thus, they have
\begin{align*}
    &\fp_{\mc B}(1)=\fp_{\mc B}(X)=\fp_{\mc B}(Y)=\fp_{\mc B}(Z)\\
    &=\fp_{\mc B}(T)=\fp_{\mc B}(U)=\fp_{\mc B}(V)=\fp_{\mc B}(W)=1,
\end{align*}
and
\[ \fp(\mc B)=8. \]
Their quantum dimensions $d_j$'s are solutions of the same multiplication rules $d_id_j=\sum_{k=1}^8{N_{ij}}^kd_k$. There are four solutions
\begin{align*}
    (d_X,d_Y,d_Z,d_T,d_U,d_V,d_W)=&(-1,-1,1,-1,-1,1,1),(-1,-1,1,1,1,-1,-1),\\
    &(1,1,1,-1,-1,-1,-1),(1,1,1,1,1,1,1).
\end{align*}
Only the last quantum dimensions give unitary MFCs. All the quantum dimensions have the same categorical dimension
\[ D^2(\mc B)=8. \]
All quantum dimensions have the same four conformal dimensions\footnote{Naively, one finds 16 consistent conformal dimensions, but the others are equivalent to one of four in the main text under permutations $(XY)(TU)$ or $(XY)(VW)$ or $(XY)(TVUW)$ of simple objects.}
\begin{align*}
    (h_X,h_Y,h_Z,h_T,h_U,h_V,h_W)=&(\frac14,\frac34,\frac12,\frac18,\frac18,\frac38,\frac38),(\frac14,\frac34,\frac12,\frac18,\frac18,\frac78,\frac78),\\
    &(\frac14,\frac34,\frac12,\frac38,\frac38,\frac58,\frac58),(\frac14,\frac34,\frac12,\frac58,\frac58,\frac78,\frac78).\quad(\mods1)
\end{align*}
The $S$-matrices are given by
\[ \widetilde S=\begin{pmatrix}1&d_X&d_Xd_Z&d_Z&d_Xd_W&d_Xd_V&d_V&d_W\\d_X&-1&-d_Z&d_Xd_Z&-d_W&-d_V&d_Xd_V&d_Xd_W\\d_Xd_Z&-d_Z&-1&d_X&-d_Zd_W&d_V&-d_Xd_V&-d_Xd_W\\d_Z&d_Xd_Z&d_X&1&-d_Xd_W&-d_Xd_V&-d_V&-d_W\\d_Xd_W&-d_W&-d_Zd_W&-d_Xd_W&\mp i&\pm i&\mp i\cdot d_X&\pm i\cdot d_X\\d_Xd_V&-d_V&d_V&-d_Xd_V&\pm i&\mp i&\pm i\cdot d_X&\mp i\cdot d_X\\d_V&d_Xd_V&-d_Xd_V&-d_V&\mp i\cdot d_X&\pm i\cdot d_X&\pm i&\mp i\\d_W&d_Xd_W&-d_Xd_W&-d_W&\pm i\cdot d_X&\mp i\cdot d_X&\mp i&\pm i\end{pmatrix}. \]
Regardless of quantum dimensions, they have additive central charges
\[ c(\mc B)=\begin{cases}2&(\text{1st }h),\\0&(\text{2nd }h),\\4&(\text{3rd }h),\\-2&(\text{4th }h).\end{cases}\quad(\mods8) \]
There are
\[ 4(\text{quantum dimensions})\times4(\text{conformal dimensions})\times2(\text{categorical dimensions})=32 \]
MFCs, among which those eight with the fourth quantum dimensions give unitary MFCs. We study connected étale algebras in all 32 MFCs simultaneously.

The most general form
\[ A\cong1\oplus n_XX\oplus n_YY\oplus n_ZZ\oplus n_TT\oplus n_UU\oplus n_VV\oplus n_WW \]
with $n_j\in\mbb N$ has
\[ \fp_{\mc B}(A)=1+n_X+n_Y+n_Z+n_T+n_U+n_V+n_W. \]
For this to obey (\ref{FPdimA2bound}), the natural numbers can take only eight values
\begin{align*}
    (n_X,n_Y,n_Z,n_T,n_U,n_V,n_W)=&(0,0,0,0,0,0,0),(1,0,0,0,0,0,0),(0,1,0,0,0,0,0),(0,0,1,0,0,0,0),\\
    &(0,0,0,1,0,0,0),(0,0,0,0,1,0,0),(0,0,0,0,0,1,0),(0,0,0,0,0,0,1).
\end{align*}
The first is nothing but the trivial connected étale algebra $A\cong1$ giving $\mc B_A^0\simeq\mc B_A\simeq\mc B$. The other six candidates except the fourth with $Z$ have nontrivial simple object, and they fail to be commutative. The fourth candidate $1\oplus Z$ does pass the necessary condition (\ref{commutativealgnecessary}), but since the simple object $Z$ has $(d_Z,h_Z)=(1,\frac12)$ (mod 1 for $h_Z$), it has $c_{Z,Z}\cong-id_1$ \cite{KK23preMFC}. Thus, it does not give commutative algebra either.

We conclude
\begin{table}[H]
\begin{center}
\begin{tabular}{c|c|c|c}
    Connected étale algebra $A$&$\mc B_A$&$\rank(\mc B_A)$&Lagrangian?\\\hline
    $1$&$\mc B$&$8$&No
\end{tabular}.
\end{center}
\caption{Connected étale algebras in rank eight MFC $\mcal B\simeq\vecG_{\mbb Z/2\mbb Z\times\mbb Z/4\mbb Z}^\alpha$}\label{rank8Z2Z4results}
\end{table}
\hspace{-17pt}All the 32 MFCs $\mc B\simeq\vecG_{\mbb Z/2\mbb Z\times\mbb Z/4\mbb Z}^\alpha$'s are completely anisotropic.

\subsubsection{$\mc B\simeq su(8)_1\simeq\vecG_{\mbb Z/8\mbb Z}^\alpha$}
The MFCs have eight simple objects $\{1,X,Y,Z,T,U,V,W\}$ obeying monoidal products
\begin{table}[H]
\begin{center}
\begin{tabular}{c|c|c|c|c|c|c|c|c}
    $\otimes$&$1$&$X$&$Y$&$Z$&$T$&$U$&$V$&$W$\\\hline
    $1$&$1$&$X$&$Y$&$Z$&$T$&$U$&$V$&$W$\\\hline
    $X$&&$1$&$U$&$T$&$Z$&$Y$&$W$&$V$\\\hline
    $Y$&&&$W$&$1$&$X$&$V$&$Z$&$T$\\\hline
    $Z$&&&&$V$&$W$&$X$&$U$&$Y$\\\hline
    $T$&&&&&$V$&$1$&$Y$&$U$\\\hline
    $U$&&&&&&$W$&$T$&$Z$\\\hline
    $V$&&&&&&&$X$&$1$\\\hline
    $W$&&&&&&&&$X$
\end{tabular}.
\end{center}
\end{table}
\hspace{-17pt}Thus, they have
\begin{align*}
    &\fp_{\mc B}(1)=\fp_{\mc B}(X)=\fp_{\mc B}(Y)=\fp_{\mc B}(Z)\\
    &=\fp_{\mc B}(T)=\fp_{\mc B}(U)=\fp_{\mc B}(V)=\fp_{\mc B}(W)=1,
\end{align*}
and
\[ \fp(\mc B)=8. \]
Their quantum dimensions $d_j$'s are solutions of the same multiplication rules $d_id_j=\sum_{k=1}^8{N_{ij}}^kd_k$. There are two solutions
\[ (d_X,d_Y,d_Z,d_T,d_U,d_V,d_W)=(1,-1,-1,-1,-1,1,1),(1,1,1,1,1,1,1). \]
Only the second gives unitary MFCs. The two quantum dimensions both have the same four conformal dimensions\footnote{Naively, one finds eight consistent conformal dimensions, but the others are equivalent to one in the main text under permutations $(YU)(ZT)$ of simple objects.}
\begin{align*}
    (h_X,h_Y,h_Z,h_T,h_U,h_V,h_W)=&(0,\frac1{16},\frac1{16},\frac9{16},\frac9{16},\frac14,\frac14),(0,\frac3{16},\frac3{16},\frac{11}{16},\frac{11}{16},\frac34,\frac34),\\
    &(0,\frac5{16},\frac5{16},\frac{13}{16},\frac{13}{16},\frac14,\frac14),(0,\frac7{16},\frac7{16},\frac{15}{16},\frac{15}{16},\frac34,\frac34).\quad(\mods1)
\end{align*}
The $S$-matrices are given by
\[ \widetilde S=\begin{pmatrix}1&d_X&d_Y&d_Z&d_T&d_U&d_V&d_W\\d_X&1&-d_Y&-d_Z&-d_T&-d_U&d_V&d_W\\d_Y&-d_Y&\pm e^{\pm\pi i/4}&\pm e^{\mp\pi i/4}&\mp e^{\mp\pi i/4}&\mp e^{\pm\pi i/4}&\mp i&\pm i\\d_Z&-d_Z&\pm e^{\mp\pi i/4}&\pm e^{\pm\pi i/4}&\mp e^{\pm\pi i/4}&\mp e^{\mp\pi i/4}&\pm i&\mp i\\d_T&-d_T&\mp e^{\mp\pi i/4}&\mp e^{\pm\pi i/4}&\pm e^{\pm\pi i/4}&\pm e^{\mp\pi i/4}&\pm i&\mp i\\d_U&-d_U&\mp e^{\pm\pi i/4}&\mp e^{\mp\pi i/4}&\pm e^{\mp\pi i/4}&\pm e^{\pm\pi i/4}&\mp i&\pm i\\d_V&d_V&\mp i&\pm i&\pm i&\mp i&-1&-1\\d_W&d_W&\pm i&\mp i&\mp i&\pm i&-1&-1\end{pmatrix}. \]
They have additive central charges
\[ c(\mc B)=
\begin{cases}1&(\text{1st\&3rd }h),\\-1&(\text{2nd\&4th }h).\end{cases}\quad(\mods8) \]
There are
\[ 2(\text{quantum dimensions})\times4(\text{conformal dimensions})\times2(\text{categorical dimensions})=16 \]
MFCs, among which those eight with the second quantum dimensions give unitary MFCs. We classify connected étale algebras in all 16 MFCs simultaneously.

We work with an ansatz
\[ A\cong1\oplus n_XX\oplus n_YY\oplus n_ZZ\oplus n_TT\oplus n_UU\oplus n_VV\oplus n_WW \]
with $n_j\in\mbb N$. It has
\[ \fp_{\mc B}(A)=1+n_X+n_Y+n_Z+n_T+n_U+n_V+n_W. \]
For this to obey (\ref{FPdimA2bound}), the natural numbers can take only eight values
\begin{align*}
    (n_X,n_Y,n_Z,n_T,n_U,n_V,n_W)=&(0,0,0,0,0,0,0),(1,0,0,0,0,0,0),(0,1,0,0,0,0,0),(0,0,1,0,0,0,0),\\
    &(0,0,0,1,0,0,0),(0,0,0,0,1,0,0),(0,0,0,0,0,1,0),(0,0,0,0,0,0,1).
\end{align*}
The first is nothing but the trivial connected étale algebra $A\cong1$ giving $\mc B_A^0\simeq\mc B_A\simeq\mc B$. Those six with $Y,Z,T,U,V,W$ fail to satisfy the necessary condition (\ref{commutativealgnecessary}), and they are ruled out. We are left with the second $n_X=1$ or $A\cong1\oplus X$. It is a $\mbb Z/2\mbb Z$ algebra by the lemma 1. It turns out that this is connected étale; it is commutative because $X$ has $(d_X,h_X)=(1,0)$ (mod 1 for $h_X$), and hence $c_{X,X}\cong id_1$ \cite{KK23preMFC}. To check the separability, we identify $\mc B_A$.

Since it has $\fp_{\mc B}(A)=2$, it demands
\[ \fp(\mc B_A^0)=2,\quad\fp(\mc B_A)=4. \]
The category $\mc B_A^0$ of dyslectic modules is identified as $\vecG_{\mbb Z/2\mbb Z}^{-1}$ because it is the only MFC with $\fp=2$. This identification also matches central charges.\footnote{The $\mbb Z/2\mbb Z$ MFC has
\[ c(\vecG_{\mbb Z/2\mbb Z}^{-1})=\begin{cases}1&(h_{\mbb Z/2\mbb Z}=\frac14),\\-1&(h_{\mbb Z/2\mbb Z}=\frac34).\end{cases} \]
We will see $\mc B_A^0=\{1\oplus X,V\oplus W\}$. The nontrivial deconfined particle has $h_V=\frac14,\frac34=h_W$, and this determines the central charge $c(\mc B_A^0)$.} The category $\mc B_A$ of right $A$-modules should contain the MFC and have $\fp=4$. It turns out that
\[ \mc B_A\simeq\vecG_{\mbb Z/4\mbb Z}^\alpha. \]
One of the easiest ways to see this fact is to perform anyon condensation. Under the procedure, we `identify' $1$ and $X$, and hence the pairs $(Y,U)$, $(Z,T)$, and $(V,W)$. Since they obey $\mbb Z/4\mbb Z$ monoidal products, we arrive the identification.\footnote{More rigorously, we should find NIM-reps. Indeed, we found a four-dimensional NIM-rep
\[ n_1=1_4=n_X,\quad n_Y=\begin{pmatrix}0&0&0&1\\0&0&1&0\\1&0&0&0\\0&1&0&0\end{pmatrix}=n_U,\quad n_Z=\begin{pmatrix}0&0&1&0\\0&0&0&1\\0&1&0&0\\1&0&0&0\end{pmatrix}=n_T,\quad n_V=\begin{pmatrix}0&1&0&0\\1&0&0&0\\0&0&0&1\\0&0&1&0\end{pmatrix}=n_W. \]
Denoting a basis of $\mc B_A$ by $\{m_1,m_2,m_3,m_4\}$, we obtain a multiplication table
\begin{table}[H]
\begin{center}
\begin{tabular}{c|c|c|c|c}
    $b_j\otimes\backslash$&$m_1$&$m_2$&$m_3$&$m_4$\\\hline
    $1,X$&$m_1$&$m_2$&$m_3$&$m_4$\\
    $Y,U$&$m_4$&$m_3$&$m_1$&$m_2$\\
    $Z,T$&$m_3$&$m_4$&$m_2$&$m_1$\\
    $V,W$&$m_2$&$m_1$&$m_4$&$m_3$
\end{tabular}.
\end{center}
\end{table}
\hspace{-10pt}In this basis, we can identify
\[ m_1\cong1\oplus X,\quad m_2\cong V\oplus W,\quad m_3\cong Z\oplus T,\quad m_4\cong Y\oplus U.\]
In the category $\mc B_A$ of right $A$-modules, they have quantum dimensions (\ref{dBAm}):
\[ d_{\mc B_A}(m_1)=1=d_{\mc B_A}(m_2),\quad d_{\mc B_A}(m_3)=d_{\mc B}(Y)=d_{\mc B_A}(m_4). \]
Furthermore, they obey monoidal products
\begin{table}[H]
\begin{center}
\begin{tabular}{c|c|c|c|c}
    $\otimes_A$&$m_1$&$m_2$&$m_3$&$m_4$\\\hline
    $m_1$&$m_1$&$m_2$&$m_3$&$m_4$\\\hline
    $m_2$&&$m_1$&$m_4$&$m_3$\\\hline
    $m_3$&&&$m_2$&$m_1$\\\hline
    $m_4$&&&&$m_2$
\end{tabular}.
\end{center}
\end{table}
\hspace{-10pt}
This shows $\mc B_A\simeq\vecG_{\mbb Z/4\mbb Z}^\alpha$.}

We conclude
\begin{table}[H]
\begin{center}
\begin{tabular}{c|c|c|c}
    Connected étale algebra $A$&$\mc B_A$&$\rank(\mc B_A)$&Lagrangian?\\\hline
    $1$&$\mc B$&$8$&No\\
    $1\oplus X$&$\vecG_{\mbb Z/4\mbb Z}^\alpha$&4&No
\end{tabular}.
\end{center}
\caption{Connected étale algebras in rank eight MFC $\mcal B\simeq\vecG_{\mbb Z/8\mbb Z}^\alpha$}\label{rank8Z8results}
\end{table}
\hspace{-17pt}That is, all the 16 MFCs $\mc B\simeq\vecG_{\mbb Z/8\mbb Z}^\alpha$'s fail to be completely anisotropic.

\subsubsection{$\mc B\simeq\fib\boxtimes\vecG_{\mbb Z/2\mbb Z\times\mbb Z/2\mbb Z}^\alpha$}
The MFCs have eight simple objects $\{1,X,Y,Z,T,U,V,W\}$ obeying monoidal products
\begin{table}[H]
\begin{center}
\begin{tabular}{c|c|c|c|c|c|c|c|c}
    $\otimes$&$1$&$X$&$Y$&$Z$&$T$&$U$&$V$&$W$\\\hline
    $1$&$1$&$X$&$Y$&$Z$&$T$&$U$&$V$&$W$\\\hline
    $X$&&$1$&$Z$&$Y$&$V$&$W$&$T$&$U$\\\hline
    $Y$&&&$1$&$X$&$W$&$V$&$U$&$T$\\\hline
    $Z$&&&&$1$&$U$&$T$&$W$&$V$\\\hline
    $T$&&&&&$1\oplus W$&$Z\oplus V$&$X\oplus U$&$Y\oplus T$\\\hline
    $U$&&&&&&$1\oplus W$&$Y\oplus T$&$X\oplus U$\\\hline
    $V$&&&&&&&$1\oplus W$&$Z\oplus V$\\\hline
    $W$&&&&&&&&$1\oplus W$
\end{tabular}.
\end{center}
\end{table}
\hspace{-17pt}(One can identify $\fib=\{1,W\},\vecG_{\mbb Z/2\mbb Z\times\mbb Z/2\mbb Z}^\alpha=\{1,X,Y,Z\}$, and $T\cong Y\otimes W,U\cong X\otimes W,V\cong Z\otimes W$.) Thus, they have
\begin{align*}
    &\fp_{\mc B}(1)=\fp_{\mc B}(X)=\fp_{\mc B}(Y)=\fp_{\mc B}(Z)=1,\\
    &\fp_{\mc B}(T)=\fp_{\mc B}(U)=\fp_{\mc B}(V)=\fp_{\mc B}(W)=\zeta:=\frac{1+\sqrt5}2,
\end{align*}
and
\[ \fp(\mc B)=10+2\sqrt5\approx14.5. \]
Their quantum dimensions $d_j$'s are solutions of the same multiplication rules $d_id_j=\sum_{k=1}^8{N_{ij}}^kd_k$. There are eight solutions
\begin{align*}
    (d_X,d_Y,d_Z,d_T,d_U,d_V,d_W)=&(-1,-1,1,-\zeta,-\zeta,\zeta,\zeta),(-1,-1,1,\zeta^{-1},\zeta^{-1},-\zeta^{-1},-\zeta^{-1}),\\
    &(-1,1,-1,-\zeta^{-1},\zeta^{-1},\zeta^{-1},-\zeta^{-1}),(-1,1,-1,\zeta,-\zeta,-\zeta,\zeta),\\
    &(1,-1,-1,-\zeta,\zeta,-\zeta,\zeta),(1,-1,-1,\zeta^{-1},-\zeta^{-1},\zeta^{-1},-\zeta^{-1}),\\
    &(1,1,1,-\zeta^{-1},-\zeta^{-1},-\zeta^{-1},-\zeta^{-1}),(1,1,1,\zeta,\zeta,\zeta,\zeta).
\end{align*}
Only the last quantum dimensions give unitary MFCs. They have categorical dimensions
\[ D^2(\mc B)=10-2\sqrt5(\approx5.5),\quad10+2\sqrt5. \]
In order to list up conformal dimensions without double-counting, we perform case analysis. (This fixes redundancies in names of simple objects, and hence superficially makes descriptions asymmetric.)

\paragraph{$\mc B\simeq\fib\boxtimes\vecG_{\mbb Z/2\mbb Z}^{-1}\boxtimes\vecG_{\mbb Z/2\mbb Z}^{-1}$.} In this case, the $\vecG_{\mbb Z/2\mbb Z}^{-1}\boxtimes\vecG_{\mbb Z/2\mbb Z}^{-1}$ factor has three classes depending on quantum dimensions $(d_X,d_Y)$ \cite{KK23MFC}. Together with $d_W=\zeta,-\zeta^{-1}$, we have six classes.
\begin{itemize}
    \item $(d_X,d_Y,d_W)=(1,1,\zeta)$. This gives unitary MFCs. Different MFCs are given by conformal dimensions
    \[ (h_X,h_Y,h_W)=(\frac14,\frac14,\frac25),(\frac14,\frac14,\frac35),(\frac14,\frac34,\frac25),(\frac14,\frac34,\frac35),(\frac34,\frac34,\frac25),(\frac34,\frac34,\frac35)\quad(\mods1). \]
    (The other conformal dimensions of $b_k\cong b_i\otimes b_j$ is given by $h_k=h_i+h_j$ mod 1.) Including two signs of categorical dimensions, we have 12 unitary MFCs.
    \item $(d_X,d_Y,d_W)=(1,1,-\zeta^{-1})$. Different MFCs are given by
    \[ (h_X,h_Y,h_W)=(\frac14,\frac14,\frac15),(\frac14,\frac14,\frac45),(\frac14,\frac34,\frac15),(\frac14,\frac34,\frac45),(\frac34,\frac34,\frac15),(\frac34,\frac34,\frac45)\quad(\mods1). \]
    With two signs of categorical dimensions, there are 12 MFCs.
    \item $(d_X,d_Y,d_W)=(1,-1,\zeta)$. Different MFCs are given by
    \[ \hspace{-60pt}(h_X,h_Y,h_W)=(\frac14,\frac14,\frac25),(\frac14,\frac14,\frac35),(\frac14,\frac34,\frac25),(\frac14,\frac34,\frac35),(\frac34,\frac14,\frac25),(\frac34,\frac14,\frac35),(\frac34,\frac34,\frac25),(\frac34,\frac34,\frac35)\quad(\mods1). \]
    Taking two signs of categorical dimensions into account, we get 16 MFCs.
    \item $(d_X,d_Y,d_W)=(1,-1,-\zeta^{-1})$. Different MFCs are given by
    \[ \hspace{-60pt}(h_X,h_Y,h_W)=(\frac14,\frac14,\frac15),(\frac14,\frac14,\frac45),(\frac14,\frac34,\frac15),(\frac14,\frac34,\frac45),(\frac34,\frac14,\frac15),(\frac34,\frac14,\frac45),(\frac34,\frac34,\frac15),(\frac34,\frac34,\frac45)\quad(\mods1). \]
    With two categorical dimensions, we get 16 MFCs.
    \item $(d_X,d_Y,d_W)=(-1,-1,\zeta)$. Different MFCs are given by
    \[ (h_X,h_Y,h_W)=(\frac14,\frac14,\frac25),(\frac14,\frac14,\frac35),(\frac14,\frac34,\frac25),(\frac14,\frac34,\frac35),(\frac34,\frac34,\frac25),(\frac34,\frac34,\frac35)\quad(\mods1). \]
    There are 12 MFCs.
    \item $(d_X,d_Y,d_W)=(-1,-1,-\zeta^{-1})$. Different MFCs are given by
    \[ (h_X,h_Y,h_W)=(\frac14,\frac14,\frac15),(\frac14,\frac14,\frac45),(\frac14,\frac34,\frac15),(\frac14,\frac34,\frac45),(\frac34,\frac34,\frac15),(\frac34,\frac34,\frac45)\quad(\mods1). \]
    There are 12 MFCs.
\end{itemize}
In total, there are
\[ 12+12+16+16+12+12=80 \]
MFCs $\mc B\simeq\fib\boxtimes\vecG_{\mbb Z/2\mbb Z}^{-1}\boxtimes\vecG_{\mbb Z/2\mbb Z}^{-1}$, among which those 12 with quantum dimensions $(d_X,d_Y,d_W)=(1,1,\zeta)$ are unitary. The $S$-matrices are given by
\[ \widetilde S=\begin{pmatrix}1&d_X&d_Y&d_Xd_Y&d_Yd_W&d_Xd_W&d_Xd_Yd_W&d_W\\d_X&-1&1&-d_X&d_W&-d_W&-d_Xd_W&d_Xd_W\\d_Y&1&-1&-d_Y&-d_W&d_W&-d_Yd_W&d_Yd_W\\d_Xd_Y&-d_X&-d_Y&1&-d_Yd_W&-d_Xd_W&d_W&d_Xd_Yd_W\\d_Yd_W&d_W&-d_W&-d_Yd_W&1&-1&d_Y&-d_Y\\d_Xd_W&-d_W&d_W&-d_Xd_W&-1&1&d_X&-d_X\\d_Xd_Yd_W&-d_Xd_W&-d_Yd_W&d_W&d_Y&d_X&-1&-d_Xd_Y\\d_W&d_Xd_W&d_Yd_W&d_Xd_Yd_W&-d_Y&-d_X&-d_Xd_Y&-1\end{pmatrix}.\]
They have additive central charges
\[ c(\mc B)=c(\fib)+c(\vecG_{\mbb Z/2\mbb Z}^{-1})+c(\vecG_{\mbb Z/2\mbb Z}^{-1})\quad(\mods8) \]
where
\[ c(\fib)=\begin{cases}\frac25&(h_W=\frac15),\\-\frac25&(h_W=\frac45),\\\frac{14}5&(h_W=\frac25),\\-\frac{14}5&(h_W=\frac35).\end{cases}\quad c(\vecG_{\mbb Z/2\mbb Z}^{-1})=\begin{cases}1&(h_{\mbb Z/2\mbb Z}=\frac14),\\-1&(h_{\mbb Z/2\mbb Z}=\frac34).\end{cases}\quad(\mods8) \]

\paragraph{$\mc B\simeq\fib\boxtimes\tc$.} In this case, simple objects $X,Y,Z$ can have conformal dimensions
\[ (h_X,h_Y,h_Z)=(0,0,\frac12),(0,\frac12,0),(\frac12,0,0),(\frac12,\frac12,\frac12)\quad(\mods1). \]
Thus, we have four classes.
\begin{itemize}
    \item $(d_X,d_Y,d_Z,d_W)=(1,1,1,\zeta)$. This gives unitary MFCs. Different MFCs are given by
    \[ (h_X,h_Y,h_Z,h_W)=(0,0,\frac12,\frac25),(0,0,\frac12,\frac35),(\frac12,\frac12,\frac12,\frac25),(\frac12,\frac12,\frac12,\frac35)\quad(\mods1). \]
    With the two signs of categorical dimensions, there are eight unitary MFCs.
    \item $(d_X,d_Y,d_Z,d_W)=(1,1,1,-\zeta^{-1})$. Different MFCs are given by
    \[ (h_X,h_Y,h_Z,h_W)=(0,0,\frac12,\frac15),(0,0,\frac12,\frac45),(\frac12,\frac12,\frac12,\frac15),(\frac12,\frac12,\frac12,\frac45)\quad(\mods1). \]
    Including the two signs of categorical dimensions, we have eight MFCs.
    \item $(d_X,d_Y,d_Z,d_W)=(1,-1,-1,\zeta)$. Different MFCs are given by
    \[ \hspace{-50pt}(h_X,h_Y,h_Z,h_W)=(0,0,\frac12,\frac25),(0,0,\frac12,\frac35),(\frac12,0,0,\frac25),(\frac12,0,0,\frac35),(\frac12,\frac12,\frac12,\frac25),(\frac12,\frac12,\frac12,\frac35)\quad(\mods1). \]
    There are 12 MFCs.
    \item $(d_X,d_Y,d_Z,d_W)=(1,-1,-1,-\zeta^{-1})$. Different MFCs are given by
    \[ \hspace{-50pt}(h_X,h_Y,h_Z,h_W)=(0,0,\frac12,\frac15),(0,0,\frac12,\frac45),(\frac12,0,0,\frac15),(\frac12,0,0,\frac45),(\frac12,\frac12,\frac12,\frac15),(\frac12,\frac12,\frac12,\frac45)\quad(\mods1). \]
    There are 12 MFCs.
\end{itemize}
In total, there are
\[ 8+8+12+12=40 \]
MFCs $\mc B\simeq\fib\boxtimes\tc$, among which those eight with the quantum dimensions $(d_X,d_Y,d_Z,d_W)=(1,1,1,\zeta)$ are unitary. The $S$-matrices are given by
\[ \widetilde S=\begin{pmatrix}1&d_X&d_Y&d_Z&d_Yd_W&d_Xd_W&d_Zd_W&d_W\\d_X&1&-d_Y&-d_Z&-d_Yd_W&d_W&-d_Zd_W&d_Xd_W\\d_Y&-d_Y&1&-1&d_W&-d_Yd_W&-d_W&d_Yd_W\\d_Z&-d_Z&-1&1&-d_W&-d_Zd_W&d_W&d_Zd_W\\d_Yd_W&-d_Yd_W&d_W&-d_W&-1&d_Y&1&-d_Y\\d_Xd_W&d_W&-d_Yd_W&-d_Zd_W&d_Y&-1&d_Z&-d_X\\d_Zd_W&-d_Zd_W&-d_W&d_W&1&d_Z&-1&-d_Z\\d_W&d_Xd_W&d_Yd_W&d_Zd_W&-d_Y&-d_X&-d_Z&-1\end{pmatrix}.\]
They have additive central charges
\[ c(\mc B)=c(\fib)+c(\tc)\quad(\mods8) \]
where
\[ c(\fib)=\begin{cases}\frac25&(h_W=\frac15),\\-\frac25&(h_W=\frac45),\\\frac{14}5&(h_W=\frac25),\\-\frac{14}5&(h_W=\frac35).\end{cases}\quad c(\tc)=\begin{cases}0&(\text{one }h=\frac12\&\text{the other }h=0),\\4&(\text{all }h=\frac12).\end{cases}\quad(\mods8). \]

Including both $\mc B\simeq\fib\boxtimes\vecG_{\mbb Z/2\mbb Z}^{-1}\boxtimes\vecG_{\mbb Z/2\mbb Z}^{-1}$ and $\mc B\simeq\fib\boxtimes\tc$, we have
\[ 80+40=120 \]
MFCs with fusion ring $K(\fib\boxtimes\vecG_{\mbb Z/2\mbb Z\times\mbb Z/2\mbb Z}^\alpha)$. Among them, 12 from the first and eight from the second MFCs are unitary, 20 in total.

Having listed MFCs, we can classify connected étale algebras in them. We work with an ansatz
\[ A\cong1\oplus n_XX\oplus n_YY\oplus n_ZZ\oplus n_TT\oplus n_UU\oplus n_VV\oplus n_WW \]
with $n_j\in\mbb N$. It has
\[ \fp_{\mc B}(A)=1+n_X+n_Y+n_Z+\zeta(n_T+n_U+n_V+n_W). \]
For this to obey (\ref{FPdimA2bound}), the natural numbers can take only 26 values
\begin{align*}
    (n_X,n_Y,n_Z,n_T,n_U,n_V,n_W)=&(0,0,0,0,0,0,0),(1,0,0,0,0,0,0),(2,0,0,0,0,0,0),\\
    &(1,1,0,0,0,0,0),(1,0,1,0,0,0,0),(1,0,0,1,0,0,0),\\
    &(1,0,0,0,1,0,0),(1,0,0,0,0,1,0),(1,0,0,0,0,0,1),\\
    &(0,1,0,0,0,0,0),(0,2,0,0,0,0,0),(0,1,1,0,0,0,0),\\
    &(0,1,0,1,0,0,0),(0,1,0,0,1,0,0),(0,1,0,0,0,1,0),\\
    &(0,1,0,0,0,0,1),(0,0,1,0,0,0,0),(0,0,2,0,0,0,0),\\
    &(0,0,1,1,0,0,0),(0,0,1,0,1,0,0),(0,0,1,0,0,1,0),\\
    &(0,0,1,0,0,0,1),(0,0,0,1,0,0,0),(0,0,0,0,1,0,0),\\
    &(0,0,0,0,0,1,0),(0,0,0,0,0,0,1).
\end{align*}
The first is nothing but the trivial connected étale algebra $A\cong1$ giving $\mc B\simeq\mc B_A^0\simeq\mc B_A$. Whether the others can be commutative depends on quantum and conformal dimensions.

For $\mc B\simeq\fib\boxtimes\vecG_{\mbb Z/2\mbb Z}^{-1}\boxtimes\vecG_{\mbb Z/2\mbb Z}^{-1}$, $X,Y,T,U,V,W$ have nontrivial conformal dimensions, and candidates with them fail to be commutative. Hence, we are left with candidates with just $Z$'s. It has $(d_Z,h_Z)=(d_Xd_Y,h_X+h_Y)$ (mod 1 for $h_Z$). Thus, those with $(d_Z,h_Z)=(1,0),(-1,\frac12)$ (mod 1 for $h_Z$) are commutative. More concretely, independent of quantum and conformal dimensions of $W$, we have
\begin{equation}
    c_{Z,Z}\cong id_1\iff(d_X,d_Y,h_X,h_Y)=(1,1,\frac14,\frac34),(-1,-1,\frac14,\frac34),(1,-1,\frac14,\frac14),(1,-1,\frac34,\frac34)\quad(\mods1\text{ for }h).\label{cZZtrivial}
\end{equation}
Therefore, we are left with two additional candidates for connected étale algebras
\[ A\cong1\oplus Z,1\oplus2Z. \]
The second candidate $A\cong1\oplus2Z$ has $\fp_{\mc B}(A)=3$, and in particular demands
\[ \fp(\mc B_A^0)=\frac{5+\sqrt5}9\approx1.61. \]
However, there is no MFC with this Frobenius-Perron dimension. Thus, the candidate is ruled out. We are left with the first candidate. It is a $\mbb Z/2\mbb Z$ algebra by the lemma 1. It also turns out separable. To show this, we identify $\mc B_A$.

It has $\fp_{\mc B}(A)=2$, and demands
\[ \fp(\mc B_A^0)=\frac{5+\sqrt5}2,\quad\fp(\mc B_A)=5+\sqrt5. \]
Since the only MFC with $\fp=\frac{5+\sqrt5}2$ is a Fibonacci MFC, we arrive
\[ \mc B_A^0\simeq\fib. \]
The matching of central charges makes the identification more precise. The Fibonacci category has
\[ c(\fib)=\begin{cases}\pm\frac25&(d_\fib=-\zeta^{-1}),\\\pm\frac{14}5&(d_\fib=\zeta).\end{cases}\quad(\mods4) \]
The $\vecG_{\mbb Z/2\mbb Z\times\mbb Z/2\mbb Z}^\alpha$ has
\[ c(\vecG_{\mbb Z/2\mbb Z\times\mbb Z/2\mbb Z}^\alpha)=\begin{cases}2&(\vecG_{\mbb Z/2\mbb Z\times\mbb Z/2\mbb Z}^\alpha\simeq\vecG_{\mbb Z/2\mbb Z}^{-1}\boxtimes\vecG_{\mbb Z/2\mbb Z}^{-1}\text{ with }h_X=h_Y),\\0&(\vecG_{\mbb Z/2\mbb Z\times\mbb Z/2\mbb Z}^\alpha\simeq\vecG_{\mbb Z/2\mbb Z}^{-1}\boxtimes\vecG_{\mbb Z/2\mbb Z}^{-1}\text{ with }h_X\neq h_Y).\end{cases}\quad(\mods4) \]
We see the second $\vecG_{\mbb Z/2\mbb Z}^{-1}\boxtimes\vecG_{\mbb Z/2\mbb Z}^{-1}$ MFCs with different conformal dimensions (mod 1) trivially match additive central charges, while the first $\vecG_{\mbb Z/2\mbb Z}^{-1}\boxtimes\vecG_{\mbb Z/2\mbb Z}^{-1}$ with the same conformal dimensions cannot. Therefore, although $(d_Z,h_Z)=(-1,\frac12)$ in (\ref{cZZtrivial}) has trivial self braiding $c_{Z,Z}\cong id_1$, they do not give connected étale algebras. We conclude
\begin{equation}
    A\cong1\oplus Z\text{ is connected étale}\iff(d_X,d_Y,h_X,h_Y)=(1,1,\frac14,\frac34),(-1,-1,\frac14,\frac34)\quad(\mods1\text{ for }h).\label{fibZ2Z21Zetale}
\end{equation}
The category $\mc B_A$ of right $A$-modules should contain this MFC as a subcategory. Together with the Frobenius-Perron dimension, it is identified as\footnote{A four-dimensional NIM-rep is given by
\[ n_1=1_4=n_Z,\quad n_X=\begin{pmatrix}0&1&0&0\\1&0&0&0\\0&0&0&1\\0&0&1&0\end{pmatrix}=n_Y,\quad n_T=\begin{pmatrix}0&0&1&0\\0&0&0&1\\1&0&0&1\\0&1&1&0\end{pmatrix}=n_U,\quad n_V=\begin{pmatrix}0&0&0&1\\0&0&1&0\\0&1&1&0\\1&0&0&1\end{pmatrix}=n_W. \]
The solution gives identifications
\[ m_1\cong1\oplus Z,\quad m_2\cong X\oplus Y,\quad m_3\cong T\oplus U,\quad m_4\cong V\oplus W. \]
The simple objects $m_{1,4}$ are deconfined, while $m_{2,3}$ are confined, i.e., $\mc B_A^0=\{m_1,m_4\}$. In $\mc B_A$, they have
\[ d_{\mc B_A}(m_1)=1,\quad d_{\mc B_A}(m_2)=\pm1,\quad d_{\mc B_A}(m_3)=\pm\zeta^{\pm1},\quad d_{\mc B_A}(m_4)=\pm\zeta^{\pm1}. \]
Signs in the last quantum dimension are correlated, but those in the third are not. Employing the free module functor $F_A$, one obtains monoidal products $\otimes_A$:
\begin{table}[H]
\begin{center}
\begin{tabular}{c|c|c|c|c}
    $\otimes_A$&$m_1$&$m_2$&$m_3$&$m_4$\\\hline
    $m_1$&$m_1$&$m_2$&$m_3$&$m_4$\\\hline
    $m_2$&&$m_1$&$m_4$&$m_3$\\\hline
    $m_3$&&&$m_1\oplus m_4$&$m_2\oplus m_3$\\\hline
    $m_4$&&&&$m_1\oplus m_4$
\end{tabular}.
\end{center}
\end{table}
\hspace{-11pt}This shows the identification.}
\[ \mc B_A\simeq su(2)_3. \]

For $\mc B\simeq\fib\boxtimes\tc$, $T,U,V,W$ have nontrivial conformal dimensions, and candidates with them fail to be commutative. On the other hand, those with just $X,Y,Z$ can have $c_{b_j,b_j}\cong id_1$ depending on their quantum and conformal dimensions. From the lists above, independent of $h_W$, we find
\begin{equation}
\begin{split}
    \hspace{-40pt}c_{X,X}\cong id_1\iff(d_X,d_Y,d_Z,d_W,h_X,h_Y,h_Z,h_W)&=(1,1,1,\zeta,0,0,\frac12,\frac25),(1,1,1,\zeta,0,0,\frac12,\frac35),\\
    &~~~~(1,1,1,-\zeta^{-1},0,0,\frac12,\frac15),(1,1,1,-\zeta^{-1},0,0,\frac12,\frac45),\\
    &~~~~(1,-1,-1,\zeta,0,0,\frac12,\frac25),(1,-1,-1,\zeta,0,0,\frac12,\frac35),\\
    &~~~~(1,-1,-1,-\zeta^{-1},0,0,\frac12,\frac15),(1,-1,-1,-\zeta^{-1},0,0,\frac12,\frac45),\\
    \hspace{-40pt}c_{Y,Y}\cong id_1\iff(d_X,d_Y,d_Z,d_W,h_X,h_Y,h_Z,h_W)&=(1,1,1,\zeta,0,0,\frac12,\frac25),(1,1,1,\zeta,0,0,\frac12,\frac35),\\
    &~~~~(1,1,1,-\zeta^{-1},0,0,\frac12,\frac15),(1,1,1,-\zeta^{-1},0,0,\frac12,\frac45),\\
    &~~~~(1,-1,-1,\zeta,\frac12,\frac12,\frac12,\frac25),(1,-1,-1,\zeta,\frac12,\frac12,\frac12,\frac35),\\
    &~~~~(1,-1,-1,-\zeta^{-1},\frac12,\frac12,\frac12,\frac15),(1,-1,-1,-\zeta^{-1},\frac12,\frac12,\frac12,\frac45),\\
    \hspace{-40pt}c_{Z,Z}\cong id_1\iff(d_X,d_Y,d_Z,d_W,h_X,h_Y,h_Z,h_W)&=(1,-1,-1,\zeta,0,0,\frac12,\frac25),(1,-1,-1,\zeta,0,0,\frac12,\frac35),\\
    &~~~~(1,-1,-1,\zeta,\frac12,\frac12,\frac12,\frac25),(1,-1,-1,\zeta,\frac12,\frac12,\frac12,\frac35),\\
    &~~~~(1,-1,-1,-\zeta^{-1},0,0,\frac12,\frac15),(1,-1,-1,-\zeta^{-1},0,0,\frac12,\frac45),\\
    &~~~~(1,-1,-1,-\zeta^{-1},\frac12,\frac12,\frac12,\frac15),(1,-1,-1,-\zeta^{-1},\frac12,\frac12,\frac12,\frac45).\\
    &\hspace{200pt}(\text{mod }1\text{ for }h)
\end{split}\label{cXYZtrivial}
\end{equation}
How about the matching of central charges? We will see the Toric Code MFC has conformal dimensions $(h_X,h_Y,h_Z)=(0,0,\frac12)$ when it admits nontrivial connected étale algebras. Such MFCs have $c(\tc)=0$ mod 8, and the central charges are trivially matched because both $c(\mc B)$ and $c(\fib)$ are determined by the same conformal dimension $h_W$.

Just as for $\mc B\simeq\fib\boxtimes\vecG_{\mbb Z/2\mbb Z}^{-1}\boxtimes\vecG_{\mbb Z/2\mbb Z}^{-1}$, the candidates with $\fp_{\mc B}=2$ can admit MFCs $\mc B_A^0$, but those with $\fp_{\mc B}=3$ do not, and are ruled out. Note that not all candidates in (\ref{cXYZtrivial}) are separable. The same observation as in section \ref{Z2Z2Z2} leads to
\begin{equation}
    A\cong\begin{cases}1&(\text{all MFCs}),\\1\oplus X&(\text{those in }(\ref{cXYZtrivial})),\\1\oplus Y&(d_X,d_Y,d_Z,h_X,h_Y,h_Z)=(1,1,1,0,0,\frac12).\end{cases}\quad(\mods1\text{ for }h)\label{fibtcetale}
\end{equation}
In order to figure out $\mc B_A$, we have to find NIM-reps. Here, since NIM-reps do not depend on conformal dimensions, we immediately learn\footnote{If one searches for NIM-reps, one realizes just names of matrices change. One finds
\[ n_1=1_4=n_X,\quad n_Y=\begin{pmatrix}0&1&0&0\\1&0&0&0\\0&0&0&1\\0&0&1&0\end{pmatrix}=n_Z,\quad n_T=\begin{pmatrix}0&0&1&0\\0&0&0&1\\1&0&0&1\\0&1&1&0\end{pmatrix}=n_V,\quad n_U=\begin{pmatrix}0&0&0&1\\0&0&1&0\\0&1&1&0\\1&0&0&1\end{pmatrix}=n_W \]
for $A\cong 1\oplus X$, and
\[ n_1=1_4=n_Y,\quad n_X=\begin{pmatrix}0&1&0&0\\1&0&0&0\\0&0&0&1\\0&0&1&0\end{pmatrix}=n_Z,\quad n_U=\begin{pmatrix}0&0&1&0\\0&0&0&1\\1&0&0&1\\0&1&1&0\end{pmatrix}=n_V,\quad n_T=\begin{pmatrix}0&0&0&1\\0&0&1&0\\0&1&1&0\\1&0&0&1\end{pmatrix}=n_W \]
for $A\cong 1\oplus Y$. They give identifications
\[ m_1\cong1\oplus X,\quad m_2\cong Y\oplus Z,\quad m_3\cong T\oplus V,\quad m_4\cong U\oplus W, \]
and
\[ m_1\cong1\oplus Y,\quad m_2\cong X\oplus Z,\quad m_3\cong U\oplus V,\quad m_4\cong T\oplus W, \]
respectively.\label{fibtcNIMrep}}
\[ \mc B_A^0\simeq\fib,\quad\mc B_A\simeq su(2)_3. \]

We conclude
\begin{table}[H]
\begin{center}
\begin{tabular}{c|c|c|c}
    Connected étale algebra $A$&$\mc B_A$&$\rank(\mc B_A)$&Lagrangian?\\\hline
    $1$&$\mc B$&$8$&No\\
    $1\oplus Z$ for (\ref{fibZ2Z21Zetale})&$su(2)_3$&4&No
\end{tabular}.
\end{center}
\caption{Connected étale algebras in rank eight MFC $\mcal B\simeq\fib\boxtimes\vecG_{\mbb Z/2\mbb Z}^{-1}\boxtimes\vecG_{\mbb Z/2\mbb Z}^{-1}$}\label{rank8fibZ2Z2results}
\end{table}
\hspace{-17pt}and
\begin{table}[H]
\begin{center}
\begin{tabular}{c|c|c|c}
    Connected étale algebra $A$&$\mc B_A$&$\rank(\mc B_A)$&Lagrangian?\\\hline
    $1$&$\mc B$&$8$&No\\
    $1\oplus X$ for (\ref{fibtcetale})&$su(2)_3$&4&No\\
    $1\oplus Y$ for (\ref{fibtcetale})&$su(2)_3$&4&No
\end{tabular}.
\end{center}
\caption{Connected étale algebras in rank eight MFC $\mcal B\simeq\fib\boxtimes\tc$}\label{rank8fibtoriccoderesults}
\end{table}
\hspace{-17pt}In particular, 64 MFCs $\mc B\simeq\fib\boxtimes\vecG_{\mbb Z/2\mbb Z}^{-1}\boxtimes\vecG_{\mbb Z/2\mbb Z}^{-1}$'s not in (\ref{fibZ2Z21Zetale}), and 24 MFCs $\mc B\simeq\fib\boxtimes\tc$'s with
\begin{align*}
    \hspace{-30pt}(d_X,d_Y,d_Z,d_W,h_X,h_Y,h_Z,h_W)&=(1,1,1,\zeta,\frac12,\frac12,\frac12,\frac25),(1,1,1,\zeta,\frac12,\frac12,\frac12,\frac35),\\
    &~~~~(1,1,1,-\zeta^{-1},\frac12,\frac12,\frac12,\frac15),(1,1,1,-\zeta^{-1},\frac12,\frac12,\frac12,\frac45),\\
    &~~~~(1,-1,-1,\zeta,\frac12,0,0,\frac25),(1,1,1,\zeta,\frac12,0,0,\frac35),\\
    &~~~~(1,-1,-1,\zeta,\frac12,
    \frac12,\frac12,\frac25),(1,1,1,\zeta,\frac12,\frac12,\frac12,\frac35),\\
    &~~~~(1,-1,-1,-\zeta^{-1},\frac12,0,0,\frac15),(1,-1,-1,-\zeta^{-1},\frac12,0,0,\frac45),\\
    &~~~~(1,-1,-1,-\zeta^{-1},\frac12,\frac12,\frac12,\frac15),(1,-1,-1,-\zeta^{-1},\frac12,\frac12,\frac12,\frac45)\quad(\mods1\text{ for }h)
\end{align*}
are completely anisotropic.

\subsubsection{$\mc B\simeq\fib\boxtimes\vecG_{\mbb Z/4\mbb Z}^\alpha$}
The MFCs have eight simple objects $\{1,X,Y,Z,T,U,V,W\}$ obeying monoidal products
\begin{table}[H]
\begin{center}
\begin{tabular}{c|c|c|c|c|c|c|c|c}
    $\otimes$&$1$&$X$&$Y$&$Z$&$T$&$U$&$V$&$W$\\\hline
    $1$&$1$&$X$&$Y$&$Z$&$T$&$U$&$V$&$W$\\\hline
    $X$&&$1$&$Z$&$Y$&$U$&$T$&$W$&$V$\\\hline
    $Y$&&&$X$&$1$&$V$&$W$&$U$&$T$\\\hline
    $Z$&&&&$X$&$W$&$V$&$T$&$U$\\\hline
    $T$&&&&&$1\oplus U$&$X\oplus T$&$Y\oplus W$&$Z\oplus V$\\\hline
    $U$&&&&&&$1\oplus U$&$Z\oplus V$&$Y\oplus W$\\\hline
    $V$&&&&&&&$X\oplus T$&$1\oplus U$\\\hline
    $W$&&&&&&&&$X\oplus T$
\end{tabular}.
\end{center}
\end{table}
\hspace{-17pt}(One can identify $\fib=\{1,U\},\vecG_{\mbb Z/4\mbb Z}^\alpha=\{1,X,Y,Z\}$, and $T\cong X\otimes U,V\cong Z\otimes U,W\cong Y\otimes U$.) Thus, they have
\begin{align*}
    &\fp_{\mc B}(1)=\fp_{\mc B}(X)=\fp_{\mc B}(Y)=\fp_{\mc B}(Z)=1,\\
    &\fp_{\mc B}(T)=\fp_{\mc B}(U)=\fp_{\mc B}(V)=\fp_{\mc B}(W)=\zeta:=\frac{1+\sqrt5}2,
\end{align*}
and
\[ \fp(\mc B)=10+2\sqrt5\approx14.5. \]
Their quantum dimensions $d_j$'s are solutions of the same multiplication rules $d_id_j=\sum_{k=1}^8{N_{ij}}^kd_k$. There are four solutions
\begin{align*}
    (d_X,d_Y,d_Z,d_T,d_U,d_V,d_W)=&(1,-1,-1,\zeta,\zeta,-\zeta,-\zeta),(1,-1,-1,-\zeta^{-1},-\zeta^{-1},\zeta^{-1},\zeta^{-1}),\\
    &(1,1,1,-\zeta^{-1},-\zeta^{-1},-\zeta^{-1},-\zeta^{-1}),(1,1,1,\zeta,\zeta,\zeta,\zeta).
\end{align*}
(Only the last quantum dimensions give unitary MFCs.) They have categorical dimensions
\[ D^2(\mc B)=10-2\sqrt5(\approx5.5),\quad10+2\sqrt5. \]
They have eight conformal dimensions.
\begin{itemize}
    \item $d_U=\zeta$. These have
    \begin{align*}
        (h_X,h_Y,h_Z,h_T,h_U,h_V,h_W)&=(\frac12,\frac18,\frac18,\frac1{10},\frac35,\frac{29}{40},\frac{29}{40}),(\frac12,\frac18,\frac18,\frac9{10},\frac25,\frac{21}{40},\frac{21}{40}),\\
        &~~~~(\frac12,\frac38,\frac38,\frac1{10},\frac35,\frac{39}{40},\frac{39}{40}),(\frac12,\frac38,\frac38,\frac9{10},\frac25,\frac{31}{40},\frac{31}{40}),\\
        &~~~~(\frac12,\frac58,\frac58,\frac1{10},\frac35,\frac{9}{40},\frac{9}{40}),(\frac12,\frac58,\frac58,\frac9{10},\frac25,\frac{1}{40},\frac{1}{40}),\\
        &~~~~(\frac12,\frac78,\frac78,\frac1{10},\frac35,\frac{19}{40},\frac{19}{40}),(\frac12,\frac78,\frac78,\frac9{10},\frac25,\frac{11}{40},\frac{11}{40})\quad(\mods1).
    \end{align*}
    \item $d_U=-\zeta^{-1}$. These have
    \begin{align*}
        (h_X,h_Y,h_Z,h_T,h_U,h_V,h_W)&=(\frac12,\frac18,\frac18,\frac3{10},\frac45,\frac{37}{40},\frac{37}{40}),(\frac12,\frac18,\frac18,\frac7{10},\frac15,\frac{13}{40},\frac{13}{40}),\\
        &~~~~(\frac12,\frac38,\frac38,\frac3{10},\frac45,\frac{7}{40},\frac{7}{40}),(\frac12,\frac38,\frac38,\frac7{10},\frac15,\frac{23}{40},\frac{23}{40}),\\
        &~~~~(\frac12,\frac58,\frac58,\frac3{10},\frac45,\frac{17}{40},\frac{17}{40}),(\frac12,\frac58,\frac58,\frac7{10},\frac15,\frac{33}{40},\frac{33}{40}),\\
        &~~~~(\frac12,\frac78,\frac78,\frac3{10},\frac45,\frac{27}{40},\frac{27}{40}),(\frac12,\frac78,\frac78,\frac7{10},\frac15,\frac{3}{40},\frac{3}{40})\quad(\mods1).
    \end{align*}
\end{itemize}
The $S$-matrices are given by
\[ \widetilde S=\begin{pmatrix}1&d_X&d_Y&d_Z&d_Xd_U&d_U&d_Zd_U&d_Yd_U\\d_X&1&-d_Y&-d_Z&d_Xd_U&d_U&-d_Zd_U&-d_Yd_U\\d_Y&-d_Y&\pm i&\mp i&-d_Yd_U&d_Yd_U&\mp i\cdot d_U&\pm i\cdot d_U\\d_Z&-d_Z&\mp i&\pm i&-d_Zd_U&d_Zd_U&\pm i\cdot d_U&\mp i\cdot d_U\\d_Xd_U&d_Xd_U&-d_Yd_U&-d_Zd_U&-1&-d_X&d_Z&d_Y\\d_U&d_U&d_Yd_U&d_Zd_U&-d_X&-1&-d_Z&-d_Y\\d_Zd_U&-d_Zd_U&\mp i\cdot d_U&\pm i\cdot d_U&d_Z&-d_Z&\mp i&\pm i\\d_Yd_U&-d_Yd_U&\pm i\cdot d_U&\mp i\cdot d_U&d_Y&-d_Y&\pm i&\mp i\end{pmatrix}. \]
They have
\[ c(\mc B)=c(\fib)+c(\vecG_{\mbb Z/4\mbb Z}^\alpha)\quad(\mods8) \]
where
\[ c(\fib)=\begin{cases}\frac25&(h_\fib=\frac15),\\-\frac25&(h_\fib=\frac45),\\\frac{14}5&(h_\fib=\frac25),\\-\frac{14}5&(h_\fib=\frac35).\end{cases}\quad c(\vecG_{\mbb Z/4\mbb Z}^\alpha)=\begin{cases}1&(h_{\mbb Z/4\mbb Z}=\frac18),\\3&(h_{\mbb Z/4\mbb Z}=\frac38),\\-3&(h_{\mbb Z/4\mbb Z}=\frac58),\\-1&(h_{\mbb Z/4\mbb Z}=\frac78),\end{cases}\quad(\mods8). \]
There are
\[ 4(\text{quantum dimensions})\times8(\text{conformal dimensions})\times2(\text{categorical dimensions})=64 \]
MFCs, among which those 16 with the fourth quantum dimensions give unitary MFCs. We classify connected étale algebras in all 64 MFCs simultaneously.

We work with an ansatz
\[ A\cong1\oplus n_XX\oplus n_YY\oplus n_ZZ\oplus n_TT\oplus n_UU\oplus n_VV\oplus n_WW \]
with $n_j\in\mbb N$. It has
\[ \fp_{\mc B}(A)=1+n_X+n_Y+n_Z+\zeta(n_T+n_U+n_V+n_W). \]
For this to obey (\ref{FPdimA2bound}), exactly the same as the previous example, the natural numbers can take only 26 values
\begin{align*}
    (n_X,n_Y,n_Z,n_T,n_U,n_V,n_W)=&(0,0,0,0,0,0,0),(1,0,0,0,0,0,0),(2,0,0,0,0,0,0),\\
    &(1,1,0,0,0,0,0),(1,0,1,0,0,0,0),(1,0,0,1,0,0,0),\\
    &(1,0,0,0,1,0,0),(1,0,0,0,0,1,0),(1,0,0,0,0,0,1),\\
    &(0,1,0,0,0,0,0),(0,2,0,0,0,0,0),(0,1,1,0,0,0,0),\\
    &(0,1,0,1,0,0,0),(0,1,0,0,1,0,0),(0,1,0,0,0,1,0),\\
    &(0,1,0,0,0,0,1),(0,0,1,0,0,0,0),(0,0,2,0,0,0,0),\\
    &(0,0,1,1,0,0,0),(0,0,1,0,1,0,0),(0,0,1,0,0,1,0),\\
    &(0,0,1,0,0,0,1),(0,0,0,1,0,0,0),(0,0,0,0,1,0,0),\\
    &(0,0,0,0,0,1,0),(0,0,0,0,0,0,1).
\end{align*}
The first is nothing but the trivial connected étale algebra $A\cong1$ giving $\mc B\simeq\mc B_A^0\simeq\mc B_A$. Those with $X$'s cannot give commutative algebras because $X$ has $(d_X,h_X)=(1,\frac12)$ (mod 1 for $h$) and $c_{X,X}\cong-id_1$ \cite{KK23preMFC}. The others also fail to be commutative because nontrivial simple objects entering them have nontrivial conformal dimensions.

We conclude
\begin{table}[H]
\begin{center}
\begin{tabular}{c|c|c|c}
    Connected étale algebra $A$&$\mc B_A$&$\rank(\mc B_A)$&Lagrangian?\\\hline
    $1$&$\mc B$&$8$&No
\end{tabular}.
\end{center}
\caption{Connected étale algebras in rank eight MFC $\mcal B\simeq\fib\boxtimes\vecG_{\mbb Z/4\mbb Z}^\alpha$}\label{rank8fibZ4results}
\end{table}
\hspace{-17pt}That is, all the 64 MFCs $\mc B\simeq\fib\boxtimes\vecG_{\mbb Z/4\mbb Z}^\alpha$'s are completely anisotropic.

\subsubsection{$\mc B\simeq\vecG_{\mbb Z/2\mbb Z}^{-1}\boxtimes\fib\boxtimes\fib$}\label{Z2fibfib}
The MFCs have eight simple objects $\{1,X,Y,Z,T,U,V,W\}$ obeying monoidal products
\begin{table}[H]
\begin{center}
\begin{tabular}{c|c|c|c|c|c|c|c|c}
    $\otimes$&$1$&$X$&$Y$&$Z$&$T$&$U$&$V$&$W$\\\hline
    $1$&$1$&$X$&$Y$&$Z$&$T$&$U$&$V$&$W$\\\hline
    $X$&&$1$&$U$&$T$&$Z$&$Y$&$W$&$V$\\\hline
    $Y$&&&$1\oplus U$&$W$&$V$&$X\oplus Y$&$T\oplus W$&$Z\oplus V$\\\hline
    $Z$&&&&$1\oplus T$&$X\oplus Z$&$V$&$U\oplus W$&$Y\oplus V$\\\hline
    $T$&&&&&$1\oplus T$&$W$&$Y\oplus V$&$U\oplus W$\\\hline
    $U$&&&&&&$1\oplus U$&$Z\oplus V$&$T\oplus W$\\\hline
    $V$&&&&&&&$1\oplus T\oplus U\oplus W$&$X\oplus Y\oplus Z\oplus V$\\\hline
    $W$&&&&&&&&$1\oplus T\oplus U\oplus W$
\end{tabular}.
\end{center}
\end{table}
\hspace{-17pt}(One can identify $\vecG_{\mbb Z/2\mbb Z}^{-1}=\{1,X\},\fib=\{1,T\},\{1,U\}$, and $Y\cong X\otimes U,Z\cong X\otimes T,V\cong X\otimes T\otimes U,W\cong T\otimes U$.) Thus, they have
\begin{align*}
    \fp_{\mc B}(1)=1=\fp_{\mc B}(X),\quad&\fp_{\mc B}(Y)=\fp_{\mc B}(Z)=\fp_{\mc B}(T)=\fp_{\mc B}(U)=\zeta,\\
    &\fp_{\mc B}(V)=\frac{3+\sqrt5}2=\fp_{\mc B}(W),
\end{align*}
and
\[ \fp(\mc B)=15+5\sqrt5\approx26.2. \]
Their quantum dimensions $d_j$'s are solutions of the same multiplication rules $d_id_j=\sum_{k=1}^8{N_{ij}}^kd_k$. There are eight solutions
\begin{align*}
    &(d_X,d_Y,d_Z,d_T,d_U,d_V,d_W)\\
    &=(-1,\zeta^{-1},\zeta^{-1},-\zeta^{-1},-\zeta^{-1},\frac{\sqrt5-3}2,\frac{3-\sqrt5}2),(1,-\zeta^{-1},-\zeta^{-1},-\zeta^{-1},-\zeta^{-1},\frac{3-\sqrt5}2,\frac{3-\sqrt5}2),\\
    &~~~~(-1,-\zeta,\zeta^{-1},-\zeta^{-1},\zeta,1,-1),(1,\zeta,-\zeta^{-1},-\zeta^{-1},\zeta,-1,-1),\\
    &~~~~(-1,\zeta^{-1},-\zeta,\zeta,-\zeta^{-1},1,-1),(1,-\zeta^{-1},\zeta,\zeta,-\zeta^{-1},-1,-1),\\
    &~~~~(-1,-\zeta,-\zeta,\zeta,\zeta,-\frac{3+\sqrt5}2,\frac{3+\sqrt5}2),(1,\zeta,\zeta,\zeta,\zeta,\frac{3+\sqrt5}2,\frac{3+\sqrt5}2)
\end{align*}
with categorical dimensions
\[ D^2(\mc B)=15-5\sqrt5(\approx3.8),\quad10,\quad15+5\sqrt5. \]
Only the last quantum dimensions give unitary MFCs.

In order to list up their conformal dimensions without double-counting, we perform case analysis.
\begin{itemize}
    \item $(d_X,d_T,d_U)=(1,\zeta,\zeta)$. This gives unitary MFCs. Different MFCs are given by
    \begin{align*}
        (h_X,h_Y,h_Z,h_T,h_U,h_V,h_W)&=(\frac14,\frac{13}{20},\frac{13}{20},\frac25,\frac25,\frac1{20},\frac45),(\frac14,\frac{17}{20},\frac{13}{20},\frac25,\frac35,\frac14,0),\\
        &~~~~(\frac14,\frac{17}{20},\frac{17}{20},\frac35,\frac35,\frac9{20},\frac15),(\frac34,\frac{3}{20},\frac{3}{20},\frac25,\frac25,\frac{11}{20},\frac45),\\
        &~~~~(\frac34,\frac{7}{20},\frac{3}{20},\frac25,\frac35,\frac34,0),(\frac34,\frac{7}{20},\frac{7}{20},\frac35,\frac35,\frac{19}{20},\frac15)\quad(\mods1).
    \end{align*}
    Including two signs of categorical dimensions, we have 12 unitary MFCs.
    \item $(d_X,d_T,d_U)=(-1,\zeta,\zeta)$. The sets of conformal dimensions are the same as the previous case. There are 12 MFCs.
    \item $(d_X,d_T,d_U)=(1,\zeta,-\zeta^{-1})$. Different MFCs are given by
    \begin{align*}
        (h_X,h_Y,h_Z,h_T,h_U,h_V,h_W)&=(\frac14,\frac{9}{20},\frac{13}{20},\frac25,\frac15,\frac{17}{20},\frac35),(\frac14,\frac{1}{20},\frac{13}{20},\frac25,\frac45,\frac9{20},\frac15),\\
        &~~~~(\frac14,\frac{9}{20},\frac{17}{20},\frac35,\frac15,\frac1{20},\frac45),(\frac14,\frac{1}{20},\frac{17}{20},\frac35,\frac45,\frac{13}{20},\frac25),\\
        &~~~~(\frac34,\frac{19}{20},\frac{3}{20},\frac25,\frac15,\frac7{20},\frac35),(\frac34,\frac{11}{20},\frac{3}{20},\frac25,\frac45,\frac{19}{20},\frac15),\\
        &~~~~(\frac34,\frac{19}{20},\frac{7}{20},\frac35,\frac15,\frac{11}{20},\frac45),(\frac34,\frac{11}{20},\frac{7}{20},\frac35,\frac45,\frac{3}{20},\frac25)\quad(\mods1).
    \end{align*}
    With two signs of categorical dimensions, there are 16 MFCs.
    \item $(d_X,d_T,d_U)=(-1,\zeta,-\zeta^{-1})$. The sets of conformal dimensions are the same as the previous case. There are 16 MFCs.
    \item $(d_X,d_T,d_U)=(1,-\zeta^{-1},-\zeta^{-1})$. Different MFCs are given by
    \begin{align*}
        (h_X,h_Y,h_Z,h_T,h_U,h_V,h_W)&=(\frac14,\frac{9}{20},\frac{9}{20},\frac15,\frac15,\frac{13}{20},\frac25),(\frac14,\frac{1}{20},\frac{9}{20},\frac15,\frac45,\frac14,0),\\
        &~~~~(\frac14,\frac{1}{20},\frac{1}{20},\frac45,\frac45,\frac{17}{20},\frac35),(\frac34,\frac{19}{20},\frac{19}{20},\frac15,\frac15,\frac{3}{20},\frac25),\\
        &~~~~(\frac34,\frac{11}{20},\frac{19}{20},\frac15,\frac45,\frac34,0),(\frac34,\frac{11}{20},\frac{11}{20},\frac45,\frac45,\frac{7}{20},\frac35)\quad(\mods1).
    \end{align*}
    There are 12 MFCs.
    \item $(d_X,d_T,d_U)=(-1,-\zeta^{-1},-\zeta^{-1})$. The sets of conformal dimensions are the same as the previous case, and there are 12 MFCs.
\end{itemize}
The $S$-matrices are given by
\[ \widetilde S=\begin{pmatrix}1&d_X&d_Xd_U&d_Xd_T&d_T&d_U&d_Xd_Td_U&d_Td_U\\d_X&-1&-d_U&-d_T&d_Xd_T&d_Xd_U&-d_Td_U&d_Xd_Td_U\\d_Xd_U&-d_U&1&-d_Td_U&d_Xd_Td_U&-d_X&d_T&-d_Xd_T\\d_Xd_T&-d_T&-d_Td_U&1&-d_X&d_Xd_Td_U&d_U&-d_Xd_U\\d_T&d_Xd_T&d_Xd_Td_U&-d_X&-1&d_Td_U&-d_Xd_U&-d_U\\d_U&d_Xd_U&-d_X&d_Xd_Td_U&d_Td_U&-1&-d_Xd_T&-d_T\\d_Xd_Td_U&-d_Td_U&d_T&d_U&-d_Xd_U&-d_Xd_T&-1&d_X\\d_Td_U&d_Xd_Td_U&-d_Xd_T&-d_Xd_U&-d_U&-d_T&d_X&1\end{pmatrix}. \]
They have additive central charges
\[ c(\mc B)=c(\vecG_{\mbb Z/2\mbb Z}^{-1})+c(\fib)+c(\fib)\quad(\mods8) \]
where
\[ c(\vecG_{\mbb Z/2\mbb Z}^{-1})=\begin{cases}1&(h_X=\frac14),\\-1&(h_X=\frac34).\end{cases},\quad c(\fib)=\begin{cases}\frac25&(h_{T,U}=\frac15),\\-\frac25&(h_{T,U}=\frac45),\\\frac{14}5&(h_{T,U}=\frac25),\\-\frac{14}5&(h_{T,U}=\frac35).\end{cases}\quad(\mods8) \]
There are
\[ 2\times12+2\times16+2\times12=80 \]
MFCs, among which those 12 with the last quantum dimensions are unitary. We classify connected étale algebras in all 80 MFCs simultaneously.

An ansatz
\[ A\cong1\oplus n_XX\oplus n_YY\oplus n_ZZ\oplus n_TT\oplus n_UU\oplus n_VV\oplus n_WW \]
with $n_j\in\mbb N$ has
\[ \fp_{\mc B}(A)=1+n_X+\zeta(n_Y+n_Z+n_T+n_U)+\frac{3+\sqrt5}2(n_V+n_W). \]
For this to obey (\ref{FPdimA2bound}), the natural numbers can take only 31 values
\begin{align*}
    (n_X,n_Y,n_Z,n_T,n_U,n_V,n_W)=&(0,0,0,0,0,0,0),(1,0,0,0,0,0,0),(2,0,0,0,0,0,0),\\
    &(3,0,0,0,0,0,0),(4,0,0,0,0,0,0),(2,1,0,0,0,0,0),\\
    &(2,0,1,0,0,0,0),(2,0,0,1,0,0,0),(2,0,0,0,1,0,0),\\
    &(1,1,0,0,0,0,0),(1,0,1,0,0,0,0),(1,0,0,1,0,0,0),\\
    &(1,0,0,0,1,0,0),(1,0,0,0,0,1,0),(1,0,0,0,0,0,1),\\
    &(0,1,0,0,0,0,0),(0,2,0,0,0,0,0),(0,1,1,0,0,0,0),\\
    &(0,1,0,1,0,0,0),(0,1,0,0,1,0,0),(0,0,1,0,0,0,0),\\
    &(0,0,2,0,0,0,0),(0,0,1,1,0,0,0),(0,0,1,0,1,0,0),\\
    &(0,0,0,1,0,0,0),(0,0,0,2,0,0,0),(0,0,0,1,1,0,0),\\
    &(0,0,0,0,1,0,0),(0,0,0,0,2,0,0),(0,0,0,0,0,1,0),\\
    &(0,0,0,0,0,0,1).
\end{align*}
The first is nothing but the trivial connected étale algebra $A\cong1$ giving $\mc B_A^0\simeq\mc B_A\simeq\mc B$. Next, those with $X$'s cannot be commutative because the simple object has $(d_X,h_X)=(\pm1,\frac14)$ (mod $1/2$ for $h$) and $c_{X,X}\cong\pm i\cdot id_1$ \cite{KK23preMFC}. Thirdly, the simple objects $Y,Z,T,U,V$ have nontrivial conformal dimensions, and candidates with them are ruled out. Thus, we are left with those with just $W$'s. Concretely, the only nontrivial candidate is $n_W=1$ or $A\cong1\oplus W$. When $h_W=0$ mod 1, it can be commutative. Indeed, it is known \cite{BD11} to be a commutative algebra in $\fib\boxtimes\fib$. It further turns out to be separable. To show this point, we identify $\mc B_A$.

It has $\fp_{\mc B}(A)=\frac{5+\sqrt5}2$, and demands
\[ \fp(\mc B_A^0)=2,\quad\fp(\mc B_A)=5+\sqrt5. \]
Calculating $b_j\otimes A$, we find simple objects
\[ 1\oplus W,\quad X\oplus V,\quad Y\oplus Z\oplus V,\quad T\oplus U\oplus W. \]
This suggests $\mc B_A$ has rank four. Indeed, we find a four-dimensional NIM-rep
\begin{align*}
    n_1=1_4,\quad n_X=\begin{pmatrix}0&1&0&0\\1&0&0&0\\0&0&0&1\\0&0&1&0\end{pmatrix},\quad n_Y=\begin{pmatrix}0&0&1&0\\0&0&0&1\\1&0&0&1\\0&1&1&0\end{pmatrix}=n_Z,\\
    n_T=\begin{pmatrix}0&0&0&1\\0&0&1&0\\0&1&1&0\\1&0&0&1\end{pmatrix}=n_U,\quad n_V=\begin{pmatrix}0&1&1&0\\1&0&0&1\\1&0&0&2\\0&1&2&0\end{pmatrix},\quad n_W=\begin{pmatrix}1&0&0&1\\0&1&1&0\\0&1&2&0\\1&0&0&2\end{pmatrix}.
\end{align*}
This gives identifications
\[ m_1\cong1\oplus W,\quad m_2\cong X\oplus V,\quad m_3\cong Y\oplus Z\oplus V,\quad m_4\cong T\oplus U\oplus W. \]
The first two have $\fp_{\mc B_A}=1$, and form $\mc B_A^0\simeq\vecG_{\mbb Z/2\mbb Z}^{-1}$. This also matches central charges because $c(\fib\boxtimes\fib)=0$ (mod 8) when $h_W=0$ (mod 1), and $h_X=h_V$ gives the conformal dimension of $m_2\in\mc B_A^0$. One can identify $\mc B_A$ by working out the monoidal products:
\begin{table}[H]
\begin{center}
\begin{tabular}{c|c|c|c|c}
    $\otimes_A$&$m_1$&$m_2$&$m_3$&$m_4$\\\hline
    $m_1$&$m_1$&$m_2$&$m_3$&$m_4$\\\hline
    $m_2$&&$m_1$&$m_4$&$m_3$\\\hline
    $m_3$&&&$m_1\oplus m_4$&$m_2\oplus m_3$\\\hline
    $m_4$&&&&$m_1\oplus m_4$
\end{tabular}.
\end{center}
\end{table}
\hspace{-17pt}We get
\[ \mc B_A\simeq su(2)_3. \]
Since this is semisimple, $A$ is separable, hence étale:
\begin{equation}
    A\cong1\oplus W\quad(d_T,d_U,h_T,h_U)=(\zeta,\zeta,\frac25,\frac35),(-\zeta^{-1},-\zeta^{-1},\frac15,\frac45).\quad(\mods1\text{ for }h)\label{Z2fibfibetale}
\end{equation}

We conclude
\begin{table}[H]
\begin{center}
\begin{tabular}{c|c|c|c}
    Connected étale algebra $A$&$\mc B_A$&$\rank(\mc B_A)$&Lagrangian?\\\hline
    $1$&$\mc B$&$8$&No\\
    $1\oplus W$ for (\ref{Z2fibfibetale})&$su(2)_3$&$4$&No
\end{tabular}.
\end{center}
\caption{Connected étale algebras in rank eight MFC $\mcal B\simeq\vecG_{\mbb Z/2\mbb Z}^{-1}\boxtimes\fib\boxtimes\fib$}\label{rank8Z2fibfibresults}
\end{table}
\hspace{-17pt}Namely, those 16 in (\ref{Z2fibfibetale}) fail to be completely anisotropic, while the other 64 MFCs $\mc B\simeq\vecG_{\mbb Z/2\mbb Z}^{-1}\boxtimes\fib\boxtimes\fib$'s are completely anisotropic.

\subsubsection{$\mc B\simeq so(9)_2$}\label{so92}
The MFCs have eight simple objects $\{1,X,Y,Z,T,U,V,W\}$ obeying monoidal products
\begin{table}[H]
\begin{center}
\makebox[1 \textwidth][c]{       
\resizebox{1.2 \textwidth}{!}{\begin{tabular}{c|c|c|c|c|c|c|c|c}
    $\otimes$&$1$&$X$&$Y$&$Z$&$T$&$U$&$V$&$W$\\\hline
    $1$&$1$&$X$&$Y$&$Z$&$T$&$U$&$V$&$W$\\\hline
    $X$&&$1$&$Y$&$Z$&$T$&$U$&$W$&$V$\\\hline
    $Y$&&&$1\oplus X\oplus U$&$T\oplus U$&$Z\oplus T$&$Y\oplus Z$&$V\oplus W$&$V\oplus W$\\\hline
    $Z$&&&&$1\oplus X\oplus Z$&$Y\oplus U$&$Y\oplus T$&$V\oplus W$&$V\oplus W$\\\hline
    $T$&&&&&$1\oplus X\oplus Y$&$Z\oplus U$&$V\oplus W$&$V\oplus W$\\\hline
    $U$&&&&&&$1\oplus X\oplus T$&$V\oplus W$&$V\oplus W$\\\hline
    $V$&&&&&&&$1\oplus Y\oplus Z\oplus T\oplus U$&$X\oplus Y\oplus Z\oplus T\oplus U$\\\hline
    $W$&&&&&&&&$1\oplus Y\oplus Z\oplus T\oplus U$
\end{tabular}.}}
\end{center}
\end{table}
\hspace{-17pt}Thus, they have
\begin{align*}
    \fp_{\mc B}(1)=1=\fp_{\mc B}(X),\quad&\fp_{\mc B}(Y)=\fp_{\mc B}(Z)=\fp_{\mc B}(T)=\fp_{\mc B}(U)=2,\\
    &\fp_{\mc B}(V)=\fp_{\mc B}(W)=3,
\end{align*}
and
\[ \fp(\mc B)=36. \]
Their quantum dimensions $d_j$'s are solutions of the same multiplication rules $d_id_j=\sum_{k=1}^8{N_{ij}}^kd_k$. There are two (nonzero) solutions
\[ (d_X,d_Y,d_Z,d_T,d_U,d_V,d_W)=(1,2,2,2,2,-3,-3),(1,2,2,2,2,3,3). \]
Only the second gives unitary MFCs. They both have the same categorical dimension
\[ D^2(\mc B)=36, \]
and four conformal dimensions\footnote{Naively, one finds 24 conformal dimensions, but the other 20 are given by permutations $(YTU)$ and $(VW)$. Thus, different MFCs are labeled by the four conformal dimensions in the main text.}
\begin{align*}
    (h_X,h_Y,h_Z,h_T,h_U,h_V,h_W)&=(0,\frac19,0,\frac79,\frac49,0,\frac12),(0,\frac19,0,\frac79,\frac49,\frac14,\frac34),\\
    &~~~~(0,\frac29,0,\frac59,\frac89,0,\frac12),(0,\frac29,0,\frac59,\frac89,\frac14,\frac34)\quad(\mods1).
\end{align*}
The $S$-matrices are given by
\[ \widetilde S=\begin{pmatrix}1&d_X&d_Y&d_Z&d_T&d_U&d_V&d_W\\d_X&d_X&d_Y&d_Z&d_T&d_U&-d_V&-d_W\\d_Y&d_Y&s&-2&s'&s''&0&0\\d_Z&d_Z&-2&4&-2&-2&0&0\\d_T&d_T&s'&-2&s''&s&0&0\\d_U&d_U&s''&-2&s&s'&0&0\\d_V&-d_V&0&0&0&0&-d_V&d_W\\d_W&-d_W&0&0&0&0&d_W&-d_W\end{pmatrix} \]
with
\[ s=4\sin\frac\pi{18},\quad s'=4\cos\frac{2\pi}9,\quad s''=-4\cos\frac\pi9, \]
or their permutations $(ss's'')$. They have additive central charges
\[ c(\mc B)=0\quad(\mods8). \]
There are
\[ 2(\text{quantum dimensions})\times4(\text{conformal dimensions})\times2(\text{categorical dimensions})=16 \]
MFCs, among which those eight with the second quantum dimensions are unitary. We classify connected étale algebras in all 16 MFCs simultaneously.

An ansatz
\[ A\cong1\oplus n_XX\oplus n_YY\oplus n_ZZ\oplus n_TT\oplus n_UU\oplus n_VV\oplus n_WW \]
with $n_j\in\mbb N$ has
\[ \fp_{\mc B}(A)=1+n_X+2(n_Y+n_Z+n_T+n_U)+3(n_V+n_W). \]
For this to obey (\ref{FPdimA2bound}), the natural numbers can take only 56 values. The sets contain one with all $n_j$'s be zero. It is nothing but the trivial connected étale algebra $A\cong1$ giving $\mc B_A^0\simeq\mc B_A\simeq\mc B$. Next, let us study candidates with nontrivial simple object(s). We can rule out those with $Y,T,U,W$ because they have nontrivial conformal dimensions and fail to meet the necessary condition (\ref{commutativealgnecessary}). Thus, we are left with candidates with just $X,Z,V$. Note that $\{1,X,Z\}$ form symmetric pre-modular fusion subcategory $\text{Rep}(S_3)$.

Solving (\ref{FPdimA2bound}) by setting $n_Y,n_T,n_U,n_W$ to zero, we get 16 sets
\begin{align*}
    (n_X,n_Z,n_V)&=(0,0,0),(1,0,0),(2,0,0),(3,0,0),\\&~~~~(4,0,0),(5,0,0),(3,1,0),(2,1,0),\\
    &~~~~(2,0,1),(1,1,0),(1,2,0),(1,0,1),\\
    &~~~~(0,2,0),(0,1,1),(0,1,0),(0,0,1).
\end{align*}
Some of these are ruled out by studying Frobenius-Perron dimensions. The three candidates
\[ (n_X,n_Z,n_V)=(3,0,0),(1,1,0),(0,0,1) \]
have $\fp_{\mc B}=4$ and demand $\fp(\mc B_A^0)=\frac{36}{16}$. However, there is no MFC with such Frobenius-Perron dimensions. Thus, the three candidates are ruled out. Similarly, four candidates
\[ (n_X,n_Z,n_V)=(4,0,0),(2,1,0),(0,2,0),(1,0,1) \]
are ruled out because there is no MFC with $\fp=\frac{36}{25}$. Thus, we are left with nine candidates
\[ (n_X,n_Z,n_V)=(0,0,0),(1,0,0),(2,0,0),(5,0,0),(3,1,0),(2,0,1),(1,2,0),(0,1,1),(0,1,0). \]
The first gives $A\cong1$, and we already found this above. Among the other eight, from the lemma 1, we know three with $(n_X,n_Z,n_V)=(1,0,0),(1,2,0),(0,1,0)$ or $A\cong1\oplus X,1\oplus X\oplus2Z,1\oplus Z$ are commutative algebras because they are so in $A\in\text{Rep}(S_3)\subset so(9)_2$ for all conformal dimensions. It turns out all of these are also separable, hence étale. Let us check this point one after another by identifying $\mc B_A$.

\paragraph{$A\cong1\oplus X$.} Since it has $\fp_{\mc B}(A)=2$, it demands
\[ \fp(\mc B_A^0)=9,\quad\fp(\mc B_A)=18. \]
In view of anyon condensation, we `identify' $1,X$, hence $V,W$. The other $Y,Z,T,U$ `split' into two each. This suggests $\mc B_A$ have rank 10. Indeed, we find a 10-dimensional NIM-rep
\begin{align*}
    \hspace{-50pt}n_1=1_{10}=n_X,&\quad n_Y=\begin{pmatrix}0&1&1&0&0&0&0&0&0&0\\1&0&0&0&0&0&0&1&0&0\\1&0&0&0&0&0&0&0&1&0\\0&0&0&0&0&0&1&0&1&0\\0&0&0&0&0&1&0&1&0&0\\0&0&0&0&1&0&1&0&0&0\\0&0&0&1&0&1&0&0&0&0\\0&1&0&0&1&0&0&0&0&0\\0&0&1&1&0&0&0&0&0&0\\0&0&0&0&0&0&0&0&0&2\end{pmatrix},\quad n_Z=\begin{pmatrix}0&0&0&1&1&0&0&0&0&0\\0&0&0&0&0&1&0&0&1&0\\0&0&0&0&0&0&1&1&0&0\\1&0&0&0&1&0&0&0&0&0\\1&0&0&1&0&0&0&0&0&0\\0&1&0&0&0&0&0&0&1&0\\0&0&1&0&0&0&0&1&0&0\\0&0&1&0&0&0&1&0&0&0\\0&1&0&0&0&1&0&0&0&0\\0&0&0&0&0&0&0&0&0&2\end{pmatrix},\\
    \hspace{-50pt}&n_T=\begin{pmatrix}0&0&0&0&0&1&1&0&0&0\\0&0&0&1&0&0&1&0&0&0\\0&0&0&0&1&1&0&0&0&0\\0&1&0&0&0&0&0&1&0&0\\0&0&1&0&0&0&0&0&1&0\\1&0&1&0&0&0&0&0&0&0\\1&1&0&0&0&0&0&0&0&0\\0&0&0&1&0&0&0&0&1&0\\0&0&0&0&1&0&0&1&0&0\\0&0&0&0&0&0&0&0&0&2\end{pmatrix},\quad n_U=\begin{pmatrix}0&0&0&0&0&0&0&1&1&0\\0&0&1&0&1&0&0&0&0&0\\0&1&0&1&0&0&0&0&0&0\\0&0&1&0&0&1&0&0&0&0\\0&1&0&0&0&0&1&0&0&0\\0&0&0&1&0&0&0&1&0&0\\0&0&0&0&1&0&0&0&1&0\\1&0&0&0&0&1&0&0&0&0\\1&0&0&0&0&0&1&0&0&0\\0&0&0&0&0&0&0&0&0&2\end{pmatrix},\\
    \hspace{-50pt}&n_V=\begin{pmatrix}0&0&0&0&0&0&0&0&0&1\\0&0&0&0&0&0&0&0&0&1\\0&0&0&0&0&0&0&0&0&1\\0&0&0&0&0&0&0&0&0&1\\0&0&0&0&0&0&0&0&0&1\\0&0&0&0&0&0&0&0&0&1\\0&0&0&0&0&0&0&0&0&1\\0&0&0&0&0&0&0&0&0&1\\0&0&0&0&0&0&0&0&0&1\\1&1&1&1&1&1&1&1&1&0\end{pmatrix}=n_W.
\end{align*}
The solution gives identifications
\[ m_1\cong1\oplus X,\quad m_2\cong Y\cong m_3,\quad m_4\cong Z\cong m_5,\quad m_6\cong T\cong m_7,\quad m_8\cong U\cong m_9,\quad m_{10}\cong V\oplus W. \]
They have quantum dimensions (\ref{dBAm})
\[ d_{\mc B_A}(m_j)=1\ (j=1,2,\dots,9),\quad d_{\mc B_A}(m_{10})=\pm3. \]
Furthermore, they obey monoidal products
\[ m_j\otimes_Am_{10}\cong m_{10}\cong m_{10}\otimes_Am_j\ (j=1,2,\dots,9),\quad m_{10}\otimes_Am_{10}\cong\bigoplus_{j=1}^9m_j. \]
This shows
\[ \mc B_A^0\simeq\vecG_{\mbb Z/9\mbb Z}^1,\quad\mc B_A\simeq\text{TY}(\mbb Z/9\mbb Z). \]
The identification also matches central charges. See section \ref{Z9}. Since $\text{TY}(\mbb Z/9\mbb Z)$ is semisimple, $A$ is separable.

\paragraph{$A\cong1\oplus X\oplus2Z$.} It has $\fp_{\mc B}(A)=6$, and demands
\[ \fp(\mc B_A^0)=1,\quad\fp(\mc B_A)=6. \]
We can identify $\mc B_A^0\simeq\vect$. This matches central charges. Calculating $b_j\otimes A$, we find candidate simple objects of $\mc B_A$
\[ 1\oplus X\oplus2Z,\quad Y\oplus T\oplus U,\quad V\oplus W \]
with Frobenius-Perron dimensions one. Logically, the latter two can have coefficients, but the possibilities are ruled out.\footnote{Here is a proof.

First, the third simple object should have coefficient one. Since it always appears in the form $3(V\oplus W)$ in $b_j\otimes A$, one can have $3(V\oplus W)$, or $V\oplus W,2(V\oplus W)$ as simple object(s). The first possibility contributes $3^2$ to $\fp(\mc B_A)$, and it exceeds six. The latter contributes $1^2+2^2$, and together with other simple objects, the contributions exceed six. Thus, the coefficient of $V\oplus W$ should be one.

Second, the second simple object should also have coefficient one. It can have coefficient two. Then, the three simple objects match Frobenius-Perron dimension $1^2+2^2+1^2=6$. Furthermore, we find a three-dimensional NIM-rep
\[ n_1=1_3=n_X,\quad n_Y=n_T=n_U=\begin{pmatrix}0&1&0\\2&1&0\\0&0&2\end{pmatrix},\quad n_Z=\begin{pmatrix}2&0&0\\0&2&0\\0&0&2\end{pmatrix},\quad n_V=\begin{pmatrix}0&0&3\\0&0&6\\1&1&0\end{pmatrix}=n_W \]
with the identifications. However, working out the monoidal products $\otimes_A$, we find the putative fusion ring is not isomorphic to any of rank three fusion categories. In particular, there is no fusion category with $m_2\otimes_Am_3\cong2m_3$. Thus, the second simple object should also have coefficient one.\label{so923dimNIMrep}} In view of anyon condensation, the second and the third simple objects can `split' into two and three, respectively. In order to match the Frobenius-Perron dimension, they should `split.' Thus, we search for six-dimensional NIM-reps. Indeed, we find a solution
\begin{align*}
    n_1=1_6=n_X,\quad n_Y=n_T&=n_U=\begin{pmatrix}0&1&1&0&0&0\\1&0&1&0&0&0\\1&1&0&0&0&0\\0&0&0&0&1&1\\0&0&0&1&0&1\\0&0&0&1&1&0\end{pmatrix},\\
    n_Z=\begin{pmatrix}2&0&0&0&0&0\\0&2&0&0&0&0\\0&0&2&0&0&0\\0&0&0&2&0&0\\0&0&0&0&2&0\\0&0&0&0&0&2\end{pmatrix},\quad&n_V=\begin{pmatrix}0&0&0&1&1&1\\0&0&0&1&1&1\\0&0&0&1&1&1\\1&1&1&0&0&0\\1&1&1&0&0&0\\1&1&1&0&0&0\end{pmatrix}=n_W.
\end{align*}
This gives identifications
\[ m_1\cong1\oplus X\oplus2Z,\quad m_2\cong Y\oplus T\oplus U\cong m_3,\quad m_4\cong m_5\cong m_6\cong V\oplus W. \]
In $\mc B_A$, they have quantum dimensions
\[ d_1=d_2=d_3=1,\quad d_4=d_5=d_6=\pm1. \]
The result tells us $\mc B_A$ should have rank six. Let us try to figure out $\mc B_A$.

As all simple objects have Frobenius-Perron dimensions one, it should be multiplicity-free. Thus, the fusion ring would be either $\text{FR}^{6,2}_1$ or $\text{FR}^{6,4}_1$.\footnote{We follow the notation of AnyonWiki \cite{anyonwiki}. A symbol $\text{FR}^{r,n}_i$ denotes a fusion ring with rank $r$ and $n$ non-self-dual simple objects. The subscript $i$ labels different fusion rings with the same $(r,n)$.} We can also work out monoidal product $\otimes_A$, but the free module functor cannot fix it uniquely. What we can fix is
\[ \hspace{-40pt}m_1\otimes_Am_j\cong m_j\cong m_j\otimes_Am_1\ (j=1,2,\dots,6),\quad m_2\otimes_Am_2\cong m_3,\quad m_2\otimes_Am_3\cong m_1\cong m_3\otimes_Am_2,\quad m_3\otimes_Am_3\cong m_2. \]
Namely, $\{m_1,m_2,m_3\}$ form the $\mbb Z/3\mbb Z$ fusion ring. We also learn
\begin{align*}
    (m_2\oplus m_3)\otimes_A(m_4\oplus m_5\oplus m_6)&\cong F_A(Y\otimes V)\\
    &\cong F_A(V\oplus W)\\
    &\cong2(m_4\oplus m_5\oplus m_6)\cong(m_4\oplus m_5\oplus m_6)\otimes_A(m_2\oplus m_3),\\
    (m_4\oplus m_5\oplus m_6)\otimes_A(m_4\oplus m_5\oplus m_6)&\cong F_A(V\otimes V)\\
    &\cong F_A(1\oplus Y\oplus Z\oplus T\oplus U)\\
    &\cong3(m_1\oplus m_2\oplus m_3).
\end{align*}
When $m_4\otimes_Am_4\cong m_5\otimes_Am_5\cong m_6\otimes_Am_6\cong m_1$, using consistency and associativity, we get two fusion rings
\begin{table}[H]
\begin{center}
\begin{tabular}{c|c|c|c|c|c|c}
    $\otimes_A$&$m_1$&$m_2$&$m_3$&$m_4$&$m_5$&$m_6$\\\hline
    $m_1$&$m_1$&$m_2$&$m_3$&$m_4$&$m_5$&$m_6$\\\hline
    $m_2$&$m_2$&$m_3$&$m_1$&$m_6$&$m_4$&$m_5$\\\hline
    $m_3$&$m_3$&$m_1$&$m_2$&$m_5$&$m_6$&$m_4$\\\hline
    $m_4$&$m_4$&$m_5$&$m_6$&$m_1$&$m_2$&$m_3$\\\hline
    $m_5$&$m_5$&$m_6$&$m_4$&$m_3$&$m_1$&$m_2$\\\hline
    $m_6$&$m_6$&$m_4$&$m_5$&$m_2$&$m_3$&$m_1$
\end{tabular},
\end{center}
\end{table}
\hspace{-17pt}or
\begin{table}[H]
\begin{center}
\begin{tabular}{c|c|c|c|c|c|c}
    $\otimes_A$&$m_1$&$m_2$&$m_3$&$m_4$&$m_5$&$m_6$\\\hline
    $m_1$&$m_1$&$m_2$&$m_3$&$m_4$&$m_5$&$m_6$\\\hline
    $m_2$&$m_2$&$m_3$&$m_1$&$m_5$&$m_6$&$m_4$\\\hline
    $m_3$&$m_3$&$m_1$&$m_2$&$m_6$&$m_4$&$m_5$\\\hline
    $m_4$&$m_4$&$m_6$&$m_5$&$m_1$&$m_3$&$m_2$\\\hline
    $m_5$&$m_5$&$m_4$&$m_6$&$m_2$&$m_1$&$m_3$\\\hline
    $m_6$&$m_6$&$m_5$&$m_4$&$m_3$&$m_2$&$m_1$
\end{tabular},
\end{center}
\end{table}
\hspace{-17pt}depending on, say, $m_4\otimes_Am_5$. Note that the non-commutativity appears automatically as a consequence of consistency and associativity. These are isomorphic to the fusion ring $\text{FR}^{6,2}_1$.\footnote{The identification with the notation in AnyonWiki \cite{anyonwiki} is given by
\[ m_1\cong1,\quad m_2\cong6,\quad m_3\cong5,\quad m_4\cong2,\quad m_5\cong3,\quad m_6\cong4, \]
or
\[ m_1\cong1,\quad m_2\cong5,\quad m_3\cong6,\quad m_4\cong2,\quad m_5\cong3,\quad m_6\cong4, \]
respectively.} On the other hand, when $m_4\otimes_Am_4\cong m_1,m_5\otimes_Am_5\cong m_2,m_6\otimes_Am_6\cong m_3$, consistency and associativity give a unique fusion ring
\begin{table}[H]
\begin{center}
\begin{tabular}{c|c|c|c|c|c|c}
    $\otimes_A$&$m_1$&$m_2$&$m_3$&$m_4$&$m_5$&$m_6$\\\hline
    $m_1$&$m_1$&$m_2$&$m_3$&$m_4$&$m_5$&$m_6$\\\hline
    $m_2$&&$m_3$&$m_1$&$m_6$&$m_4$&$m_5$\\\hline
    $m_3$&&&$m_2$&$m_5$&$m_6$&$m_4$\\\hline
    $m_4$&&&&$m_1$&$m_3$&$m_2$\\\hline
    $m_5$&&&&&$m_2$&$m_1$\\\hline
    $m_6$&&&&&&$m_3$
\end{tabular}.
\end{center}
\end{table}
\hspace{-17pt}Note that the commutativity also appears automatically as a consequence of consistency and associativity. This is isomorphic to the fusion ring $\text{FR}^{6,4}_1$.\footnote{The identification with the notation in AnyonWiki \cite{anyonwiki} is given by
\[ m_1\cong1,\quad m_2\cong5,\quad m_3\cong6,\quad m_4\cong2,\quad m_5\cong3,\quad m_6\cong4. \]} (Since $m_{4,5,6}$ enter symmetrically, the other monoidal products are given by permutations of simple objects.) Therefore, we find
\[ \mc B_A\simeq\mc C(\text{FR}^{6,2}_1)\text{ or }\mc C(\text{FR}^{6,4}_1). \]
(It is unsatisfactory that we cannot fix $\mc B_A$ uniquely. However, since the two have the same ranks, the ambiguity does not affect our physical applications.) Since this is semisimple, $A$ is separable, hence étale.

\paragraph{$A\cong1\oplus Z$.} It has $\fp_{\mc B}(A)=3$, and demands
\[ \fp(\mc B_A^0)=4,\quad\fp(\mc B_A)=12. \]
Matching of additive central charges, Frobenius-Perron dimensions, and the invariance of topological twists only allow
\[ \mc B_A^0\simeq\begin{cases}\vecG_{\mbb Z/2\mbb Z}^{-1}\boxtimes\vecG_{\mbb Z/2\mbb Z}^{-1}&(h_V,h_W)=(\frac14,\frac34),\\\tc&(h_V,h_W)=(0,\frac12).\end{cases}\quad(\mods1) \]
The category $\mc B_A$ of right $A$-modules has additional simple objects. A natural scenario has two more simple objects with Frobenius-Perron dimensions two so that $1+1+1+1+2^2+2^2=12$. Indeed, we find a six-dimensional NIM-rep
\begin{align*}
    n_1=1_6,\quad n_X=\begin{pmatrix}0&1&0&0&0&0\\1&0&0&0&0&0\\0&0&1&0&0&0\\0&0&0&1&0&0\\0&0&0&0&0&1\\0&0&0&0&1&0\end{pmatrix},\quad n_Y=n_T=n_U=\begin{pmatrix}0&0&1&0&0&0\\0&0&1&0&0&0\\1&1&1&0&0&0\\0&0&0&1&1&1\\0&0&0&1&0&0\\0&0&0&1&0&0\end{pmatrix},\\
    n_Z=\begin{pmatrix}1&1&0&0&0&0\\1&1&0&0&0&0\\0&0&2&0&0&0\\0&0&0&2&0&0\\0&0&0&0&1&1\\0&0&0&0&1&1\end{pmatrix},\quad n_V=\begin{pmatrix}0&0&0&1&1&0\\0&0&0&1&0&1\\0&0&0&2&1&1\\1&1&2&0&0&0\\1&0&1&0&0&0\\0&1&1&0&0&0\end{pmatrix},\quad n_W=\begin{pmatrix}0&0&0&1&0&1\\0&0&0&1&1&0\\0&0&0&2&1&1\\1&1&2&0&0&0\\0&1&1&0&0&0\\1&0&1&0&0&0\end{pmatrix}.
\end{align*}
The solution gives identifications
\[ m_1\cong1\oplus Z,\quad m_2\cong X\oplus Z,\quad m_3\cong Y\oplus T\oplus U,\quad m_4\cong V\oplus W,\quad m_5\cong V,\quad m_6\cong W. \]
In view of anyon condensation, $m_{1,2,5,6}$ are deconfined, and form $\mc B_A^0$. The other two simple objects $m_{3,4}$ are confined. In $\mc B_A$, they have quantum dimensions
\[ d_{\mc B_A}(m_1)=1=d_{\mc B_A}(m_2),\quad d_{\mc B_A}(m_3)=2,\quad d_{\mc B_A}(m_4)=\pm2,\quad d_{\mc B_A}(m_5)=\pm1=d_{\mc B_A}(m_6). \]
Working out the monoidal products $\otimes_A$, we find
\begin{table}[H]
\begin{center}
\begin{tabular}{c|c|c|c|c|c|c}
    $\otimes_A$&$m_1$&$m_2$&$m_3$&$m_4$&$m_5$&$m_6$\\\hline
    $m_1$&$m_1$&$m_2$&$m_3$&$m_4$&$m_5$&$m_6$\\\hline
    $m_2$&&$m_1$&$m_3$&$m_4$&$m_6$&$m_5$\\\hline
    $m_3$&&&$m_1\oplus m_2\oplus m_3$&$m_4\oplus m_5\oplus m_6$&$m_4$&$m_4$\\\hline
    $m_4$&&&&$m_1\oplus m_2\oplus m_3$&$m_3$&$m_3$\\\hline
    $m_5$&&&&&$m_1$&$m_2$\\\hline
    $m_6$&&&&&&$m_1$
\end{tabular}.
\end{center}
\end{table}
\hspace{-17pt}The fusion ring is isomorphic to $\text{FR}^{6,0}_2$.\footnote{The identification of simple objects with the notation in AnyonWiki is given by
\[ m_1\cong1,\quad m_2\cong4,\quad m_3\cong6,\quad m_4\cong5,\quad m_5\cong2,\quad m_6\cong3, \]
or its permutation $(23)$.} This in particular implies $1\oplus Z$ is separable, and étale.\newline

Now, we are left with $(n_X,n_Z,n_V)=(2,0,0),(5,0,0),(3,1,0),(2,0,1),(0,1,1)$. We study these in turn. It turns out only the last, $A\cong1\oplus Z\oplus V$, is connected étale.

\paragraph{$A\cong1\oplus2X$.} It has $\fp_{\mc B}(A)=3$, and demands
\[ \fp(\mc B_A^0)=4,\quad\fp(\mc B_A)=12. \]
Calculating $b_j\otimes A$, we find candidate simple objects of $\mc B_A$:
\[ 1\oplus2X,\quad2\oplus X,\quad3Y,\quad3Z,\quad3T,\quad3U,\quad V,\quad W. \]
(The putative $\mc B_A$ may also have $V\oplus W$, but our discussion below does not depend on the presence/absence.) In a putative $\mc B_A$, they have Frobenius-Perron dimensions
\[ 1,\quad1,\quad2,\quad2,\quad2,\quad2,\quad1,\quad1, \]
respectively. In order to match additive central charges and Frobenius-Perron dimensions, $\mc B_A^0$ should consist of $\{1\oplus2X,2\oplus X,V,W\}$. The four additional objects contribute $4\times2^2=16$ to $\fp(\mc B_A)$, and exceeds 12.\footnote{One may wonder why it is not allowed to pick just two candidate simple objects from those four with Frobenius-Perron dimensions two. This can match Frobenius-Perron dimensions, but one finds the actions $b_j\otimes-$ cannot be closed.} Thus, the candidate is ruled out.

\paragraph{$A\cong1\oplus5X$.} It has $\fp_{\mc B}(A)=6$, and demands
\[ \fp(\mc B_A^0)=1,\quad\fp(\mc B_A)=6. \]
Computing $b_j\otimes A$, we find candidate simple objects of $\mc B_A$  with the smallest Frobenius-Perron dimensions:
\[ 1\oplus5X,\quad5\oplus X,\quad3Y,\quad3Z,\quad3T,\quad3U,\quad V\oplus W,\quad\dots\ . \]
They have Frobenius-Perron dimensions
\[ 1,\quad1,\quad1,\quad1,\quad1,\quad1,\quad1,\quad\dots\ , \]
respectively. Again, their contributions to the Frobenius-Perron dimension exceed six, and the candidate is ruled out.

\paragraph{$A\cong1\oplus3X\oplus Z$.} It has $\fp_{\mc B}(A)=6$, and demands
\[ \fp(\mc B_A^0)=1,\quad\fp(\mc B_A)=6. \]
Calculating $b_j\otimes A$, we find candidate simple objects of $\mc B_A$:
\[ 1\oplus3X\oplus Z,\quad3\oplus X\oplus Z,\quad Y\oplus T\oplus U,\quad3Y,\quad3T,\quad3U,\quad\dots\ . \]
Note that there should exist candidate simple objects made of $V,W$. They have Frobenius-Perron dimensions
\[ 1,\quad1,\quad1,\quad1,\quad1,\quad1,\quad\dots\ , \]
respectively. Again, their contributions to Frobenius-Perron dimension exceed six, and the candidate is ruled out.

\paragraph{$A\cong1\oplus2X\oplus V$.} It has $\fp_{\mc B}(A)=6$, and demands
\[ \fp(\mc B_A^0)=1,\quad\fp(\mc B_A)=6. \]
With the free module functor $F_A(b_j)=b_j\otimes A$, we find candidate simple objects of $\mc B_A$:
\[ 1\oplus2X\oplus V,\quad2\oplus X\oplus W,\quad3Y,\quad3Z,\quad3T,\quad3U,\quad V\oplus W,\quad\dots\ . \]
There should exist more simple object, but they do not affect the following discussion. They have Frobenius-Perron dimensions
\[ 1,\quad1,\quad1,\quad1,\quad1,\quad1,\quad1,\quad\dots\ , \]
respectively. Their contributions to $\fp(\mc B_A)$ exceed six, and the candidate is ruled out.

\paragraph{$A\cong1\oplus Z\oplus V$.} It has $\fp_{\mc B}(A)=6$, and demands
\[ \fp(\mc B_A^0)=1,\quad\fp(\mc B_A)=6. \]
Computing $b_j\otimes A$, we find candidate simple objects of $\mc B_A$:
\[ 1\oplus Z\oplus V,\quad X\oplus Z\oplus W,\quad Y\oplus T\oplus U\oplus V\oplus W. \]
They have Frobenius-Perron dimensions
\[ 1,\quad1,\quad2, \]
respectively. Their contributions to $\fp(\mc B_A)$ matches, $1^2+1^2+2^2=6$. This suggests $\mc B_A$ has rank three. Indeed, we find a three-dimensional NIM-rep
\begin{align*}
    &n_1=1_3,\quad n_X=\begin{pmatrix}0&1&0\\1&0&0\\0&0&1\end{pmatrix},\quad n_Y=n_T=n_U=\begin{pmatrix}0&0&1\\0&0&1\\1&1&1\end{pmatrix},\\
    &n_Z=\begin{pmatrix}1&1&0\\1&1&0\\0&0&2\end{pmatrix},\quad n_V=\begin{pmatrix}1&0&1\\0&1&1\\1&1&2\end{pmatrix},\quad n_W=\begin{pmatrix}0&1&1\\1&0&1\\1&1&2\end{pmatrix}.
\end{align*}
The solution gives identifications
\[ m_1\cong1\oplus Z\oplus V,\quad m_2\cong X\oplus Z\oplus W,\quad m_3\cong Y\oplus T\oplus U\oplus V\oplus W. \]
Computing their quantum dimensions (\ref{dBAm}), we find the candidate can be separable only when $\mc B$ is unitary. Then, the right $A$-modules have quantum dimensions
\[ d_{\mc B_A}(m_1)=1=d_{\mc B_A}(m_2),\quad d_{\mc B_A}(m_3)=2. \]
Working out the monoidal products $\otimes_A$, we find
\begin{table}[H]
\begin{center}
\begin{tabular}{c|c|c|c}
    $\otimes_A$&$m_1$&$m_2$&$m_3$\\\hline
    $m_1$&$m_1$&$m_2$&$m_3$\\\hline
    $m_2$&&$m_1$&$m_3$\\\hline
    $m_3$&&&$m_1\oplus m_2\oplus m_3$
\end{tabular}.
\end{center}
\end{table}
\hspace{-17pt}We can identify $\mc B_A\simeq\text{Rep}(S_3)$. Since $\text{Rep}(S_3)$ is semisimple, $A$ is separable.

Here, a careful reader would notice that we have not shown commutativity of the candidate. For $V$ to give commutative algebras, we need to limit our MFCs to those with the first and third conformal dimensions. Then, since these are conformal dimensions of $\mc C(B_4,2)$ or its conjugate, the commutativity can in principle be proven using the $F$- and $R$-symbols obtained in \cite{AFT16}, and computing multiplication morphism $\mu$. However, without these computations, we can conclude this should be connected étale from the lemma 2. We have already seen $1\oplus Z$ gives $so(9)_2\to\tc$ (for the first and third conformal dimensions). The latter MFC is known \cite{KK23MFC} to have two (for unitary $\tc$) or one (for non-unitary $\tc$) nontrivial connected étale algebra(s). This implies, for unitary $\tc$, there are two inequivalent operations
\[ \tc\to\vect. \]
Therefore, for $\mc B\simeq so(9)_2$ which admits unitary $\mc B_A^0\simeq\tc$, there should exist two inequivalent operations
\[ so(9)_2\to\vect. \]
We have already found one nontrivial connected étale algebra, $A\cong1\oplus X\oplus2Z$, giving one such operation.\footnote{Note that
\[ (1\oplus X)\otimes(1\oplus Z)\cong1\oplus X\oplus2Z. \]} From our considerations so far, the only candidate which could give another operation is $A\cong1\oplus Z\oplus V$. Thus, we learn the candidate is connected étale. Indeed, this matches known results in $\mc C(B_4,2)$.\footnote{We thank Terry Gannon for teaching this fact to us.} On the other hand, when $A\cong1\oplus Z$ gives non-unitary $\mc B_A^0\simeq\tc$, the Toric Code MFC admits \cite{KK23MFC} only one nontrivial connected étale algebra. Thus, $A\cong1\oplus Z\oplus V$ cannot be connected étale. This conclusion is consistent with our result above that the candidate is not separable in non-unitary $so(9)_2$. We found a connected étale algebra
\begin{equation}
    A\cong1\oplus Z\oplus V\quad(d_V,d_W,h_V,h_W)=(3,3,0,\frac12).\quad(\mods1\text{ for }h)\label{so921ZVetale}
\end{equation}

We conclude
\begin{table}[H]
\begin{center}
\begin{tabular}{c|c|c|c}
    Connected étale algebra $A$&$\mc B_A$&$\rank(\mc B_A)$&Lagrangian?\\\hline
    $1$&$\mc B$&$8$&No\\
    $1\oplus X$&$\text{TY}(\mbb Z/9\mbb Z)$&10&No\\
    $1\oplus X\oplus2Z$&$\mc C(\text{FR}^{6,2}_1)\text{ or }\mc C(\text{FR}^{6,4}_1)$&6&Yes\\
    $1\oplus Z$&$\mc C(\text{FR}^{6,0}_2)$&6&No\\
    $1\oplus Z\oplus V$ for (\ref{so921ZVetale})&$\text{Rep}(S_3)$&3&Yes
\end{tabular}.
\end{center}
\caption{Connected étale algebras in rank eight MFC $\mcal B\simeq so(9)_2$}\label{rank8so92results}
\end{table}
\hspace{-17pt}All the 16 MFCs $\mc B\simeq so(9)_2$'s fail to be completely anisotropic.

\subsubsection{$\mc B\simeq\text{Rep}(D(D_3))$}
The MFCs have eight simple objects $\{1,X,Y,Z,T,U,V,W\}$ obeying monoidal products
\begin{table}[H]
\begin{center}
\makebox[1 \textwidth][c]{       
\resizebox{1.2 \textwidth}{!}{\begin{tabular}{c|c|c|c|c|c|c|c|c}
    $\otimes$&$1$&$X$&$Y$&$Z$&$T$&$U$&$V$&$W$\\\hline
    $1$&$1$&$X$&$Y$&$Z$&$T$&$U$&$V$&$W$\\\hline
    $X$&&$1$&$Y$&$Z$&$T$&$U$&$W$&$V$\\\hline
    $Y$&&&$1\oplus X\oplus Y$&$T\oplus U$&$Z\oplus U$&$Z\oplus T$&$V\oplus W$&$V\oplus W$\\\hline
    $Z$&&&&$1\oplus X\oplus Z$&$Y\oplus U$&$Y\oplus T$&$V\oplus W$&$V\oplus W$\\\hline
    $T$&&&&&$1\oplus X\oplus T$&$Y\oplus Z$&$V\oplus W$&$V\oplus W$\\\hline
    $U$&&&&&&$1\oplus X\oplus U$&$V\oplus W$&$V\oplus W$\\\hline
    $V$&&&&&&&$1\oplus Y\oplus Z\oplus T\oplus U$&$X\oplus Y\oplus Z\oplus T\oplus U$\\\hline
    $W$&&&&&&&&$1\oplus Y\oplus Z\oplus T\oplus U$
\end{tabular}.}}
\end{center}
\end{table}
\hspace{-17pt}Thus, they have
\begin{align*}
    \fp_{\mc B}(1)=1=\fp_{\mc B}(X),\quad&\fp_{\mc B}(Y)=\fp_{\mc B}(Z)=\fp_{\mc B}(T)=\fp_{\mc B}(U)=2,\\
    &\fp_{\mc B}(V)=3=\fp_{\mc B}(W),
\end{align*}
and
\[ \fp(\mc B)=36. \]
Their quantum dimensions $d_j$'s are solutions of the same multiplication rules $d_id_j=\sum_{k=1}^8{N_{ij}}^kd_k$. There are two (nonzero) solutions
\[ (d_X,d_Y,d_Z,d_T,d_U,d_V,d_W)=(1,2,2,2,2,-3,-3),(1,2,2,2,2,3,3) \]
with the same categorical dimension
\[ D^2(\mc B)=36. \]
They both have the same four conformal dimensions\footnote{Naively, one finds 72 consistent conformal dimensions, but the others give the same MFCs with the one in the main text under permutations of simple objects $(YZ),(YT),(YU),(ZT),(ZU),(TU),(VW)$.}
\begin{align*}
    (h_X,h_Y,h_Z,h_T,h_U,h_V,h_W)&=(0,0,0,\frac13,\frac23,0,\frac12),(0,0,0,\frac13,\frac23,\frac14,\frac34),\\
    &~~~~(0,\frac13,\frac13,\frac23,\frac23,0,\frac12),(0,\frac13,\frac13,\frac23,\frac23,\frac14,\frac34)\quad(\mods1).
\end{align*}
They have $S$-matrices
\[ \widetilde S=\begin{pmatrix}1&1&2&2&2&2&d_V&d_W\\1&1&2&2&2&2&-d_V&-d_W\\2&2&s&s'&-2&-2&0&0\\2&2&s'&s&-2&-2&0&0\\2&2&-2&-2&-2&4&0&0\\2&2&-2&-2&4&-2&0&0\\d_V&-d_V&0&0&0&0&\pm3&\mp3\\d_W&-d_W&0&0&0&0&\mp3&\pm3\end{pmatrix}, \]
with
\[ (s,s')=\begin{cases}(4,-2)&(\text{1st and 2nd }h),\\(-2,4)&(\text{3rd and 4th }h).\end{cases} \]
They have additive central charges
\[ c(\mc B)=\begin{cases}0&(\text{1st\&2nd }h),\\4&(\text{3rd\&4th }h).\end{cases}\quad(\mods8). \]
There are
\[ 2(\text{quantum dimensions})\times4(\text{conformal dimensions})\times2(\text{categorical dimensions})=16 \]
MFCs, among which those eight with the second quantum dimensions are unitary. We classify connected étale algebras in all 16 MFCs simultaneously.

An ansatz
\[ A\cong1\oplus n_XX\oplus n_YY\oplus n_ZZ\oplus n_TT\oplus n_UU\oplus n_VV\oplus n_WW \]
with $n_j\in\mbb N$ has
\[ \fp_{\mc B}(A)=1+n_X+2(n_Y+n_Z+n_T+n_U)+3(n_V+n_W). \]
For this to obey (\ref{FPdimA2bound}), the natural numbers can take only 56 values just as in the previous example. The sets contain the one with all $n_j$'s be zero. It is the trivial connected étale algebra $A\cong1$ giving $\mc B_A^0\simeq\mc B_A\simeq\mc B$. Other candidates with $T,U,W$ have nontrivial conformal dimensions, and they fail to be commutative. Thus, we are left with candidates with just $X,Y,Z,V$. Note that $\{1,X,Y\},\{1,X,Z\}$ form $\text{Rep}(S_3)\subset\text{Rep}(D(D_3))$. Apart from the trivial one, there are 24 candidates
\begin{align*}
    (n_X,n_Y,n_Z,n_V)&=(1,0,0,0),(0,1,0,0),(0,0,1,0),(0,0,0,1),\\
    &~~~~(2,0,0,0),(1,1,0,0),(1,0,1,0),(1,0,0,1),\\
    &~~~~(0,2,0,0),(0,1,1,0),(0,1,0,1),(0,0,2,0),\\
    &~~~~(0,0,1,1),(3,0,0,0),(2,1,0,0),(2,0,1,0),\\
    &~~~~(2,0,0,1),(1,2,0,0),(1,1,1,0),(1,0,2,0),\\
    &~~~~(4,0,0,0),(3,1,0,0),(3,0,1,0),(5,0,0,0).
\end{align*}
Some of these are ruled out by studying Frobenius-Perron dimensions. The four candidates
\[ (n_X,n_Y,n_Z,n_V)=(0,0,0,1),(1,1,0,0),(1,0,1,0),(3,0,0,0) \]
have $\fp_{\mc B}=4$ and demand $\fp(\mc B_A^0)=\frac94$. However, there is no MFC with such Frobenius-Perron dimension. Thus, the candidates are ruled out. Similarly, seven candidates
\[ (n_X,n_Y,n_Z,n_V)=(1,0,0,1),(0,1,1,0),(0,2,0,0),(0,0,2,0),(2,1,0,0),(2,0,1,0),(4,0,0,0) \]
are ruled out because there is no MFC $\mc B_A^0$ with Frobenius-Perron dimension $\fp(\mc B_A^0)=\frac{36}{25}$. Thus, we are left with 13 nontrivial candidates
\begin{align*}
    (n_X,n_Y,n_Z,n_V)&=(1,0,0,0),(0,1,0,0),(0,0,1,0),(2,0,0,0),(0,1,0,1),(0,0,1,1),(2,0,0,1),\\
    &~~~~(1,2,0,0),(1,1,1,0),(1,0,2,0),(3,1,0,0),(3,0,1,0),(5,0,0,0).
\end{align*}
Thanks to the lemma 1, we know five of them $(n_X,n_Y,n_Z,n_V)=(1,0,0,0),(0,1,0,0)$, $(0,0,1,0),(1,2,0,0),(1,0,2,0)$ or $A\cong1\oplus X,1\oplus Y,1\oplus Z,1\oplus X\oplus2Y,1\oplus X\oplus2Z$ are commutative algebras (for certain conformal dimensions) because they are so in (symmetric) braided fusion subcategory $\text{Rep}(S_3)\subset\mc B$ \cite{KK23preMFC}. Let us see whether they are separable or not by identifying $\mc B_A$. As a result, we find all of them are separable, hence étale (for certain conformal dimensions).

\paragraph{$A\cong1\oplus X$.} The simple object $X$ has $(d_X,h_X)=(1,0)$ (mod 1 for $h_X$), and it has trivial braiding $c_{X,X}\cong id_1$ \cite{KK23preMFC}. Thus it is a commutative algebra for all four conformal dimensions. Furthermore, it turns out separable, hence étale. Let us check this point by identifying $\mc B_A$.

Since it has $\fp_{\mc B}(A)=2$, it demands
\[ \fp(\mc B_A^0)=9,\quad\fp(\mc B_A)=18. \]
Employing the matching of central charges (\ref{propBA0}) and the invariance of topological twists (\ref{invtopologicaltwist}), we identify
\begin{equation}
    \mc B_A^0\simeq\vecG_{\mbb Z/3\mbb Z}^1\boxtimes\vecG_{\mbb Z/3\mbb Z}^1.\label{repDD31XBA0}
\end{equation}
(See section \ref{Z3Z3}.) What is the category $\mc B_A$ of right $A$-modules? We find
\begin{equation}
    \mc B_A\simeq\text{TY}(\mbb Z/3\mbb Z\times\mbb Z/3\mbb Z).\label{repDD31XBA}
\end{equation}
One of the easiest ways to find this fact is to perform anyon condensation. It `identifies' $1,X$, and hence $V,W$. The resulting $1\oplus X$ and $V\oplus W$ have Frobenius-Perron dimensions one, three in $\mc B_A$. The other simple objects $Y,Z,T,U$ split into two each. They all have Frobenius-Perron dimensions one in $\mc B_A$. In total, we have nine invertible simple objects forming $\mbb Z/3\times\mbb Z/3$ MFC (\ref{repDD31XBA0}), and one additional non-invertible simple object with $\fp_{\mc B_A}=3$. The fusion category is identified with a $\mbb Z/3\mbb Z\times\mbb Z/3\mbb Z$ Tambara-Yamagami category (\ref{repDD31XBA}).\footnote{More rigorously, we should find NIM-reps. Indeed, we find four 10-dimensional solutions. One of them is given by
\begin{align*} n_1=1_{10}=n_X,\quad n_Y=&\begin{pmatrix}0&1&1&0&0&0&0&0&0&0\\1&0&1&0&0&0&0&0&0&0\\1&1&0&0&0&0&0&0&0&0\\0&0&0&0&0&0&1&0&1&0\\0&0&0&0&0&1&0&1&0&0\\0&0&0&0&1&0&0&1&0&0\\0&0&0&1&0&0&0&0&1&0\\0&0&0&0&1&1&0&0&0&0\\0&0&0&1&0&0&1&0&0&0\\0&0&0&0&0&0&0&0&0&2\end{pmatrix},\quad n_Z=\begin{pmatrix}0&0&0&1&1&0&0&0&0&0\\0&0&0&0&0&0&1&1&0&0\\0&0&0&0&0&1&0&0&1&0\\1&0&0&0&1&0&0&0&0&0\\1&0&0&1&0&0&0&0&0&0\\0&0&1&0&0&0&0&0&1&0\\0&1&0&0&0&0&0&1&0&0\\0&1&0&0&0&0&1&0&0&0\\0&0&1&0&0&1&0&0&0&0\\0&0&0&0&0&0&0&0&0&2\end{pmatrix},\\
n_T=&\begin{pmatrix}0&0&0&0&0&1&1&0&0&0\\0&0&0&0&1&0&0&0&1&0\\0&0&0&1&0&0&0&1&0&0\\0&0&1&0&0&0&0&1&0&0\\0&1&0&0&0&0&0&0&1&0\\1&0&0&0&0&0&1&0&0&0\\1&0&0&0&0&1&0&0&0&0\\0&0&1&1&0&0&0&0&0&0\\0&1&0&0&1&0&0&0&0&0\\0&0&0&0&0&0&0&0&0&2\end{pmatrix},\quad n_U=\begin{pmatrix}0&0&0&0&0&0&0&1&1&0\\0&0&0&1&0&1&0&0&0&0\\0&0&0&0&1&0&1&0&0&0\\0&1&0&0&0&1&0&0&0&0\\0&0&1&0&0&0&1&0&0&0\\0&1&0&1&0&0&0&0&0&0\\0&0&1&0&1&0&0&0&0&0\\1&0&0&0&0&0&0&0&1&0\\1&0&0&0&0&0&0&1&0&0\\0&0&0&0&0&0&0&0&0&2\end{pmatrix},\\n_V=&\begin{pmatrix}0&0&0&0&0&0&0&0&0&1\\0&0&0&0&0&0&0&0&0&1\\0&0&0&0&0&0&0&0&0&1\\0&0&0&0&0&0&0&0&0&1\\0&0&0&0&0&0&0&0&0&1\\0&0&0&0&0&0&0&0&0&1\\0&0&0&0&0&0&0&0&0&1\\0&0&0&0&0&0&0&0&0&1\\0&0&0&0&0&0&0&0&0&1\\1&1&1&1&1&1&1&1&1&0\end{pmatrix}=n_W.
\end{align*}
The NIM-rep gives identifications
\[ m_1\cong1\oplus X,\quad m_2\cong Y\cong m_3,\quad m_4\cong Z\cong m_5,\quad m_6\cong T\cong m_7,\quad m_8\cong U\cong m_9,\quad m_{10}\cong V\oplus W. \]
(The other three solutions give the same identifications.) They have
\[ d_{\mc B_A}(m_j)=1\ (j=1,2,\dots,9),\quad d_{\mc B_A}(m_{10})=\pm3. \]
Working out the monoidal products $\otimes_A$, we find
\[ m_j\otimes_Am_{10}\cong m_{10}\cong m_{10}\otimes_Am_j\ (j=1,2,\dots,9),\quad m_{10}\otimes_Am_{10}\cong\bigoplus_{j=1}^9m_j, \]
showing (\ref{repDD31XBA}).} Since $\text{TY}(\mbb Z/3\mbb Z\times\mbb Z/3\mbb Z)$ is a fusion category, it is semisimple. Therefore, $A$ is separable, and étale.

\paragraph{$A\cong1\oplus Y$.} For the first and second conformal dimensions, $1\oplus Y$ gives commutative algebra. What is left is to check its separability. It has $\fp_{\mc B}(A)=3$, and demands
\[ \fp(\mc B_A^0)=4,\quad\fp(\mc B_A)=12. \]
Computing $b_j\otimes A$, we find candidate simple objects of $\mc B_A$:
\[ 1\oplus Y,\quad X\oplus Y,\quad Z\oplus T\oplus U,\quad V\oplus W,\quad V,\quad W \]
with some multiplicities. In $\mc B_A$, they have Frobenius-Perron dimensions
\[ 1,\quad1,\quad2,\quad2,\quad1,\quad1, \]
respectively. The invariance of topological twists, matching of central charges, and Frobenius-Perron dimensions demand
\[ \mc B_A^0=\{1\oplus Y,X\oplus Y,V,W\}\simeq\begin{cases}\tc&(h_V,h_W)=(0,\frac12),\\\vecG_{\mbb Z/2\mbb Z}^{-1}\boxtimes\vecG_{\mbb Z/2\mbb Z}^{-1}&(h_V,h_W)=(\frac14,\frac34).\end{cases}\quad(\mods1) \]
The two additional simple objects give correct Frobenius-Perron dimensions for them to match $\fp(\mc B_A)=12$. Indeed, we find four six-dimensional NIM-reps. One of them is given by
\begin{align*}
    n_1=1_6,\quad n_X=\begin{pmatrix}0&1&0&0&0&0\\1&0&0&0&0&0\\0&0&1&0&0&0\\0&0&0&1&0&0\\0&0&0&0&0&1\\0&0&0&0&1&0\end{pmatrix}&,\quad n_Y=\begin{pmatrix}1&1&0&0&0&0\\1&1&0&0&0&0\\0&0&2&0&0&0\\0&0&0&2&0&0\\0&0&0&0&1&1\\0&0&0&0&1&1\end{pmatrix},\\
    n_Z=n_T=n_U=\begin{pmatrix}0&0&1&0&0&0\\0&0&1&0&0&0\\1&1&1&0&0&0\\0&0&0&1&1&1\\0&0&0&1&0&0\\0&0&0&1&0&0\end{pmatrix},\quad n_V=&\begin{pmatrix}0&0&0&1&1&0\\0&0&0&1&0&1\\0&0&0&2&1&1\\1&1&2&0&0&0\\1&0&1&0&0&0\\0&1&1&0&0&0\end{pmatrix},\quad n_W=\begin{pmatrix}0&0&0&1&0&1\\0&0&0&1&1&0\\0&0&0&2&1&1\\1&1&2&0&0&0\\0&1&1&0&0&0\\1&0&1&0&0&0\end{pmatrix}.
\end{align*}
The solution gives the identifications
\[ m_1\cong1\oplus Y,\quad m_2\cong X\oplus Y,\quad m_3\cong Z\oplus T\oplus U,\quad m_4\cong V\oplus W,\quad m_5\cong V,\quad m_6\cong W. \]
(The other three solutions give the same identifications.) In $\mc B_A$, they have quantum dimensions
\[ d_{\mc B_A}(m_1)=d_{\mc B_A}(m_2)=d_{\mc B_A}(m_3)=1,\quad d_{\mc B_A}(m_4)=\pm2,\quad d_{\mc B_A}(m_5)=\pm1=d_{\mc B_A}(m_6), \]
where signs are correlated. Working out the monoidal product $\otimes_A$, we find
\begin{table}[H]
\begin{center}
\begin{tabular}{c|c|c|c|c|c|c}
    $\otimes_A$&$m_1$&$m_2$&$m_3$&$m_4$&$m_5$&$m_6$\\\hline
    $m_1$&$m_1$&$m_2$&$m_3$&$m_4$&$m_5$&$m_6$\\\hline
    $m_2$&&$m_1$&$m_3$&$m_4$&$m_6$&$m_5$\\\hline
    $m_3$&&&$m_1\oplus m_2\oplus m_3$&$m_4\oplus m_5\oplus m_6$&$m_4$&$m_4$\\\hline
    $m_4$&&&&$m_1\oplus m_2\oplus m_3$&$m_3$&$m_3$\\\hline
    $m_5$&&&&&$m_1$&$m_2$\\\hline
    $m_6$&&&&&&$m_1$
\end{tabular}.
\end{center}
\end{table}
\hspace{-17pt}Just as in the previous MFC, the category $\mc B_A$ of right $A$-modules is isomorphic to a fusion category $\mc C(\text{FR}^{6,0}_2)$. Thus, $1\oplus Y$ is separable, hence étale.

\paragraph{$A\cong1\oplus Z$.} Since the ambient category $\mc B\simeq\text{Rep}(D(D_3))$ is invariant under the change of names $Y\leftrightarrow Z$, our previous analysis implies $A\cong1\oplus Z$ is a connected étale algebra for the first and second conformal dimensions. It gives $\mc B_A\simeq\mc C(\text{FR}^{6,0}_2)$.

\paragraph{$A\cong1\oplus X\oplus2Y$.} For the first and second conformal dimensions, the algebra is commutative. Furthermore, it turns out separable, hence étale. So as to check this point, we identify $\mc B_A$.

Calculating $b_j\otimes A$, we find candidate simple objects:
\[ 1\oplus X\oplus2Y,\quad Z\oplus T\oplus U,\quad V\oplus W, \]
with Frobenius-Perron dimensions
\[ 1,\quad1,\quad1, \]
respectively. Logically, the latter two can have coefficients, but the possibilities are ruled out.\footnote{The proof is the same as in the footnote \ref{so923dimNIMrep}.} For them to saturate $\fp(\mc B_A)=6$, the second and third should `split' into two and three, respectively. Therefore, we search for six-dimensional NIM-reps. Indeed, we find four solutions. One of them is given by
\begin{align*}
    n_1=1_6=n_X,\quad n_Y=\begin{pmatrix}2&0&0&0&0&0\\0&2&0&0&0&0\\0&0&2&0&0&0\\0&0&0&0&1&1\\0&0&0&1&0&1\\0&0&0&1&1&0\end{pmatrix},\quad n_Z=\begin{pmatrix}0&1&1&0&0&0\\1&0&1&0&0&0\\1&1&0&0&0&0\\0&0&0&0&1&1\\0&0&0&1&0&1\\0&0&0&1&1&0\end{pmatrix}=n_T,\\
    n_U=\begin{pmatrix}0&1&1&0&0&0\\1&0&1&0&0&0\\1&1&0&0&0&0\\0&0&0&2&0&0\\0&0&0&0&2&0\\0&0&0&0&0&2\end{pmatrix},\quad n_V=\begin{pmatrix}0&0&0&1&1&1\\0&0&0&1&1&1\\0&0&0&1&1&1\\1&1&1&0&0&0\\1&1&1&0&0&0\\1&1&1&0&0&0\end{pmatrix}=n_W.
\end{align*}
The NIM-rep gives identifications
\[ m_1\cong1\oplus X\oplus2Y,\quad m_2\cong Z\oplus T\oplus U\cong m_3,\quad m_4\cong m_5\cong m_6\cong V\oplus W. \]
(The other three solutions give the same identifications.) They have
\[ d_{\mc B_A}(m_1)=d_{\mc B_A}(m_2)=d_{\mc B_A}(m_3)=1,\quad d_{\mc B_A}(m_4)=d_{\mc B_A}(m_5)=d_{\mc B_A}(m_6)=\pm1. \]
The result tells us $\mc B_A$ should have rank six. Since $F_A(b)\otimes_AF_A(b')$ are the same as in $1\oplus X\oplus2Z\in so(9)_2$, we immediately learn $\mc B_A\simeq\mc C(\text{FR}^{6,2}_1)\text{ or }\mc C(\text{FR}^{6,4}_1)$. Since these are semisimple, $A$ is separable, hence étale.

\paragraph{$A\cong1\oplus X\oplus2Z$.} The symmetry $Y\leftrightarrow Z$ and the previous analysis imply $A\cong1\oplus X\oplus2Z$ is a connected étale algebra for the first and second conformal dimensions. It gives $\mc B_A\simeq\mc C(\text{FR}^{6,2}_1)\text{ or }\mc C(\text{FR}^{6,4}_1)$.\newline

Now, we are left with eight candidates
\[ (n_X,n_Y,n_Z,n_V)=(2,0,0,0),(0,1,0,1),(0,0,1,1),(2,0,0,1),(1,1,1,0),(3,1,0,0),(3,0,1,0),(5,0,0,0) \]
or $A\cong1\oplus2X,1\oplus Y\oplus V,1\oplus Z\oplus V,1\oplus2X\oplus V,1\oplus X\oplus Y\oplus Z,1\oplus3X\oplus Y,1\oplus3X\oplus Z,1\oplus5X$. We study these in turn. It turns out that all but two fail to be étale.

\paragraph{$A\cong1\oplus2X$.} It has $\fp_{\mc B}(A)=3$, and demands
\[ \fp(\mc B_A^0)=4,\quad\fp(\mc B_A)=12. \]
Calculating $b_j\otimes A$, we find candidate simple objects of $\mc B_A$:
\[ 1\oplus2X,\quad2\oplus X,\quad 3Y,\quad3Z,\quad3T,\quad3U,\quad V,\quad W,\quad\dots\ . \]
They have Frobenius-Perron dimensions
\[ 1,\quad1,\quad2,\quad2,\quad2,\quad2,\quad1,\quad1,\quad\dots\ , \]
respectively. In order to match Frobenius-Perron dimensions, $\mc B_A^0$ should consist of $\{1\oplus2X,2\oplus X,V,W\}$. The four additional candidate simple objects contribute $4\times2^2=16$ to Frobenius-Perron dimension, and exceeds 12. Thus, the candidate $1\oplus2X$ is ruled out.

\paragraph{$A\cong1\oplus Y\oplus V$.} For the first and third conformal dimensions, $V$ has $h_V=0$ (mod 1), and the candidate can be commutative. Taking $Y$ into account, we learn this candidate can be commutative only for the first conformal dimension. It has $\fp_{\mc B}(A)=6$, and demands
\[ \fp(\mc B_A^0)=1,\quad\fp(\mc B_A)=6. \]
Computing $b_j\otimes A$, we find candidate simple objects:
\[ 1\oplus Y\oplus V,\quad X\oplus Y\oplus W,\quad Z\oplus T\oplus U\oplus V\oplus W. \]
They have Frobenius-Perron dimensions
\[ 1,\quad1,\quad2, \]
respectively. Their contributions to $\fp(\mc B_A)$ match, $1^2+1^2+2^2=6$. The only possibility for the category of dyslectic (right) $A$-modules is
\[ \mc B_A^0\simeq\vect. \]
This identification also matches central charges. For $\mc B_A$, the analysis above suggests it has rank three. Indeed, we find a three-dimensional NIM-rep
\begin{align*}
    n_1=1_3,\quad n_X=&\begin{pmatrix}0&1&0\\1&0&0\\0&0&1\end{pmatrix},\quad n_Y=\begin{pmatrix}1&1&0\\1&1&0\\0&0&2\end{pmatrix},\\
    n_Z=n_T=n_U=\begin{pmatrix}0&0&1\\0&0&1\\1&1&1\end{pmatrix},&\quad n_V=\begin{pmatrix}1&0&1\\0&1&1\\1&1&2\end{pmatrix},\quad n_W=\begin{pmatrix}0&1&1\\1&0&1\\1&1&2\end{pmatrix}.
\end{align*}
The solution gives identifications
\[ m_1\cong1\oplus Y\oplus V,\quad m_2\cong X\oplus Y\oplus W,\quad m_3\cong Z\oplus T\oplus U\oplus V\oplus W. \]
By computing quantum dimensions, we find the candidate can be separable only when $d_V=3=d_W$. Then, the right $A$-modules have
\[ d_{\mc B_A}(m_1)=1=d_{\mc B_A}(m_2),\quad d_{\mc B_A}(m_3)=2. \]
Working out the monoidal products, we find
\begin{table}[H]
\begin{center}
\begin{tabular}{c|c|c|c}
    $\otimes_A$&$m_1$&$m_2$&$m_3$\\\hline
    $m_1$&$m_1$&$m_2$&$m_3$\\\hline
    $m_2$&&$m_1$&$m_3$\\\hline
    $m_3$&&&$m_1\oplus m_2\oplus m_3$
\end{tabular}.
\end{center}
\end{table}
\hspace{-17pt}We can identify
\[ \mc B_A\simeq\text{Rep}(S_3). \]
Since $\text{Rep}(S_3)$ is semisimple, $A$ is separable.

We still need to show commutativity. However, since $F$- and $R$-symbols of the MFCs are unknown (to the best of our knowledge), we take the indirect path through the lemma 2. We have already found $A\cong1\oplus Y,1\oplus Z$ are connected étale, giving $\mc B_A^0\simeq\vecG_{\mbb Z/2\mbb Z\times\mbb Z/2\mbb Z}^\alpha$. The MFCs further admit connected étale algebras $A'\in\mc B_A^0$ giving $(\mc B_A^0)_{A'}^0\simeq\vect$. More precisely, for the first conformal dimensions, we get
\[ \mc B_A^0\simeq\tc, \]
and it admits one (when $\tc$ is non-unitary) or two (when $\tc$ is unitary) connected étale algebra(s). (Taking both $A\cong1\oplus Y,1\oplus Z$ into account, we have two or four operations $\text{Rep}(D(D_3))\to\vect$, respectively.) One such operation $\text{Rep}(D(D_3))\to\vect$ is given by $A\cong1\oplus X\oplus2Y$ or $A\cong1\oplus X\oplus2Z$ we found above. These exist regardless of quantum dimensions, and we can identify them as connected étale algebras giving the composition
\[ \text{Rep}(D(D_3))\to\text{ non-unitary }\tc\to\vect. \]
When $A\cong1\oplus Y,1\oplus Z$ give unitary $\tc$, i.e., $d_V=3=d_W$, we need two more connected étale algebras. The only possibility is this and its symmetric partner. Therefore, we conclude the candidate $A\cong1\oplus Z\oplus V$ is connected étale for the first conformal dimension. It gives
\[ \text{Rep}(D(D_3))\to\vect. \]
We found
\begin{equation}
    A\cong1\oplus Y\oplus V\quad(d_V=3=d_W\&\text{1st }h).\label{RepDD31YVetale}
\end{equation}

\paragraph{$A\cong1\oplus Z\oplus V$.} The symmetry $Y\leftrightarrow Z$ and the previous analysis imply this candidate is connected étale
\begin{equation}
    A\cong1\oplus Z\oplus V\quad(d_V=3=d_W\&\text{1st }h),\label{RepDD31ZVetale}
\end{equation}
giving
\[ \mc B_A^0\simeq\vect,\quad\mc B_A\simeq\text{Rep}(S_3). \]

\paragraph{$A\cong1\oplus2X\oplus V$.} For the first and third conformal dimensions, this can be commutative. It has $\fp_{\mc B}(A)=6$, and demands
\[ \fp(\mc B_A^0)=1,\quad\fp(\mc B_A)=6. \]
With the free module functor $F_A(b_j)=b_j\otimes A$, we find candidate simple objects
\[ 1\oplus2X\oplus V,\quad2\oplus X\oplus W,\quad3Y,\quad3Z,\quad3T,\quad3U,\quad V\oplus W,\quad\dots\ , \]
with additional simple objects which do not affect our discussion below. In $\mc B_A$, they have Frobenius-Perron dimensions
\[ 1,\quad1,\quad1,\quad1,\quad1,\quad1,\quad1,\quad\dots, \]
respectively. Their contributions to $\fp(\mc B_A)$ exceed six, and the candidate is ruled out.

\paragraph{$A\cong1\oplus X\oplus Y\oplus Z$.}
It has $\fp_{\mc B}(A)=6$, and demands
\[ \fp(\mc B_A^0)=1,\quad\fp(\mc B_A)=6. \]
Computing $b_j\otimes A$, we find candidate simple objects
\[ 1\oplus X\oplus Y\oplus Z,\quad1\oplus X\oplus2Y,\quad1\oplus X\oplus2Z,\quad Y\oplus T\oplus U,\quad Z\oplus T\oplus U,\quad V\oplus W, \]
with Frobenius-Perron dimensions
\[ 1,\quad1,\quad1,\quad1,\quad1,\quad1, \]
respectively. Their contributions to Frobenius-Perron dimension match. We proceed to search for a rank six fusion category. Employing the free module functor, one can work out monoidal product $\otimes_A$. As a result, we find
\[ m_j\otimes_Am_6\cong m_6\cong m_6\otimes_Am_j\ (j=1,2,\dots,5). \]
We conclude $1\oplus X\oplus Y\oplus Z$ cannot be étale. The reason is as follows. If it were étale, $\mc B_A$ should be a fusion category. As all simple objects $m_j$ have $\fp_{\mc B_A}(m_j)=1$, the fusion ring should be multiplicity-free and listed in \cite{LPR20} or AnyonWiki. However, there is no fusion category with such monoidal product, a contradiction. Thus, the candidate is ruled out.

\paragraph{$A\cong1\oplus3X\oplus Y$.} It has $\fp_{\mc B}(A)=6$, and demands
\[ \fp(\mc B_A^0)=1,\quad\fp(\mc B_A)=6 \]
as in the previous candidate. Calculating $b_j\otimes A$, we find candidate simple objects (assuming the smallest possible Frobenius-Perron dimensions)
\[ 1\oplus3X\oplus Y,\quad3\oplus X\oplus Y,\quad1\oplus X\oplus2Y,\quad3Y,\quad Z\oplus T\oplus U,\quad3Z,\quad3T,\quad3U,\quad2V,\quad2W,\quad\dots \]
with Frobenius-Perron dimensions
\[ 1,\quad1,\quad1,\quad1,\quad1,\quad1,\quad1,\quad1,\quad1,\quad1,\quad\dots\ , \]
respectively. Their contributions to Frobenius-Perron dimension exceed six, and the candidate is ruled out.

\paragraph{$A\cong1\oplus3X\oplus Z$.} Symmetry $(YZ)$ of our ambient category $\mc B\simeq\text{Rep}(D(D_3))$ and the previous analysis rule out the candidate.

\paragraph{$A\cong1\oplus5X$.} It has $\fp_{\mc B}(A)=6$, and demands
\[ \fp(\mc B_A^0)=1,\quad\fp(\mc B_A)=6. \]
Computing $b_j\otimes A$, we find candidate simple objects (assuming the smallest possible Frobenius-Perron dimensions)
\[ 1\oplus5X,\quad5\oplus X,\quad3Y,\quad3Z,\quad3T,\quad3U,\quad2V,\quad V\oplus W,\quad2W, \]
with Frobenius-Perron dimensions
\[ 1,\quad1,\quad1,\quad1,\quad1,\quad1,\quad1,\quad1,\quad1, \]
respectively. Their contributions to Frobenius-Perron dimension exceed six, and the candidate is ruled out.\newline

We conclude
\begin{table}[H]
\begin{center}
\begin{tabular}{c|c|c|c}
    Connected étale algebra $A$&$\mc B_A$&$\rank(\mc B_A)$&Lagrangian?\\\hline
    $1$&$\mc B$&$8$&No\\
    $1\oplus X$&$\text{TY}(\mbb Z/3\mbb Z\times\mbb Z/3\mbb Z)$&10&No\\
    $1\oplus Y$ for 1st\&2nd $h$&$\mc C(\text{FR}^{6,0}_2)$&6&No\\
    $1\oplus Z$ for 1st\&2nd $h$&$\mc C(\text{FR}^{6,0}_2)$&6&No\\
    $1\oplus X\oplus2Y$ for 1st\&2nd $h$&$\mc C(\text{FR}^{6,2}_1)\text{ or }\mc C(\text{FR}^{6,4}_1)$&6&Yes\\
    $1\oplus X\oplus2Z$ for 1st\&2nd $h$&$\mc C(\text{FR}^{6,2}_1)\text{ or }\mc C(\text{FR}^{6,4}_1)$&6&Yes\\
    $1\oplus Y\oplus V$ for (\ref{RepDD31YVetale})&$\text{Rep}(S_3)$&3&Yes\\
    $1\oplus Z\oplus V$ for (\ref{RepDD31ZVetale})&$\text{Rep}(S_3)$&3&Yes
\end{tabular}.
\end{center}
\caption{Connected étale algebras in rank eight MFC $\mcal B\simeq\text{Rep}(D(D_3))$}\label{rank8RepDD3results}
\end{table}
\hspace{-17pt}All the 16 MFCs $\mc B\simeq\text{Rep}(D(D_3))$'s fail to be completely anisotropic. It may be useful to summarize these results in a ``cascade'' of conformal embeddings:
\begin{figure}[H]
\begin{center}
\includegraphics[keepaspectratio, scale=0.65]{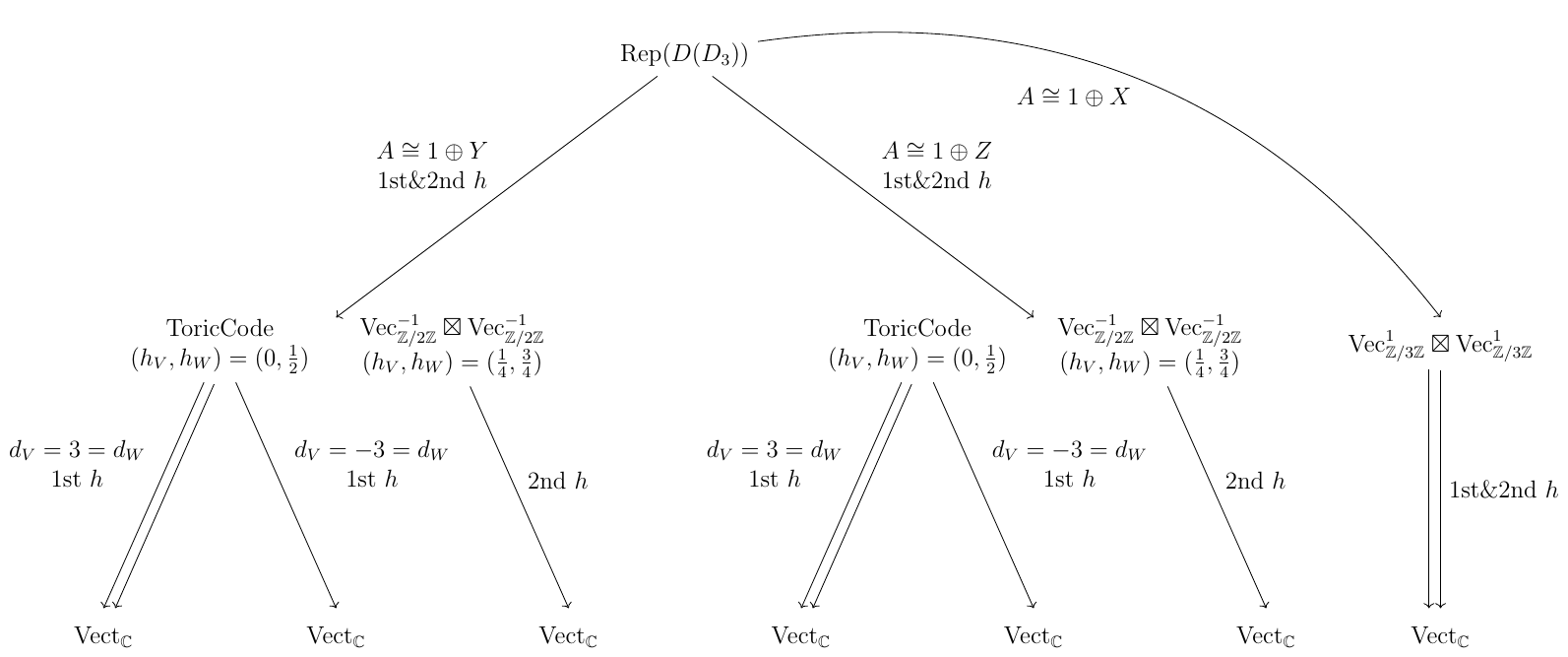}
\caption{``Cascades'' of conformal embeddings}\label{cascades}
\end{center}
\end{figure}\hspace{-22pt}
Let us see the consistency of our results. 

\paragraph{3rd and 4th conformal dimensions.} For these conformal dimensions, the only connected étale algebra is $A\cong1\oplus X$. Correspondingly, the only path is $\text{Rep}(D(D_3))\to\vecG_{\mbb Z/3\mbb Z}^1\boxtimes\vecG_{\mbb Z/3\mbb Z}^1$. Since the resulting MFC has central charge (``anomaly'') different from $\vect$, it cannot descend to $\vect$ and the path terminates there.

\paragraph{2nd conformal dimension.} For the second conformal dimensions, there are two additional paths given by $A\cong1\oplus Y,1\oplus Z$. The compositions of the two operations
\[ \text{Rep}(D(D_3))\to\vecG_{\mbb Z/3\mbb Z}^1\boxtimes\vecG_{\mbb Z/3\mbb Z}^1,\quad\text{Rep}(D(D_3))\to\vecG_{\mbb Z/2\mbb Z}^{-1}\boxtimes\vecG_{\mbb Z/2\mbb Z}^{-1} \]
give two inequivalent paths
\[ \text{Rep}(D(D_3))\to\vecG_{\mbb Z/3\mbb Z}^1\boxtimes\vecG_{\mbb Z/3\mbb Z}^1\to\vect. \]
This is consistent with the result that $\vecG_{\mbb Z/3\mbb Z}^1\boxtimes\vecG_{\mbb Z/3\mbb Z}^1$ with such conformal dimensions admits two connected étale algebras. (See section \ref{Z3Z3}.) The other paths via $\vecG_{\mbb Z/2\mbb Z}^{-1}\boxtimes\vecG_{\mbb Z/2\mbb Z}^{-1}$ also lead to $\vect$ because the MFC admits a Lagrangian algebra. The two inequivalent compositions
\[ \text{Rep}(D(D_3))\to\vecG_{\mbb Z/2\mbb Z}^{-1}\boxtimes\vecG_{\mbb Z/2\mbb Z}^{-1}\to\vect \]
are given by two additional connected étale algebras $A\cong1\oplus X\oplus2Y,1\oplus X\oplus2Z$.\footnote{Note that
\[ (1\oplus Y)\otimes(1\oplus X)\cong1\oplus X\oplus2Y,\quad(1\oplus Z)\otimes(1\oplus X)\cong1\oplus X\oplus2Z. \]}

\paragraph{1st conformal dimension.} For the first conformal dimensions, we still have $A\cong1\oplus X,1\oplus Y,1\oplus Z,1\oplus X\oplus2Y,1\oplus X\oplus2Z$, and the story about the path
\[ \text{Rep}(D(D_3))\to\vecG_{\mbb Z/3\mbb Z}^1\boxtimes\vecG_{\mbb Z/3\mbb Z}^1\to\vect \]
is the same as the last case. However, a difference appears in the other paths. For non-unitary $\text{Rep}(D(D_3))$ with the first conformal dimensions, $A\cong1\oplus Y,1\oplus Z$ lead to $\tc$. It admits one Lagrangian algebra giving $\tc\to\vect$. The two inequivalent compositions
\[ \text{Rep}(D(D_3))\to\tc\to\vect \]
are given by $A\cong1\oplus X\oplus2Y,1\oplus X\oplus2Z$. On the other hand, for a unitary $\text{Rep}(D(D_3))$, $A\cong1\oplus Y,1\oplus Z$ lead to $\tc$ with two Lagrangian algebras. The other inequivalent paths are given by additional connected étale algebras $A\cong1\oplus Y\oplus V,1\oplus Y\oplus V$. In total, there are four inequivalent paths
\[ \text{Rep}(D(D_3))\to\tc\to\vect \]
given by four Lagrangian algebras $A\cong1\oplus X\oplus2Y,1\oplus X\oplus2Z,1\oplus Y\oplus V,1\oplus Z\oplus V$.

\subsubsection{$\mc B\simeq su(2)_7$}
The MFCs have eight simple objects $\{1,X,Y,Z,T,U,V,W\}$ obeying monoidal products
\begin{table}[H]
\begin{center}
\begin{tabular}{c|c|c|c|c|c|c|c|c}
    $\otimes$&$1$&$X$&$Y$&$Z$&$T$&$U$&$V$&$W$\\\hline
    $1$&$1$&$X$&$Y$&$Z$&$T$&$U$&$V$&$W$\\\hline
    $X$&&$1$&$Z$&$Y$&$U$&$T$&$W$&$V$\\\hline
    $Y$&&&$1\oplus T$&$X\oplus U$&$Y\oplus V$&$Z\oplus W$&$T\oplus V$&$U\oplus W$\\\hline
    $Z$&&&&$1\oplus T$&$Z\oplus W$&$Y\oplus V$&$U\oplus W$&$T\oplus V$\\\hline
    $T$&&&&&$1\oplus T\oplus V$&$X\oplus U\oplus W$&$Y\oplus T\oplus V$&$Z\oplus U\oplus W$\\\hline
    $U$&&&&&&$1\oplus T\oplus V$&$Z\oplus U\oplus W$&$Y\oplus T\oplus V$\\\hline
    $V$&&&&&&&$1\oplus Y\oplus T\oplus V$&$X\oplus Z\oplus U\oplus W$\\\hline
    $W$&&&&&&&&$1\oplus Y\oplus T\oplus V$
\end{tabular}.
\end{center}
\end{table}
\hspace{-17pt}(One can identify $\vecG_{\mbb Z/2\mbb Z}^{-1}=\{1,X\},psu(2)_7=\{1,Y,T,V\}$, and $Z\cong X\otimes Y,U\cong X\otimes T,W\cong X\otimes V$.) Thus, they have
\begin{align*}
    \fp_{\mc B}(1)=1=\fp_{\mc B}(X),\quad\fp_{\mc B}(Y)=\frac{\sin\frac{2\pi}9}{\sin\frac\pi9}=\fp_{\mc B}(Z),\\
    \fp_{\mc B}(T)=\frac{\sin\frac{3\pi}9}{\sin\frac\pi9}=\fp_{\mc B}(U),\quad\fp_{\mc B}(V)=\frac{\sin\frac{4\pi}9}{\sin\frac\pi9}=\fp_{\mc B}(W),
\end{align*}
and
\[ \fp(\mc B)=\frac9{2\sin^2\frac\pi9}\approx38.5. \]
Their quantum dimensions $d_j$'s are solutions of the same multiplication rules $d_id_j=\sum_{k=1}^8{N_{ij}}^kd_k$. There are six (nonzero) solutions
\begin{align*}
    (d_X,d_Y,d_Z,d_T,d_U,d_V,d_W)&=(-1,-\frac{\sin\frac{\pi}9}{\cos\frac\pi{18}},\frac{\sin\frac{\pi}9}{\cos\frac\pi{18}},-\frac{\sin\frac{3\pi}9}{\cos\frac\pi{18}},\frac{\sin\frac{3\pi}9}{\cos\frac\pi{18}},\frac{\sin\frac{2\pi}9}{\cos\frac\pi{18}},-\frac{\sin\frac{2\pi}9}{\cos\frac\pi{18}}),\\
    &~~~~(1,-\frac{\sin\frac{\pi}9}{\cos\frac\pi{18}},-\frac{\sin\frac{\pi}9}{\cos\frac\pi{18}},-\frac{\sin\frac{3\pi}9}{\cos\frac\pi{18}},-\frac{\sin\frac{3\pi}9}{\cos\frac\pi{18}},\frac{\sin\frac{2\pi}9}{\cos\frac\pi{18}},\frac{\sin\frac{2\pi}9}{\cos\frac\pi{18}}),\\
    &~~~~(-1,-\frac{\sin\frac{4\pi}9}{\cos\frac{5\pi}{18}},\frac{\sin\frac{4\pi}9}{\cos\frac{5\pi}{18}},\frac{\sin\frac{3\pi}9}{\cos\frac{5\pi}{18}},-\frac{\sin\frac{3\pi}9}{\cos\frac{5\pi}{18}},-\frac{\sin\frac{\pi}9}{\cos\frac{5\pi}{18}},\frac{\sin\frac{\pi}9}{\cos\frac{5\pi}{18}}),\\
    &~~~~(1,-\frac{\sin\frac{4\pi}9}{\cos\frac{5\pi}{18}},-\frac{\sin\frac{4\pi}9}{\cos\frac{5\pi}{18}},\frac{\sin\frac{3\pi}9}{\cos\frac{5\pi}{18}},\frac{\sin\frac{3\pi}9}{\cos\frac{5\pi}{18}},-\frac{\sin\frac{\pi}9}{\cos\frac{5\pi}{18}},-\frac{\sin\frac{\pi}9}{\cos\frac{5\pi}{18}}),\\
    &~~~~(-1,\frac{\sin\frac{2\pi}9}{\sin\frac\pi9},-\frac{\sin\frac{2\pi}9}{\sin\frac\pi9},\frac{\sin\frac{3\pi}9}{\sin\frac\pi9},-\frac{\sin\frac{3\pi}9}{\sin\frac\pi9},\frac{\sin\frac{4\pi}9}{\sin\frac\pi9},-\frac{\sin\frac{4\pi}9}{\sin\frac\pi9}),\\
    &~~~~(1,\frac{\sin\frac{2\pi}9}{\sin\frac\pi9},\frac{\sin\frac{2\pi}9}{\sin\frac\pi9},\frac{\sin\frac{3\pi}9}{\sin\frac\pi9},\frac{\sin\frac{3\pi}9}{\sin\frac\pi9},\frac{\sin\frac{4\pi}9}{\sin\frac\pi9},\frac{\sin\frac{4\pi}9}{\sin\frac\pi9}),
\end{align*}
with categorical dimensions
\[ D^2(\mc B)=\frac9{2\cos^2\frac\pi{18}}(\approx4.6),\quad\frac9{2\cos^2\frac{5\pi}{18}}(\approx10.9),\quad\frac9{2\sin^2\frac\pi9}, \]
respectively. Only the last quantum dimensions give unitary MFCs. They each have four conformal dimensions
\begin{align*}
    \hspace{-40pt}&(h_X,h_Y,h_Z,h_T,h_U,h_V,h_W)\\
    \hspace{-40pt}&=\begin{cases}(\frac14,\frac13,\frac7{12},\frac89,\frac{5}{36},\frac23,\frac{11}{12}),(\frac14,\frac23,\frac{11}{12},\frac19,\frac{13}{36},\frac13,\frac{7}{12}),(\frac34,\frac13,\frac1{12},\frac89,\frac{23}{36},\frac23,\frac{5}{12}),(\frac34,\frac23,\frac5{12},\frac19,\frac{31}{36},\frac13,\frac{1}{12})&(\text{1st\&2nd}),\\(\frac14,\frac13,\frac7{12},\frac59,\frac{29}{36},\frac23,\frac{11}{12}),(\frac14,\frac23,\frac{11}{12},\frac49,\frac{25}{36},\frac13,\frac{7}{12}),(\frac34,\frac13,\frac1{12},\frac59,\frac{11}{36},\frac23,\frac{5}{12}),(\frac34,\frac23,\frac5{12},\frac49,\frac{7}{36},\frac13,\frac{1}{12})&(\text{3rd\&4th}),\\(\frac14,\frac13,\frac7{12},\frac29,\frac{17}{36},\frac23,\frac{11}{12}),(\frac14,\frac23,\frac{11}{12},\frac79,\frac{1}{36},\frac13,\frac{7}{12}),(\frac34,\frac13,\frac1{12},\frac29,\frac{35}{36},\frac23,\frac{5}{12}),(\frac34,\frac23,\frac5{12},\frac79,\frac{19}{36},\frac13,\frac{1}{12})&(\text{5th\&6th}).\end{cases}\quad(\mods1)
\end{align*}
The $S$-matrices are given by
\[ \widetilde S=\begin{pmatrix}1&d_X&d_Y&d_Xd_Y&d_T&d_Xd_T&d_V&d_Xd_V\\d_X&-1&d_Xd_Y&-d_Y&d_Xd_T&-d_T&d_Xd_V&-d_V\\d_Y&d_Xd_Y&-d_V&-d_Xd_V&d_T&d_Xd_T&-1&-d_X\\d_Xd_Y&-d_Y&-d_Xd_V&d_V&d_Xd_T&-d_T&-d_X&1\\d_T&d_Xd_T&d_T&d_Xd_T&0&0&-d_T&-d_Xd_T\\d_Xd_T&-d_T&d_Xd_T&-d_T&0&0&-d_Xd_T&d_T\\d_V&d_Xd_V&-1&-d_X&-d_T&-d_Xd_T&d_Y&d_Xd_Y\\d_Xd_V&-d_V&-d_X&1&-d_Xd_T&d_T&d_Xd_Y&-d_Y\end{pmatrix}. \]
There are
\[ 6(\text{quantum dimensions})\times4(\text{conformal dimensions})\times2(\text{categorical dimensions})=48 \]
MFCs, among which those eight with the last quantum dimensions give unitary MFCs. We classify connected étale algebras in all 48 MFCs simultaneously.

We work with an ansatz
\[ A\cong1\oplus n_XX\oplus n_YY\oplus n_ZZ\oplus n_TT\oplus n_UU\oplus n_VV\oplus n_WW \]
with $n_j\in\mbb N$. It has
\[ \fp_{\mc B}(A)=1+n_X+\frac{\sin\frac{2\pi}9}{\sin\frac\pi9}(n_Y+n_Z)+\frac{\sin\frac{3\pi}9}{\sin\frac\pi9}(n_T+n_U)+\frac{\sin\frac{4\pi}9}{\sin\frac\pi9}(n_V+n_W). \]
For this to obey (\ref{FPdimA2bound}), the natural numbers can take only 43 values. The sets contain the one with all $n_j$'s be zero. This is the trivial connected étale algebra $A\cong1$ giving $\mc B_A^0\simeq\mc B_A\simeq\mc B$. The other 42 candidates contain nontrivial simple object(s) $b_j\not\cong1$. They all fail to be commutative because the nontrivial simple objects have nontrivial conformal dimensions.

We conclude
\begin{table}[H]
\begin{center}
\begin{tabular}{c|c|c|c}
    Connected étale algebra $A$&$\mc B_A$&$\rank(\mc B_A)$&Lagrangian?\\\hline
    $1$&$\mc B$&$8$&No
\end{tabular}.
\end{center}
\caption{Connected étale algebras in rank eight MFC $\mcal B\simeq su(2)_7$}\label{rank8su27results}
\end{table}
\hspace{-17pt}All the 48 MFCs $\mc B\simeq su(2)_7$'s are completely anisotropic.

\subsubsection{$\mc B\simeq\fib\boxtimes\fib\boxtimes\fib$}
The MFCs have eight simple objects $\{1,X,Y,Z,T,U,V,W\}$ obeying monoidal products
\begin{table}[H]
\begin{center}
\makebox[1 \textwidth][c]{       
\resizebox{1.2 \textwidth}{!}{\begin{tabular}{c|c|c|c|c|c|c|c|c}
    $\otimes$&$1$&$X$&$Y$&$Z$&$T$&$U$&$V$&$W$\\\hline
    $1$&$1$&$X$&$Y$&$Z$&$T$&$U$&$V$&$W$\\\hline
    $X$&&$1\oplus X$&$V$&$T$&$Z\oplus T$&$W$&$Y\oplus V$&$U\oplus W$\\\hline
    $Y$&&&$1\oplus Y$&$U$&$W$&$Z\oplus U$&$X\oplus V$&$T\oplus W$\\\hline
    $Z$&&&&$1\oplus Z$&$X\oplus T$&$Y\oplus U$&$W$&$V\oplus W$\\\hline
    $T$&&&&&$1\oplus X\oplus Z\oplus T$&$V\oplus W$&$U\oplus W$&$Y\oplus U\oplus V\oplus W$\\\hline
    $U$&&&&&&$1\oplus Y\oplus Z\oplus U$&$T\oplus W$&$X\oplus T\oplus V\oplus W$\\\hline
    $V$&&&&&&&$1\oplus X\oplus Y\oplus V$&$Z\oplus T\oplus U\oplus W$\\\hline
    $W$&&&&&&&&$1\oplus X\oplus Y\oplus Z\oplus T\oplus U\oplus V\oplus W$
\end{tabular}.}}
\end{center}
\end{table}
\hspace{-17pt}(One can identify $\fib=\{1,X\},\{1,Y\},\{1,Z\}$, and $T\cong X\otimes Z,U\cong Y\otimes Z,V\cong X\otimes Y,W\cong X\otimes Y\otimes Z$.) Thus, they have
\begin{align*}
    \fp_{\mc B}(1)=1,\quad\fp_{\mc B}(X)=\fp_{\mc B}(Y)=\fp_{\mc B}(Z)=\zeta,\\
    \fp_{\mc B}(T)=\fp_{\mc B}(U)=\fp_{\mc B}(V)=\frac{3+\sqrt5}2,\quad\fp_{\mc B}(W)=2+\sqrt5,
\end{align*}
and
\[ \fp(\mc B)=25+10\sqrt5\approx47.4. \]
Their quantum dimensions $d_j$'s are solutions of the same multiplication rules $d_id_j=\sum_{k=1}^8{N_{ij}}^kd_k$. There are eight solutions
\begin{align*}
    &(d_X,d_Y,d_Z,d_T,d_U,d_V,d_W)\\
    &=(-\zeta^{-1},-\zeta^{-1},-\zeta^{-1},\frac{3-\sqrt5}2,\frac{3-\sqrt5}2,\frac{3-\sqrt5}2,2-\sqrt5),(-\zeta^{-1},-\zeta^{-1},\zeta,-1,-1,\frac{3-\sqrt5}2,\zeta^{-1}),\\
    &~~~~(-\zeta^{-1},\zeta,-\zeta^{-1},\frac{3-\sqrt5}2,-1,-1,\zeta^{-1}),(\zeta,-\zeta^{-1},-\zeta^{-1},-1,\frac{3-\sqrt5}2,-1,\zeta^{-1}),\\
    &~~~~(-\zeta^{-1},\zeta,\zeta,-1,\frac{3+\sqrt5}2,-1,-\zeta),(\zeta,-\zeta^{-1},\zeta,\frac{3+\sqrt5}2,-1,-1,-\zeta),\\
    &~~~~(\zeta,\zeta,-\zeta^{-1},-1,-1,\frac{3+\sqrt5}2,-\zeta),(\zeta,\zeta,\zeta,\frac{3+\sqrt5}2,\frac{3+\sqrt5}2,\frac{3+\sqrt5}2,2+\sqrt5).
\end{align*}
Only the last quantum dimensions give unitary MFCs. The first, second to fourth, fifth to seventh, and the eighth have categorical dimensions
\[ D^2(\mc B)=25-10\sqrt5(\approx2.6),\quad\frac{25-5\sqrt5}2(\approx6.9),\quad\frac{25+5\sqrt5}2(\approx18.1),\quad25+10\sqrt5, \]
respectively.

In order to list up their conformal dimensions without double-counting, we perform case analysis. Depending on quantum dimensions $\zeta,-\zeta^{-1}$ of Fibonacci objects, we have four classes.
\begin{itemize}
    \item $(d_X,d_Y,d_Z)=(\zeta,\zeta,\zeta)$. This class gives unitary MFCs. Different MFCs are given by
    \begin{align*}
        (h_X,h_Y,h_Z,h_T,h_U,h_V,h_W)&=(\frac25,\frac25,\frac25,\frac45,\frac45,\frac45,\frac15),(\frac25,\frac25,\frac35,0,0,\frac45,\frac25),\\
        &~~~~(\frac25,\frac35,\frac35,0,\frac15,0,\frac35),(\frac35,\frac35,\frac35,\frac15,\frac15,\frac15,\frac45).\quad(\mods1)
    \end{align*}
    Including the two signs of categorical dimensions, we have eight unitary MFCs.
    \item $(d_X,d_Y,d_Z)=(\zeta,\zeta,-\zeta^{-1})$. Different MFCs are given by
    \begin{align*}
        (h_X,h_Y,h_Z,h_T,h_U,h_V,h_W)&=(\frac25,\frac25,\frac15,\frac35,\frac35,\frac45,0),(\frac25,\frac25,\frac45,\frac15,\frac15,\frac45,\frac35),\\
        &~~~~(\frac25,\frac35,\frac15,\frac35,\frac45,0,\frac15),(\frac25,\frac35,\frac45,\frac15,\frac25,0,\frac45),\\
        &~~~~(\frac35,\frac35,\frac15,\frac45,\frac45,\frac15,\frac25),(\frac35,\frac35,\frac45,\frac25,\frac25,\frac15,0).\quad(\mods1)
    \end{align*}
    With two categorical dimensions, there are 12 MFCs.
    \item $(d_X,d_Y,d_Z)=(\zeta,-\zeta^{-1},-\zeta^{-1})$. Different MFCs are given by
    \begin{align*}
        (h_X,h_Y,h_Z,h_T,h_U,h_V,h_W)&=(\frac25,\frac15,\frac15,\frac35,\frac25,\frac35,\frac45),(\frac25,\frac15,\frac45,\frac15,0,\frac35,\frac25),\\
        &~~~~(\frac25,\frac45,\frac45,\frac15,\frac35,\frac15,0),(\frac35,\frac15,\frac15,\frac45,\frac25,\frac45,0),\\
        &~~~~(\frac35,\frac15,\frac45,\frac25,0,\frac45,\frac35),(\frac35,\frac45,\frac45,\frac25,\frac35,\frac25,\frac15).\quad(\mods1)
    \end{align*}
    There are 12 MFCs.
    \item $(d_X,d_Y,d_Z)=(-\zeta^{-1},-\zeta^{-1},-\zeta^{-1})$. Different MFCs are given by
    \begin{align*}
        (h_X,h_Y,h_Z,h_T,h_U,h_V,h_W)&=(\frac15,\frac15,\frac15,\frac25,\frac25,\frac25,\frac35),(\frac15,\frac15,\frac45,0,0,\frac25,\frac15),\\
        &~~~~(\frac15,\frac45,\frac45,0,\frac35,0,\frac45),(\frac45,\frac45,\frac45,\frac35,\frac35,\frac35,\frac25).\quad(\mods1)
    \end{align*}
    There are eight MFCs.
\end{itemize}
The $S$-matrices are given by
\[ \widetilde S=\begin{pmatrix}1&d_X&d_Y&d_Z&d_Xd_Z&d_Yd_Z&d_Xd_Y&d_Xd_Yd_Z\\d_X&-1&d_Xd_Y&d_Xd_Z&-d_Z&d_Xd_Yd_Z&-d_Y&-d_Yd_Z\\d_Y&d_Xd_Y&-1&d_Yd_Z&d_Xd_Yd_Z&-d_Z&-d_X&-d_Xd_Z\\d_Z&d_Xd_Z&d_Yd_Z&-1&-d_Z&-d_Y&d_Xd_Yd_Z&-d_Xd_Y\\d_Xd_Z&-d_Z&d_Xd_Yd_Z&-d_Z&1&-d_Xd_Y&-d_Yd_Z&d_Y\\d_Yd_Z&d_Xd_Yd_Z&-d_Z&-d_Y&-d_Xd_Y&1&-d_Xd_Z&d_X\\d_Xd_Y&-d_Y&-d_X&d_Xd_Yd_Z&-d_Yd_Z&-d_Xd_Z&1&d_Z\\d_Xd_Yd_Z&-d_Yd_Z&-d_Xd_Z&-d_Xd_Y&d_Y&d_X&d_Z&-1\end{pmatrix}. \]
They have additive central charge
\[ c(\mc B)=c(\fib)+c(\fib)+c(\fib)\quad(\mods8) \]
where
\[ c(\fib)=\begin{cases}\frac25&(h_\fib=\frac15),\\-\frac25&(h_\fib=\frac45),\\\frac{14}5&(h_\fib=\frac25),\\-\frac{14}5&(h_\fib=\frac35).\end{cases}\quad(\mods8) \]
There are
\[ 8+12+12+8=40 \]
MFCs, among which those eight in the first class are unitary. We classify connected étale algebras in all 40 MFCs simultaneously.

An ansatz
\[ A\cong1\oplus n_XX\oplus n_YY\oplus n_ZZ\oplus n_TT\oplus n_UU\oplus n_VV\oplus n_WW \]
with $n_j\in\mbb N$ has
\[ \fp_{\mc B}(A)=1+\zeta(n_X+n_Y+n_Z)+\frac{3+\sqrt5}2(n_T+n_U+n_V)+(2+\sqrt5)n_W. \]
For this to obey (\ref{FPdimA2bound}), the natural numbers can take only 60 values. It contains the one with all $n_j$'s be zero. It is the trivial connected étale algebra $A\cong1$ giving $\mc B_A^0\simeq\mc B_A\simeq\mc B$. Those with $X,Y,Z$ cannot be commutative because they have nontrivial conformal dimensions. Setting $n_X=n_Y=n_Z=0$, apart from the trivial one $A\cong1$, we get 10 sets
\begin{align*}
    (n_T,n_U,n_V,n_W)&=(1,0,0,0),(0,1,0,0),(0,0,1,0),(0,0,0,1),(2,0,0,0),\\
    &~~~~(1,1,0,0),(1,0,1,0),(0,2,0,0),(0,1,1,0),(0,0,2,0).
\end{align*}
Some of them are ruled out due to their Frobenius-Perron dimensions. The six candidates $(n_T,n_U,n_V,n_W)=(2,0,0,0),(1,1,0,0),(1,0,1,0),(0,2,0,0),(0,1,1,0),(0,0,2,0)$ have $\fp_{\mc B}=4+\sqrt5$, and demands $\fp(\mc B_A^0)\approx1.2$. However, there is no MFC with such Frobenius-Perron dimension. Thus, the six candidates are ruled out. Also, a candidate $(n_T,n_U,n_V,n_W)=(0,0,0,1)$ or $A\cong1\oplus W$ is ruled out because there is no MFC with $\fp(\mc B_A^0)\approx1.7$. Thus, we are left with three candidates
\[ (n_T,n_U,n_V,n_W)=(1,0,0,0),(0,1,0,0),(0,0,1,0). \]
They all have $\fp_{\mc B}(A)=\frac{5+\sqrt5}2$, and demands
\[ \fp(\mc B_A^0)=\frac{5+\sqrt5}2,\quad\fp(\mc B_A)=\frac{15+5\sqrt5}2. \]
As we recalled in section \ref{Z2fibfib}, these are commutative algebras when $h=0$ mod 1 \cite{BD11}. What is left is to check separability of them. In order to judge this point, we identify $\mc B_A$. We start from $A\cong1\oplus T$.

The Frobenius-Perron dimensions only allow
\[ \mc B_A^0\simeq\fib,\quad\mc B_A\simeq\fib\boxtimes\fib. \]
This in particular suggests $\mc B_A$ has rank four. Indeed, we find a four-dimensional NIM-rep
\begin{align*}
    n_1=1_4,\quad n_X&=\begin{pmatrix}0&1&0&0\\1&1&0&0\\0&0&0&1\\0&0&1&1\end{pmatrix}=n_Z,\quad n_Y=\begin{pmatrix}0&0&1&0\\0&0&0&1\\1&0&1&0\\0&1&0&1\end{pmatrix},\\
    n_T=\begin{pmatrix}1&1&0&0\\1&2&0&0\\0&0&1&1\\0&0&1&2\end{pmatrix},&\quad n_U=\begin{pmatrix}0&0&0&1\\0&0&1&1\\0&1&0&1\\1&1&1&1\end{pmatrix}=n_V,\quad n_W=\begin{pmatrix}0&0&1&1\\0&0&1&2\\1&1&1&1\\1&2&1&2\end{pmatrix}.
\end{align*}
The solution gives identifications
\[ m_1\cong1\oplus T,\quad m_2\cong X\oplus Z\oplus T,\quad m_3\cong Y\oplus W,\quad m_4\cong U\oplus V\oplus W. \]
For the candidate to be commutative, we need $h_T=0$ mod 1. This happens when $d_X=d_Z$. Then, the right $A$-modules have
\[ d_{\mc B_A}(m_1)=1,\quad d_{\mc B_A}(m_2)=d_X,\quad d_{\mc B_A}(m_3)=d_Y,\quad d_{\mc B_A}(m_4)=d_Xd_Y. \]
For MFCs with $h_T=0$ mod 1, $m_{1,3}$ are deconfined. They obey the monoidal products $\otimes_A$:
\begin{table}[H]
\begin{center}
\begin{tabular}{c|c|c|c|c}
    $\otimes_A$&$m_1$&$m_2$&$m_3$&$m_4$\\\hline
    $m_1$&$m_1$&$m_2$&$m_3$&$m_4$\\\hline
    $m_2$&&$m_1\oplus m_2$&$m_4$&$m_3\oplus m_4$\\\hline
    $m_3$&&&$m_1\oplus m_3$&$m_2\oplus m_4$\\\hline
    $m_4$&&&&$m_1\oplus m_2\oplus m_3\oplus m_4$
\end{tabular}.
\end{center}
\end{table}
\hspace{-17pt}
We see
\[ \mc B_A^0=\{m_1,m_3\}\simeq\fib,\quad\mc B_A=\{m_1,m_2,m_3,m_4\}\simeq\fib\boxtimes\fib. \]
The identification also matches central charges because two $\fib$ generated by $X,Z$ have opposite central charges, and both $c(\mc B)$ and $c(\mc B_A^0)$ are determined by $h_Y$. The identification of $\mc B_A$ shows $A$ is separable (hence étale) because $\fib\boxtimes\fib$ is semisimple:
\begin{equation}
    A\cong1\oplus T\quad(d_X,d_Z,h_X,h_Z)=(\zeta,\zeta,\frac25,\frac35),(-\zeta^{-1},-\zeta^{-1}.\frac15,\frac45)\quad(\mods1\text{ for }h)\label{FibFibFib1Tetale}
\end{equation}

How about the other two candidates $A\cong 1\oplus U,1\oplus V$? Since $T,U,V$ enter symmetrically, we immediately learn they give connected étale algebras:
\begin{equation}
\begin{split}
    A\cong1\oplus U&\quad(d_Y,d_Z,h_Y,h_Z)=(\zeta,\zeta,\frac25,\frac35),(-\zeta^{-1},-\zeta^{-1},\frac15,\frac45),\\
    A\cong1\oplus V&\quad(d_X,d_Y,h_X,h_Y)=(\zeta,\zeta,\frac25,\frac35),(-\zeta^{-1},-\zeta^{-1},\frac15,\frac45).
\end{split}\quad(\mods1\text{ for }h)\label{FibFibFib1UVetale}
\end{equation}

We conclude
\begin{table}[H]
\begin{center}
\begin{tabular}{c|c|c|c}
    Connected étale algebra $A$&$\mc B_A$&$\rank(\mc B_A)$&Lagrangian?\\\hline
    $1$&$\mc B$&$8$&No\\
    $1\oplus T$ for (\ref{FibFibFib1Tetale})&$\fib\boxtimes\fib$&4&No\\
    $1\oplus U$ for (\ref{FibFibFib1UVetale})&$\fib\boxtimes\fib$&4&No\\
    $1\oplus V$ for (\ref{FibFibFib1UVetale})&$\fib\boxtimes\fib$&4&No
\end{tabular}.
\end{center}
\caption{Connected étale algebras in rank eight MFC $\mcal B\simeq\fib\boxtimes\fib\boxtimes\fib$}\label{rank8fibfibfibresults}
\end{table}
\hspace{-17pt}Namely, those 16 MFCs $\mc B\simeq\fib\boxtimes\fib\boxtimes\fib$'s in (\ref{FibFibFib1Tetale},\ref{FibFibFib1UVetale}) fail to be completely anisotropic, while the other 24 are completely anisotropic.

\subsubsection{$\mc B\simeq\fib\boxtimes psu(2)_7$}
The MFCs have eight simple objects $\{1,X,Y,Z,T,U,V,W\}$ obeying monoidal products
\begin{table}[H]
\begin{center}
\makebox[1 \textwidth][c]{       
\resizebox{1.2 \textwidth}{!}{\begin{tabular}{c|c|c|c|c|c|c|c|c}
    $\otimes$&$1$&$X$&$Y$&$Z$&$T$&$U$&$V$&$W$\\\hline
    $1$&$1$&$X$&$Y$&$Z$&$T$&$U$&$V$&$W$\\\hline
    $X$&&$1\oplus X$&$U$&$V$&$W$&$Y\oplus U$&$Z\oplus V$&$T\oplus W$\\\hline
    $Y$&&&$1\oplus Z$&$Y\oplus T$&$Z\oplus T$&$X\oplus V$&$U\oplus W$&$V\oplus W$\\\hline
    $Z$&&&&$1\oplus Z\oplus T$&$Y\oplus Z\oplus T$&$U\oplus W$&$X\oplus V\oplus W$&$U\oplus V\oplus W$\\\hline
    $T$&&&&&$1\oplus Y\oplus Z\oplus T$&$V\oplus W$&$U\oplus V\oplus W$&$X\oplus U\oplus V\oplus W$\\\hline
    $U$&&&&&&$1\oplus X\oplus Z\oplus V$&$Y\oplus T\oplus U\oplus W$&$Z\oplus T\oplus V\oplus W$\\\hline
    $V$&&&&&&&$1\oplus X\oplus Z\oplus T\oplus V\oplus W$&$Y\oplus Z\oplus T\oplus U\oplus V\oplus W$\\\hline
    $W$&&&&&&&&$1\oplus X\oplus Y\oplus Z\oplus T\oplus U\oplus V\oplus W$
\end{tabular}.}}
\end{center}
\end{table}
\hspace{-17pt}(One can identify $\fib=\{1,X\},psu(2)_7=\{1,Y,Z,T\}$, and $U\cong X\otimes Y,V\cong X\otimes Z,W\cong X\otimes T$.) Thus, they have
\begin{align*}
    \fp_{\mc B}(1)=1,\quad\fp_{\mc B}(X)=\zeta,\quad&\fp_{\mc B}(Y)=\frac{\sin\frac{2\pi}9}{\sin\frac\pi9},\quad\fp_{\mc B}(Z)=\frac{\sin\frac{3\pi}9}{\sin\frac\pi9},\\
    \fp_{\mc B}(T)=\frac{\sin\frac{4\pi}9}{\sin\frac\pi9},\quad\fp_{\mc B}(U)=\zeta\frac{\sin\frac{2\pi}9}{\sin\frac\pi9},\quad&\fp_{\mc B}(V)=\zeta\frac{\sin\frac{3\pi}9}{\sin\frac\pi9},\quad\fp_{\mc B}(W)=\zeta\frac{\sin\frac{4\pi}9}{\sin\frac\pi9},
\end{align*}
and
\[ \fp(\mc B)=\frac{45+9\sqrt5}{8\sin^2\frac\pi9}\approx69.6. \]
Their quantum dimensions $d_j$'s are solutions of the same multiplication rules $d_id_j=\sum_{k=1}^8{N_{ij}}^kd_k$. There are six (nonzero) solutions
\begin{align*}
    \hspace{-30pt}(d_X,d_Y,d_Z,d_T,d_U,d_V,d_W)&=(-\zeta^{-1},-\frac{\sin\frac\pi9}{\cos\frac\pi{18}},-\frac{\sin\frac\pi3}{\cos\frac\pi{18}},1-\frac{\sin\frac\pi9}{\cos\frac\pi{18}},\zeta^{-1}\frac{\sin\frac\pi9}{\cos\frac\pi{18}},\zeta^{-1}\frac{\sin\frac\pi3}{\cos\frac\pi{18}}-\zeta^{-1}+\zeta^{-1}\frac{\sin\frac\pi9}{\cos\frac\pi{18}}),\\
    &~~~~(-\zeta^{-1},-\frac{\cos\frac\pi{18}}{\sin\frac{2\pi}9},\frac{\sin\frac\pi3}{\sin\frac{2\pi}9},1-\frac{\cos\frac\pi{18}}{\sin\frac{2\pi}9},\zeta^{-1}\frac{\cos\frac\pi{18}}{\sin\frac{2\pi}9},-\zeta^{-1}\frac{\sin\frac\pi3}{\sin\frac{2\pi}9},-\zeta^{-1}+\zeta^{-1}\frac{\cos\frac\pi{18}}{\sin\frac{2\pi}9}),\\
    &~~~~(\zeta,-\frac{\sin\frac\pi9}{\cos\frac\pi{18}},-\frac{\sin\frac\pi3}{\cos\frac\pi{18}},1-\frac{\sin\frac\pi9}{\cos\frac\pi{18}},-\zeta\frac{\sin\frac\pi9}{\cos\frac\pi{18}},-\zeta\frac{\sin\frac\pi3}{\cos\frac\pi{18}},\zeta-\zeta\frac{\sin\frac\pi9}{\cos\frac\pi{18}}),\\
    &~~~~(\zeta,-\frac{\cos\frac\pi{18}}{\sin\frac{2\pi}9},\frac{\sin\frac\pi3}{\sin\frac{2\pi}9},1-\frac{\cos\frac\pi{18}}{\sin\frac{2\pi}9},-\zeta\frac{\cos\frac\pi{18}}{\sin\frac{2\pi}9},\zeta\frac{\sin\frac\pi3}{\sin\frac{2\pi}9},\zeta-\zeta\frac{\cos\frac\pi{18}}{\sin\frac{2\pi}9}),\\
    &~~~~(-\zeta^{-1},\frac{\sin\frac{2\pi}9}{\sin\frac\pi9},\frac{\sin\frac{3\pi}9}{\sin\frac\pi9},\frac{\sin\frac{4\pi}9}{\sin\frac\pi9},-\zeta^{-1}\frac{\sin\frac{2\pi}9}{\sin\frac\pi9},-\zeta^{-1}\frac{\sin\frac{3\pi}9}{\sin\frac\pi9},-\zeta^{-1}\frac{\sin\frac{4\pi}9}{\sin\frac\pi9}),\\
    &~~~~(\zeta,\frac{\sin\frac{2\pi}9}{\sin\frac\pi9},\frac{\sin\frac{3\pi}9}{\sin\frac\pi9},\frac{\sin\frac{4\pi}9}{\sin\frac\pi9},\zeta\frac{\sin\frac{2\pi}9}{\sin\frac\pi9},\zeta\frac{\sin\frac{3\pi}9}{\sin\frac\pi9},\zeta\frac{\sin\frac{4\pi}9}{\sin\frac\pi9}),
\end{align*}
with categorical dimensions
\[ D^2(\mc B)=\begin{cases}\frac{45-9\sqrt5}{8\cos^2\frac\pi{18}}(\approx3.2)&(\text{1st}),\\\frac{45-9\sqrt5}{8\sin^2\frac{2\pi}9}(\approx7.5)&(\text{2nd}),\\\frac{45+9\sqrt5}{8\cos^2\frac\pi{18}}(\approx8.4)&(\text{3rd}),\\\frac{45+9\sqrt5}{8\sin^2\frac{2\pi}9}(\approx19.7)&(\text{4th}),\\\frac{45-9\sqrt5}{8\sin^2\frac\pi9}(\approx26.6)&(\text{5th}),\\\frac{45+9\sqrt5}{8\sin^2\frac\pi9}&(\text{6th}),\end{cases} \]
respectively. Each of them has four conformal dimensions:
\begin{align*}
    \hspace{-30pt}&(h_X,h_Y,h_Z,h_T,h_U,h_V,h_W)\\
    \hspace{-30pt}&=\begin{cases}(\frac15,\frac13,\frac89,\frac23,\frac8{15},\frac4{45},\frac{13}{15}),(\frac15,\frac23,\frac19,\frac13,\frac{13}{15},\frac{14}{45},\frac{8}{15}),(\frac45,\frac13,\frac89,\frac23,\frac{2}{15},\frac{31}{45},\frac{7}{15}),(\frac45,\frac23,\frac19,\frac13,\frac{7}{15},\frac{41}{45},\frac{2}{15})&(\text{1st}),\\(\frac15,\frac13,\frac59,\frac23,\frac8{15},\frac{34}{45},\frac{13}{15}),(\frac15,\frac23,\frac49,\frac13,\frac{13}{15},\frac{29}{45},\frac{8}{15}),(\frac45,\frac13,\frac59,\frac23,\frac{2}{15},\frac{16}{45},\frac{7}{15}),(\frac45,\frac23,\frac49,\frac13,\frac{7}{15},\frac{11}{45},\frac{2}{15})&(\text{2nd}),\\(\frac25,\frac13,\frac89,\frac23,\frac{11}{15},\frac{13}{45},\frac{1}{15}),(\frac25,\frac23,\frac19,\frac13,\frac{1}{15},\frac{23}{45},\frac{11}{15}),(\frac35,\frac13,\frac89,\frac23,\frac{14}{15},\frac{22}{45},\frac{4}{15}),(\frac35,\frac23,\frac19,\frac13,\frac{4}{15},\frac{32}{45},\frac{14}{15})&(\text{3rd}),\\(\frac25,\frac13,\frac59,\frac23,\frac{11}{15},\frac{43}{45},\frac{1}{15}),(\frac25,\frac23,\frac49,\frac13,\frac{1}{15},\frac{38}{45},\frac{11}{15}),(\frac35,\frac13,\frac59,\frac23,\frac{14}{15},\frac{7}{45},\frac{4}{15}),(\frac35,\frac23,\frac49,\frac13,\frac{4}{15},\frac{2}{45},\frac{14}{15})&(\text{4th}),\\(\frac15,\frac23,\frac79,\frac13,\frac{13}{15},\frac{44}{45},\frac{8}{15}),(\frac15,\frac13,\frac29,\frac23,\frac{8}{15},\frac{19}{45},\frac{13}{15}),(\frac45,\frac23,\frac79,\frac13,\frac{7}{15},\frac{26}{45},\frac{2}{15}),(\frac45,\frac13,\frac29,\frac23,\frac{2}{15},\frac{1}{45},\frac{7}{15})&(\text{5th}),\\(\frac25,\frac23,\frac79,\frac13,\frac1{15},\frac8{45},\frac{11}{15}),(\frac25,\frac13,\frac29,\frac23,\frac{11}{15},\frac{28}{45},\frac{1}{15}),(\frac35,\frac23,\frac79,\frac13,\frac{4}{15},\frac{17}{45},\frac{14}{15}),(\frac35,\frac13,\frac29,\frac23,\frac{14}{15},\frac{37}{45},\frac{4}{15})&(\text{6th}).\end{cases}\quad(\mods1)
\end{align*}
The $S$-matrices are given by
\[ \widetilde S=\begin{pmatrix}1&d_X&d_Y&d_Z&d_T&d_Xd_Y&d_Xd_Z&d_Xd_T\\d_X&-1&d_Xd_Y&d_Xd_Y&d_Xd_T&-d_Y&-d_Z&-d_T\\d_Y&d_Xd_Y&-d_T&d_Z&-1&-d_Xd_T&d_Xd_Z&-d_X\\d_Z&d_Xd_Z&d_Z&0&-d_Z&d_Xd_Z&0&-d_Xd_Z\\d_T&d_Xd_T&-1&-d_Z&d_Y&-d_X&-d_Xd_Z&d_Xd_Y\\d_Xd_Y&-d_Y&-d_Xd_T&d_Xd_Z&-d_X&d_T&-d_Z&1\\d_Xd_Z&-d_Z&d_Xd_Z&0&-d_Xd_Z&-d_Z&0&d_Z\\d_Xd_T&-d_T&-d_X&-d_Xd_Z&d_Xd_Y&1&d_Z&-d_Y\end{pmatrix}. \]
There are
\[ 6(\text{quantum dimensions})\times4(\text{conformal dimensions})\times2(\text{categorical dimensions})=48 \]
MFCs, among which those eight with the sixth quantum dimensions are unitary. We classify connected étale algebras in all 48 MFCs simultaneously.

An ansatz
\[ A\cong1\oplus n_XX\oplus n_YY\oplus n_ZZ\oplus n_TT\oplus n_UU\oplus n_VV\oplus n_WW \]
with $n_j\in\mbb N$ has
\[ \fp_{\mc B}(A)=1+\zeta n_X+\frac{\sin\frac{2\pi}9}{\sin\frac\pi9}n_Y+\frac{\sin\frac{3\pi}9}{\sin\frac\pi9}n_Z+\frac{\sin\frac{4\pi}9}{\sin\frac\pi9}n_T+\zeta\frac{\sin\frac{2\pi}9}{\sin\frac\pi9}n_U+\zeta\frac{\sin\frac{3\pi}9}{\sin\frac\pi9}n_V+\zeta\frac{\sin\frac{4\pi}9}{\sin\frac\pi9}n_W. \]
For this to obey (\ref{FPdimA2bound}), the natural numbers can take only 54 values. The sets contain the one with all $n_j$'s zero. This corresponds to the trivial connected étale algebra $A\cong1$ giving $\mc B_A^0\simeq\mc B_A\simeq\mc B$. The other 53 candidates fail to be commutative because they contain nontrivial simple objects with nontrivial conformal dimensions.

We conclude
\begin{table}[H]
\begin{center}
\begin{tabular}{c|c|c|c}
    Connected étale algebra $A$&$\mc B_A$&$\rank(\mc B_A)$&Lagrangian?\\\hline
    $1$&$\mc B$&$8$&No
\end{tabular}.
\end{center}
\caption{Connected étale algebras in rank eight MFC $\mcal B\simeq\fib\boxtimes psu(2)_7$}\label{rank8fibpsu27results}
\end{table}
\hspace{-17pt}All the 48 MFCs $\mc B\simeq\fib\boxtimes psu(2)_7$'s are completely anisotropic.

\subsubsection{$\mc B\simeq psu(2)_{15}$}
The MFCs have eight simple objects $\{1,X,Y,Z,T,U,V,W\}$ obeying monoidal products
\begin{table}[H]
\begin{center}
\makebox[1 \textwidth][c]{       
\resizebox{1.2 \textwidth}{!}{\begin{tabular}{c|c|c|c|c|c|c|c|c}
    $\otimes$&$1$&$X$&$Y$&$Z$&$T$&$U$&$V$&$W$\\\hline
    $1$&$1$&$X$&$Y$&$Z$&$T$&$U$&$V$&$W$\\\hline
    $X$&&$1\oplus Y$&$X\oplus Z$&$Y\oplus T$&$Z\oplus U$&$T\oplus V$&$U\oplus W$&$V\oplus W$\\\hline
    $Y$&&&$1\oplus Y\oplus T$&$X\oplus Z\oplus U$&$Y\oplus T\oplus V$&$Z\oplus U\oplus W$&$T\oplus V\oplus W$&$U\oplus V\oplus W$\\\hline
    $Z$&&&&$1\oplus Y\oplus T\oplus V$&$X\oplus Z\oplus U\oplus W$&$Y\oplus T\oplus V\oplus W$&$Z\oplus U\oplus V\oplus W$&$T\oplus U\oplus V\oplus W$\\\hline
    $T$&&&&&$1\oplus Y\oplus T\oplus V\oplus W$&$X\oplus Z\oplus U\oplus V\oplus W$&$Y\oplus T\oplus U\oplus V\oplus W$&$Z\oplus T\oplus U\oplus V\oplus W$\\\hline
    $U$&&&&&&$1\oplus Y\oplus T\oplus U\oplus V\oplus W$&$X\oplus Z\oplus T\oplus U\oplus V\oplus W$&$Y\oplus Z\oplus T\oplus U\oplus V\oplus W$\\\hline
    $V$&&&&&&&$1\oplus Y\oplus Z\oplus T\oplus U\oplus V\oplus W$&$X\oplus Y\oplus Z\oplus T\oplus U\oplus V\oplus W$\\\hline
    $W$&&&&&&&&$1\oplus X\oplus Y\oplus Z\oplus T\oplus U\oplus V\oplus W$
\end{tabular}.}}
\end{center}
\end{table}
\hspace{-17pt}Thus, they have
\begin{align*}
    \fp_{\mc B}(1)=1,\quad\fp_{\mc B}(X)=\frac{\sin\frac{2\pi}{17}}{\sin\frac\pi{17}},\quad&\fp_{\mc B}(Y)=\frac{\sin\frac{3\pi}{17}}{\sin\frac\pi{17}},\quad\fp_{\mc B}(Z)=\frac{\sin\frac{4\pi}{17}}{\sin\frac\pi{17}},\\
    \fp_{\mc B}(T)=\frac{\sin\frac{5\pi}{17}}{\sin\frac\pi{17}},\quad\fp_{\mc B}(U)=\frac{\sin\frac{6\pi}{17}}{\sin\frac\pi{17}},\quad&\fp_{\mc B}(V)=\frac{\sin\frac{7\pi}{17}}{\sin\frac\pi{17}},\quad\fp_{\mc B}(W)=\frac{\sin\frac{8\pi}{17}}{\sin\frac\pi{17}},
\end{align*}
and
\[ \fp(\mc B)=\frac{17}{4\sin^2\frac\pi{17}}\approx125.9. \]
Their quantum dimensions $d_j$'s are solutions of the same multiplication rules $d_id_j=\sum_{k=1}^8{N_{ij}}^kd_k$. There are eight solutions
\begin{align*}
    (d_X,d_Y,d_Z,d_T,d_U,d_V,d_W)&=(-\frac{\sin\frac{\pi}{17}}{\cos\frac{\pi}{34}},-\frac{\sin\frac{7\pi}{17}}{\cos\frac{\pi}{34}},\frac{\sin\frac{2\pi}{17}}{\cos\frac{\pi}{34}},\frac{\sin\frac{6\pi}{17}}{\cos\frac{\pi}{34}},-\frac{\sin\frac{3\pi}{17}}{\cos\frac{\pi}{34}},-\frac{\sin\frac{5\pi}{17}}{\cos\frac{\pi}{34}},\frac{\sin\frac{4\pi}{17}}{\cos\frac{\pi}{34}}),\\
    &~~~~(\frac{\sin\frac{3\pi}{17}}{\cos\frac{3\pi}{34}},-\frac{\sin\frac{4\pi}{17}}{\cos\frac{3\pi}{34}},-\frac{\sin\frac{6\pi}{17}}{\cos\frac{3\pi}{34}},\frac{\sin\frac{\pi}{17}}{\cos\frac{3\pi}{34}},\frac{\sin\frac{8\pi}{17}}{\cos\frac{3\pi}{34}},\frac{\sin\frac{2\pi}{17}}{\cos\frac{3\pi}{34}},-\frac{\sin\frac{5\pi}{17}}{\cos\frac{3\pi}{34}}),\\
    &~~~~(-\frac{\sin\frac{5\pi}{17}}{\cos\frac{5\pi}{34}},-\frac{\sin\frac{\pi}{17}}{\cos\frac{5\pi}{34}},\frac{\sin\frac{7\pi}{17}}{\cos\frac{5\pi}{34}},-\frac{\sin\frac{4\pi}{17}}{\cos\frac{5\pi}{34}},-\frac{\sin\frac{2\pi}{17}}{\cos\frac{5\pi}{34}},\frac{\sin\frac{8\pi}{17}}{\cos\frac{5\pi}{34}},-\frac{\sin\frac{3\pi}{17}}{\cos\frac{5\pi}{34}}),\\
    &~~~~(\frac{\sin\frac{7\pi}{17}}{\cos\frac{7\pi}{34}},\frac{\sin\frac{2\pi}{17}}{\cos\frac{7\pi}{34}},-\frac{\sin\frac{3\pi}{17}}{\cos\frac{7\pi}{34}},-\frac{\sin\frac{8\pi}{17}}{\cos\frac{7\pi}{34}},-\frac{\sin\frac{4\pi}{17}}{\cos\frac{7\pi}{34}},\frac{\sin\frac{\pi}{17}}{\cos\frac{7\pi}{34}},\frac{\sin\frac{6\pi}{17}}{\cos\frac{7\pi}{34}}),\\
    &~~~~(-\frac{\sin\frac{8\pi}{17}}{\cos\frac{9\pi}{34}},\frac{\sin\frac{5\pi}{17}}{\cos\frac{9\pi}{34}},-\frac{\sin\frac{\pi}{17}}{\cos\frac{9\pi}{34}},-\frac{\sin\frac{3\pi}{17}}{\cos\frac{9\pi}{34}},\frac{\sin\frac{7\pi}{17}}{\cos\frac{9\pi}{34}},-\frac{\sin\frac{6\pi}{17}}{\cos\frac{9\pi}{34}},\frac{\sin\frac{2\pi}{17}}{\cos\frac{9\pi}{34}}),\\
    &~~~~(\frac{\sin\frac{6\pi}{17}}{\cos\frac{11\pi}{34}},\frac{\sin\frac{8\pi}{17}}{\cos\frac{11\pi}{34}},\frac{\sin\frac{5\pi}{17}}{\cos\frac{11\pi}{34}},\frac{\sin\frac{2\pi}{17}}{\cos\frac{11\pi}{34}},-\frac{\sin\frac{\pi}{17}}{\cos\frac{11\pi}{34}},-\frac{\sin\frac{4\pi}{17}}{\cos\frac{11\pi}{34}},-\frac{\sin\frac{7\pi}{17}}{\cos\frac{11\pi}{34}}),\\
    &~~~~(-\frac{\sin\frac{4\pi}{17}}{\cos\frac{13\pi}{34}},\frac{\sin\frac{6\pi}{17}}{\cos\frac{13\pi}{34}},-\frac{\sin\frac{8\pi}{17}}{\cos\frac{13\pi}{34}},\frac{\sin\frac{7\pi}{17}}{\cos\frac{13\pi}{34}},-\frac{\sin\frac{5\pi}{17}}{\cos\frac{13\pi}{34}},\frac{\sin\frac{3\pi}{17}}{\cos\frac{13\pi}{34}},-\frac{\sin\frac{\pi}{17}}{\cos\frac{13\pi}{34}}),\\
    &~~~~(\frac{\sin\frac{2\pi}{17}}{\sin\frac{\pi}{17}},\frac{\sin\frac{3\pi}{17}}{\sin\frac{\pi}{17}},\frac{\sin\frac{4\pi}{17}}{\sin\frac{\pi}{17}},\frac{\sin\frac{5\pi}{17}}{\sin\frac{\pi}{17}},\frac{\sin\frac{6\pi}{17}}{\sin\frac{\pi}{17}},\frac{\sin\frac{7\pi}{17}}{\sin\frac{\pi}{17}},\frac{\sin\frac{8\pi}{17}}{\sin\frac{\pi}{17}})
\end{align*}
with categorical dimensions
\begin{align*}
    D^2(\mc B)&=\frac{17}{4\cos^2\frac{\pi}{34}}(\approx4.3),\quad\frac{17}{4\cos^2\frac{3\pi}{34}}(\approx4.6),\quad\frac{17}{4\cos^2\frac{5\pi}{34}}(\approx5.3),\quad\frac{17}{4\cos^2\frac{7\pi}{34}}(\approx6.7),\\
    &~~~~\frac{17}{4\cos^2\frac{9\pi}{34}}(\approx9.4),\quad\frac{17}{4\cos^2\frac{11\pi}{34}}(\approx15.3),\quad\frac{17}{4\cos^2\frac{13\pi}{34}}(\approx32.6),\quad\frac{17}{4\sin^2\frac{\pi}{17}},
\end{align*}
respectively. They have two conformal dimensions each:
\[ \hspace{-30pt}(h_X,h_Y,h_Z,h_T,h_U,h_V,h_W)=\begin{cases}(\frac{6}{17},\frac{16}{17},\frac{13}{17},\frac{14}{17},\frac{2}{17},\frac{11}{17},\frac{7}{17}),(\frac{11}{17},\frac{1}{17},\frac{4}{17},\frac{3}{17},\frac{15}{17},\frac{6}{17},\frac{10}{17})&(\text{1st}),\\(\frac{1}{17},\frac{14}{17},\frac{5}{17},\frac{8}{17},\frac{6}{17},\frac{16}{17},\frac{4}{17}),(\frac{16}{17},\frac{3}{17},\frac{12}{17},\frac{9}{17},\frac{11}{17},\frac{1}{17},\frac{13}{17})&(\text{2nd}),\\(\frac{4}{17},\frac{5}{17},\frac{3}{17},\frac{15}{17},\frac{7}{17},\frac{13}{17},\frac{16}{17}),(\frac{13}{17},\frac{12}{17},\frac{14}{17},\frac{2}{17},\frac{10}{17},\frac{4}{17},\frac{1}{17})&(\text{3rd}),\\(\frac{8}{17},\frac{10}{17},\frac{6}{17},\frac{13}{17},\frac{14}{17},\frac{9}{17},\frac{15}{17}),(\frac{9}{17},\frac{7}{17},\frac{11}{17},\frac{4}{17},\frac{3}{17},\frac{8}{17},\frac{2}{17})&(\text{4th}),\\(\frac{3}{17},\frac{8}{17},\frac{15}{17},\frac{7}{17},\frac{1}{17},\frac{14}{17},\frac{12}{17}),(\frac{14}{17},\frac{9}{17},\frac{2}{17},\frac{10}{17},\frac{16}{17},\frac{3}{17},\frac{5}{17})&(\text{5th}),\\(\frac{2}{17},\frac{11}{17},\frac{10}{17},\frac{16}{17},\frac{12}{17},\frac{15}{17},\frac{8}{17}),(\frac{15}{17},\frac{6}{17},\frac{7}{17},\frac{1}{17},\frac{5}{17},\frac{2}{17},\frac{9}{17})&(\text{6th}),\\(\frac{7}{17},\frac{13}{17},\frac{1}{17},\frac{5}{17},\frac{8}{17},\frac{10}{17},\frac{11}{17}),(\frac{10}{17},\frac{4}{17},\frac{16}{17},\frac{12}{17},\frac{9}{17},\frac{7}{17},\frac{6}{17})&(\text{7th}),\\(\frac5{17},\frac2{17},\frac8{17},\frac{6}{17},\frac{13}{17},\frac{12}{17},\frac{3}{17}),(\frac{12}{17},\frac{15}{17},\frac9{17},\frac{11}{17},\frac{4}{17},\frac{5}{17},\frac{14}{17})&(\text{8th}).\end{cases}\quad(\mods1) \]
The $S$-matrices are given by
\[ \widetilde S=\begin{pmatrix}1&d_X&d_Y&d_Z&d_T&d_U&d_V&d_W\\d_X&-d_Z&d_U&-d_W&d_V&-d_T&d_Y&-1\\d_Y&d_U&d_W&d_T&d_X&-1&-d_T&-d_V\\d_Z&-d_W&d_T&-1&-d_Y&d_V&-d_U&d_X\\d_T&d_V&d_X&-d_Y&-d_W&-d_Z&1&d_U\\d_U&-d_T&-1&d_V&-d_Z&-d_X&d_W&-d_Y\\d_V&d_Y&-d_T&-d_U&1&d_W&d_X&-d_T\\d_W&-1&-d_V&d_X&d_U&-d_Y&-d_T&d_Z\end{pmatrix}. \]
There are
\[ 8(\text{quantum dimensions})\times2(\text{conformal dimensions})\times2(\text{categorical dimensions})=32 \]
MFCs, among which those four with the eighth quantum dimensions are unitary. We classify connected étale algebras in all 32 MFCs simultaneously.

An ansatz
\[ A\cong1\oplus n_XX\oplus n_YY\oplus n_ZZ\oplus n_TT\oplus n_UU\oplus n_VV\oplus n_WW \]
with $n_j\in\mbb N$ has
\[ \fp_{\mc B}(A)=1+\frac{\sin\frac{2\pi}{17}}{\sin\frac\pi{17}}n_X+\frac{\sin\frac{3\pi}{17}}{\sin\frac\pi{17}}n_Y+\frac{\sin\frac{4\pi}{17}}{\sin\frac\pi{17}}n_Z+\frac{\sin\frac{5\pi}{17}}{\sin\frac\pi{17}}n_T+\frac{\sin\frac{6\pi}{17}}{\sin\frac\pi{17}}n_U+\frac{\sin\frac{7\pi}{17}}{\sin\frac\pi{17}}n_V+\frac{\sin\frac{8\pi}{17}}{\sin\frac\pi{17}}n_W. \]
For this to obey (\ref{FPdimA2bound}), the natural numbers can take only 55 values. The sets contain the one with all $n_j$'s zero. This corresponds to the trivial connected étale algebra $A\cong1$ giving $\mc B_A^0\simeq\mc B_A\simeq\mc B$. The other 54 candidates fail to be commutative because they contain nontrivial simple objects with nontrivial conformal dimensions.

We conclude
\begin{table}[H]
\begin{center}
\begin{tabular}{c|c|c|c}
    Connected étale algebra $A$&$\mc B_A$&$\rank(\mc B_A)$&Lagrangian?\\\hline
    $1$&$\mc B$&$8$&No
\end{tabular}.
\end{center}
\caption{Connected étale algebras in rank eight MFC $\mcal B\simeq psu(2)_{15}$}\label{rank8psu215results}
\end{table}
\hspace{-17pt}All the 32 MFCs $\mc B\simeq psu(2)_{15}$'s are completely anisotropic.

\subsection{Rank nine}
\subsubsection{$\mc B\simeq su(9)_1\simeq\vecG_{\mbb Z/9\mbb Z}^1$}\label{Z9}
The MFCs have nine simple objects $\{1,X,Y,Z,S,T,U,V,W\}$ obeying monoidal products
\begin{table}[H]
\begin{center}
\begin{tabular}{c|c|c|c|c|c|c|c|c|c}
    $\otimes$&$1$&$X$&$Y$&$Z$&$S$&$T$&$U$&$V$&$W$\\\hline
    $1$&$1$&$X$&$Y$&$Z$&$S$&$T$&$U$&$V$&$W$\\\hline
    $X$&&$Y$&$1$&$T$&$V$&$W$&$S$&$U$&$Z$\\\hline
    $Y$&&&$X$&$W$&$U$&$Z$&$V$&$S$&$T$\\\hline
    $Z$&&&&$V$&$1$&$U$&$Y$&$X$&$S$\\\hline
    $S$&&&&&$W$&$X$&$T$&$Z$&$Y$\\\hline
    $T$&&&&&&$S$&$1$&$Y$&$V$\\\hline
    $U$&&&&&&&$Z$&$W$&$X$\\\hline
    $V$&&&&&&&&$T$&$1$\\\hline
    $W$&&&&&&&&&$U$
\end{tabular}.
\end{center}
\end{table}
\hspace{-17pt}(One can identify $\vecG_{\mbb Z/3\mbb Z}^1=\{1,X,Y\}$.) Thus, they have
\begin{align*}
    &\fp_{\mc B}(1)=\fp_{\mc B}(X)=\fp_{\mc B}(Y)=\fp_{\mc B}(Z)=\fp_{\mc B}(S)\\
    &=\fp_{\mc B}(T)=\fp_{\mc B}(U)=\fp_{\mc B}(V)=\fp_{\mc B}(W)=1,
\end{align*}
and
\[ \fp(\mc B)=9. \]
Their quantum dimensions $d_j$'s are solutions of the same multiplication rules $d_id_j=\sum_{k=1}^9{N_{ij}}^kd_k$. The only solution is
\[ (d_X,d_Y,d_Z,d_S,d_T,d_U,d_V,d_W)=(1,1,1,1,1,1,1,1) \]
with categorical dimension
\[ D^2(\mc B)=9. \]
Thus, all MFCs are unitary. They have two\footnote{Naively, one finds six consistent conformal dimensions. However, the other four are related to one of the two in the main text by permutations $(XY)(ZUWSTV)$ of simple objects.} conformal dimensions
\[ (h_X,h_Y,h_Z,h_S,h_T,h_U,h_V,h_W)=(0,0,\frac19,\frac19,\frac79,\frac79,\frac49,\frac49),(0,0,\frac89,\frac89,\frac29,\frac29,\frac59,\frac59)\quad(\mods1). \]
The $S$-matrices are given by
\[ \widetilde S=\begin{pmatrix}1&d_X&d_Y&d_Z&d_S&d_T&d_U&d_V&d_W\\d_X&1&1&e^{\mp2\pi i/3}&e^{\pm2\pi i/3}&e^{\mp2\pi i/3}&e^{\pm2\pi i/3}&e^{\pm2\pi i/3}&e^{\mp2\pi i/3}\\d_Y&1&1&e^{\pm2\pi i/3}&e^{\mp2\pi i/3}&e^{\pm2\pi i/3}&e^{\mp2\pi i/3}&e^{\mp2\pi i/3}&e^{\pm2\pi i/3}\\d_Z&e^{\mp2\pi i/3}&e^{\pm2\pi i/3}&e^{\pm4\pi i/9}&e^{\mp4\pi i/9}&e^{\mp2\pi i/9}&e^{\pm2\pi i/9}&-e^{\mp\pi i/9}&-e^{\pm\pi i/9}\\d_S&e^{\pm2\pi i/3}&e^{\mp2\pi i/3}&e^{\mp4\pi i/9}&e^{\pm4\pi i/9}&e^{\pm2\pi i/9}&e^{\mp2\pi i/9}&-e^{\pm\pi i/9}&-e^{\mp\pi i/9}\\d_T&e^{\mp2\pi i/3}&e^{\pm2\pi i/3}&e^{\mp2\pi i/9}&e^{\pm2\pi i/9}&-e^{\pm\pi i/9}&-e^{\mp\pi i/9}&e^{\mp4\pi i/9}&e^{\pm4\pi i/9}\\d_U&e^{\pm2\pi i/3}&e^{\mp2\pi i/3}&e^{\pm2\pi i/9}&e^{\mp2\pi i/9}&-e^{\mp\pi i/9}&-e^{\pm\pi i/9}&e^{\pm4\pi i/9}&e^{\mp4\pi i/9}\\d_V&e^{\pm2\pi i/3}&e^{\mp2\pi i/3}&-e^{\mp\pi i/9}&-e^{\pm\pi i/9}&e^{\mp4\pi i/9}&e^{\pm4\pi i/9}&e^{\mp2\pi i/9}&e^{\pm2\pi i/9}\\d_W&e^{\mp2\pi i/3}&e^{\pm2\pi i/3}&-e^{\pm\pi i/9}&-e^{\mp\pi i/9}&e^{\pm4\pi i/9}&e^{\mp4\pi i/9}&e^{\pm2\pi i/9}&e^{\mp2\pi i/9}\end{pmatrix}. \]
(All signs are correlated. In other words, one $S$-matrix is given by choosing upper signs, and the other is its complex conjugate.) They have additive central charges
\[ c(\mc B)=0\quad(\mods8). \]
There are
\[ 1(\text{quantum dimension})\times2(\text{conformal dimensions})\times2(\text{categorical dimensions})=4 \]
MFCs, which are all unitary. We classify connected étale algebras in them.

We work with an ansatz
\[ A\cong1\oplus n_XX\oplus n_YY\oplus n_ZZ\oplus n_SS\oplus n_TT\oplus n_UU\oplus n_VV\oplus n_WW \]
with $n_j\in\mbb N$. It has
\[ \fp_{\mc B}(A)=1+n_X+n_Y+n_Z+n_S+n_T+n_U+n_V+n_W. \]
For this to obey (\ref{FPdimA2bound}), the natural numbers can take only 45 values. The sets contain the one with all $n_j$'s be zero. It corresponds to the trivial connected étale algebra $A\cong1$ giving $\mc B_A^0\simeq\mc B_A\simeq\mc B$.

The candidates with simple object(s) $Z,S,T,U,V,W$ fail to be commutative because they have nontrivial conformal dimensions. Thus, we are left with nontrivial candidates with simple objects $X,Y$. Since any étale algebra $A\in\mc B$ is self-dual $A^*\cong A$ \cite{DMNO10}, the only candidate is
\[ A\cong1\oplus X\oplus Y. \]
In fact, this is a commutative algebra; the braidings were computed in \cite{KK23preMFC}, and it was found they are trivial. Thanks to the lemma 1, $A\in\vecG_{\mbb Z/3\mbb Z}^1\subset\vecG_{\mbb Z/9\mbb Z}^1$ is a commutative algebra. It also turns out separable, hence étale. Let us check this point. The algebra has $\fp_{\mc B}=3$, and demands
\[ \fp(\mc B_A^0)=1,\quad\fp(\mc B_A)=3. \]
One identifies $\mc B_A^0\simeq\vect$. This identification also matches central charges. The category $\mc B_A$ of right $A$-modules is identified as
\[ \mc B_A\simeq\vecG_{\mbb Z/3\mbb Z}^1. \]
To find this, we search for NIM-reps. Since there is no fusion category with Frobenius-Perron dimension three up to rank two, we start from three-dimensional NIM-reps. Indeed, we find a three-dimensional NIM-rep
\[ n_1=n_X=n_Y=1_3,\quad n_Z=n_T=n_W=\begin{pmatrix}0&0&1\\1&0&0\\0&1&0\end{pmatrix},\quad n_S=n_U=n_V=\begin{pmatrix}0&1&0\\0&0&1\\1&0&0\end{pmatrix}. \]
Denoting a basis of $\mc B_A$ by $\{m_1,m_2,m_3\}$, we get a multiplication table
\begin{table}[H]
\begin{center}
\begin{tabular}{c|c|c|c}
    $b_j\otimes\backslash$&$m_1$&$m_2$&$m_3$\\\hline
    $1,X,Y$&$m_1$&$m_2$&$m_3$\\
    $Z,T,W$&$m_3$&$m_1$&$m_2$\\
    $S,U,V$&$m_2$&$m_3$&$m_1$
\end{tabular}.
\end{center}
\end{table}
\hspace{-17pt}In this basis, we can identify
\[ m_1\cong1\oplus X\oplus Y,\quad m_2\cong S\oplus U\oplus V,\quad m_3\cong Z\oplus T\oplus W. \]
In $\mc B_A$, they have quantum dimensions (\ref{dBAm})
\[ d_{\mc B_A}(m_1)=d_{\mc B_A}(m_2)=d_{\mc B_A}(m_3)=1, \]
and obey the $\mbb Z/3\mbb Z$ monoidal products
\begin{table}[H]
\begin{center}
\begin{tabular}{c|c|c|c}
    $\otimes_A$&$m_1$&$m_2$&$m_3$\\\hline
    $m_1$&$m_1$&$m_2$&$m_3$\\\hline
    $m_2$&&$m_3$&$m_1$\\\hline
    $m_3$&&&$m_2$
\end{tabular}.
\end{center}
\end{table}
\hspace{-17pt}This shows $\mc B_A$ is a fusion category. Thus, $A$ is separable, showing the fact.

We conclude
\begin{table}[H]
\begin{center}
\begin{tabular}{c|c|c|c}
    Connected étale algebra $A$&$\mc B_A$&$\rank(\mc B_A)$&Lagrangian?\\\hline
    $1$&$\mc B$&$9$&No\\
    $1\oplus X\oplus Y$&$\vecG_{\mbb Z/3\mbb Z}^1$&3&Yes
\end{tabular}.
\end{center}
\caption{Connected étale algebras in rank nine MFC $\mcal B\simeq\vecG_{\mbb Z/9\mbb Z}^1$}\label{rank9Z9results}
\end{table}
\hspace{-17pt}All the four MFCs $\mc B\simeq\vecG_{\mbb Z/9\mbb Z}^1$'s fail to be completely anisotropic.

\subsubsection{$\mc B\simeq\vecG_{\mbb Z/3\mbb Z\times\mbb Z/3\mbb Z}^\alpha$}\label{Z3Z3}
The MFCs have nine simple objects $\{1,X,Y,Z,S,T,U,V,W\}$ obeying monoidal products
\begin{table}[H]
\begin{center}
\begin{tabular}{c|c|c|c|c|c|c|c|c|c}
    $\otimes$&$1$&$X$&$Y$&$Z$&$S$&$T$&$U$&$V$&$W$\\\hline
    $1$&$1$&$X$&$Y$&$Z$&$S$&$T$&$U$&$V$&$W$\\\hline
    $X$&&$Y$&$1$&$T$&$V$&$W$&$S$&$U$&$Z$\\\hline
    $Y$&&&$X$&$W$&$U$&$Z$&$V$&$S$&$T$\\\hline
    $Z$&&&&$S$&$1$&$V$&$Y$&$X$&$U$\\\hline
    $S$&&&&&$Z$&$X$&$W$&$T$&$Y$\\\hline
    $T$&&&&&&$U$&$1$&$Y$&$S$\\\hline
    $U$&&&&&&&$T$&$Z$&$X$\\\hline
    $V$&&&&&&&&$W$&$1$\\\hline
    $W$&&&&&&&&&$V$
\end{tabular}.
\end{center}
\end{table}
\hspace{-17pt}(We have $X^*\cong Y,Z^*\cong S,T^*\cong U,V^*\cong W$. One can identify $\vecG_{\mbb Z/3\mbb Z}^1=\{1,T,U\},\{1,V,W\}$, and $X\cong U\otimes W,Y\cong T\otimes V,Z\cong U\otimes V,S\cong T\otimes W$.) Thus, they have
\begin{align*}
    &\fp_{\mc B}(1)=\fp_{\mc B}(X)=\fp_{\mc B}(Y)=\fp_{\mc B}(Z)=\fp_{\mc B}(S)\\
    &=\fp_{\mc B}(T)=\fp_{\mc B}(U)=\fp_{\mc B}(V)=\fp_{\mc B}(W)=1,
\end{align*}
and
\[ \fp(\mc B)=9. \]
Their quantum dimensions $d_j$'s are solutions of the same multiplication rules $d_id_j=\sum_{k=1}^9{N_{ij}}^kd_k$. There is only one solution
\[ (d_X,d_Y,d_Z,d_S,d_T,d_U,d_V,d_W)=(1,1,1,1,1,1,1,1) \]
with categorical dimension
\[ D^2(\mc B)=9. \]
Thus, all MFCs are unitary. They have two conformal dimensions\footnote{Naively, one finds 18 consistent conformal dimensions, but the other 16 are equivalent to one of the two in the main text under permutations of simple objects.}
\[ (h_X,h_Y,h_Z,h_S,h_T,h_U,h_V,h_W)=(0,0,0,0,\frac13,\frac13,\frac23,\frac23),(\frac13,\frac13,\frac13,\frac13,\frac23,\frac23,\frac23,\frac23)\quad(\mods1). \]
The $S$-matrices are given by
\[ \hspace{-57pt}\widetilde S=\begin{pmatrix}1&d_X&d_Y&d_Xd_W&d_Yd_V&d_Yd_W&d_Xd_V&d_V&d_W\\d_X&e^{\mp2\pi i/3}&e^{\pm2\pi i/3}&e^{\mp2\pi i/3}d_W&e^{\pm2\pi i/3}d_V&e^{\pm2\pi i/3}d_W&e^{\mp2\pi i/3}d_V&d_Xd_V&d_Xd_W\\d_Y&e^{\pm2\pi i/3}&e^{\mp2\pi i/3}&e^{\pm2\pi i/3}d_W&e^{\mp2\pi i/3}d_V&e^{\mp2\pi i/3}d_W&e^{\pm2\pi i/3}d_V&d_Yd_V&d_Yd_W\\d_Xd_W&e^{\mp2\pi i/3}d_W&e^{\pm2\pi i/3}d_W&e^{\mp2\pi i/3\mp2\pi i/3}&e^{\pm2\pi i/3\pm2\pi i/3}&e^{\pm2\pi i/3\mp2\pi i/3}&e^{\mp2\pi i/3\pm2\pi i/3}&e^{\pm2\pi i/3}d_X&e^{\mp2\pi i/3}d_X\\d_Yd_V&e^{\pm2\pi i/3}d_V&e^{\mp2\pi i/3}d_V&e^{\pm2\pi i/3\pm2\pi i/3}&e^{\mp2\pi i/3\mp2\pi i/3}&e^{\mp2\pi i/3\pm2\pi i/3}&e^{\pm2\pi i/3\mp2\pi i/3}&e^{\mp2\pi i/3}d_Y&e^{\pm2\pi i/3}d_Y\\d_Yd_W&e^{\pm2\pi i/3}d_W&e^{\mp2\pi i/3}d_W&e^{\pm2\pi i/3\mp2\pi i/3}&e^{\mp2\pi i/3\pm2\pi i/3}&e^{\mp2\pi i/3\mp2\pi i/3}&e^{\pm2\pi i/3\pm2\pi i/3}&e^{\pm2\pi i/3}d_Y&e^{\mp2\pi i/3}d_Y\\d_Xd_V&e^{\mp2\pi i/3}d_V&e^{\pm2\pi i/3}d_V&e^{\mp2\pi i/3\pm2\pi i/3}&e^{\pm2\pi i/3\mp2\pi i/3}&e^{\pm2\pi i/3\pm2\pi i/3}&e^{\mp2\pi i/3\mp2\pi i/3}&e^{\mp2\pi i/3}d_X&e^{\pm2\pi i/3}d_X\\d_V&d_Xd_V&d_Yd_V&e^{\pm2\pi i/3}d_X&e^{\mp2\pi i/3}d_Y&e^{\pm2\pi i/3}d_Y&e^{\mp2\pi i/3}d_X&e^{\mp2\pi i/3}&e^{\pm2\pi i/3}\\d_W&d_Xd_W&d_Yd_W&e^{\mp2\pi i/3}d_X&e^{\pm2\pi i/3}d_Y&e^{\mp2\pi i/3}d_Y&e^{\pm2\pi i/3}d_X&e^{\pm2\pi i/3}&e^{\mp2\pi i/3}\end{pmatrix}. \]
They have additive central charges
\[ c(\mc B)=c(\vecG_{\mbb Z/3\mbb Z}^1)+c(\vecG_{\mbb Z/3\mbb Z}^1)\quad(\mods8) \]
where
\[ c(\vecG_{\mbb Z/3\mbb Z}^1)=\begin{cases}2&(h_{\mbb Z/3\mbb Z}=\frac13),\\-2&(h_{\mbb Z/3\mbb Z}=\frac23).\end{cases}\quad(\mods8) \]
There are
\[ 1(\text{quantum dimension})\times2(\text{conformal dimensions})\times2(\text{categorical dimensions})=4 \]
MFCs, and all of them are unitary. We classify connected étale algebras in them.

An ansatz
\[ A\cong1\oplus n_XX\oplus n_YY\oplus n_ZZ\oplus n_SS\oplus n_TT\oplus n_UU\oplus n_VV\oplus n_WW \]
with $n_j\in\mbb N$ has
\[ \fp_{\mc B}(A)=1+n_X+n_Y+n_Z+n_S+n_T+n_U+n_V+n_W. \]
For this to obey (\ref{FPdimA2bound}), the natural numbers can take only 45 values just as in the previous example. The sets contain the one with all $n_j$'s be zero. It is the trivial connected étale algebra $A\cong1$ giving $\mc B_A^0\simeq\mc B_A\simeq\mc B$. The others contain nontrivial simple object(s). Those 30 with $T,U,V,W$ fail to be commutative because they have nontrivial conformal dimensions. Together with the self-duality \cite{DMNO10} of étale algebras $A^*\cong A$, we are left with two nontrivial candidates
\[ A\cong1\oplus X\oplus Y,1\oplus Z\oplus S. \]
Since $\{1,X,Y\}$ or $\{1,Z,S\}$ form $\mbb Z/3\mbb Z$ pre-modular fusion subcategories, they were studied in \cite{KK23preMFC}. As a result, the candidates fail to be commutative for the second conformal dimension, while they are commutative for the first conformal dimension.

Let us check their separability by identifying $\mc B_A$. Since the commutative algebras have $\fp_{\mc B}(A)=3$, they demand
\[ \fp(\mc B_A^0)=1,\quad\fp(\mc B_A)=3. \]
The MFC $\mc B_A^0$ is identified as
\[ \mc B_A^0\simeq\vect. \]
This identification also matches central charges. The category $\mc B_A$ of right $A$-modules is identified as
\[ \mc B_A\simeq\vecG_{\mbb Z/3\mbb Z}^1. \]
To show this, we search for NIM-reps. We start from $A\cong1\oplus X\oplus Y$. We find a three-dimensional NIM-rep
\[ n_1=n_X=n_Y=1_3,\quad n_Z=n_T=n_W=\begin{pmatrix}0&0&1\\1&0&0\\0&1&0\end{pmatrix},\quad n_S=n_U=n_V=\begin{pmatrix}0&1&0\\0&0&1\\1&0&0\end{pmatrix}. \]
Denoting a basis of $\mc B_A$ by $\{m_1,m_2,m_3\}$, we obtain a multiplication table
\begin{table}[H]
\begin{center}
\begin{tabular}{c|c|c|c}
    $b_j\otimes\backslash$&$m_1$&$m_2$&$m_3$\\\hline
    $1,X,Y$&$m_1$&$m_2$&$m_3$\\
    $Z,T,W$&$m_3$&$m_1$&$m_2$\\
    $S,U,V$&$m_2$&$m_3$&$m_1$
\end{tabular}.
\end{center}
\end{table}
\hspace{-17pt}In this basis, we can identify
\[ m_1\cong1\oplus X\oplus Y,\quad m_2\cong S\oplus U\oplus V,\quad m_3\cong Z\oplus T\oplus W. \]
They have
\[ d_{\mc B_A}(m_1)=d_{\mc B_A}(m_2)=d_{\mc B_A}(m_3)=1. \]
Similarly, for the other commutative algebra $A\cong 1\oplus Z\oplus S$, just names of matrices change. We find a three-dimensional NIM-rep
\[ n_1=n_Z=n_S=1_3,\quad n_X=n_T=n_V=\begin{pmatrix}0&0&1\\1&0&0\\0&1&0\end{pmatrix},\quad n_Y=n_U=n_W=\begin{pmatrix}0&1&0\\0&0&1\\1&0&0\end{pmatrix}. \]
We obtain identifications
\[ m_1\cong1\oplus Z\oplus S,\quad m_2\cong X\oplus T\oplus V,\quad m_3\cong Y\oplus U\oplus W. \]
They both obey the $\mbb Z/3\mbb Z$ monoidal products
\begin{table}[H]
\begin{center}
\begin{tabular}{c|c|c|c}
    $\otimes_A$&$m_1$&$m_2$&$m_3$\\\hline
    $m_1$&$m_1$&$m_2$&$m_3$\\\hline
    $m_2$&&$m_3$&$m_1$\\\hline
    $m_3$&&&$m_2$
\end{tabular}.
\end{center}
\end{table}
\hspace{-17pt}This shows $\mc B_A\simeq\vecG_{\mbb Z/3\mbb Z}^1$ for the two commutative algebras. Since $\mc B_A$ is semisimple, this means the $A$'s are separable, hence étale.

We conclude
\begin{table}[H]
\begin{center}
\begin{tabular}{c|c|c|c}
    Connected étale algebra $A$&$\mc B_A$&$\rank(\mc B_A)$&Lagrangian?\\\hline
    $1$&$\mc B$&$9$&No\\
    $1\oplus X\oplus Y$ for the 1st $h$&$\vecG_{\mbb Z/3\mbb Z}^1$&3&Yes\\
    $1\oplus Z\oplus S$ for the 1st $h$&$\vecG_{\mbb Z/3\mbb Z}^1$&3&Yes
\end{tabular}.
\end{center}
\caption{Connected étale algebras in rank nine MFC $\mcal B\simeq\vecG_{\mbb Z/3\mbb Z\times\mbb Z/3\mbb Z}^\alpha$}\label{rank9Z3Z3results}
\end{table}
\hspace{-17pt}That is, two MFCs $\mc B\simeq\vecG_{\mbb Z/3\mbb Z\times\mbb Z/3\mbb Z}^\alpha$'s with the second conformal dimensions are completely anisotropic, while the other two are not.

\subsubsection{$\mc B\simeq\vecG_{\mbb Z/3\mbb Z}^1\boxtimes\ising$}
The MFCs have nine simple objects $\{1,X,Y,Z,S,T,U,V,W\}$ obeying monoidal products
\begin{table}[H]
\begin{center}
\begin{tabular}{c|c|c|c|c|c|c|c|c|c}
    $\otimes$&$1$&$X$&$Y$&$Z$&$S$&$T$&$U$&$V$&$W$\\\hline
    $1$&$1$&$X$&$Y$&$Z$&$S$&$T$&$U$&$V$&$W$\\\hline
    $X$&&$1$&$T$&$S$&$Z$&$Y$&$U$&$V$&$W$\\\hline
    $Y$&&&$Z$&$1$&$X$&$S$&$V$&$W$&$U$\\\hline
    $Z$&&&&$Y$&$T$&$X$&$W$&$U$&$V$\\\hline
    $S$&&&&&$Y$&$1$&$W$&$U$&$V$\\\hline
    $T$&&&&&&$Z$&$V$&$W$&$U$\\\hline
    $U$&&&&&&&$1\oplus X$&$Y\oplus T$&$Z\oplus S$\\\hline
    $V$&&&&&&&&$Z\oplus S$&$1\oplus X$\\\hline
    $W$&&&&&&&&&$Y\oplus T$
\end{tabular}.
\end{center}
\end{table}
\hspace{-17pt}(One can identify $\vecG_{\mbb Z/3\mbb Z}^1=\{1,Y,Z\},\ising=\{1,X,U\}$, and $S\cong X\otimes Z,T\cong X\otimes Y,V\cong Y\otimes U,W\cong Z\otimes U$.) Thus, they have
\begin{align*}
    \fp_{\mc B}(1)=\fp_{\mc B}(X)&=\fp_{\mc B}(Y)=\fp_{\mc B}(Z)=\fp_{\mc B}(S)=\fp_{\mc B}(T)=1,\\
    &\fp_{\mc B}(U)=\fp_{\mc B}(V)=\fp_{\mc B}(W)=\sqrt2,
\end{align*}
and
\[ \fp(\mc B)=12. \]
Their quantum dimensions $d_j$'s are solutions of the same multiplication rules $d_id_j=\sum_{k=1}^9{N_{ij}}^kd_k$. There are two (nonzero) solutions
\[ (d_X,d_Y,d_Z,d_S,d_T,d_U,d_V,d_W)=(1,1,1,1,1,-\sqrt2,-\sqrt2,-\sqrt2),(1,1,1,1,1,\sqrt2,\sqrt2,\sqrt2) \]
with categorical dimension
\[ D^2(\mc B)=12. \]
Only the second quantum dimension gives unitary MFCs. They both have the same 16 conformal dimensions
\begin{align*}
    (h_X,h_Y,h_Z,h_S,h_T,h_U,h_V,h_W)&=(\frac12,\frac13,\frac13,\frac56,\frac56,\frac1{16},\frac{19}{48},\frac{19}{48}),(\frac12,\frac13,\frac13,\frac56,\frac56,\frac3{16},\frac{25}{48},\frac{25}{48}),\\
    &~~~~(\frac12,\frac13,\frac13,\frac56,\frac56,\frac5{16},\frac{31}{48},\frac{31}{48}),(\frac12,\frac13,\frac13,\frac56,\frac56,\frac7{16},\frac{37}{48},\frac{37}{48}),\\
    &~~~~(\frac12,\frac13,\frac13,\frac56,\frac56,\frac9{16},\frac{43}{48},\frac{43}{48}),(\frac12,\frac13,\frac13,\frac56,\frac56,\frac{11}{16},\frac{1}{48},\frac{1}{48}),\\
    &~~~~(\frac12,\frac13,\frac13,\frac56,\frac56,\frac{13}{16},\frac{7}{48},\frac{7}{48}),(\frac12,\frac13,\frac13,\frac56,\frac56,\frac{15}{16},\frac{13}{48},\frac{13}{48}),\\
    &~~~~(\frac12,\frac23,\frac23,\frac16,\frac16,\frac1{16},\frac{35}{48},\frac{35}{48}),(\frac12,\frac23,\frac23,\frac16,\frac16,\frac3{16},\frac{41}{48},\frac{41}{48}),\\
    &~~~~(\frac12,\frac23,\frac23,\frac16,\frac16,\frac5{16},\frac{47}{48},\frac{47}{48}),(\frac12,\frac23,\frac23,\frac16,\frac16,\frac7{16},\frac{5}{48},\frac{5}{48}),\\
    &~~~~(\frac12,\frac23,\frac23,\frac16,\frac16,\frac9{16},\frac{11}{48},\frac{11}{48}),(\frac12,\frac23,\frac23,\frac16,\frac16,\frac{11}{16},\frac{17}{48},\frac{17}{48}),\\
    &~~~~(\frac12,\frac23,\frac23,\frac16,\frac16,\frac{13}{16},\frac{23}{48},\frac{23}{48}),(\frac12,\frac23,\frac23,\frac16,\frac16,\frac{15}{16},\frac{29}{48},\frac{29}{48})\quad(\mods1).
\end{align*}
The $S$-matrices are given by
\[ \hspace{-30pt}\widetilde S=\begin{pmatrix}1&d_X&d_Y&d_Z&d_Xd_Z&d_Xd_Y&d_U&d_Yd_U&d_Zd_U\\d_X&1&d_Xd_Y&d_Xd_Z&d_Z&d_Y&-d_U&-d_Y&-d_Z\\d_Y&d_Xd_Y&e^{\mp2\pi i/3}&e^{\pm2\pi i/3}&e^{\pm2\pi i/3} d_X&e^{\mp2\pi i/3}d_X&d_Yd_U&e^{\mp2\pi i/3}d_U&e^{\pm2\pi i/3}d_U\\d_Z&d_Xd_Z&e^{\pm2\pi i/3}&e^{\mp2\pi i/3}&e^{\mp2\pi i/3}d_X&e^{\pm2\pi i/3}d_X&d_Zd_U&e^{\pm2\pi i/3}d_U&e^{\mp2\pi i/3}d_U\\d_Xd_Z&d_Z&e^{\pm2\pi i/3}d_X&e^{\mp2\pi i/3}d_X&e^{\mp2\pi i/3}&e^{\pm2\pi i/3}&-d_Zd_U&-e^{\pm2\pi i/3}d_U&-e^{\mp2\pi i/3}d_U\\d_Xd_Y&d_Y&e^{\mp2\pi i/3}d_X&e^{\pm2\pi i/3}d_X&e^{\pm2\pi i/3}&e^{\mp2\pi i/3}&-d_Yd_U&-e^{\mp2\pi i/3}d_U&-e^{\pm2\pi i/3}d_U\\d_U&-d_U&d_Yd_U&d_Zd_U&-d_Zd_U&-d_Yd_U&0&0&0\\d_Yd_U&-d_Y&e^{\mp2\pi i/3}d_U&e^{\pm2\pi i/3}d_U&-e^{\pm2\pi i/3}d_U&-e^{\mp2\pi i/3}d_U&0&0&0\\d_Zd_U&-d_Z&e^{\pm2\pi i/3}d_U&e^{\mp2\pi i/3}d_U&-e^{\mp2\pi i/3}d_U&-e^{\pm2\pi i/3}d_U&0&0&0\end{pmatrix}. \]
(All signs are correlated. In other words, one $S$-matrix is given by choosing upper signs, and the other is its complex conjugate.) There are
\[ 2(\text{quantum dimensions})\times16(\text{conformal dimensions})\times2(\text{categorical dimensions})=64 \]
MFCs, among which those 32 with the second quantum dimensions are unitary. We classify connected étale algebras in all 64 MFCs simultaneously.

An ansatz
\[ A\cong1\oplus n_XX\oplus n_YY\oplus n_ZZ\oplus n_SS\oplus n_TT\oplus n_UU\oplus n_VV\oplus n_WW \]
with $n_j\in\mbb N$ has
\[ \fp_{\mc B}(A)=1+n_X+n_Y+n_Z+n_S+n_T+\sqrt2(n_U+n_V+n_W). \]
For this to obey (\ref{FPdimA2bound}), the natural numbers can take only 39 values. The sets contain the one with all $n_j$'s be zero. It is the trivial connected étale algebra $A\cong1$ giving $\mc B_A^0\simeq\mc B_A\simeq\mc B$. The other 38 candidates contain nontrivial simple object(s) with nontrivial conformal dimensions, and they all fail to be commutative. (Those with just $X$'s do satisfy the necessary condition (\ref{commutativealgnecessary}), but it has $c_{X,X}\cong-id_1$ \cite{KK23preMFC} and fail to be commutative.)

We conclude
\begin{table}[H]
\begin{center}
\begin{tabular}{c|c|c|c}
    Connected étale algebra $A$&$\mc B_A$&$\rank(\mc B_A)$&Lagrangian?\\\hline
    $1$&$\mc B$&$9$&No
\end{tabular}.
\end{center}
\caption{Connected étale algebras in rank nine MFC $\mcal B\simeq\vecG_{\mbb Z/3\mbb Z}^1\boxtimes\ising$}\label{rank9Z3isingresults}
\end{table}
\hspace{-17pt}All the 64 MFCs $\mc B\simeq\vecG_{\mbb Z/3\mbb Z}^1\boxtimes\ising$'s are completely anisotropic.

\subsubsection{$\mc B\simeq\ising\boxtimes\ising$}
The MFCs have nine simple objects $\{1,X,Y,Z,S,T,U,V,W\}$ obeying monoidal products
\begin{table}[H]
\begin{center}
\begin{tabular}{c|c|c|c|c|c|c|c|c|c}
    $\otimes$&$1$&$X$&$Y$&$Z$&$S$&$T$&$U$&$V$&$W$\\\hline
    $1$&$1$&$X$&$Y$&$Z$&$S$&$T$&$U$&$V$&$W$\\\hline
    $X$&&$1$&$Z$&$Y$&$V$&$T$&$U$&$S$&$W$\\\hline
    $Y$&&&$1$&$X$&$V$&$U$&$T$&$S$&$W$\\\hline
    $Z$&&&&$1$&$S$&$U$&$T$&$V$&$W$\\\hline
    $S$&&&&&$1\oplus Z$&$W$&$W$&$X\oplus Y$&$T\oplus U$\\\hline
    $T$&&&&&&$1\oplus X$&$Y\oplus Z$&$W$&$S\oplus V$\\\hline
    $U$&&&&&&&$1\oplus X$&$W$&$S\oplus V$\\\hline
    $V$&&&&&&&&$1\oplus Z$&$T\oplus U$\\\hline
    $W$&&&&&&&&&$1\oplus X\oplus Y\oplus Z$
\end{tabular}.
\end{center}
\end{table}
\hspace{-17pt}(One can identify $\ising=\{1,X,T\},\{1,Z,S\}$, and $Y\cong X\otimes Z, U\cong Z\otimes T,V\cong X\otimes S,W\cong S\otimes T$.) Thus, they have
\begin{align*}
    &\fp_{\mc B}(1)=\fp_{\mc B}(X)=\fp_{\mc B}(Y)=\fp_{\mc B}(Z)=1,\\
    \fp_{\mc B}(S)=&\fp_{\mc B}(T)=\fp_{\mc B}(U)=\fp_{\mc B}(V)=\sqrt2,\quad\fp_{\mc B}(W)=2,
\end{align*}
and
\[ \fp(\mc B)=16. \]
Their quantum dimensions $d_j$'s are solutions of the same multiplication rules $d_id_j=\sum_{k=1}^9{N_{ij}}^kd_k$. There are four (nonzero) solutions
\begin{align*}
    (d_X,d_Y,d_Z,d_S,d_T,d_U,d_V,d_W)&=(1,1,1,-\sqrt2,-\sqrt2,-\sqrt2,-\sqrt2,2),(1,1,1,-\sqrt2,\sqrt2,\sqrt2,-\sqrt2,-2),\\
    &~~~~(1,1,1,\sqrt2,-\sqrt2,-\sqrt2,\sqrt2,-2),(1,1,1,\sqrt2,\sqrt2,\sqrt2,\sqrt2,2)
\end{align*}
with categorical dimension
\[ D^2(\mc B)=16. \]
Only the last quantum dimension gives unitary MFCs. They have many conformal dimensions. In order to avoid double-counting, we perform case analysis.

\paragraph{$(d_S,d_T)=(\sqrt2,\sqrt2)$.} These give unitary MFCs. Different MFCs are labeled by 20 conformal dimensions\footnote{Naively, one has 64 conformal dimensions, but two conformal dimensions related by permutations $(SV),(TU),(XZ)(STVU),(XZ)(SUVT)$ give the same MFC. By studying some examples, one finds conformal dimensions are determined by two pairs $(h_S,h_V),(h_T,h_U)$. In the 64 sets, one conformal dimension $(h_S,h_T,h_U,h_V)$ appears in eight different orders. The number eight is given by
\[ 8=2(\text{1st pair in/outside})\times2(\text{order within 1st pair})\times2(\text{order within 2nd pair}). \]
If two combinations give the same $h_W$, they give the same MFC. Then, the 20 conformal dimensions are given by 10 pairs
\begin{align*}
    &(\frac1{16},\frac9{16})(\frac1{16},\frac9{16}),(\frac1{16},\frac9{16})(\frac3{16},\frac{11}{16}),(\frac1{16},\frac9{16})(\frac5{16},\frac{13}{16}),(\frac1{16},\frac9{16})(\frac7{16},\frac{15}{16}),(\frac3{16},\frac{11}{16})(\frac3{16},\frac{11}{16}),\\
    &(\frac3{16},\frac{11}{16})(\frac5{16},\frac{13}{16}),(\frac3{16},\frac{11}{16})(\frac7{16},\frac{15}{16}),(\frac5{16},\frac{13}{16})(\frac5{16},\frac{13}{16}),(\frac5{16},\frac{13}{16})(\frac7{16},\frac{15}{16}),(\frac7{16},\frac{15}{16})(\frac7{16},\frac{15}{16}),
\end{align*}
each pair giving two different $h_W$'s and hence two different MFCs.}
\begin{align*}
    \hspace{-10pt}(h_X,h_Y,h_Z,h_S,h_T,h_U,h_V,h_W)&=(\frac12,0,\frac12,\frac1{16},\frac1{16},\frac9{16},\frac9{16},\frac18),(\frac12,0,\frac12,\frac1{16},\frac3{16},\frac{11}{16},\frac9{16},\frac14),\\
    &~~~~(\frac12,0,\frac12,\frac1{16},\frac5{16},\frac{13}{16},\frac9{16},\frac38),(\frac12,0,\frac12,\frac1{16},\frac7{16},\frac{15}{16},\frac9{16},\frac12),\\
    &~~~~(\frac12,0,\frac12,\frac1{16},\frac9{16},\frac1{16},\frac9{16},\frac58),(\frac12,0,\frac12,\frac1{16},\frac{11}{16},\frac3{16},\frac9{16},\frac34),\\
    &~~~~(\frac12,0,\frac12,\frac1{16},\frac{13}{16},\frac5{16},\frac9{16},\frac78),(\frac12,0,\frac12,\frac1{16},\frac{15}{16},\frac7{16},\frac9{16},0),\\
    &~~~~(\frac12,0,\frac12,\frac3{16},\frac3{16},\frac{11}{16},\frac{11}{16},\frac38),(\frac12,0,\frac12,\frac3{16},\frac5{16},\frac{13}{16},\frac{11}{16},\frac12),\\
    &~~~~(\frac12,0,\frac12,\frac3{16},\frac7{16},\frac{15}{16},\frac{11}{16},\frac58),(\frac12,0,\frac12,\frac3{16},\frac{11}{16},\frac3{16},\frac{11}{16},\frac78),\\
    &~~~~(\frac12,0,\frac12,\frac3{16},\frac{13}{16},\frac5{16},\frac{11}{16},0),(\frac12,0,\frac12,\frac3{16},\frac{15}{16},\frac7{16},\frac{11}{16},\frac18),\\
    &~~~~(\frac12,0,\frac12,\frac5{16},\frac5{16},\frac{13}{16},\frac{13}{16},\frac58),(\frac12,0,\frac12,\frac5{16},\frac7{16},\frac{15}{16},\frac{13}{16},\frac34),\\
    &~~~~(\frac12,0,\frac12,\frac5{16},\frac{13}{16},\frac5{16},\frac{13}{16},\frac18),(\frac12,0,\frac12,\frac5{16},\frac{15}{16},\frac7{16},\frac{13}{16},\frac14),\\
    &~~~~(\frac12,0,\frac12,\frac7{16},\frac7{16},\frac{15}{16},\frac{15}{16},\frac78),(\frac12,0,\frac12,\frac7{16},\frac{15}{16},\frac7{16},\frac{15}{16},\frac38)\quad(\mods1).
\end{align*}
Including the two signs of categorical dimensions, we have 40 unitary MFCs.

\paragraph{$(d_S,d_T)=(\sqrt2,-\sqrt2)$.} In this case, we have smaller symmetries; just exchange of orders in two outer and inner pairs, $(h_S,h_V),(h_T,h_U)$. Different MFCs are given by 32 conformal dimensions
\begin{align*}
    \hspace{-10pt}(h_X,h_Y,h_Z,h_S,h_T,h_U,h_V,h_W)&=(\frac12,0,\frac12,\frac1{16},\frac1{16},\frac9{16},\frac9{16},\frac18),(\frac12,0,\frac12,\frac1{16},\frac3{16},\frac{11}{16},\frac9{16},\frac14),\\
    &~~~~(\frac12,0,\frac12,\frac1{16},\frac5{16},\frac{13}{16},\frac9{16},\frac38),(\frac12,0,\frac12,\frac1{16},\frac7{16},\frac{15}{16},\frac9{16},\frac12),\\
    &~~~~(\frac12,0,\frac12,\frac1{16},\frac9{16},\frac1{16},\frac9{16},\frac58),(\frac12,0,\frac12,\frac1{16},\frac{11}{16},\frac3{16},\frac9{16},\frac34),\\
    &~~~~(\frac12,0,\frac12,\frac1{16},\frac{13}{16},\frac5{16},\frac9{16},\frac78),(\frac12,0,\frac12,\frac1{16},\frac{15}{16},\frac7{16},\frac9{16},0),\\
    &~~~~(\frac12,0,\frac12,\frac3{16},\frac1{16},\frac9{16},\frac{11}{16},\frac14),(\frac12,0,\frac12,\frac3{16},\frac3{16},\frac{11}{16},\frac{11}{16},\frac38),\\
    &~~~~(\frac12,0,\frac12,\frac3{16},\frac5{16},\frac{13}{16},\frac{11}{16},\frac12),(\frac12,0,\frac12,\frac3{16},\frac7{16},\frac{15}{16},\frac{11}{16},\frac58),\\
    &~~~~(\frac12,0,\frac12,\frac3{16},\frac9{16},\frac1{16},\frac{11}{16},\frac34),(\frac12,0,\frac12,\frac3{16},\frac{11}{16},\frac3{16},\frac{11}{16},\frac78),\\
    &~~~~(\frac12,0,\frac12,\frac3{16},\frac{13}{16},\frac5{16},\frac{11}{16},0),(\frac12,0,\frac12,\frac3{16},\frac{15}{16},\frac7{16},\frac{11}{16},\frac18),\\
    &~~~~(\frac12,0,\frac12,\frac5{16},\frac1{16},\frac9{16},\frac{13}{16},\frac38),(\frac12,0,\frac12,\frac5{16},\frac3{16},\frac{11}{16},\frac{13}{16},\frac12),\\
    &~~~~(\frac12,0,\frac12,\frac5{16},\frac5{16},\frac{13}{16},\frac{13}{16},\frac58),(\frac12,0,\frac12,\frac5{16},\frac7{16},\frac{15}{16},\frac{13}{16},\frac34),\\
    &~~~~(\frac12,0,\frac12,\frac5{16},\frac9{16},\frac1{16},\frac{13}{16},\frac78),(\frac12,0,\frac12,\frac5{16},\frac{11}{16},\frac3{16},\frac{13}{16},0),\\
    &~~~~(\frac12,0,\frac12,\frac5{16},\frac{13}{16},\frac5{16},\frac{13}{16},\frac18),(\frac12,0,\frac12,\frac5{16},\frac{15}{16},\frac7{16},\frac{13}{16},\frac14),\\
    &~~~~(\frac12,0,\frac12,\frac7{16},\frac1{16},\frac9{16},\frac{15}{16},\frac12),(\frac12,0,\frac12,\frac7{16},\frac3{16},\frac{11}{16},\frac{15}{16},\frac58),\\
    &~~~~(\frac12,0,\frac12,\frac7{16},\frac5{16},\frac{13}{16},\frac{15}{16},\frac34),(\frac12,0,\frac12,\frac7{16},\frac7{16},\frac{15}{16},\frac{15}{16},\frac78),\\
    &~~~~(\frac12,0,\frac12,\frac7{16},\frac9{16},\frac1{16},\frac{15}{16},0),(\frac12,0,\frac12,\frac7{16},\frac{11}{16},\frac3{16},\frac{15}{16},\frac18),\\
    &~~~~(\frac12,0,\frac12,\frac7{16},\frac{13}{16},\frac5{16},\frac{15}{16},\frac14),(\frac12,0,\frac12,\frac7{16},\frac{15}{16},\frac7{16},\frac{15}{16},\frac38)\quad(\mods1).
\end{align*}
With two signs of categorical dimensions, there are 64 MFCs.

\paragraph{$(d_S,d_T)=(-\sqrt2,-\sqrt2)$.} In this case, symmetries are the same as $(d_S,d_T)=(\sqrt2,\sqrt2)$. Thus, the same 20 conformal dimensions give 40 MFCs.\newline

Therefore, there are
\[ 40+64+40=144 \]
MFCs, among which those 40 in the first case are unitary. The $S$-matrices are given by
\[ \widetilde S=\begin{pmatrix}1&d_X&d_Xd_Z&d_Z&d_S&d_T&d_Zd_T&d_Xd_S&d_Sd_T\\d_X&1&d_Z&d_Xd_Z&d_Xd_S&-d_T&-d_Zd_T&d_S&-d_Sd_T\\d_Xd_Z&d_Z&1&d_X&-d_Xd_S&-d_Zd_T&-d_T&-d_S&d_Sd_T\\d_Z&d_Xd_Z&d_X&1&-d_S&d_Zd_T&d_T&-d_Xd_S&-d_Sd_T\\d_S&d_Xd_S&-d_Xd_S&-d_S&0&d_Sd_T&-d_Sd_T&0&0\\d_T&-d_T&-d_Zd_T&d_Zd_T&d_Sd_T&0&0&-d_Sd_T&0\\d_Zd_T&-d_Zd_T&-d_T&d_T&-d_Sd_T&0&0&d_Sd_T&0\\d_Xd_S&d_S&-d_S&-d_Xd_S&0&-d_Sd_T&d_Sd_T&0&0\\d_Sd_T&-d_Sd_T&d_Sd_T&-d_Sd_T&0&0&0&0&0\end{pmatrix}. \]
They have additive central charges
\[ c(\mc B)=c(\ising)+c(\ising)\quad(\mods8) \]
where
\[ c(\ising)=\begin{cases}\frac12&(h_{S,T}=\frac1{16}),\\\frac32&(h_{S,T}=\frac3{16}),\\\frac52&(h_{S,T}=\frac5{16}),\\\frac72&(h_{S,T}=\frac7{16}),\\-\frac72&(h_{S,T}=\frac9{16}),\\-\frac52&(h_{S,T}=\frac{11}{16}),\\-\frac32&(h_{S,T}=\frac{13}{16}),\\-\frac12&(h_{S,T}=\frac{15}{16}).\end{cases}\quad(\mods8) \]

In order to classify connected étale algebras, we set an ansatz
\[ A\cong1\oplus n_XX\oplus n_YY\oplus n_ZZ\oplus n_SS\oplus n_TT\oplus n_UU\oplus n_VV\oplus n_WW \]
with $n_j\in\mbb N$. It has
\[ \fp_{\mc B}(A)=1+n_X+n_Y+n_Z+\sqrt2(n_S+n_T+n_U+n_V)+2n_W. \]
For this to obey (\ref{FPdimA2bound}), the natural numbers can take only 50 values. The sets contain the one with all the $n_j$'s be zero. It is the trivial connected étale algebra $A\cong1$ giving $\mc B_A^0\simeq\mc B_A\simeq\mc B$. All other candidates with nontrivial simple objects except $Y$ or $W$ fail to be commutative because they have nontrivial conformal dimensions. Thus, we are left with nontrivial candidates just with $Y,W$. Setting $n_X,n_Z,n_S,n_T,n_U,n_V$ to zero, apart from the trivial one, we find five natural numbers
\[ (n_Y,n_W)=(1,0),(0,1),(2,0),(1,1),(3,0). \]
Some of them are ruled out by studying Frobenius-Perron dimensions. The two candidates $(n_Y,n_W)=(0,1),(2,0)$ have $\fp_{\mc B}(A)=3$, and demands $\fp(\mc B_A^0)=\frac{16}9$, but there is no MFC with such Frobenius-Perron dimension. Thus, the two candidates are ruled out. Therefore, we are left with three candidates
\[ (n_Y,n_W)=(1,0),(1,1),(3,0). \]
We study these one after another.

\paragraph{$A\cong1\oplus Y$.} This is a $\mbb Z/2\mbb Z$ algebra, and it is also commutative because it has $(d_Y,h_Y)=(1,0)$ (mod 1 for $h_Y$) and $c_{Y,Y}\cong id_1$ \cite{KK23preMFC}. It further turns out separable (hence étale). Let us check this point by identifying $\mc B_A$.

First, $\fp_{\mc B}(A)=2$ demands
\[ \fp(\mc B_A^0)=4,\quad\fp(\mc B_A)=8. \]
The MFC $\mc B_A^0$ is identified as
\begin{equation}
    \mc B_A^0\simeq\begin{cases}\vecG_{\mbb Z/2\mbb Z}^{-1}\boxtimes\vecG_{\mbb Z/2\mbb Z}^{-1}&(h_W=\frac{2n+1}4\text{ with }n\in\mbb N\quad(\mods1)),\\\tc&(h_W=0\quad(\mods\frac12)),\\\vecG_{\mbb Z/4\mbb Z}^\alpha&(h_W=\frac{2n+1}8\text{ with }n\in\mbb N\quad(\mods1)).\end{cases}\label{isingisingBA0}
\end{equation}
This can be seen in two ways. First, let us perform anyon condensation. This `identifies' $1,Y$, and hence pairs $(X,Z),(S,V),(T,U)$, and $W$ splits into two. Since $S$ and $V$ have different conformal dimensions, the resulting simple object in $\mc B_A$ is confined. Similarly, the simple object from $T,U$ is confined. The simple object coming from $X,Z$ is deconfined, and is an object of $\mc B_A^0$ with conformal dimension $\frac12$. Finally, $W$ splits into two invertible simple objects with conformal dimensions $h_W$ mod 1. Thus, if $h_W=\frac n8$ (mod 1) with an odd natural number $n$, the four simple objects form $\mbb Z/4\mbb Z$ MFC, while if $h_W=\frac n8$ (mod 1) with an even natural number $n$, the four simple objects form $\mbb Z/2\mbb Z\times\mbb Z/2\mbb Z$ MFC. More precisely, if $h_W=\frac n4$ (mod 1) with $n=1,3$, it is $\vecG_{\mbb Z/2\mbb Z}^{-1}\boxtimes\vecG_{\mbb Z/2\mbb Z}^{-1}$, and if $h_W=\frac n4$ (mod 1) with $n=0,2$, then it is $\tc$. This gives the identifications.

Another way to see the identification is to compute the (additive) central charge $c(\mc B)$. The Ising MFC has central charge $\frac{2n+1}2$ (mod 8) with an integer $n=-4,-3,\dots,3$. Let us denote the central charges of the two factors $\frac{2n_1+1}2,\frac{2n_2+1}2$, respectively. Then, the Deligne tensor product $\ising\boxtimes\ising$ has (additive) central charge
\[ c(\mc B)=n_1+n_2+1\quad(\mods8). \]
Therefore, the MFC $\mc B_A^0$ has (additive) central charge
\[ c(\mc B_A^0)=\begin{cases}0&(n_1+n_2\in2\mbb Z+1),\\1&(n_1+n_2\in2\mbb Z).\end{cases}.\quad(\mods2) \]
An MFC with $\fp=4$ and this (additive) central charge is either $\mbb Z/2\mbb Z\times\mbb Z/2\mbb Z$ MFC or $\mbb Z/4\mbb Z$ MFC, respectively.\footnote{They have additive central charges
\[ c(\vecG_{\mbb Z/2\mbb Z}^{-1}\boxtimes\vecG_{\mbb Z/2\mbb Z}^{-1})=\begin{cases}0&(h_x\neq h_y),\\2&(h_x=\frac14=h_y),\\-2&(h_x=\frac34=h_y),\end{cases}\quad(\mods8) \]
\[ c(\tc)=\begin{cases}0&(h_x,h_y,h_z)=(\frac12,0,0),\\4&(h_x,h_y,h_z)=(\frac12,\frac12,\frac12),\end{cases}\quad(\mods8) \]
or
\[ c(\vecG_{\mbb Z/4\mbb Z}^\alpha)=\begin{cases}1&(h_{\mbb Z/4\mbb Z}=\frac18),\\3&(h_{\mbb Z/4\mbb Z}=\frac38),\\-3&(h_{\mbb Z/4\mbb Z}=\frac58),\\-1&(h_{\mbb Z/4\mbb Z}=\frac78),\end{cases}\quad(\mods8) \]
respectively. Here, conformal dimensions mean those of two nontrivial simple objects in two factors of $\mbb Z/2\mbb Z$ MFC ($h_x,h_y$), those of three nontrivial simple objects ($\tc$), or those of $\mbb Z/4\mbb Z$ simple objects ($\mbb Z/4\mbb Z$ MFC).} (One may naively think $\ising$ MFCs also match the Frobenius-Perron dimensions, but they are ruled out because they cannot match the additive central charges.) More precise matching is given as follows.

\paragraph{Even $c(\mc B)$.} When $n_1+n_2\in2\mbb Z+1$, the ambient category has even additive central charge. It is matched as follows:
\[ c(\mc B)=\begin{cases}0&(\mc B_A^0\simeq\vecG_{\mbb Z/2\mbb Z}^{-1}\boxtimes\vecG_{\mbb Z/2\mbb Z}^{-1}\text{ with }h_x\neq h_y\text{, or }\tc\text{ with }(h_x,h_y,h_z)=(\frac12,0,0)),\\2&(\vecG_{\mbb Z/2\mbb Z}^{-1}\boxtimes\vecG_{\mbb Z/2\mbb Z}^{-1}\text{ with }h_x=\frac14=h_y),\\4&(\tc\text{ with }(h_x,h_y,h_z)=(\frac12,\frac12,\frac12)),\\6&(\vecG_{\mbb Z/2\mbb Z}^{-1}\boxtimes\vecG_{\mbb Z/2\mbb Z}^{-1}\text{ with }h_x=\frac34=h_y).\end{cases} \]

\paragraph{Odd $c(\mc B)$.} When $n_1+n_2\in2\mbb Z$, the central charge of the ambient category can only be matched by $\mbb Z/4\mbb Z$ MFCs:
\[ c(\mc B)=\begin{cases}1&(\vecG_{\mbb Z/4\mbb Z}^\alpha\text{ with }h_{\mbb Z/4\mbb Z}=\frac18),\\3&(\vecG_{\mbb Z/4\mbb Z}^\alpha\text{ with }h_{\mbb Z/4\mbb Z}=\frac38),\\5&(\vecG_{\mbb Z/4\mbb Z}^\alpha\text{ with }h_{\mbb Z/4\mbb Z}=\frac58),\\7&(\vecG_{\mbb Z/4\mbb Z}^\alpha\text{ with }h_{\mbb Z/4\mbb Z}=\frac78).\end{cases}\quad(\mods8) \]

Together with the invariance of topological twist (\ref{invtopologicaltwist}), we can uniquely identify $\mc B_A^0$.

Having specified the category $\mc B_A^0$ of dyslectic right $A$-modules, let us next figure out the category $\mc B_A$ of right $A$-modules. In view of anyon condensation, it has two more simple objects with Frobenius-Perron dimensions $\sqrt2$. Thus, we learn
\[ \fp(\mc B_A)=8,\quad\rank(\mc B_A)=6. \]
With this information, we search for six-dimensional NIM-reps. Indeed, we find a solution
\begin{align*}
    n_1=1_6=n_Y,&\quad n_X=\begin{pmatrix}0&1&0&0&0&0\\1&0&0&0&0&0\\0&0&1&0&0&0\\0&0&0&1&0&0\\0&0&0&0&0&1\\0&0&0&0&1&0\end{pmatrix}=n_Z,\quad n_S=\begin{pmatrix}0&0&1&0&0&0\\0&0&1&0&0&0\\1&1&0&0&0&0\\0&0&0&0&1&1\\0&0&0&1&0&0\\0&0&0&1&0&0\end{pmatrix}=n_V,\\
    &n_T=\begin{pmatrix}0&0&0&1&0&0\\0&0&0&1&0&0\\0&0&0&0&1&1\\1&1&0&0&0&0\\0&0&1&0&0&0\\0&0&1&0&0&0\end{pmatrix}=n_U,\quad n_W=\begin{pmatrix}0&0&0&0&1&1\\0&0&0&0&1&1\\0&0&0&2&0&0\\0&0&2&0&0&0\\1&1&0&0&0&0\\1&1&0&0&0&0\end{pmatrix}.
\end{align*}
The NIM-rep gives a multiplication table
\begin{table}[H]
\begin{center}
\begin{tabular}{c|c|c|c|c|c|c}
    $b_j\otimes\backslash$&$m_1$&$m_2$&$m_3$&$m_4$&$m_5$&$m_6$\\\hline
    $1,Y$&$m_1$&$m_2$&$m_3$&$m_4$&$m_5$&$m_6$\\
    $X,Z$&$m_2$&$m_1$&$m_3$&$m_4$&$m_6$&$m_5$\\
    $S,V$&$m_3$&$m_3$&$m_1\oplus m_2$&$m_5\oplus m_6$&$m_4$&$m_4$\\
    $T,U$&$m_4$&$m_4$&$m_5\oplus m_6$&$m_1\oplus m_2$&$m_3$&$m_3$\\
    $W$&$m_5\oplus m_6$&$m_5\oplus m_6$&$2m_4$&$2m_3$&$m_1\oplus m_2$&$m_1\oplus m_2$
\end{tabular}.
\end{center}
\end{table}
\hspace{-17pt}From this, we can identify
\[ m_1\cong1\oplus Y,\quad m_2\cong X\oplus Z,\quad m_3\cong S\oplus V,\quad m_4\cong T\oplus U,\quad m_5\cong W\cong m_6. \]
In $\mc B_A$, they have quantum dimensions
\[ d_{\mc B_A}(m_1)=1=d_{\mc B_A}(m_2),\quad d_{\mc B_A}(m_3)=\pm\sqrt2,\quad d_{\mc B_A}(m_4)=\pm\sqrt2,\quad d_{\mc B_A}(m_5)=\pm1=d_{\mc B_A}(m_6). \]
Working out the monoidal products $\otimes_A$, we find
\begin{table}[H]
\begin{center}
\begin{tabular}{c|c|c|c|c|c|c}
    $\otimes_A$&$m_1$&$m_2$&$m_3$&$m_4$&$m_5$&$m_6$\\\hline
    $m_1$&$m_1$&$m_2$&$m_3$&$m_4$&$m_5$&$m_6$\\\hline
    $m_2$&&$m_1$&$m_3$&$m_4$&$m_6$&$m_5$\\\hline
    $m_3$&&&$m_1\oplus m_2$&$m_5\oplus m_6$&$m_4$&$m_4$\\\hline
    $m_4$&&&&$m_1\oplus m_2$&$m_3$&$m_3$\\\hline
    $m_5$&&&&&$m_1$&$m_2$\\\hline
    $m_6$&&&&&&$m_1$
\end{tabular}
\end{center}
\end{table}\hspace{-17pt}for $\mc B_A^0\simeq\vecG_{\mbb Z/2\mbb Z\times\mbb Z/2\mbb Z}^\alpha$ (i.e., $h_W=0,\frac14,\frac12,\frac34$ mod 1), and
\begin{table}[H]
\begin{center}
\begin{tabular}{c|c|c|c|c|c|c}
    $\otimes_A$&$m_1$&$m_2$&$m_3$&$m_4$&$m_5$&$m_6$\\\hline
    $m_1$&$m_1$&$m_2$&$m_3$&$m_4$&$m_5$&$m_6$\\\hline
    $m_2$&&$m_1$&$m_3$&$m_4$&$m_6$&$m_5$\\\hline
    $m_3$&&&$m_1\oplus m_2$&$m_5\oplus m_6$&$m_4$&$m_4$\\\hline
    $m_4$&&&&$m_1\oplus m_2$&$m_3$&$m_3$\\\hline
    $m_5$&&&&&$m_2$&$m_1$\\\hline
    $m_6$&&&&&&$m_2$
\end{tabular}
\end{center}
\end{table}\hspace{-17pt}for $\mc B_A^0\simeq\vecG_{\mbb Z/4\mbb Z}^\alpha$ (i.e., $h_W=\frac18,\frac38,\frac58,\frac78$ mod 1). These give identifications\footnote{The identification with the notation in AnyonWiki is given by
\[ m_1\cong1,\quad m_2\cong4,\quad m_3\cong5,\quad m_4\cong6,\quad m_5\cong2,\quad m_6\cong3, \]
or its permutations $(56),(23)$, and
\[ m_1\cong1,\quad m_2\cong2,\quad m_3\cong5,\quad m_4\cong6,\quad m_5\cong3,\quad m_6\cong4, \]
or its permutations $(34),(56)$, respectively.}
\[ \mc B_A\simeq\begin{cases}\mc C(\text{FR}^{6,0}_1)&(c(\mc B)=0\text{ mod }2),\\\mc C(\text{FR}^{6,0}_2)&(c(\mc B)=1\text{ mod }2).\end{cases} \]
This shows $\mc B_A$ is semisimple, hence $A$ is separable and étale.

\paragraph{$A\cong1\oplus Y\oplus W$.} When $h_W=0$ (mod 1), $W$ can give commutative algebra. It has $\fp_{\mc B}(A)=4$, and demands
\[ \fp(\mc B_A^0)=1,\quad\fp(\mc B_A)=4. \]
Computing $b_j\otimes A$, we find candidate simple objects
\[ 1\oplus Y\oplus W,\quad X\oplus Z\oplus W,\quad S\oplus T\oplus U\oplus V. \]
They have Frobenius-Perron dimensions
\[ 1,\quad1,\quad\sqrt2, \]
and their contributions to $\fp(\mc B_A)$ match, $1^2+1^2+\sqrt2^2=4$. This suggests $\mc B_A$ have rank three. Indeed, we find a three-dimensional NIM-rep
\[ \hspace{-30pt}n_1=1_3=n_Y,\quad n_X=\begin{pmatrix}0&1&0\\1&0&0\\0&0&1\end{pmatrix}=n_Z,\quad n_S=n_T=n_U=n_V=\begin{pmatrix}0&0&1\\0&0&1\\1&1&0\end{pmatrix},\quad n_W=\begin{pmatrix}1&1&0\\1&1&0\\0&0&2\end{pmatrix}. \]
The solution gives identifications
\[ m_1\cong1\oplus Y\oplus W,\quad m_2\cong X\oplus Z\oplus W,\quad m_3\cong S\oplus T\oplus U\oplus V. \]
By computing quantum dimensions (\ref{dBAm}), we find the candidate can be separable only when $d_S=d_T$. Then, in $\mc B_A$, they have
\[ d_{\mc B_A}(m_1)=1=d_{\mc B_A}(m_2),\quad d_{\mc B_A}(m_3)=\pm\sqrt2. \]
Working out the monoidal products $\otimes_A$, we find
\begin{table}[H]
\begin{center}
\begin{tabular}{c|c|c|c}
    $\otimes_A$&$m_1$&$m_2$&$m_3$\\\hline
    $m_1$&$m_1$&$m_2$&$m_3$\\\hline
    $m_2$&&$m_1$&$m_3$\\\hline
    $m_3$&&&$m_1\oplus m_2$
\end{tabular}.
\end{center}
\end{table}
\hspace{-17pt}We arrive the identification
\[ \mc B_A\simeq\ising. \]
Since this is semisimple, $A$ is separable and étale. We found\footnote{The existence of the Lagrangian algebra can also be understood from the lemma 3 because when $h_W=0$ (mod 1), our ambient MFC is a Drinfeld center $\mc B\simeq Z(\ising)$.}
\begin{equation}
\begin{split}
    A\cong1\oplus Y\oplus W\quad(d_S,d_T,h_S,h_T)=&(\sqrt2,\sqrt2,\frac1{16},\frac{15}{16}),(\sqrt2,\sqrt2,\frac3{16},\frac{13}{16}),\\
    &(-\sqrt2,-\sqrt2,\frac1{16},\frac{15}{16}),(-\sqrt2,-\sqrt2,\frac3{16},\frac{13}{16}).
\end{split}\quad(\mods1\text{ for }h)\label{isingising1YWetale}
\end{equation}

\paragraph{$A\cong1\oplus3Y$.} It has $\fp_{\mc B}(A)=4$, and demands
\[ \fp(\mc B_A^0)=1,\quad\fp(\mc B_A)=4. \]
Calculating $b_j\otimes A$, we find candidate simple objects (assuming the smallest Frobenius-Perron dimensions)
\[ 1\oplus3Y,\quad X\oplus3Z,\quad3\oplus Y,\quad3X\oplus Z,\quad S\oplus3V,\quad T\oplus3U,\quad3T\oplus U,\quad3S\oplus V,\quad2W. \]
In a putative $\mc B_A$, they have Frobenius-Perron dimensions
\[ 1,\quad1,\quad1,\quad1,\quad\sqrt2,\quad\sqrt2,\quad\sqrt2,\quad\sqrt2,\quad1. \]
Their contributions to $\fp(\mc B_A)$ exceed four, and the candidate is ruled out.

We conclude
\begin{table}[H]
\begin{center}
\begin{tabular}{c|c|c|c}
    Connected étale algebra $A$&$\mc B_A$&$\rank(\mc B_A)$&Lagrangian?\\\hline
    $1$&$\mc B$&$9$&No\\
    $1\oplus Y$&$\begin{cases}\mc C(\text{FR}^{6,0}_1)&(c(\mc B)=0),\\\mc C(\text{FR}^{6,0}_2)&(c(\mc B)=1).\end{cases}$ mod 2&6&No\\
    $1\oplus Y\oplus W$ for (\ref{isingising1YWetale})&$\ising$&3&Yes
\end{tabular}.
\end{center}
\caption{Connected étale algebras in rank nine MFC $\mcal B\simeq\ising\boxtimes\ising$}\label{rank9isingisingresults}
\end{table}
\hspace{-17pt}All the 144 MFCs $\mc B\simeq\ising\boxtimes\ising$'s fail to be completely anisotropic.

Let us comment on the consistency of our results. The nontrivial connected étale algebra leads to the category of dyslelctic right $A$-modules, $\vecG_{\mbb Z/4\mbb Z}^\alpha,\vecG_{\mbb Z/2\mbb Z}^{-1}\boxtimes\vecG_{\mbb Z/2\mbb Z}^{-1}\ (\text{with }h_x=h_y)$, or $\tc$ with $(h_x,h_y,h_z)=(\frac12,\frac12,\frac12)$ or $(h_x,h_y,h_z)=(\frac12,0,0)$. The first three are completely anisotropic, and the ``cascades'' of conformal embeddings terminate there. On the other hand, for the last $\mc B_A^0$, it can admit nontrivial connected étale algebras. More precisely, the number of additional nontrivial connected étale algebras are given by
\[ \begin{cases}0&(d_W,h_W)=(-2,0),\\2&(d_W,h_W)=(2,0).\end{cases} \]
The additional connected étale algebras lead to $\vect$. One of the two paths is given by composition with itself, $1\oplus Y$. The other path
\[ \ising\boxtimes\ising\to\vect \]
is given by the additional connected étale algebra $1\oplus Y\oplus W$. This appears exactly when $(d_W,h_W)=(2,0)$. See (\ref{isingising1YWetale}).

\subsubsection{$\mc B\simeq\vecG_{\mbb Z/3\mbb Z}^1\boxtimes psu(2)_5$}
The MFCs have nine simple objects $\{1,X,Y,Z,S,T,U,V,W\}$ obeying monoidal products
\begin{table}[H]
\begin{center}
\begin{tabular}{c|c|c|c|c|c|c|c|c|c}
    $\otimes$&$1$&$X$&$Y$&$Z$&$S$&$T$&$U$&$V$&$W$\\\hline
    $1$&$1$&$X$&$Y$&$Z$&$S$&$T$&$U$&$V$&$W$\\\hline
    $X$&&$Y$&$1$&$T$&$Z$&$S$&$V$&$W$&$U$\\\hline
    $Y$&&&$X$&$S$&$T$&$Z$&$W$&$U$&$V$\\\hline
    $Z$&&&&$1\oplus U$&$Y\oplus W$&$X\oplus V$&$Z\oplus U$&$T\oplus V$&$S\oplus W$\\\hline
    $S$&&&&&$X\oplus V$&$1\oplus U$&$S\oplus W$&$Z\oplus U$&$T\oplus V$\\\hline
    $T$&&&&&&$Y\oplus W$&$T\oplus V$&$S\oplus W$&$Z\oplus U$\\\hline
    $U$&&&&&&&$1\oplus Z\oplus U$&$X\oplus T\oplus V$&$Y\oplus S\oplus W$\\\hline
    $V$&&&&&&&&$Y\oplus S\oplus W$&$1\oplus Z\oplus U$\\\hline
    $W$&&&&&&&&&$X\oplus T\oplus V$
\end{tabular}.
\end{center}
\end{table}
\hspace{-17pt}(One can identify $\vecG_{\mbb Z/3\mbb Z}^1=\{1,X,Y\},psu(2)_5=\{1,Z,U\}$, and $S\cong Y\otimes Z,T\cong X\otimes Z,V\cong X\otimes U,W\cong Y\otimes U$.) Thus, they have
\begin{align*}
    \fp_{\mc B}(1)=\fp_{\mc B}(X)&=\fp_{\mc B}(Y)=1,\quad\fp_{\mc B}(Z)=\fp_{\mc B}(S)=\fp_{\mc B}(T)=\frac{\sin\frac{2\pi}7}{\sin\frac\pi7},\\
    &\fp_{\mc B}(U)=\fp_{\mc B}(V)=\fp_{\mc B}(W)=\frac{\sin\frac{3\pi}7}{\sin\frac\pi7},
\end{align*}
and
\[ \fp(\mc B)=\frac{21}{4\sin^2\frac\pi7}\approx27.9. \]
Their quantum dimensions $d_j$'s are solutions of the same multiplication rules $d_id_j=\sum_{k=1}^9{N_{ij}}^kd_k$. There are three solutions
\begin{align*}
    (d_X,d_Y,d_Z,d_S,d_T,d_U,d_V,d_W)&=(1,1,\frac{\sin\frac\pi7}{\cos\frac\pi{14}},\frac{\sin\frac\pi7}{\cos\frac\pi{14}},\frac{\sin\frac\pi7}{\cos\frac\pi{14}},-\frac{\sin\frac{2\pi}7}{\cos\frac\pi{14}},-\frac{\sin\frac{2\pi}7}{\cos\frac\pi{14}},-\frac{\sin\frac{2\pi}7}{\cos\frac\pi{14}}),\\
    &~~~~(1,1,-\frac{\sin\frac{3\pi}7}{\cos\frac{3\pi}{14}},-\frac{\sin\frac{3\pi}7}{\cos\frac{3\pi}{14}},-\frac{\sin\frac{3\pi}7}{\cos\frac{3\pi}{14}},\frac{\sin\frac{\pi}7}{\cos\frac{3\pi}{14}},\frac{\sin\frac{\pi}7}{\cos\frac{3\pi}{14}},\frac{\sin\frac{\pi}7}{\cos\frac{3\pi}{14}}),\\
    &~~~~(1,1,\frac{\sin\frac{2\pi}7}{\sin\frac\pi7},\frac{\sin\frac{2\pi}7}{\sin\frac\pi7},\frac{\sin\frac{2\pi}7}{\sin\frac\pi7},\frac{\sin\frac{3\pi}7}{\sin\frac\pi7},\frac{\sin\frac{3\pi}7}{\sin\frac\pi7},\frac{\sin\frac{3\pi}7}{\sin\frac\pi7})
\end{align*}
with categorical dimensions
\[ D^2(\mc B)=\frac{21}{4\cos^2\frac\pi{14}}(\approx5.5),\quad\frac{21}{4\cos^2\frac{3\pi}{14}}(\approx8.6),\quad\frac{21}{4\sin^2\frac\pi7}, \]
respectively. They have four conformal dimensions each:
\begin{align*}
    \hspace{-30pt}&(h_X,h_Y,h_Z,h_S,h_T,h_U,h_V,h_W)\\
    \hspace{-30pt}&=\begin{cases}(\frac13,\frac13,\frac37,\frac{16}{21},\frac{16}{21},\frac17,\frac{10}{21},\frac{10}{21}),(\frac13,\frac13,\frac47,\frac{19}{21},\frac{19}{21},\frac67,\frac{4}{21},\frac{4}{21}),(\frac23,\frac23,\frac37,\frac{2}{21},\frac{2}{21},\frac17,\frac{17}{21},\frac{17}{21}),(\frac23,\frac23,\frac47,\frac{5}{21},\frac{5}{21},\frac67,\frac{11}{21},\frac{11}{21})&(\text{1st}),\\(\frac13,\frac13,\frac27,\frac{13}{21},\frac{13}{21},\frac37,\frac{16}{21},\frac{16}{21}),(\frac13,\frac13,\frac57,\frac{1}{21},\frac{1}{21},\frac47,\frac{19}{21},\frac{19}{21}),(\frac23,\frac23,\frac27,\frac{20}{21},\frac{20}{21},\frac37,\frac2{21},\frac2{21}),(\frac23,\frac23,\frac57,\frac{8}{21},\frac{8}{21},\frac47,\frac{5}{21},\frac{5}{21})&(\text{2nd}),\\(\frac13,\frac13,\frac17,\frac{10}{21},\frac{10}{21},\frac57,\frac1{21},\frac1{21}),(\frac13,\frac13,\frac67,\frac{4}{21},\frac{4}{21},\frac27,\frac{13}{21},\frac{13}{21}),(\frac23,\frac23,\frac17,\frac{17}{21},\frac{17}{21},\frac57,\frac8{21},\frac8{21}),(\frac23,\frac23,\frac67,\frac{11}{21},\frac{11}{21},\frac27,\frac{20}{21},\frac{20}{21})&(\text{3rd}).\end{cases}\\
    \hspace{-30pt}&\hspace{500pt}(\mods1)
\end{align*}
The $S$-matrices are given by
\[ \hspace{-30pt}\widetilde S=\begin{pmatrix}1&d_X&d_Y&d_Z&d_Yd_Z&d_Xd_Z&d_U&d_Xd_U&d_Yd_U\\d_X&e^{\mp2\pi i/3}&e^{\pm2\pi i/3}&d_Xd_Z&e^{\pm2\pi i/3}d_Z&e^{\mp2\pi i/3}d_Z&d_Xd_U&e^{\mp2\pi i/3}d_U&e^{\pm2\pi i/3}d_U\\d_Y&e^{\pm2\pi i/3}&e^{\mp2\pi i/3}&d_Yd_Z&e^{\mp2\pi i/3}d_Z&e^{\pm2\pi i/3}d_Z&d_Yd_U&e^{\pm2\pi i/3}d_U&e^{\mp2\pi i/3}d_U\\d_Z&d_Xd_Z&d_Yd_Z&-d_U&-d_Yd_U&-d_Xd_U&1&d_X&d_Y\\d_Yd_Z&e^{\pm2\pi i/3}d_Z&e^{\mp2\pi i/3}d_Z&-d_Yd_U&-e^{\mp2\pi i/3}d_U&-e^{\pm2\pi i/3}d_U&d_Y&e^{\pm2\pi i/3}&e^{\mp2\pi i/3}\\d_Xd_Z&e^{\mp2\pi i/3}d_Z&e^{\pm2\pi i/3}d_Z&-d_Xd_U&-e^{\pm2\pi i/3}d_U&-e^{\mp2\pi i/3}d_U&d_X&e^{\mp2\pi i/3}&e^{\pm2\pi i/3}\\d_U&d_Xd_U&d_Yd_U&1&d_Y&d_X&-d_Z&-d_Xd_Z&-d_Yd_Z\\d_Xd_U&e^{\mp2\pi i/3}d_U&e^{\pm2\pi i/3}d_U&d_X&e^{\pm2\pi i/3}&e^{\mp2\pi i/3}&-d_Xd_Z&-e^{\mp2\pi i/3}d_Z&-e^{\pm2\pi i/3}d_Z\\d_Yd_U&e^{\pm2\pi i/3}d_U&e^{\mp2\pi i/3}d_U&d_Y&e^{\mp2\pi i/3}&e^{\pm2\pi i/3}&-d_Yd_Z&-e^{\pm2\pi i/3}d_Z&-e^{\mp2\pi i/3}d_Z\end{pmatrix}. \]
(All signs are correlated. In other words, one $S$-matrix is given by choosing upper signs, and the other is its complex conjugate.) There are
\[ 3(\text{quantum dimensions})\times4(\text{conformal dimensions})\times2(\text{categorical dimensions})=24 \]
MFCs, among which those eight with the third quantum dimensions give unitary MFCs. We classify connected étale algebras in all 24 MFCs simultaneously.

An ansatz
\[ A\cong1\oplus n_XX\oplus n_YY\oplus n_ZZ\oplus n_SS\oplus n_TT\oplus n_UU\oplus n_VV\oplus n_WW \]
with $n_j\in\mbb N$ has
\[ \fp_{\mc B}(A)=1+n_X+n_Y+\frac{\sin\frac{2\pi}7}{\sin\frac\pi7}(n_Z+n_S+n_T)+\frac{\sin\frac{3\pi}7}{\sin\frac\pi7}(n_U+n_V+n_W). \]
For this to obey (\ref{FPdimA2bound}), the natural number can take only 66 values. The sets contain the one with all $n_j$'s be zero. It is the trivial connected étale algebra $A\cong1$ giving $\mc B_A^0\simeq\mc B_A\simeq\mc B$. The other 65 sets contain simple object(s) with nontrivial conformal dimensions, and they all fail to be commutative.

We conclude
\begin{table}[H]
\begin{center}
\begin{tabular}{c|c|c|c}
    Connected étale algebra $A$&$\mc B_A$&$\rank(\mc B_A)$&Lagrangian?\\\hline
    $1$&$\mc B$&$9$&No
\end{tabular}.
\end{center}
\caption{Connected étale algebras in rank nine MFC $\mcal B\simeq\vecG_{\mbb Z/3\mbb Z}^1\boxtimes psu(2)_5$}\label{rank9Z3psu25results}
\end{table}
\hspace{-17pt}All the 24 MFCs $\mc B\simeq\vecG_{\mbb Z/3\mbb Z}^1\boxtimes psu(2)_5$'s are completely anisotropic.

\subsubsection{$\mc B\simeq\ising\boxtimes psu(2)_5$}
The MFCs have nine simple objects $\{1,X,Y,Z,S,T,U,V,W\}$ obeying monoidal products
\begin{table}[H]
\begin{center}
\makebox[1 \textwidth][c]{       
\resizebox{1.2 \textwidth}{!}{\begin{tabular}{c|c|c|c|c|c|c|c|c|c}
    $\otimes$&$1$&$X$&$Y$&$Z$&$S$&$T$&$U$&$V$&$W$\\\hline
    $1$&$1$&$X$&$Y$&$Z$&$S$&$T$&$U$&$V$&$W$\\\hline
    $X$&&$1$&$Y$&$S$&$Z$&$U$&$T$&$V$&$W$\\\hline
    $Y$&&&$1\oplus X$&$V$&$V$&$W$&$W$&$Z\oplus S$&$T\oplus U$\\\hline
    $Z$&&&&$1\oplus T$&$X\oplus U$&$Z\oplus T$&$S\oplus U$&$Y\oplus W$&$V\oplus W$\\\hline
    $S$&&&&&$1\oplus T$&$S\oplus U$&$Z\oplus T$&$Y\oplus W$&$V\oplus W$\\\hline
    $T$&&&&&&$1\oplus Z\oplus T$&$X\oplus S\oplus U$&$V\oplus W$&$Y\oplus V\oplus W$\\\hline
    $U$&&&&&&&$1\oplus Z\oplus T$&$V\oplus W$&$Y\oplus V\oplus W$\\\hline
    $V$&&&&&&&&$1\oplus X\oplus T\oplus U$&$Z\oplus S\oplus T\oplus U$\\\hline
    $W$&&&&&&&&&$1\oplus X\oplus Z\oplus S\oplus T\oplus U$
\end{tabular}.}}
\end{center}
\end{table}
\hspace{-17pt}(One can identify $\ising=\{1,X,Y\},psu(2)_5=\{1,Z,T\}$, and $S\cong X\otimes Z,U\cong X\otimes T,V\cong Y\otimes Z,W\cong Y\otimes T$.) Thus, they have
\begin{align*}
    \fp_{\mc B}(1)=1=\fp_{\mc B}(X),\quad\fp_{\mc B}(Y)=\sqrt2,\quad\fp_{\mc B}(Z)=\frac{\sin\frac{2\pi}7}{\sin\frac\pi7}=\fp_{\mc B}(S),\\
    \fp_{\mc B}(T)=\frac{\sin\frac{3\pi}7}{\sin\frac\pi7}=\fp_{\mc B}(U),\quad\fp_{\mc B}(V)=\sqrt2\frac{\sin\frac{2\pi}7}{\sin\frac\pi7},\quad\fp_{\mc B}(W)=\sqrt2\frac{\sin\frac{3\pi}7}{\sin\frac\pi7},
\end{align*}
and
\[ \fp(\mc B)=\frac7{\sin^2\frac\pi7}\approx37.2. \]
Their quantum dimensions $d_j$'s are solutions of the same multiplication rules $d_id_j=\sum_{k=1}^9{N_{ij}}^kd_k$. There are six (nonzero) solutions
\begin{align*}
    (d_X,d_Y,d_Z,d_S,d_T,d_U,d_V,d_W)&=(1,-\sqrt2,\frac{\sin\frac\pi7}{\cos\frac\pi{14}},\frac{\sin\frac\pi7}{\cos\frac\pi{14}},-\frac{\sin\frac{2\pi}7}{\cos\frac\pi{14}},-\frac{\sin\frac{2\pi}7}{\cos\frac\pi{14}},-\frac{\sqrt2\sin\frac\pi7}{\cos\frac\pi{14}},\frac{\sqrt2\sin\frac{2\pi}7}{\cos\frac\pi{14}}),\\
    &~~~~(1,\sqrt2,\frac{\sin\frac\pi7}{\cos\frac\pi{14}},\frac{\sin\frac\pi7}{\cos\frac\pi{14}},-\frac{\sin\frac{2\pi}7}{\cos\frac\pi{14}},-\frac{\sin\frac{2\pi}7}{\cos\frac\pi{14}},\frac{\sqrt2\sin\frac\pi7}{\cos\frac\pi{14}},-\frac{\sqrt2\sin\frac{2\pi}7}{\cos\frac\pi{14}}),\\
    &~~~~(1,-\sqrt2,-\frac{\sin\frac{3\pi}7}{\cos\frac{3\pi}{14}},-\frac{\sin\frac{3\pi}7}{\cos\frac{3\pi}{14}},\frac{\sin\frac{\pi}7}{\cos\frac{3\pi}{14}},\frac{\sin\frac{\pi}7}{\cos\frac{3\pi}{14}},\frac{\sqrt2\sin\frac{3\pi}7}{\cos\frac{3\pi}{14}},-\frac{\sqrt2\sin\frac{\pi}7}{\cos\frac{3\pi}{14}}),\\
    &~~~~(1,\sqrt2,-\frac{\sin\frac{3\pi}7}{\cos\frac{3\pi}{14}},-\frac{\sin\frac{3\pi}7}{\cos\frac{3\pi}{14}},\frac{\sin\frac{\pi}7}{\cos\frac{3\pi}{14}},\frac{\sin\frac{\pi}7}{\cos\frac{3\pi}{14}},-\frac{\sqrt2\sin\frac{3\pi}7}{\cos\frac{3\pi}{14}},\frac{\sqrt2\sin\frac{\pi}7}{\cos\frac{3\pi}{14}}),\\
    &~~~~(1,-\sqrt2,\frac{\sin\frac{2\pi}7}{\sin\frac{\pi}{7}},\frac{\sin\frac{2\pi}7}{\sin\frac{\pi}{7}},\frac{\sin\frac{3\pi}7}{\sin\frac{\pi}{7}},\frac{\sin\frac{3\pi}7}{\sin\frac{\pi}{7}},-\frac{\sqrt2\sin\frac{2\pi}7}{\sin\frac{\pi}{7}},-\frac{\sqrt2\sin\frac{3\pi}7}{\cos\frac{\pi}{7}}),\\
    &~~~~(1,\sqrt2,\frac{\sin\frac{2\pi}7}{\sin\frac{\pi}{7}},\frac{\sin\frac{2\pi}7}{\sin\frac{\pi}{7}},\frac{\sin\frac{3\pi}7}{\sin\frac{\pi}{7}},\frac{\sin\frac{3\pi}7}{\sin\frac{\pi}{7}},\frac{\sqrt2\sin\frac{2\pi}7}{\sin\frac{\pi}{7}},\frac{\sqrt2\sin\frac{3\pi}7}{\cos\frac{\pi}{7}})
\end{align*}
with categorical dimensions
\[ D^2(\mc B)=\frac7{\cos^2\frac\pi{14}}(\approx7.4),\quad\frac7{\cos^2\frac{3\pi}{14}}(\approx11.5),\quad\frac7{\sin^2\frac\pi7}, \]
respectively for each pair. Each pair has 16 conformal dimensions
\[ \hspace{-50pt}(h_X,h_Y,h_Z,h_S,h_T,h_U,h_V,h_W)=\begin{cases}(\frac12,\frac1{16},\frac37,\frac{13}{14},\frac17,\frac9{14},\frac{55}{112},\frac{23}{112}),(\frac12,\frac1{16},\frac47,\frac1{14},\frac67,\frac5{14},\frac{71}{112},\frac{103}{112}),\\(\frac12,\frac3{16},\frac37,\frac{13}{14},\frac17,\frac9{14},\frac{69}{112},\frac{37}{112}),(\frac12,\frac3{16},\frac47,\frac1{14},\frac67,\frac5{14},\frac{85}{112},\frac5{112}),\\
    (\frac12,\frac5{16},\frac37,\frac{13}{14},\frac17,\frac9{14},\frac{83}{112},\frac{51}{112}),(\frac12,\frac5{16},\frac47,\frac1{14},\frac67,\frac5{14},\frac{99}{112},\frac{19}{112}),\\(\frac12,\frac7{16},\frac37,\frac{13}{14},\frac17,\frac9{14},\frac{97}{112},\frac{65}{112}),(\frac12,\frac7{16},\frac47,\frac1{14},\frac67,\frac5{14},\frac1{112},\frac{33}{112}),\\
    (\frac12,\frac9{16},\frac37,\frac{13}{14},\frac17,\frac9{14},\frac{111}{112},\frac{79}{112}),(\frac12,\frac9{16},\frac47,\frac1{14},\frac67,\frac5{14},\frac{15}{112},\frac{47}{112}),\\(\frac12,\frac{11}{16},\frac37,\frac{13}{14},\frac17,\frac9{14},\frac{13}{112},\frac{93}{112}),(\frac12,\frac{11}{16},\frac47,\frac1{14},\frac67,\frac5{14},\frac{29}{112},\frac{61}{112}),\\(\frac12,\frac{13}{16},\frac37,\frac{13}{14},\frac17,\frac9{14},\frac{27}{112},\frac{107}{112}),(\frac12,\frac{13}{16},\frac47,\frac1{14},\frac67,\frac5{14},\frac{43}{112},\frac{75}{112}),\\(\frac12,\frac{15}{16},\frac37,\frac{13}{14},\frac17,\frac9{14},\frac{41}{112},\frac9{112}),(\frac12,\frac{15}{16},\frac47,\frac1{14},\frac67,\frac5{14},\frac{57}{112},\frac{89}{112}),&(\text{1st\&2nd})\\
    \\
    (\frac12,\frac1{16},\frac27,\frac{11}{14},\frac37,\frac{13}{14},\frac{39}{112},\frac{55}{112}),(\frac12,\frac1{16},\frac57,\frac3{14},\frac47,\frac1{14},\frac{87}{112},\frac{71}{112}),\\(\frac12,\frac3{16},\frac27,\frac{11}{14},\frac37,\frac{13}{14},\frac{53}{112},\frac{69}{112}),(\frac12,\frac3{16},\frac57,\frac3{14},\frac47,\frac1{14},\frac{101}{112},\frac{85}{112}),\\
    (\frac12,\frac5{16},\frac27,\frac{11}{14},\frac37,\frac{13}{14},\frac{67}{112},\frac{83}{112}),(\frac12,\frac5{16},\frac57,\frac3{14},\frac47,\frac1{14},\frac3{112},\frac{99}{112}),\\(\frac12,\frac7{16},\frac27,\frac{11}{14},\frac37,\frac{13}{14},\frac{81}{112},\frac{97}{112}),(\frac12,\frac7{16},\frac57,\frac3{14},\frac47,\frac1{14},\frac{17}{112},\frac1{112}),\\
    (\frac12,\frac9{16},\frac27,\frac{11}{14},\frac37,\frac{13}{14},\frac{95}{112},\frac{111}{112}),(\frac12,\frac9{16},\frac57,\frac3{14},\frac47,\frac1{14},\frac{31}{112},\frac{15}{112}),\\(\frac12,\frac{11}{16},\frac27,\frac{11}{14},\frac37,\frac{13}{14},\frac{109}{112},\frac{13}{112}),(\frac12,\frac{11}{16},\frac57,\frac3{14},\frac47,\frac1{14},\frac{45}{112},\frac{29}{112}),\\
    (\frac12,\frac{13}{16},\frac27,\frac{11}{14},\frac37,\frac{13}{14},\frac{11}{112},\frac{27}{112}),(\frac12,\frac{13}{16},\frac57,\frac3{14},\frac47,\frac1{14},\frac{59}{112},\frac{43}{112}),\\(\frac12,\frac{15}{16},\frac27,\frac{11}{14},\frac37,\frac{13}{14},\frac{25}{112},\frac{41}{112}),(\frac12,\frac{15}{16},\frac57,\frac3{14},\frac47,\frac1{14},\frac{73}{112},\frac{57}{112}),&(\text{3rd\&4th})\\
    \\
    (\frac12,\frac1{16},\frac17,\frac9{14},\frac57,\frac3{14},\frac{23}{112},\frac{87}{112}),(\frac12,\frac1{16},\frac67,\frac5{14},\frac27,\frac{11}{14},\frac{103}{112},\frac{39}{112}),\\(\frac12,\frac3{16},\frac17,\frac9{14},\frac57,\frac3{14},\frac{37}{112},\frac{101}{112}),(\frac12,\frac3{16},\frac67,\frac5{14},\frac27,\frac{11}{14},\frac5{112},\frac{53}{112}),\\
    (\frac12,\frac5{16},\frac17,\frac9{14},\frac57,\frac3{14},\frac{51}{112},\frac3{112}),(\frac12,\frac5{16},\frac67,\frac5{14},\frac27,\frac{11}{14},\frac{19}{112},\frac{67}{112}),\\(\frac12,\frac7{16},\frac17,\frac9{14},\frac57,\frac3{14},\frac{65}{112},\frac{17}{112}),(\frac12,\frac7{16},\frac67,\frac5{14},\frac27,\frac{11}{14},\frac{33}{112},\frac{81}{112}),\\
    (\frac12,\frac9{16},\frac17,\frac9{14},\frac57,\frac3{14},\frac{79}{112},\frac{31}{112}),(\frac12,\frac9{16},\frac67,\frac5{14},\frac27,\frac{11}{14},\frac{47}{112},\frac{95}{112}),\\(\frac12,\frac{11}{16},\frac17,\frac9{14},\frac57,\frac3{14},\frac{93}{112},\frac{45}{112}),(\frac12,\frac{11}{16},\frac67,\frac5{14},\frac27,\frac{11}{14},\frac{61}{112},\frac{109}{112}),\\
    (\frac12,\frac{13}{16},\frac17,\frac9{14},\frac57,\frac3{14},\frac{107}{112},\frac{59}{112}),(\frac12,\frac{13}{16},\frac67,\frac5{14},\frac27,\frac{11}{14},\frac{75}{112},\frac{11}{112}),\\(\frac12,\frac{15}{16},\frac17,\frac9{14},\frac57,\frac3{14},\frac9{112},\frac{73}{112}),(\frac12,\frac{15}{16},\frac67,\frac5{14},\frac27,\frac{11}{14},\frac{89}{112},\frac{25}{112}).&(\text{5th\&6th}).\end{cases}\quad(\mods1) \]
The $S$-matrices are given by
\[ \widetilde S=\begin{pmatrix}1&d_X&d_Y&d_Z&d_Xd_Z&d_T&d_Xd_T&d_Yd_Z&d_Yd_T\\d_X&1&-d_Y&d_Xd_Z&d_Z&d_Xd_T&d_T&-d_Yd_Z&-d_Yd_T\\d_Y&-d_Y&0&d_Yd_Z&-d_Yd_Z&d_Yd_T&-d_Yd_T&0&0\\d_Z&d_Xd_Z&d_Yd_Z&-d_T&-d_Xd_T&1&d_X&-d_Yd_T&d_Y\\d_Xd_Z&d_Z&-d_Yd_Z&-d_Xd_T&-d_T&d_X&1&d_Yd_T&-d_Y\\d_T&d_Xd_T&d_Yd_T&1&d_X&-d_Z&-d_Xd_Z&d_Y&-d_Yd_Z\\d_Xd_T&d_T&-d_Yd_T&d_X&1&-d_Xd_Z&-d_Z&-d_Y&d_Yd_Z\\d_Yd_Z&-d_Yd_Z&0&-d_Yd_T&d_Yd_T&d_Y&-d_Y&0&0\\d_Yd_T&-d_Yd_T&0&d_Y&-d_Y&-d_Yd_Z&d_Yd_Z&0&0\end{pmatrix}. \]
There are
\[ 6(\text{quantum dimensions})\times16(\text{conformal dimensions})\times2(\text{categorical dimensions})=192 \]
MFCs, among which those 32 with the sixth quantum dimensions are unitary. We classify connected étale algebras in all 192 MFCs simultaneously.

An ansatz
\[ A\cong1\oplus n_XX\oplus n_YY\oplus n_ZZ\oplus n_SS\oplus n_TT\oplus n_UU\oplus n_VV\oplus n_WW \]
with $n_j\in\mbb N$ has
\[ \fp_{\mc B}(A)=1+n_X+\sqrt2n_Y+\frac{\sin\frac{2\pi}7}{\sin\frac\pi7}(n_Z+n_S)+\frac{\sin\frac{3\pi}7}{\sin\frac\pi7}(n_T+n_U)+\sqrt2\frac{\sin\frac{2\pi}7}{\sin\frac\pi7}n_V+\sqrt2\frac{\sin\frac{3\pi}7}{\sin\frac\pi7}n_W. \]
For this to obey (\ref{FPdimA2bound}), the natural number can take only 75 values. The sets contain the one with all $n_j$'s be zero. It is the trivial connected étale algebra $A\cong1$ giving $\mc B_A^0\simeq\mc B_A\simeq\mc B$. The other 74 candidates contain nontrivial simple object(s) $b_j\not\cong1$ with nontrivial conformal dimensions. Thus, they all fail to be commutative.

We conclude
\begin{table}[H]
\begin{center}
\begin{tabular}{c|c|c|c}
    Connected étale algebra $A$&$\mc B_A$&$\rank(\mc B_A)$&Lagrangian?\\\hline
    $1$&$\mc B$&$9$&No
\end{tabular}.
\end{center}
\caption{Connected étale algebras in rank nine MFC $\mcal B\simeq\ising\boxtimes psu(2)_5$}\label{rank9isingpsu25results}
\end{table}
\hspace{-17pt}All the 192 MFCs $\mc B\simeq\ising\boxtimes psu(2)_5$'s are completely anisotropic.

\subsubsection{$\mc B\simeq so(11)_2$}
The MFCs have nine simple objects $\{1,X,Y,Z,S,T,U,V,W\}$ obeying monoidal products
\begin{table}[H]
\begin{center}
\makebox[1 \textwidth][c]{       
\resizebox{1.2 \textwidth}{!}{\begin{tabular}{c|c|c|c|c|c|c|c|c|c}
    $\otimes$&$1$&$X$&$Y$&$Z$&$S$&$T$&$U$&$V$&$W$\\\hline
    $1$&$1$&$X$&$Y$&$Z$&$S$&$T$&$U$&$V$&$W$\\\hline
    $X$&&$1$&$Y$&$Z$&$S$&$T$&$U$&$W$&$V$\\\hline
    $Y$&&&$1\oplus X\oplus U$&$T\oplus U$&$S\oplus T$&$Z\oplus S$&$Y\oplus Z$&$V\oplus W$&$V\oplus W$\\\hline
    $Z$&&&&$1\oplus X\oplus S$&$Z\oplus U$&$Y\oplus T$&$Y\oplus S$&$V\oplus W$&$V\oplus W$\\\hline
    $S$&&&&&$1\oplus X\oplus Y$&$Y\oplus U$&$Z\oplus T$&$V\oplus W$&$V\oplus W$\\\hline
    $T$&&&&&&$1\oplus X\oplus Z$&$S\oplus U$&$V\oplus W$&$V\oplus W$\\\hline
    $U$&&&&&&&$1\oplus X\oplus T$&$V\oplus W$&$V\oplus W$\\\hline
    $V$&&&&&&&&$1\oplus Y\oplus Z\oplus S\oplus T\oplus U$&$X\oplus Y\oplus Z\oplus S\oplus T\oplus U$\\\hline
    $W$&&&&&&&&&$1\oplus Y\oplus Z\oplus S\oplus T\oplus U$
\end{tabular}.}}
\end{center}
\end{table}
\hspace{-17pt}Thus, they have
\begin{align*}
    \hspace{-40pt}\fp_{\mc B}(1)=1=\fp_{\mc B}(X),\quad\fp_{\mc B}(Y)=&\fp_{\mc B}(Z)=\fp_{\mc B}(S)=\fp_{\mc B}(T)=\fp_{\mc B}(U)=2,\\
    \fp_{\mc B}(V)=&\sqrt{11}=\fp_{\mc B}(W),
\end{align*}
and
\[ \fp(\mc B)=44. \]
Their quantum dimensions $d_j$'s are solutions of the same multiplication rules $d_id_j=\sum_{k=1}^9{N_{ij}}^kd_k$. There are two (nonzero) solutions
\[ (d_X,d_Y,d_Z,d_S,d_T,d_U,d_V,d_W)=(1,2,2,2,2,2,-\sqrt{11},-\sqrt{11}),(1,2,2,2,2,2,\sqrt{11},\sqrt{11}) \]
with the same categorical dimension
\[ D^2(\mc B)=44. \]
They have four conformal dimensions\footnote{Naively one finds 40 consistent conformal dimensions, however, the other 36 are equivalent to one of the four in the main text under permutations $(VW)$ or $(YZUST)$ of simple objects.}
\begin{align*}
    (h_X,h_Y,h_Z,h_S,h_T,h_U,h_V,h_W)&=(0,\frac1{11},\frac9{11},\frac3{11},\frac5{11},\frac4{11},\frac18,\frac58),(0,\frac1{11},\frac9{11},\frac3{11},\frac5{11},\frac4{11},\frac38,\frac78),\\
    &~~~~(0,\frac{10}{11},\frac2{11},\frac8{11},\frac{6}{11},\frac7{11},\frac18,\frac58),(0,\frac{10}{11},\frac2{11},\frac8{11},\frac{6}{11},\frac7{11},\frac38,\frac78)\quad(\mods1).
\end{align*}
The $S$-matrices are given by
\[ \widetilde S=\begin{pmatrix}1&d_X&d_Y&d_Z&d_S&d_T&d_U&d_V&d_W\\d_X&d_X&d_Y&d_Z&d_S&d_T&d_U&-d_V&-d_W\\d_Y&d_Y&s_1&s_2&s_3&s_4&s_5&0&0\\d_Z&d_Z&s_2&s_5&s_4&s_1&s_3&0&0\\d_S&d_S&s_3&s_4&s_2&s_5&s_1&0&0\\d_T&d_T&s_4&s_1&s_5&s_3&s_2&0&0\\d_U&d_U&s_5&s_3&s_1&s_2&s_4&0&0\\d_V&-d_V&0&0&0&0&0&\pm d_V&\mp d_W\\d_W&-d_W&0&0&0&0&0&\mp d_W&\pm d_W\end{pmatrix} \]
with
\[ s_1=4\sin\frac{3\pi}{22},\quad s_2=-4\cos\frac\pi{11},\quad s_3=4\cos\frac{2\pi}{11},\quad s_4=-4\sin\frac\pi{22},\quad s_5=-4\sin\frac{5\pi}{22}. \]
They have additive central charges
\[ c(\mc B)=\begin{cases}2&(\text{1st\&2nd }h),\\-2&(\text{3rd\&4th }h).\end{cases}\quad(\mods8) \]
There are
\[ 2(\text{quantum dimensions})\times4(\text{conformal dimensions})\times2(\text{categorical dimensions})=16 \]
MFCs, among which those eight with the second quantum dimensions are unitary. We classify connected étale algebras in all 16 MFCs simultaneously.

An ansatz
\[ A\cong1\oplus n_XX\oplus n_YY\oplus n_ZZ\oplus n_SS\oplus n_TT\oplus n_UU\oplus n_VV\oplus n_WW \]
with $n_j\in\mbb N$ has
\[ \fp_{\mc B}(A)=1+n_X+2(n_Y+n_Z+n_S+n_T+n_U)+\sqrt{11}(n_V+n_W). \]
For this to obey (\ref{FPdimA2bound}), the natural numbers can take only 72 values. The sets contain the one with all $n_j$'s be zero. It is the trivial connected étale algebra giving $\mc B_A^0\simeq\mc B_A\simeq\mc B$. Other 32 candidates without $X$ all fail to be commutative because nontrivial simple object(s) entering them have nontrivial conformal dimensions. Thus, we are left with five nontrivial candidates with $n_X=1,2,3,4,5$. Since $X$ has $(d_X,h_X)=(1,0)$ (mod 1 for $h_X$), it has $c_{X,X}\cong id_1$ \cite{KK23preMFC}, and can give commutative algebras. An obstruction is the Frobenius-Perron dimensions. The Frobenius-Perron dimensions $2,3,4,5,6$ of the candidates demand $\fp(\mc B_A^0)=11,\frac{44}9,\frac{11}4,\frac{44}{25},\frac{11}9$, but there are no MFCs with such Frobenius-Perron dimensions for all of them but the first. Therefore, the only nontrivial candidate is
\[ A\cong1\oplus X. \]
This is just a $\mbb Z/2\mbb Z$ algebra, and it is commutative. It further turns out separable. Let us check this point by identifying $\mc B_A$.

Its Frobenius-Perron dimension $\fp=2$ demands
\[ \fp(\mc B_A^0)=11,\quad\fp(\mc B_A)=22. \]
It turns out
\[ \mc B_A^0\simeq\vecG_{\mbb Z/11\mbb Z}^1,\quad\mc B_A\simeq\text{TY}(\mbb Z/11\mbb Z). \]
This identification also matches central charges because
\[ \hspace{-30pt}c(\vecG_{\mbb Z/11\mbb Z}^1)=\begin{cases}2&(h_X,h_Y,h_Z,h_Q,h_R,h_S,h_T,h_U,h_V,h_W)=(\frac1{11},\frac4{11},\frac9{11},\frac5{11},\frac3{11},\frac3{11},\frac5{11},\frac9{11},\frac4{11},\frac1{11}),\\-2&(h_X,h_Y,h_Z,h_Q,h_R,h_S,h_T,h_U,h_V,h_W)=(\frac{10}{11},\frac7{11},\frac2{11},\frac7{11},\frac8{11},\frac8{11},\frac6{11},\frac2{11},\frac7{11},\frac{10}{11}).\end{cases}\quad(\mods8) \]
One of the easiest ways to see the identifications is to perform anyon condensation. Since we `identify' $1$ and $X$, $V$ and $W$ are also identified. The other simple objects `split' into two each. Since $V$ and $W$ have different conformal dimensions, they are confined. Thus, the category of dyslectic right $A$-modules consists of 11 invertible simple objects, and it is identified as a $\mbb Z/11\mbb Z$ MFC. In the category $\mc B_A$ of right $A$-modules, $V\oplus W$ has Frobenius-Perron dimension $\sqrt{11}$. Thus, the category is identified as a $\mbb Z/11\mbb Z$ Tambara-Yamagami category $\text{TY}(\mbb Z/11\mbb Z)$.\footnote{More rigorously, we should search for NIM-reps. Indeed, we find a 12-dimensional NIM-rep
\begin{align*}
    \hspace{-20pt}&n_1=1_{12}=n_X,\quad n_Y=\begin{pmatrix}0&1&1&0&0&0&0&0&0&0&0&0\\1&0&0&0&0&0&0&0&0&1&0&0\\1&0&0&0&0&0&0&0&0&0&1&0\\0&0&0&0&0&0&0&1&0&1&0&0\\0&0&0&0&0&0&0&0&1&0&1&0\\0&0&0&0&0&0&1&1&0&0&0&0\\0&0&0&0&0&1&0&0&1&0&0&0\\0&0&0&1&0&1&0&0&0&0&0&0\\0&0&0&0&1&0&1&0&0&0&0&0\\0&1&0&1&0&0&0&0&0&0&0&0\\0&0&1&0&1&0&0&0&0&0&0&0\\0&0&0&0&0&0&0&0&0&0&0&2\end{pmatrix},\quad n_Z=\begin{pmatrix}0&0&0&1&1&0&0&0&0&0&0&0\\0&0&0&0&0&0&0&1&0&0&1&0\\0&0&0&0&0&0&0&0&1&1&0&0\\1&0&0&0&0&0&1&0&0&0&0&0\\1&0&0&0&0&1&0&0&0&0&0&0\\0&0&0&0&1&0&0&0&0&1&0&0\\0&0&0&1&0&0&0&0&0&0&1&0\\0&1&0&0&0&0&0&0&1&0&0&0\\0&0&1&0&0&0&0&1&0&0&0&0\\0&0&1&0&0&1&0&0&0&0&0&0\\0&1&0&0&0&0&1&0&0&0&0&0\\0&0&0&0&0&0&0&0&0&0&0&2\end{pmatrix},\\
    \hspace{-20pt}&n_S=\begin{pmatrix}0&0&0&0&0&1&1&0&0&0&0&0\\0&0&0&0&0&0&1&0&1&0&0&0\\0&0&0&0&0&1&0&1&0&0&0&0\\0&0&0&0&1&0&0&0&0&0&1&0\\0&0&0&1&0&0&0&0&0&1&0&0\\1&0&1&0&0&0&0&0&0&0&0&0\\1&1&0&0&0&0&0&0&0&0&0&0\\0&0&1&0&0&0&0&0&0&0&1&0\\0&1&0&0&0&0&0&0&0&1&0&0\\0&0&0&0&1&0&0&0&1&0&0&0\\0&0&0&1&0&0&0&1&0&0&0&0\\0&0&0&0&0&0&0&0&0&0&0&2\end{pmatrix},\quad n_T=\begin{pmatrix}0&0&0&0&0&0&0&1&1&0&0&0\\0&0&0&0&1&1&0&0&0&0&0&0\\0&0&0&1&0&0&1&0&0&0&0&0\\0&0&1&0&0&0&0&0&1&0&0&0\\0&1&0&0&0&0&0&1&0&0&0&0\\0&1&0&0&0&0&0&0&0&0&1&0\\0&0&1&0&0&0&0&0&0&1&0&0\\1&0&0&0&1&0&0&0&0&0&0&0\\1&0&0&1&0&0&0&0&0&0&0&0\\0&0&0&0&0&0&1&0&0&0&1&0\\0&0&0&0&0&1&0&0&0&1&0&0\\0&0&0&0&0&0&0&0&0&0&0&2\end{pmatrix},\\
    \hspace{-20pt}&n_U=\begin{pmatrix}0&0&0&0&0&0&0&0&0&1&1&0\\0&0&1&1&0&0&0&0&0&0&0&0\\0&1&0&0&1&0&0&0&0&0&0&0\\0&1&0&0&0&1&0&0&0&0&0&0\\0&0&1&0&0&0&1&0&0&0&0&0\\0&0&0&1&0&0&0&0&1&0&0&0\\0&0&0&0&1&0&0&1&0&0&0&0\\0&0&0&0&0&0&1&0&0&1&0&0\\0&0&0&0&0&1&0&0&0&0&1&0\\1&0&0&0&0&0&0&1&0&0&0&0\\1&0&0&0&0&0&0&0&1&0&0&0\\0&0&0&0&0&0&0&0&0&0&0&2\end{pmatrix},\quad n_V=\begin{pmatrix}0&0&0&0&0&0&0&0&0&0&0&1\\0&0&0&0&0&0&0&0&0&0&0&1\\0&0&0&0&0&0&0&0&0&0&0&1\\0&0&0&0&0&0&0&0&0&0&0&1\\0&0&0&0&0&0&0&0&0&0&0&1\\0&0&0&0&0&0&0&0&0&0&0&1\\0&0&0&0&0&0&0&0&0&0&0&1\\0&0&0&0&0&0&0&0&0&0&0&1\\0&0&0&0&0&0&0&0&0&0&0&1\\0&0&0&0&0&0&0&0&0&0&0&1\\0&0&0&0&0&0&0&0&0&0&0&1\\1&1&1&1&1&1&1&1&1&1&1&0\end{pmatrix}=n_W.
\end{align*}
The NIM-rep gives identifications
\begin{align*}
    m_1\cong1\oplus X,\quad m_2\cong&Y\cong m_3,\quad m_4\cong Z\cong m_5,\\
    m_6\cong S\cong m_7,\quad m_8\cong T\cong m_9,&\quad m_{10}\cong U\cong m_{11},\quad m_{12}\cong V\oplus W.
\end{align*}
Since $m_1,m_2,\dots,m_{11}$ have quantum dimensions
\[ d_{\mc B_A}(m_j)=\frac{d_{\mc B}(m_j)}{d_{\mc B}(A)}=+1, \]
and conformal dimensions $(\frac1{11},\frac9{11},\frac3{11},\frac5{11},\frac4{11})$ or its conjugate, $\mc B_A^0$ is identified as a $\mbb Z/11\mbb Z$ MFC. The non-invertible simple object $m_{12}$ has $d_{\mc B_A}(m_{12})=\pm\sqrt{11}$ and obeys monoidal products
\[ m_j\otimes_A m_{12}\cong m_{12}\cong m_{12}\otimes_Am_j\ (j=1,2,\dots,11),\quad m_{12}\otimes_Am_{12}\cong\bigoplus_{j=1}^{11}m_j. \]
This shows $\mc B_A\simeq\text{TY}(\mbb Z/11\mbb Z)$.} Since $\text{TY}(\mbb Z/11\mbb Z)$ is semisimple, $A$ is separable and étale. Note that this example also has
\[ \rank(\mc B_A)>\rank(\mc B). \]

We conclude
\begin{table}[H]
\begin{center}
\begin{tabular}{c|c|c|c}
    Connected étale algebra $A$&$\mc B_A$&$\rank(\mc B_A)$&Lagrangian?\\\hline
    $1$&$\mc B$&$9$&No\\
    $1\oplus X$&$\text{TY}(\mbb Z/11\mbb Z)$&12&No
\end{tabular}.
\end{center}
\caption{Connected étale algebras in rank nine MFC $\mcal B\simeq so(11)_2$}\label{rank9so112results}
\end{table}
\hspace{-17pt}All the 16 MFCs $\mc B\simeq so(11)_2$'s fail to be completely anisotropic.

\subsubsection{$\mc B\simeq su(2)_8$}
The MFCs have nine simple objects $\{1,X,Y,Z,S,T,U,V,W\}$ obeying monoidal products
\begin{table}[H]
\begin{center}
\makebox[1 \textwidth][c]{       
\resizebox{1.2 \textwidth}{!}{\begin{tabular}{c|c|c|c|c|c|c|c|c|c}
    $\otimes$&$1$&$X$&$Y$&$Z$&$S$&$T$&$U$&$V$&$W$\\\hline
    $1$&$1$&$X$&$Y$&$Z$&$S$&$T$&$U$&$V$&$W$\\\hline
    $X$&&$1$&$Z$&$Y$&$T$&$S$&$V$&$U$&$W$\\\hline
    $Y$&&&$1\oplus T$&$X\oplus S$&$Z\oplus V$&$Y\oplus U$&$T\oplus W$&$S\oplus W$&$U\oplus V$\\\hline
    $Z$&&&&$1\oplus T$&$Y\oplus U$&$Z\oplus V$&$S\oplus W$&$T\oplus W$&$U\oplus V$\\\hline
    $S$&&&&&$1\oplus T\oplus W$&$X\oplus S\oplus W$&$Z\oplus U\oplus V$&$Y\oplus U\oplus V$&$S\oplus T\oplus W$\\\hline
    $T$&&&&&&$1\oplus T\oplus W$&$Y\oplus U\oplus V$&$Z\oplus U\oplus V$&$S\oplus T\oplus W$\\\hline
    $U$&&&&&&&$1\oplus S\oplus T\oplus W$&$X\oplus S\oplus T\oplus W$&$Y\oplus Z\oplus U\oplus V$\\\hline
    $V$&&&&&&&&$1\oplus S\oplus T\oplus W$&$Y\oplus Z\oplus U\oplus V$\\\hline
    $W$&&&&&&&&&$1\oplus X\oplus S\oplus T\oplus W$
\end{tabular}.}}
\end{center}
\end{table}
\hspace{-17pt}Thus, they have
\begin{align*}
    \hspace{-50pt}\fp_{\mc B}(1)=1=\fp_{\mc B}(X),\quad\fp_{\mc B}(Y)=\sqrt{\frac{5+\sqrt5}2}=&\fp_{\mc B}(Z),\quad\fp_{\mc B}(S)=\frac{3+\sqrt5}2=\fp_{\mc B}(T),\\
    \fp_{\mc B}(U)=\sqrt{5+2\sqrt5}=
    \fp_{\mc B}(V),\quad&\fp_{\mc B}(W)=1+\sqrt5,
\end{align*}
and
\[ \fp(\mc B)=30+10\sqrt5\approx52.4. \]
Their quantum dimensions $d_j$'s are solutions of the same multiplication rules $d_id_j=\sum_{k=1}^9{N_{ij}}^kd_k$. There are four (nonzero) solutions
\begin{align*}
    \hspace{-50pt}(d_X,d_Y,d_Z,d_S,d_T,d_U,d_V,d_W)&=(1,-\sqrt{\frac{5-\sqrt5}2},-\sqrt{\frac{5-\sqrt5}2},\frac{3-\sqrt5}2,\frac{3-\sqrt5}2,\sqrt{5-2\sqrt5},\sqrt{5-2\sqrt5},1-\sqrt5),\\
    \hspace{-50pt}&~~~~(1,\sqrt{\frac{5-\sqrt5}2},\sqrt{\frac{5-\sqrt5}2},\frac{3-\sqrt5}2,\frac{3-\sqrt5}2,-\sqrt{5-2\sqrt5},-\sqrt{5-2\sqrt5},1-\sqrt5),\\
    \hspace{-50pt}&~~~~(1,-\sqrt{\frac{5+\sqrt5}2},-\sqrt{\frac{5+\sqrt5}2},\frac{3+\sqrt5}2,\frac{3+\sqrt5}2,-\sqrt{5+2\sqrt5},-\sqrt{5+2\sqrt5},1+\sqrt5),\\
    \hspace{-50pt}&~~~~(1,\sqrt{\frac{5+\sqrt5}2},\sqrt{\frac{5+\sqrt5}2},\frac{3+\sqrt5}2,\frac{3+\sqrt5}2,\sqrt{5+2\sqrt5},\sqrt{5+2\sqrt5},1+\sqrt5),
\end{align*}
with categorical dimensions
\[ D^2(\mc B)=30-10\sqrt5(\approx7.6),\quad30+10\sqrt5, \]
respectively for each pair. They have four conformal dimensions\footnote{Naively, one finds eight consistent conformal dimensions, but the other four are equivalent to one in the main text under permutations $(YZ)(UV)$ of simple objects.}
\begin{align*}
    (h_X,h_Y,h_Z,h_S,h_T,h_U,h_V,h_W)&=(0,\frac1{40},\frac{21}{40},\frac25,\frac25,\frac18,\frac58,\frac15),(0,\frac9{40},\frac{29}{40},\frac35,\frac35,\frac18,\frac58,\frac45),\\
    &~~~~(0,\frac{11}{40},\frac{31}{40},\frac25,\frac25,\frac38,\frac78,\frac15),(0,\frac{19}{40},\frac{39}{40},\frac35,\frac35,\frac38,\frac78,\frac45)\quad(\mods1)
\end{align*}
for the first two quantum dimensions and
\begin{align*}
    (h_X,h_Y,h_Z,h_S,h_T,h_U,h_V,h_W)&=(0,\frac3{40},\frac{23}{40},\frac15,\frac15,\frac38,\frac78,\frac35),(0,\frac7{40},\frac{27}{40},\frac45,\frac45,\frac78,\frac38,\frac25),\\
    &~~~~(0,\frac{13}{40},\frac{33}{40},\frac15,\frac15,\frac58,\frac18,\frac35),(0,\frac{17}{40},\frac{37}{40},\frac45,\frac45,\frac18,\frac58,\frac25)\quad(\mods1)
\end{align*}
for the last two quantum dimensions. The $S$-matrices are given by
\[ \widetilde S=\begin{pmatrix}1&d_X&d_Y&d_Z&d_S&d_T&d_U&d_V&d_W\\d_X&d_X&-d_Y&-d_Z&d_S&d_T&-d_U&-d_V&d_W\\d_Y&-d_Y&\pm d_U&\mp d_U&-d_U&d_U&\pm d_Y&\mp d_Y&0\\d_Z&-d_Z&\mp d_U&\pm d_U&-d_U&d_U&\mp d_Y&\pm d_Y&0\\d_S&d_S&-d_U&-d_U&1&1&d_Y&d_Y&-d_W\\d_T&d_T&d_U&d_U&1&1&-d_Y&-d_Y&-d_W\\d_U&-d_U&\pm d_Y&\mp d_Y&d_Y&-d_Y&\pm d_U&\mp d_U&0\\d_V&-d_V&\mp d_Y&\pm d_Y&d_Y&-d_Y&\mp d_U&\pm d_U&0\\d_W&d_W&0&0&-d_W&-d_W&0&0&d_W\end{pmatrix}. \]
They have additive central charges
\[ c(\mc B)=\begin{cases}\begin{cases}\frac45&(\text{1st\&3rd }h),\\-\frac45&(\text{2nd\&4th }h),\end{cases}&(\text{1st\&2nd quantum dimensions})\\\begin{cases}\frac{12}5&(\text{1st\&3rd }h),\\-\frac{12}5&(\text{2nd\&4th }h).\end{cases}&(\text{3rd\&4th quantum dimensions})\end{cases}\quad(\mods8) \]
There are
\[ 4(\text{quantum dimensions})\times4(\text{conformal dimensions})\times2(\text{categorical dimensions})=32 \]
MFCs, among which those eight with the fourth quantum dimensions are unitary. We classify connected étale algebras in all 32 MFCs simultaneously.

An ansatz
\[ A\cong1\oplus n_XX\oplus n_YY\oplus n_ZZ\oplus n_SS\oplus n_TT\oplus n_UU\oplus n_VV\oplus n_WW \]
with $n_j\in\mbb N$ has
\[ \fp_{\mc B}(A)=1+n_X+\sqrt{\frac{5+\sqrt5}2}(n_Y+n_Z)+\frac{3+\sqrt5}2(n_S+n_T)+\sqrt{5+2\sqrt5}(n_U+n_V)+(1+\sqrt5)n_W. \]
For this to obey (\ref{FPdimA2bound}), the natural numbers can take only 85 values. The sets contain the one with all $n_j$'s be zero. It is the trivial connected étale algebra $A\cong1$ giving $\mc B_A^0\simeq\mc B_A\simeq\mc B$. Other candidates without $X$ contain simple object(s) with nontrivial conformal dimensions, and fail to be commutative. Thus, we are left with candidates $1\oplus n_XX$ with $n_X=1,2,3,4,5,6$. Since $X$ has $(d_X,h_X)=(1,0)$, it has trivial braiding $c_{X,X}\cong id_1$ \cite{KK23preMFC}, and the candidates can be commutative. The obstruction is the Frobenius-Perron dimensions. The Frobenius-Perron dimensions $2,3,4,5,6,7$ of the candidates demand $\fp(\mc B_A^0)=\frac{15+5\sqrt5}2,\frac{30+10\sqrt5}9,\frac{30+10\sqrt5}{16},\frac{30+10\sqrt5}{25},\frac{30+10\sqrt5}{36},\frac{30+10\sqrt5}{49}$, but there is no MFC with such Frobenius-Perron dimensions for $n_X>1$. Therefore, the only candidate is $A\cong1\oplus X$. This is indeed a commutative $\mbb Z/2\mbb Z$ algebra, and it turns out separable. To see this, we identify $\mc B_A$.

Its Frobenius-Perron dimension demands
\[ \fp(\mc B_A^0)=\frac{15+5\sqrt5}2,\quad\fp(\mc B_A)=15+5\sqrt5. \]
From the Frobenius-Perron dimension, the MFC $\mc B_A^0$ is identified as
\[ \mc B_A^0\simeq\fib\boxtimes\fib. \]
This can also be seen via anyon condensation. Condensation of $A\cong1\oplus X$ `identifies' a pair $(1,X)$, and hence $(Y,Z),(S,T),(U,V)$, and it `splits' $W$ into two simple objects. Since $Y$ and $Z$ have different conformal dimensions, the resulting simple object in $\mc B_A$ is confined. Similarly for the simple object descending from $U,V$. Thus, the deconfined phase $\mc B_A^0$ consists of four simple objects with quantum dimensions $1,\frac{3\pm\sqrt5}2,\frac{1\pm\sqrt5}2,\frac{1\pm\sqrt5}2$. This is identified as $\fib\boxtimes\fib$ where two Fibonacci objects have the same quantum and conformal dimensions. The identification also matches central charges. To describe the matching more precisely, let us denote two Fibonacci objects $x,y$. The central charges are matched as
\[ c(\mc B)=c(\mc B_A^0)=\begin{cases}\pm\frac45&(d_W,h_W,d_x,d_y,h_x,h_y)=(1-\sqrt5,\pm\frac15,\frac{1-\sqrt5}2,\frac{1-\sqrt5}2,\pm\frac15,\pm\frac15),\\\pm\frac{12}5&(d_W,h_W,d_x,d_y,h_x,h_y)=(1+\sqrt5,\pm\frac25,\frac{1+\sqrt5}2,\frac{1+\sqrt5}2,\pm\frac25,\pm\frac25),\end{cases}\quad(\mods1\text{ for }h) \]
where signs are correlated. The category $\mc B_A$ of right $A$-modules have two more simple objects with quantum dimensions $\pm\sqrt{\frac{5\pm\sqrt5}2}$ and $\sqrt{5\pm2\sqrt5}$. Such a rank six fusion category is not listed in the AnyonWiki, but we can figure out their properties.

More rigorously, we should find NIM-reps. Indeed, we find a six-dimensional solution
\begin{align*}
    n_1=1_6=n_X,&\quad n_Y=\begin{pmatrix}0&1&0&0&0&0\\1&0&1&0&0&0\\0&1&0&1&0&0\\0&0&1&0&1&1\\0&0&0&1&0&0\\0&0&0&1&0&0\end{pmatrix}=n_Z,\quad n_S=\begin{pmatrix}0&0&1&0&0&0\\0&1&0&1&0&0\\1&0&1&0&1&1\\0&1&0&2&0&0\\0&0&1&0&0&1\\0&0&1&0&1&0\end{pmatrix}=n_T,\\
    &n_U=\begin{pmatrix}0&0&0&1&0&0\\0&0&1&0&1&1\\0&1&0&2&0&0\\1&0&2&0&1&1\\0&1&0&1&0&0\\0&1&0&1&0&0\end{pmatrix}=n_V,\quad n_W=\begin{pmatrix}0&0&0&0&1&1\\0&0&0&2&0&0\\0&0&2&0&1&1\\0&2&0&2&0&0\\1&0&1&0&1&0\\1&0&1&0&0&1\end{pmatrix}.
\end{align*}
This gives a multiplication table
\begin{table}[H]
\begin{center}
\makebox[1 \textwidth][c]{       
\resizebox{1.2 \textwidth}{!}{
\begin{tabular}{c|c|c|c|c|c|c}
    $b_j\otimes\backslash$&$m_1$&$m_2$&$m_3$&$m_4$&$m_5$&$m_6$\\\hline
    $1,X$&$m_1$&$m_2$&$m_3$&$m_4$&$m_5$&$m_6$\\
    $Y,Z$&$m_2$&$m_1\oplus m_3$&$m_2\oplus m_4$&$m_3\oplus m_5\oplus m_6$&$m_4$&$m_4$\\
    $S,T$&$m_3$&$m_2\oplus m_4$&$m_1\oplus m_3\oplus m_5\oplus m_6$&$m_2\oplus2m_4$&$m_3\oplus m_6$&$m_3\oplus m_5$\\
    $U,V$&$m_4$&$m_3\oplus m_5\oplus m_6$&$m_2\oplus2m_4$&$m_1\oplus2m_3\oplus m_5\oplus m_6$&$m_2\oplus m_4$&$m_2\oplus m_4$\\
    $W$&$m_5\oplus m_6$&$2m_4$&$2m_3\oplus m_5\oplus m_6$&$2m_2\oplus2m_4$&$m_1\oplus m_3\oplus m_5$&$m_1\oplus m_3\oplus m_6$
\end{tabular}.}}
\end{center}
\end{table}
\hspace{-17pt}We obtain identifications
\[ m_1\cong1\oplus X,\quad m_2\cong Y\oplus Z,\quad m_3\cong S\oplus T,\quad m_4\cong U\oplus V,\quad m_5\cong W\cong m_6. \]
They have
\[ d_{\mc B_A}(m_1)=1,\quad d_{\mc B_A}(m_2)=d_Y,\quad d_{\mc B_A}(m_3)=d_S,\quad d_{\mc B_A}(m_4)=d_U,\quad dd_{\mc B_A}(m_5)=\frac{d_W}2=d_{\mc B_A}(m_6). \]
They obey monoidal products
\begin{table}[H]
\begin{center}
\begin{tabular}{c|c|c|c|c|c|c}
    $\otimes_A$&$m_1$&$m_2$&$m_3$&$m_4$&$m_5$&$m_6$\\\hline
    $m_1$&$m_1$&$m_2$&$m_3$&$m_4$&$m_5$&$m_6$\\\hline
    $m_2$&&$m_1\oplus m_3$&$m_2\oplus m_4$&$m_3\oplus m_5\oplus m_6$&$m_4$&$m_4$\\\hline
    $m_3$&&&$m_1\oplus m_3\oplus m_5\oplus m_6$&$m_2\oplus2m_4$&$m_3\oplus m_6$&$m_3\oplus m_5$\\\hline
    $m_4$&&&&$m_1\oplus2m_3\oplus m_5\oplus m_6$&$m_2\oplus m_4$&$m_2\oplus m_4$\\\hline
    $m_5$&&&&&$m_1\oplus m_5$&$m_3$\\\hline
    $m_6$&&&&&&$m_1\oplus m_6$
\end{tabular}.
\end{center}
\end{table}
\hspace{-17pt}Note that the fusion ring has multiplicity two. Thus, it is natural that the ring is not listed in \cite{LPR20} and AnyonWiki. Luckily, however, we know this is a fusion category as follows. Let us choose the specific $su(2)_8$ MFC describing $su(2)_8$ WZW model. (The MFC is also denoted as $\mc C(A_1,8)$.) It is known \cite{KO01} that $\mc C(A_1,8)$ has two connected étale algebras, $1,1\oplus X$. The second thus gives a fusion category $\mc C(A_1,8)_A$. The fusion category corresponds to the $D_6$ Dynkin diagram. The correspondence especially implies $\rank(\mc C(A_1,8)_A)=6$. Therefore, we learn the fusion ring with multiplicity two gives rank six fusion category. Let us collectively denote the fusion categories with the fusion ring as $\mc C(D_6)$.

We now relax the assumption on conformal dimensions and consider all $su(2)_8$ MFCs. Our analysis above showed $A\cong1\oplus X$ is a connected commutative algebra giving a fusion category $\mc B_A\simeq\mc C(D_6)$. Since $\mc C(D_6)$ is semisimple, $1\oplus X$ is separable, hence étale.

We conclude
\begin{table}[H]
\begin{center}
\begin{tabular}{c|c|c|c}
    Connected étale algebra $A$&$\mc B_A$&$\rank(\mc B_A)$&Lagrangian?\\\hline
    $1$&$\mc B$&$9$&No\\
    $1\oplus X$&$\mc C(D_6)$&6&No
\end{tabular}.
\end{center}
\caption{Connected étale algebras in rank nine MFC $\mcal B\simeq su(2)_8$}\label{rank9su28results}
\end{table}
\hspace{-17pt}All the 32 MFCs $\mc B\simeq su(2)_8$'s fail to be completely anisotropic.

\subsubsection{$\mc B\simeq psu(2)_5\boxtimes psu(2)_5$}\label{psu25psu25}
The MFCs have nine simple objects $\{1,X,Y,Z,S,T,U,V,W\}$ obeying monoidal products
\begin{table}[H]
\begin{center}
\makebox[1 \textwidth][c]{       
\resizebox{1.2 \textwidth}{!}{\begin{tabular}{c|c|c|c|c|c|c|c|c|c}
    $\otimes$&$1$&$X$&$Y$&$Z$&$S$&$T$&$U$&$V$&$W$\\\hline
    $1$&$1$&$X$&$Y$&$Z$&$S$&$T$&$U$&$V$&$W$\\\hline
    $X$&&$1\oplus S$&$T$&$V$&$X\oplus S$&$Y\oplus U$&$T\oplus U$&$Z\oplus W$&$V\oplus W$\\\hline
    $Y$&&&$1\oplus Z$&$Y\oplus Z$&$U$&$X\oplus V$&$S\oplus W$&$T\oplus V$&$U\oplus W$\\\hline
    $Z$&&&&$1\oplus Y\oplus Z$&$W$&$T\oplus V$&$U\oplus W$&$X\oplus T\oplus V$&$S\oplus U\oplus W$\\\hline
    $S$&&&&&$1\oplus X\oplus S$&$T\oplus U$&$Y\oplus T\oplus U$&$V\oplus W$&$Z\oplus V\oplus W$\\\hline
    $T$&&&&&&$1\oplus Z\oplus S\oplus W$&$X\oplus S\oplus V\oplus W$&$Y\oplus Z\oplus U\oplus W$&$T\oplus U\oplus V\oplus W$\\\hline
    $U$&&&&&&&$1\oplus X\oplus Z\oplus S\oplus V\oplus W$&$T\oplus U\oplus V\oplus W$&$Y\oplus Z\oplus T\oplus U\oplus V\oplus W$\\\hline
    $V$&&&&&&&&$1\oplus Y\oplus Z\oplus S\oplus U\oplus W$&$X\oplus S\oplus T\oplus U\oplus V\oplus W$\\\hline
    $W$&&&&&&&&&$1\oplus X\oplus Y\oplus Z\oplus S\oplus T\oplus U\oplus V\oplus W$
\end{tabular}.}}
\end{center}
\end{table}
\hspace{-17pt}(One can identify $psu(2)_5=\{1,X,S\},\{1,Y,Z\}$, and $T\cong X\otimes Y,U\cong Y\otimes S,V\cong X\otimes Z,W\cong Z\otimes S$.) Thus, they have
\begin{align*}
    \hspace{-00pt}\fp_{\mc B}(1)=1,\quad\fp_{\mc B}(X)=\frac{\sin\frac{2\pi}7}{\sin\frac\pi7}=\fp_{\mc B}(Y),\quad\fp_{\mc B}(Z)=\frac{\sin\frac{3\pi}7}{\sin\frac\pi7}=\fp_{\mc B}(S),\\
    \fp_{\mc B}(T)=\left(\frac{\sin\frac{2\pi}7}{\sin\frac\pi7}\right)^2,\quad
    \fp_{\mc B}(U)=\frac{\sin\frac{2\pi}7\sin\frac{3\pi}7}{\sin^2\frac\pi7}=
    \fp_{\mc B}(V),\quad\fp_{\mc B}(W)=\left(\frac{\sin\frac{3\pi}7}{\sin\frac\pi7}\right)^2,
\end{align*}
and
\[ \fp(\mc B)=\frac{49}{16\sin^4\frac\pi7}\approx86.4. \]
Their quantum dimensions $d_j$'s are solutions of the same multiplication rules $d_id_j=\sum_{k=1}^9{N_{ij}}^kd_k$. There are nine solutions
\begin{align*}
    &(d_X,d_Y,d_Z,d_S,d_T,d_U,d_V,d_W)\\
    &=(\frac{\sin\frac\pi7}{\cos\frac\pi{14}},\frac{\sin\frac\pi7}{\cos\frac\pi{14}},-\frac{\sin\frac{2\pi}7}{\cos\frac\pi{14}},-\frac{\sin\frac{2\pi}7}{\cos\frac\pi{14}},\left(\frac{\sin\frac\pi7}{\cos\frac\pi{14}}\right)^2,-\frac{\sin\frac\pi7\sin\frac{2\pi}7}{\cos^2\frac\pi{14}},-\frac{\sin\frac\pi7\sin\frac{2\pi}7}{\cos^2\frac\pi{14}},\left(\frac{\sin\frac{2\pi}7}{\cos\frac\pi{14}}\right)^2,),\\
    &~~~~(-\frac{\sin\frac{3\pi}7}{\cos\frac{3\pi}{14}},\frac{\sin\frac\pi7}{\cos\frac\pi{14}},-\frac{\sin\frac{2\pi}7}{\cos\frac\pi{14}},\frac{\sin\frac\pi7}{\cos\frac{3\pi}{14}},-\frac{\sin\frac\pi7}{\cos\frac{3\pi}{14}},\frac{\sin^2\frac\pi7}{\cos\frac\pi{14}\cos\frac{3\pi}{14}},1,-\frac{\sin\frac\pi7\sin\frac{2\pi}7}{\cos\frac\pi{14}\cos\frac{3\pi}{14}}),\\
    &~~~~(\frac{\sin\frac\pi7}{\cos\frac\pi{14}},-\frac{\sin\frac{3\pi}7}{\cos\frac{3\pi}{14}},\frac{\sin\frac\pi7}{\cos\frac{3\pi}{14}},-\frac{\sin\frac{2\pi}7}{\cos\frac\pi{14}},-\frac{\sin\frac\pi7}{\cos\frac{3\pi}{14}},1,\frac{\sin^2\frac\pi7}{\cos\frac\pi{14}\cos\frac{3\pi}{14}},-\frac{\sin\frac\pi7\sin\frac{2\pi}7}{\cos\frac\pi{14}\cos\frac{3\pi}{14}}),\\
    &~~~~(-\frac{\sin\frac{3\pi}7}{\cos\frac{3\pi}{14}},-\frac{\sin\frac{3\pi}7}{\cos\frac{3\pi}{14}},\frac{\sin\frac\pi7}{\cos\frac{3\pi}{14}},\frac{\sin\frac\pi7}{\cos\frac{3\pi}{14}},\left(\frac{\sin\frac{3\pi}7}{\cos\frac{3\pi}{14}}\right)^2,-\frac{\sin\frac\pi7\sin\frac{3\pi}7}{\cos^2\frac{3\pi}{14}},-\frac{\sin\frac\pi7\sin\frac{3\pi}7}{\cos^2\frac{3\pi}{14}},\left(\frac{\sin\frac\pi7}{\cos\frac{3\pi}{14}}\right)^2),\\
    &~~~~(\frac{\sin\frac\pi7}{\cos\frac\pi{14}},\frac{\sin\frac{2\pi}7}{\sin\frac\pi7},\frac{\sin\frac{3\pi}7}{\sin\frac\pi7},-\frac{\sin\frac{2\pi}7}{\cos\frac\pi{14}},\frac{\sin\frac\pi7\sin\frac{2\pi}7}{\cos\frac\pi{14}\sin\frac\pi7},-\frac{\sin^2\frac{2\pi}7}{\cos\frac\pi{14}\sin\frac\pi7},1,-\frac{\sin\frac{2\pi}7\sin\frac{3\pi}7}{\cos\frac\pi{14}\sin\frac\pi7}),\\
    &~~~~(\frac{\sin\frac{2\pi}7}{\sin\frac\pi7},\frac{\sin\frac\pi7}{\cos\frac\pi{14}},-\frac{\sin\frac{2\pi}7}{\cos\frac\pi{14}},\frac{\sin\frac{3\pi}7}{\sin\frac\pi7},\frac{\sin\frac{2\pi}7\sin\frac\pi7}{\sin\frac\pi7\cos\frac\pi{14}},1,-\frac{\sin^2\frac{2\pi}7}{\sin\frac\pi7\cos\frac\pi{14}},-\frac{\sin\frac{2\pi}7\sin\frac{3\pi}7}{\sin\frac\pi7\cos\frac\pi{14}}),\\
    &~~~~(-\frac{\sin\frac{3\pi}7}{\cos\frac{3\pi}{14}},\frac{\sin\frac{2\pi}7}{\sin\frac\pi7},\frac{\sin\frac{3\pi}7}{\sin\frac\pi7},\frac{\sin\frac{\pi}7}{\cos\frac{3\pi}{14}},-\frac{\sin\frac{3\pi}7\sin\frac{2\pi}7}{\cos\frac{3\pi}{14}\sin\frac\pi7},1,-\frac{\sin^2\frac{3\pi}7}{\cos\frac{3\pi}{14}\sin\frac\pi7},\frac{\sin\frac{3\pi}7}{\cos\frac{3\pi}{14}}),\\
    &~~~~(\frac{\sin\frac{2\pi}7}{\sin\frac\pi7},-\frac{\sin\frac{3\pi}7}{\cos\frac{3\pi}{14}},\frac{\sin\frac{\pi}7}{\cos\frac{3\pi}{14}},\frac{\sin\frac{3\pi}7}{\sin\frac\pi7},-\frac{\sin\frac{3\pi}7\sin\frac{2\pi}7}{\cos\frac{3\pi}{14}\sin\frac\pi7},-\frac{\sin^2\frac{3\pi}7}{\cos\frac{3\pi}{14}\sin\frac\pi7},1,\frac{\sin\frac{3\pi}7}{\cos\frac{3\pi}{14}}),\\
    &~~~~(\frac{\sin\frac{2\pi}7}{\sin\frac\pi7},\frac{\sin\frac{2\pi}7}{\sin\frac\pi7},\frac{\sin\frac{3\pi}7}{\sin\frac\pi7},\frac{\sin\frac{3\pi}7}{\sin\frac\pi7},\left(\frac{\sin\frac{2\pi}7}{\sin\frac\pi7}\right)^2,\frac{\sin\frac{2\pi}7\sin\frac{3\pi}7}{\sin^2\frac\pi7},\frac{\sin\frac{2\pi}7\sin\frac{3\pi}7}{\sin^2\frac\pi7},\left(\frac{\sin\frac{3\pi}7}{\sin\frac\pi7}\right)^2)
\end{align*}
with categorical dimensions
\[ D^2(\mc B)=\begin{cases}\frac{49}{16\cos^2\frac\pi{14}}(\approx3.4)&(\text{1st}),\\\frac{49}{16\cos^2\frac\pi{14}\cos^2\frac{3\pi}{14}}(\approx5.3)&(\text{2nd\&3rd}),\\\frac{49}{16\cos^4\frac{3\pi}{14}}(\approx8.2)&(\text{4th}),\\\frac{49}{16\cos^2\frac\pi{14}\sin^2\frac{\pi}7}(\approx17.1)&(\text{5th\&6th}),\\\frac{49}{16\cos^2\frac{3\pi}{14}\sin^2\frac{\pi}7}(\approx26.6)&(\text{7th\&8th}),\\\frac{49}{16\sin^4\frac\pi7}&(\text{9th}),\end{cases} \]
respectively. Since the two quantum dimensions giving the same categorical dimension are equivalent under permutations $(XY)(ZS)(UV)$ of simple objects, we choose the first quantum dimensions from each pair and fix the order of two factors to remove redundancies. They have conformal dimensions
\begin{align*}
    \hspace{-30pt}&(h_X,h_Y,h_Z,h_S,h_T,h_U,h_V,h_W)\\
    \hspace{-40pt}&=\begin{cases}(\frac47,\frac37,\frac17,\frac67,0,\frac27,\frac57,0),(\frac47,\frac47,\frac67,\frac67,\frac17,\frac37,\frac37,\frac57),(\frac37,\frac37,\frac17,\frac17,\frac67,\frac47,\frac47,\frac27)&(\text{1st}),\\(\frac57,\frac37,\frac17,\frac47,\frac17,0,\frac67,\frac57),(\frac57,\frac47,\frac67,\frac47,\frac27,\frac17,\frac47,\frac37),(\frac27,\frac37,\frac17,\frac37,\frac57,\frac67,\frac37,\frac47),(\frac27,\frac47,\frac67,\frac37,\frac67,0,\frac17,\frac27)&(\text{2nd}),\\(\frac27,\frac57,\frac47,\frac37,0,\frac17,\frac67,0),(\frac57,\frac57,\frac47,\frac47,\frac37,\frac27,\frac27,\frac17),(\frac27,\frac27,\frac37,\frac37,\frac47,\frac57,\frac57,\frac67)&(\text{4th}),\\(\frac37,\frac67,\frac27,\frac17,\frac27,0,\frac57,\frac37),(\frac47,\frac67,\frac27,\frac67,\frac37,\frac57,\frac67,\frac17),(\frac37,\frac17,\frac57,\frac17,\frac47,\frac27,\frac17,\frac67),(\frac47,\frac17,\frac57,\frac67,\frac57,0,\frac27,\frac47)&(\text{5th}),\\(\frac27,\frac67,\frac27,\frac37,\frac17,\frac27,\frac47,\frac57),(\frac27,\frac17,\frac57,\frac37,\frac37,\frac47,0,\frac17),(\frac57,\frac67,\frac27,\frac47,\frac47,\frac37,0,\frac67),(\frac57,\frac17,\frac57,\frac47,\frac67,\frac57,\frac37,\frac27)&(\text{7th}),\\(\frac67,\frac17,\frac57,\frac27,0,\frac37,\frac47,0),(\frac17,\frac17,\frac57,\frac57,\frac27,\frac67,\frac67,\frac37),(\frac67,\frac67,\frac27,\frac27,\frac57,\frac17,\frac17,\frac47)&(\text{9th}).\end{cases}\quad(\mods1)
\end{align*}
The $S$-matrices are given by
\[ \widetilde S=\begin{pmatrix}1&d_X&d_Y&d_Z&d_S&d_Xd_Y&d_Yd_S&d_Xd_Z&d_Zd_S\\d_X&-d_S&d_Xd_Y&d_Xd_Z&1&-d_Yd_S&d_Y&-d_Zd_S&d_Z\\d_Y&d_Xd_Y&-d_Z&1&d_Yd_S&-d_Xd_Z&-d_Zd_S&d_X&d_S\\d_Z&d_Xd_Z&1&-d_Y&d_Zd_S&d_X&d_S&-d_Xd_Y&-d_Yd_S\\d_S&1&d_Yd_S&d_Yd_S&-d_X&d_Y&-d_Xd_Y&d_Z&-d_Xd_Z\\d_Xd_Y&-d_Yd_S&-d_Xd_Z&d_X&d_Y&d_Zd_S&-d_Z&-d_S&1\\d_Yd_S&d_Y&-d_Zd_S&d_S&-d_Xd_Y&-d_Z&d_Xd_Z&1&-d_X\\d_Xd_Z&-d_Zd_S&d_X&-d_Xd_Y&d_Z&-d_S&1&d_Yd_S&-d_Y\\d_Zd_S&d_Z&d_S&-d_Yd_S&-d_Xd_Z&1&-d_X&-d_Y&d_Xd_Y\end{pmatrix}. \]
They have additive central charges
\[ c(\mc B)=c(psu(2)_5)+c(psu(2)_5)\quad(\mods8) \]
where
\[ c(psu(2)_5)=\begin{cases}\begin{cases}\frac47&(h_{X,Y},h_{S,Z})=(\frac37,\frac17),\\-\frac47&(h_{X,Y},h_{S,Z})=(\frac47,\frac67),\end{cases}&(d_{X,Y},d_{S,Z})=(\frac{\sin\frac\pi7}{\cos\frac\pi{14}},-\frac{\sin\frac{2\pi}7}{\cos\frac\pi{14}}),\\\begin{cases}\frac{12}7&(h_{X,Y},h_{S,Z})=(\frac27,\frac37),\\-\frac{12}7&(h_{X,Y},h_{S,Z})=(\frac57,\frac47),\end{cases}&(d_{X,Y},d_{S,Z})=(-\frac{\sin\frac{3\pi}7}{\cos\frac{3\pi}{14}},\frac{\sin\frac{\pi}7}{\cos\frac{3\pi}{14}}),\\\begin{cases}-\frac87&(h_{X,Y},h_{S,Z})=(\frac17,\frac57),\\\frac87&(h_{X,Y},h_{S,Z})=(\frac67,\frac27),\end{cases}&(d_{X,Y},d_{S,Z})=(\frac{\sin\frac{2\pi}7}{\sin\frac\pi7},\frac{\sin\frac{3\pi}7}{\sin\frac\pi7}).\end{cases}\quad(\mods8) \]
Taking the two signs of categorical dimensions, each class of quantum dimensions has $6,8,6,8,8,6$ MFCs, respectively. Therefore, there are
\[ 6+8+6+8+8+6=42 \]
MFCs, among which those six with the ninth quantum dimensions are unitary. We classify connected étale algebras in all 42 MFCs simultaneously.

An ansatz
\[ A\cong1\oplus n_XX\oplus n_YY\oplus n_ZZ\oplus n_SS\oplus n_TT\oplus n_UU\oplus n_VV\oplus n_WW \]
with $n_j\in\mbb N$ has
\[ \hspace{-30pt}\fp_{\mc B}(A)=1+\frac{\sin\frac{2\pi}7}{\sin\frac\pi7}(n_X+n_Y)+\frac{\sin\frac{3\pi}7}{\sin\frac\pi7}(n_Z+n_S)+\left(\frac{\sin\frac{2\pi}7}{\sin\frac\pi7}\right)^2n_T+\frac{\sin\frac{2\pi}7\sin\frac{3\pi}7}{\sin^2\frac\pi7}(n_U+n_V)+\left(\frac{\sin\frac{3\pi}7}{\sin\frac\pi7}\right)^2n_W. \]
For this to obey (\ref{FPdimA2bound}), the natural numbers can take only 110 values. The sets contain the one with all $n_j$'s be zero. It is the trivial connected étale algebra $A\cong1$ giving $\mc B_A^0\simeq\mc B_A\simeq\mc B_A$. Other candidates with $X,Y,Z,S$ contain nontrivial simple object(s) with nontrivial conformal dimensions, and they all fail to be commutative. We are left with just $T,U,V,W$. Setting $n_X,n_Y,n_Z,n_S$ to zero, apart from the trivial one, we find 11 sets
\begin{align*}
    (n_T,n_U,n_V,n_W)&=(1,0,0,0),(0,1,0,0),(0,0,1,0),(0,0,0,1),\\
    &~~~~(2,0,0,0),(1,1,0,0),(1,0,1,0),(1,0,0,1),\\
    &~~~~(0,2,0,0),(0,1,1,0),(0,0,2,0).
\end{align*}
All but the eighth are ruled out because they demand
\[ \fp(\mc B_A^0)\approx4.8,\ 3.4,\ 3.4,\ 2.3,\ 1.5,\ 1.3,\ 1.3,\ 1,\ 1.04,\ 1.04,\ 1.04, \]
respectively. There is no MFC with these Frobenius-Perron dimensions except $\fp(\mc B_A^0)=1$. Thus, we are only left with
\[ A\cong1\oplus T\oplus W. \]
This candidate can be commutative only for the first, fourth, and ninth quantum dimensions and the first conformal dimensions in each quantum dimension. Namely, when the two factors have the same quantum dimensions and opposite conformal dimensions. In other words, when the ambient MFC is a Drinfeld center, $\mc B\simeq Z(psu(2)_5)$. In this case, from the lemma 3, we know there exists a Lagrangian algebra. From our analysis above, the only candidate is $A\cong1\oplus T\oplus W$. Therefore, we learn it is connected étale without any computation. In order to see the consistency, let us identify the category of right $A$-modules.

Its Frobenius-Perron dimension $\fp_{\mc B}(A)=\frac7{4\sin^2\frac\pi7}$ requires
\[ \fp(\mc B_A^0)=1,\quad\fp(\mc B_A)=\frac7{4\sin^2\frac\pi7}. \]
The only possible category of dyslectic right $A$-modules is the trivial MFC
\[ \mc B_A^0\simeq\vect. \]
This identification also matches central charges because the Drinfeld center has additive central charge zero.

So as to identify $\mc B_A$, we follow our routine method. Calculating $b_j\otimes A$, we find simple objects:
\[ 1\oplus T\oplus W,\quad X\oplus Y\oplus U\oplus V\oplus W,\quad Z\oplus S\oplus T\oplus U\oplus V\oplus W. \]
They have Frobenius-Perron dimensions
\[ 1,\quad\frac{\sin\frac{2\pi}7}{\sin\frac\pi7},\quad\frac{\sin\frac{3\pi}7}{\sin\frac\pi7}, \]
respectively, and their contributions to $\fp(\mc B_A)$ match as it should be. This suggests $\mc B_A$ have rank three. Indeed, we find a three-dimensional NIM-rep
\begin{align*}
    n_1=1_3,\quad n_X=\begin{pmatrix}0&1&0\\1&0&1\\0&1&1\end{pmatrix}=n_Y,\quad n_Z=\begin{pmatrix}0&0&1\\0&1&1\\1&1&1\end{pmatrix}=n_S,\\
    n_T=\begin{pmatrix}1&0&1\\0&2&1\\1&1&2\end{pmatrix},\quad n_U=\begin{pmatrix}0&1&1\\1&1&2\\1&2&2\end{pmatrix}=n_V,\quad n_W=\begin{pmatrix}1&1&1\\1&2&2\\2&2&3\end{pmatrix}.
\end{align*}
The solution gives identifications
\[ m_1\cong1\oplus T\oplus W,\quad m_2\cong X\oplus Y\oplus U\oplus V\oplus W,\quad m_3\cong Z\oplus S\oplus T\oplus U\oplus V\oplus W. \]
In $\mc B_A$, they have quantum dimensions
\[ (d_{\mc B_A}(m_1),d_{\mc B_A}(m_2),d_{\mc B_A}(m_3))=\begin{cases}(1,\frac{\sin\frac\pi7}{\cos\frac\pi{14}},-\frac{\sin\frac{2\pi}7}{\cos\frac\pi{14}})&(\text{1st quantum dimension}),\\(1,-\frac{\sin\frac{3\pi}7}{\cos\frac{3\pi}{14}},\frac{\sin\frac\pi7}{\cos\frac{3\pi}{14}})&(\text{4th quantum dimension}),\\(1,\frac{\sin\frac{2\pi}7}{\sin\frac\pi7},\frac{\sin\frac{3\pi}7}{\sin\frac\pi7})&(\text{9th quantum dimension}).\end{cases}. \]
These are nothing but the quantum dimensions of $psu(2)_5$. Furthermore, working out the monoidal products, we find
\begin{table}[H]
\begin{center}
\begin{tabular}{c|c|c|c}
    $\otimes_A$&$m_1$&$m_2$&$m_3$\\\hline
    $m_1$&$m_1$&$m_2$&$m_3$\\\hline
    $m_2$&&$m_1\oplus m_3$&$m_2\oplus m_3$\\\hline
    $m_3$&&&$m_1\oplus m_2\oplus m_3$
\end{tabular}.
\end{center}
\end{table}
\hspace{-17pt}This also shows $\mc B_A\simeq psu(2)_5$. Therefore, we found a connected étale algebra
\begin{equation}
    A\cong1\oplus T\oplus W\quad(\text{1st or 4th or 9th quantum dimensions with 1st }h\text{'s}).\label{psu25psu25etale}
\end{equation}

We conclude
\begin{table}[H]
\begin{center}
\begin{tabular}{c|c|c|c}
    Connected étale algebra $A$&$\mc B_A$&$\rank(\mc B_A)$&Lagrangian?\\\hline
    $1$&$\mc B$&$9$&No\\
    $1\oplus T\oplus W$ for (\ref{psu25psu25etale})&$psu(2)_5$&3&Yes
\end{tabular}.
\end{center}
\caption{Connected étale algebras in rank nine MFC $\mcal B\simeq psu(2)_5\boxtimes psu(2)_5$}\label{rank9psu25psu25results}
\end{table}
\hspace{-17pt}Namely, those six in (\ref{psu25psu25etale}) fail to be completely anisotropic, while the other 36 MFCs $\mc B\simeq psu(2)_5\boxtimes psu(2)_5$'s are completely anisotropic.

\subsubsection{$\mc B\simeq psu(2)_{17}$}
The MFCs have nine simple objects $\{1,X,Y,Z,S,T,U,V,W\}$ obeying monoidal products
\begin{table}[H]
\begin{center}
\makebox[1 \textwidth][c]{       
\resizebox{1.2 \textwidth}{!}{\begin{tabular}{c|c|c|c|c|c|c|c|c|c}
    $\otimes$&$1$&$X$&$Y$&$Z$&$S$&$T$&$U$&$V$&$W$\\\hline
    $1$&$1$&$X$&$Y$&$Z$&$S$&$T$&$U$&$V$&$W$\\\hline
    $X$&&$1\oplus Y$&$X\oplus Z$&$Y\oplus S$&$Z\oplus T$&$S\oplus U$&$T\oplus V$&$U\oplus W$&$V\oplus W$\\\hline
    $Y$&&&$1\oplus Y\oplus S$&$X\oplus Z\oplus T$&$Y\oplus S\oplus U$&$Z\oplus T\oplus V$&$S\oplus U\oplus W$&$T\oplus V\oplus W$&$U\oplus V\oplus W$\\\hline
    $Z$&&&&$1\oplus Y\oplus S\oplus U$&$X\oplus Z\oplus T\oplus V$&$Y\oplus S\oplus U\oplus W$&$Z\oplus T\oplus V\oplus W$&$S\oplus U\oplus V\oplus W$&$T\oplus U\oplus V\oplus W$\\\hline
    $S$&&&&&$1\oplus Y\oplus S\oplus U\oplus W$&$X\oplus Z\oplus T\oplus V\oplus W$&$Y\oplus S\oplus U\oplus V\oplus W$&$Z\oplus T\oplus U\oplus V\oplus W$&$S\oplus T\oplus U\oplus V\oplus W$\\\hline
    $T$&&&&&&$1\oplus Y\oplus S\oplus U\oplus V\oplus W$&$X\oplus Z\oplus T\oplus U\oplus V\oplus W$&$Y\oplus S\oplus T\oplus U\oplus V\oplus W$&$Z\oplus S\oplus T\oplus U\oplus V\oplus W$\\\hline
    $U$&&&&&&&$1\oplus Y\oplus S\oplus T\oplus U\oplus V\oplus W$&$X\oplus Z\oplus S\oplus T\oplus U\oplus V\oplus W$&$Y\oplus Z\oplus S\oplus T\oplus U\oplus V\oplus W$\\\hline
    $V$&&&&&&&&$1\oplus Y\oplus Z\oplus S\oplus T\oplus U\oplus V\oplus W$&$X\oplus Y\oplus Z\oplus S\oplus T\oplus U\oplus V\oplus W$\\\hline
    $W$&&&&&&&&&$1\oplus X\oplus Y\oplus Z\oplus S\oplus T\oplus U\oplus V\oplus W$
\end{tabular}.}}
\end{center}
\end{table}
\hspace{-17pt}Thus, they have
\begin{align*}
    \hspace{-50pt}\fp_{\mc B}(1)=1,\quad\fp_{\mc B}(X)=\frac{\sin\frac{2\pi}{19}}{\sin\frac\pi{19}},\quad\fp_{\mc B}(Y)=&\frac{\sin\frac{3\pi}{19}}{\sin\frac\pi{19}},\quad\fp_{\mc B}(Z)=\frac{\sin\frac{4\pi}{19}}{\sin\frac\pi{19}},\quad\fp_{\mc B}(S)=\frac{\sin\frac{5\pi}{19}}{\sin\frac\pi{19}},\\
    \fp_{\mc B}(T)=\frac{\sin\frac{6\pi}{19}}{\sin\frac\pi{19}},\quad
    \fp_{\mc B}(U)=\frac{\sin\frac{7\pi}{19}}{\sin\frac\pi{19}},&\quad
    \fp_{\mc B}(V)=\frac{\sin\frac{8\pi}{19}}{\sin\frac\pi{19}},\quad\fp_{\mc B}(W)=\frac{\sin\frac{9\pi}{19}}{\sin\frac\pi{19}},
\end{align*}
and
\[ \fp(\mc B)=\frac{19}{4\sin^2\frac\pi{19}}\approx175.3. \]
Their quantum dimensions $d_j$'s are solutions of the same multiplication rules $d_id_j=\sum_{k=1}^9{N_{ij}}^kd_k$. There are nine solutions
\begin{align*}
    &(d_X,d_Y,d_Z,d_S,d_T,d_U,d_V,d_W)\\
    &=(\frac{\sin\frac\pi{19}}{\cos\frac{\pi}{38}},-\frac{\sin\frac{8\pi}{19}}{\cos\frac{\pi}{38}},-\frac{\sin\frac{2\pi}{19}}{\cos\frac{\pi}{38}},\frac{\sin\frac{7\pi}{19}}{\cos\frac{\pi}{38}},\frac{\sin\frac{3\pi}{19}}{\cos\frac{\pi}{38}},-\frac{\sin\frac{6\pi}{19}}{\cos\frac{\pi}{38}},-\frac{\sin\frac{4\pi}{19}}{\cos\frac{\pi}{38}},\frac{\sin\frac{5\pi}{19}}{\cos\frac{\pi}{38}}),\\
    &~~~~(-\frac{\sin\frac{3\pi}{19}}{\cos\frac{3\pi}{38}},-\frac{\sin\frac{5\pi}{19}}{\cos\frac{3\pi}{38}},\frac{\sin\frac{6\pi}{19}}{\cos\frac{3\pi}{38}},\frac{\sin\frac{2\pi}{19}}{\cos\frac{3\pi}{38}},-\frac{\sin\frac{9\pi}{19}}{\cos\frac{3\pi}{38}},\frac{\sin\frac{\pi}{19}}{\cos\frac{3\pi}{38}},\frac{\sin\frac{7\pi}{19}}{\cos\frac{3\pi}{38}},-\frac{\sin\frac{4\pi}{19}}{\cos\frac{3\pi}{38}}),\\
    &~~~~(\frac{\sin\frac{5\pi}{19}}{\cos\frac{5\pi}{38}},-\frac{\sin\frac{2\pi}{19}}{\cos\frac{5\pi}{38}},-\frac{\sin\frac{9\pi}{19}}{\cos\frac{5\pi}{38}},-\frac{\sin\frac{3\pi}{19}}{\cos\frac{5\pi}{38}},\frac{\sin\frac{4\pi}{19}}{\cos\frac{5\pi}{38}},\frac{\sin\frac{8\pi}{19}}{\cos\frac{5\pi}{38}},\frac{\sin\frac{\pi}{19}}{\cos\frac{5\pi}{38}},-\frac{\sin\frac{6\pi}{19}}{\cos\frac{5\pi}{38}}),\\
    &~~~~(-\frac{\sin\frac{7\pi}{19}}{\cos\frac{7\pi}{38}},\frac{\sin\frac{\pi}{19}}{\cos\frac{7\pi}{38}},\frac{\sin\frac{5\pi}{19}}{\cos\frac{7\pi}{38}},-\frac{\sin\frac{8\pi}{19}}{\cos\frac{7\pi}{38}},\frac{\sin\frac{2\pi}{19}}{\cos\frac{7\pi}{38}},\frac{\sin\frac{4\pi}{19}}{\cos\frac{7\pi}{38}},-\frac{\sin\frac{9\pi}{19}}{\cos\frac{7\pi}{38}},\frac{\sin\frac{3\pi}{19}}{\cos\frac{7\pi}{38}}),\\
    &~~~~(\frac{\sin\frac{9\pi}{19}}{\cos\frac{9\pi}{38}},\frac{\sin\frac{4\pi}{19}}{\cos\frac{9\pi}{38}},-\frac{\sin\frac{\pi}{19}}{\cos\frac{9\pi}{38}},-\frac{\sin\frac{6\pi}{19}}{\cos\frac{9\pi}{38}},-\frac{\sin\frac{8\pi}{19}}{\cos\frac{9\pi}{38}},-\frac{\sin\frac{3\pi}{19}}{\cos\frac{9\pi}{38}},\frac{\sin\frac{2\pi}{19}}{\cos\frac{9\pi}{38}},\frac{\sin\frac{7\pi}{19}}{\cos\frac{9\pi}{38}}),\\
    &~~~~(-\frac{\sin\frac{8\pi}{19}}{\cos\frac{11\pi}{38}},\frac{\sin\frac{7\pi}{19}}{\cos\frac{11\pi}{38}},-\frac{\sin\frac{3\pi}{19}}{\cos\frac{11\pi}{38}},-\frac{\sin\frac{\pi}{19}}{\cos\frac{11\pi}{38}},\frac{\sin\frac{5\pi}{19}}{\cos\frac{11\pi}{38}},-\frac{\sin\frac{9\pi}{19}}{\cos\frac{11\pi}{38}},\frac{\sin\frac{6\pi}{19}}{\cos\frac{11\pi}{38}},-\frac{\sin\frac{2\pi}{19}}{\cos\frac{11\pi}{38}}),\\
    &~~~~(\frac{\sin\frac{6\pi}{19}}{\cos\frac{13\pi}{38}},\frac{\sin\frac{9\pi}{19}}{\cos\frac{13\pi}{38}},\frac{\sin\frac{7\pi}{19}}{\cos\frac{13\pi}{38}},\frac{\sin\frac{4\pi}{19}}{\cos\frac{13\pi}{38}},\frac{\sin\frac{\pi}{19}}{\cos\frac{13\pi}{38}},-\frac{\sin\frac{2\pi}{19}}{\cos\frac{13\pi}{38}},-\frac{\sin\frac{5\pi}{19}}{\cos\frac{13\pi}{38}},-\frac{\sin\frac{8\pi}{19}}{\cos\frac{13\pi}{38}}),\\
    &~~~~(-\frac{\sin\frac{4\pi}{19}}{\cos\frac{15\pi}{38}},\frac{\sin\frac{6\pi}{19}}{\cos\frac{15\pi}{38}},-\frac{\sin\frac{8\pi}{19}}{\cos\frac{15\pi}{38}},\frac{\sin\frac{9\pi}{19}}{\cos\frac{15\pi}{38}},-\frac{\sin\frac{7\pi}{19}}{\cos\frac{15\pi}{38}},\frac{\sin\frac{5\pi}{19}}{\cos\frac{15\pi}{38}},-\frac{\sin\frac{3\pi}{19}}{\cos\frac{15\pi}{38}},\frac{\sin\frac{\pi}{19}}{\cos\frac{15\pi}{38}}),\\
    &~~~~(\frac{\sin\frac{2\pi}{19}}{\sin\frac\pi{19}},\frac{\sin\frac{3\pi}{19}}{\sin\frac\pi{19}},\frac{\sin\frac{4\pi}{19}}{\sin\frac\pi{19}},\frac{\sin\frac{5\pi}{19}}{\sin\frac\pi{19}},\frac{\sin\frac{6\pi}{19}}{\sin\frac\pi{19}},\frac{\sin\frac{7\pi}{19}}{\sin\frac\pi{19}},\frac{\sin\frac{8\pi}{19}}{\sin\frac\pi{19}},\frac{\sin\frac{9\pi}{19}}{\sin\frac\pi{19}})
\end{align*}
with categorical dimensions
\[ D^2(\mc B)=\begin{cases}\frac{19}{4\cos^2\frac\pi{38}}(\approx4.8)&(\text{1st}),\\\frac{19}{4\cos^2\frac{3\pi}{38}}(\approx5.1)&(\text{2nd}),\\\frac{19}{4\cos^2\frac{5\pi}{38}}(\approx5.7)&(\text{3rd}),\\\frac{19}{4\cos^2\frac{7\pi}{38}}(\approx6.8)&(\text{4th}),\\\frac{19}{4\cos^2\frac{9\pi}{38}}(\approx8.8)&(\text{5th}),\\\frac{19}{4\cos^2\frac{11\pi}{38}}(\approx12.6)&(\text{6th}),\\\frac{19}{4\cos^2\frac{13\pi}{38}}(\approx21.0)&(\text{7th}),\\\frac{19}{4\cos^2\frac{15\pi}{38}}(\approx45.1)&(\text{8th}),\\\frac{19}{4\sin^2\frac\pi{19}}&(\text{9th}),\end{cases} \]
respectively. They have two conformal dimensions each:
\[ \hspace{-40pt}(h_X,h_Y,h_Z,h_S,h_T,h_U,h_V,h_W)=\begin{cases}(\frac2{19},\frac{18}{19},\frac{10}{19},\frac{16}{19},\frac{17}{19},\frac{13}{19},\frac4{19},\frac9{19}),(\frac{17}{19},\frac1{19},\frac9{19},\frac3{19},\frac2{19},\frac6{19},\frac{15}{19},\frac{10}{19})&(\text{1st}),\\(\frac6{19},\frac{16}{19},\frac{11}{19},\frac{10}{19},\frac{13}{19},\frac1{19},\frac{12}{19},\frac8{19}),(\frac{13}{19},\frac3{19},\frac8{19},\frac9{19},\frac6{19},\frac{18}{19},\frac7{19},\frac{11}{19})&(\text{2nd}),\\(\frac9{19},\frac5{19},\frac7{19},\frac{15}{19},\frac{10}{19},\frac{11}{19},\frac{18}{19},\frac{12}{19}),(\frac{10}{19},\frac{14}{19},\frac{12}{19},\frac4{19},\frac9{19},\frac8{19},\frac1{19},\frac7{19})&(\text{3rd}),\\(\frac5{19},\frac7{19},\frac6{19},\frac2{19},\frac{14}{19},\frac4{19},\frac{10}{19},\frac{13}{19}),(\frac{14}{19},\frac{12}{19},\frac{13}{19},\frac{17}{19},\frac5{19},\frac{15}{19},\frac9{19},\frac6{19})&(\text{4th}),\\(\frac1{19},\frac9{19},\frac5{19},\frac8{19},\frac{18}{19},\frac{16}{19},\frac2{19},\frac{14}{19}),(\frac{18}{19},\frac{10}{19},\frac{14}{19},\frac{11}{19},\frac1{19},\frac3{19},\frac{17}{19},\frac5{19})&(\text{5th}),\\(\frac3{19},\frac8{19},\frac{15}{19},\frac5{19},\frac{16}{19},\frac{10}{19},\frac6{19},\frac4{19}),(\frac{16}{19},\frac{11}{19},\frac4{19},\frac{14}{19},\frac3{19},\frac9{19},\frac{13}{19},\frac{15}{19})&(\text{6th}),\\(\frac{7}{19},\frac{6}{19},\frac{16}{19},\frac{18}{19},\frac{12}{19},\frac{17}{19},\frac{14}{19},\frac{3}{19}),(\frac{12}{19},\frac{13}{19},\frac{3}{19},\frac{1}{19},\frac{7}{19},\frac{2}{19},\frac{5}{19},\frac{16}{19})&(\text{7th}),\\(\frac{8}{19},\frac{15}{19},\frac{2}{19},\frac{7}{19},\frac{11}{19},\frac{14}{19},\frac{16}{19},\frac{17}{19}),(\frac{11}{19},\frac{4}{19},\frac{17}{19},\frac{12}{19},\frac{8}{19},\frac{5}{19},\frac{3}{19},\frac{2}{19})&(\text{8th}),\\(\frac4{19},\frac{17}{19},\frac{1}{19},\frac{13}{19},\frac{15}{19},\frac{7}{19},\frac{8}{19},\frac{18}{19}),(\frac{15}{19},\frac{2}{19},\frac{18}{19},\frac{6}{19},\frac{4}{19},\frac{12}{19},\frac{11}{19},\frac{1}{19})&(\text{9th}).\end{cases}\quad(\mods1) \]
The $S$-matrices are given by
\[\widetilde S=\begin{pmatrix}1&d_X&d_Y&d_Z&d_S&d_T&d_U&d_V&d_W\\d_X&-d_Z&d_T&-d_V&d_W&-d_U&d_S&-d_Y&1\\d_Y&d_T&d_W&d_U&d_Z&1&-d_X&-d_S&-d_V\\d_Z&-d_V&d_U&-d_Y&-1&d_S&-d_W&d_T&-d_X\\d_S&d_W&d_Z&-1&-d_T&-d_V&-d_Y&d_X&d_U\\d_T&-d_U&1&d_S&-d_V&d_X&d_Z&-d_W&d_Y\\d_U&d_S&-d_X&-d_W&-d_Y&d_Z&d_V&1&-d_T\\d_V&-d_Y&-d_S&d_T&d_X&-d_W&1&d_U&-d_Z\\d_W&1&-d_V&-d_X&d_U&d_Y&-d_T&-d_Z&d_S\end{pmatrix}. \]
There are
\[ 9(\text{quantum dimensions})\times2(\text{conformal dimensions})\times2(\text{categorical dimensions})=36 \]
MFCs, among which those four with the ninth quantum dimensions are unitary. We classify connected étale algebras in all 36 MFCs simultaneously.

An ansatz
\[ A\cong1\oplus n_XX\oplus n_YY\oplus n_ZZ\oplus n_SS\oplus n_TT\oplus n_UU\oplus n_VV\oplus n_WW \]
with $n_j\in\mbb N$ has
\[ \hspace{-20pt}\fp_{\mc B}(A)=1+\frac{\sin\frac{2\pi}{19}}{\sin\frac\pi{19}}n_X+\frac{\sin\frac{3\pi}{19}}{\sin\frac\pi{19}}n_Y+\frac{\sin\frac{4\pi}{19}}{\sin\frac\pi{19}}n_Z+\frac{\sin\frac{5\pi}{19}}{\sin\frac\pi{19}}n_S+\frac{\sin\frac{6\pi}{19}}{\sin\frac\pi{19}}n_T+\frac{\sin\frac{7\pi}{19}}{\sin\frac\pi{19}}n_U+\frac{\sin\frac{8\pi}{19}}{\sin\frac\pi{19}}n_V+\frac{\sin\frac{9\pi}{19}}{\sin\frac\pi{19}}n_W. \]
For this to obey (\ref{FPdimA2bound}), the natural numbers can take only 107 values. The sets contain the one with all $n_j$'s be zero. It is the trivial connected étale algebra $A\cong1$ giving $\mc B_A^0\simeq\mc B_A\simeq\mc B$. The other 106 candidates contain nontrivial simple object(s) with nontrivial conformal dimensions, and they all fail to be commutative.

We conclude
\begin{table}[H]
\begin{center}
\begin{tabular}{c|c|c|c}
    Connected étale algebra $A$&$\mc B_A$&$\rank(\mc B_A)$&Lagrangian?\\\hline
    $1$&$\mc B$&$9$&No
\end{tabular}.
\end{center}
\caption{Connected étale algebras in rank nine MFC $\mcal B\simeq psu(2)_{17}$}\label{rank9psu217results}
\end{table}
\hspace{-17pt}All the 36 MFCs $\mc B\simeq psu(2)_{17}$'s are completely anisotropic.

\section{Physical applications}
Our classification results have physical implications in constraining renormalization group (RG) flows. Concretely, in massive RG flows, they constrain ground state degeneracies (GSDs) and prove spontaneous symmetry breaking (SSB). In this section, we discuss the physical applications.

\subsection{Theorems}
Let $\mc C$ be a fusion category. Two-dimensional gapped phases with $\mc C$ symmetry stand in bijection with $\mc C$-module categories \cite{TW19,HLS21}
\begin{equation}
    \{\text{2d }\mcal C\text{-symmetric gapped phases}\}\cong\{\mcal C\text{-module categories }\mcal M\}.\label{1to1}
\end{equation}
The correspondence in particular implies GSD of a $\mc C$-symmetric gapped phase in the LHS is equal to $\rank(\mc M)$ of a $\mc C$-module category in the RHS:
\[ \text{GSD}=\rank(\mc M). \]
This leads to constraints on GSDs:\newline

\textbf{Theorem.} \textit{Let $\mcal B$ be a multiplicity-free modular fusion category up to rank nine and $A\in\mcal B$ a connected étale algebra. Two-dimensional $\mcal B$-symmetric gapped phases described by indecomposable $\mcal B_A$'s have}
\[ \text{GSD}\in\begin{cases}\{7\}&(\mc B\simeq su(7)_1),\\\{7\}&(\mc B\simeq su(2)_6),\\\{7,8\}&(\mc B\simeq so(7)_2),\\\{7\}&(\mc B\simeq psu(2)_{13})\end{cases} \]
\textit{for rank seven,}
\[ \text{GSD}\in\begin{cases}\{4,8\}&(\mc B\simeq\vecG_{\mbb Z/2\mbb Z\times\mbb Z/2\mbb Z\times\mbb Z/2\mbb Z}^\alpha\simeq\vecG_{\mbb Z/2\mbb Z}^{-1}\boxtimes\tc\ (\text{16 with (\ref{Z2Z2Z2etale})})),\\\{8\}&(\mc B\simeq\vecG_{\mbb Z/2\mbb Z\times\mbb Z/2\mbb Z\times\mbb Z/2\mbb Z}^\alpha\simeq\vecG_{\mbb Z/2\mbb Z}^{-1}\boxtimes\tc\ (\text{the other 24})),\\\{8\}&(\mc B\simeq\vecG_{\mbb Z/2\mbb Z\times\mbb Z/4\mbb Z}^\alpha),\\\{4,8\}&(\mc B\simeq su(8)_1),\\\{4,8\}&(\mc B\simeq\fib\boxtimes\vecG_{\mbb Z/2\mbb Z}^{-1}\boxtimes\vecG_{\mbb Z/2\mbb Z}^{-1}\ (\text{16 with (\ref{fibZ2Z21Zetale})})),\\\{8\}&(\mc B\simeq\fib\boxtimes\vecG_{\mbb Z/2\mbb Z}^{-1}\boxtimes\vecG_{\mbb Z/2\mbb Z}^{-1}\ (\text{the other 64})),\\\{4,8\}&(\mc B\simeq\fib\boxtimes\tc\ (\text{16 with (\ref{fibtcetale})})),\\\{8\}&(\mc B\simeq\fib\boxtimes\tc\ (\text{the other 24})),\\\{8\}&(\mc B\simeq\fib\boxtimes\vecG_{\mbb Z/4\mbb Z}^\alpha),\\\{4,8\}&(\mc B\simeq\vecG_{\mbb Z/2\mbb Z}^{-1}\boxtimes\fib\boxtimes\fib\ (\text{16 with (\ref{Z2fibfibetale})})),\\\{8\}&(\mc B\simeq\vecG_{\mbb Z/2\mbb Z}^{-1}\boxtimes\fib\boxtimes\fib\ (\text{the other 64})),\\\{3,6,8,10\}&(\mc B\simeq so(9)_2\ (\text{four with }(\ref{so921ZVetale}))),\\\{6,8,10\}&(\mc B\simeq so(9)_2\ (\text{the other 12})),\\\{3,6,8,10\}&(\mc B\simeq\text{Rep}(D(D_3))\ (\text{two with }(\ref{RepDD31YVetale},\ref{RepDD31ZVetale}))),\\\{6,8,10\}&(\mc B\simeq\text{Rep}(D(D_3))\ (\text{the other six with 1st\&2nd }h)),\\\{8,10\}&(\mc B\simeq\text{Rep}(D(D_3))\ (\text{eight with 3rd\&4th }h)),\\\{8\}&(\mc B\simeq su(2)_7),\\\{4,8\}&(\mc B\simeq\fib\boxtimes\fib\boxtimes\fib\ (\text{16 with }(\ref{FibFibFib1Tetale},\ref{FibFibFib1UVetale}))),\\\{8\}&(\mc B\simeq\fib\boxtimes\fib\boxtimes\fib\ (\text{the other 24})),\\\{8\}&(\mc B\simeq\fib\boxtimes psu(2)_7),\\\{8\}&(\mc B\simeq psu(2)_{15})\end{cases} \]
\textit{for rank eight, and}
\[ \text{GSD}\in\begin{cases}\{3,9\}&(\mc B\simeq su(9)_1),\\\{3,9\}&(\mc B\simeq\vecG_{\mbb Z/3\mbb Z\times\mbb Z/3\mbb Z}^\alpha\ (\text{two with 1st }h)),\\\{9\}&(\mc B\simeq\vecG_{\mbb Z/3\mbb Z\times\mbb Z/3\mbb Z}^\alpha\ (\text{two with 2nd }h)),\\\{9\}&(\mc B\simeq\vecG_{\mbb Z/3\mbb Z}^1\boxtimes\ising),\\\{3,6,9\}&(\mc B\simeq\ising\boxtimes\ising\ (\text{eight with (\ref{isingising1YWetale})})),\\\{6,9\}&(\mc B\simeq\ising\boxtimes\ising\ (\text{the other 136})),\\\{9\}&(\mc B\simeq\vecG_{\mbb Z/3\mbb Z}^1\boxtimes psu(2)_5),\\\{9\}&(\mc B\simeq\ising\boxtimes psu(2)_5),\\\{9,12\}&(\mc B\simeq so(11)_2),\\\{6,9\}&(\mc B\simeq su(2)_8),\\\{3,9\}&(\mc B\simeq psu(2)_5\boxtimes psu(2)_5\ (\text{six with }(\ref{psu25psu25etale}))),\\\{9\}&(\mc B\simeq psu(2)_5\boxtimes psu(2)_5\ (\text{the other 36})),\\\{9\}&(\mc B\simeq psu(2)_{17})\end{cases} \]
\textit{for rank nine.}\newline

As demonstrated \cite{KK23MFC} for general MFC (without assumption on multiplicity), we see all of them are spontaneously broken. Here, we have the\newline

\textbf{Definition.} \cite{KK23GSD} Let $\mc C$ be a fusion category and $\mc M$ a (left) $\mc C$-module category describing a $\mc C$-symmetric gapped phase. A symmetry $c\in\mc C$ is called \textit{spontaneously broken} if $\exists m\in\mc M$ such that $c\triangleright m\not\cong m$. We also say $\mc C$ \textit{is spontaneously broken} if there exists a spontaneously broken object $c\in\mc C$. A categorical symmetry $\mc C$ is \textit{preserved} (i.e., not spontaneously broken) if all objects act trivially.\newline

With the definition, one can show a\newline

\textbf{Lemma 3.} \cite{KK23GSD} \textit{Let $\mc C$ be a fusion category and $\mc M$ be an indecomposable (left) $\mc C$-module category. Then, $\rank(\mc M)>1$ implies SSB of $\mc C$ (i.e., $\mc C$ is spontaneously broken.)}\newline

Therefore, we obtain a\newline

\textbf{Corollary.} \textit{Let $\mc B$ be a nontrivial (i.e., $\rank(\mc B)>1$) multiplicity-free modular fusion category up to rank nine and $A\in\mc B$ a connected étale algebra. In two-dimensional $\mc B$-symmetric gapped phases described by indecomposable $\mc B_A$'s, $\mc B$ symmetries are spontaneously broken.}\newline

\textbf{Remark.} As noted in \cite{KK23GSD}, commutativity of an algebra seems too strong; numerical computation suggests an existence of $\mc B$-symmetric gapped phase described by $\mc B_A$ with non-commutative connected separable algebra.

\subsection{Examples}
In this section, we discuss concrete examples and \textit{predict} GSDs. Since unitary cases are relatively well-understood, we study less-understood non-unitary examples.

Pick a non-unitary minimal model\footnote{We basically follow the notations of \cite{FMS}.} $M(p,2p+1)$ with $p\ge2$ as an ultraviolet (UV) theory. It was proved \cite{KK22free} that its relevant $\phi_{5,1}$-deformation preserves rank $(p-1)$ modular fusion subcategory formed by symmetry operators $\{\mc L_{1,1},\mc L_{1,2},\dots,\mc L_{1,p-1}\}$. If the relevant deformation triggers massless renormalization group (RG) flow, it is known \cite{Z90,Z91,M91,RST94,DDT00} that the infrared (IR) theory is another non-unitary minimal model $M(p,2p-1)$. Its further relevant $\phi_{1,2}$-deformation also preserves \cite{KK22free} rank $(p-1)$ modular fusion subcategory $\{\mc L_{1,1},\mc L_{3,1},\dots,\mc L_{2p-3,1}\}$, and IR theory can be another non-unitary minimal model $M(p-1,2p-1)$. However, massless RG flows typically require fine-tuning. Therefore, generic relevant deformations are expected to trigger massive RG flows. The IR theories are $\mc B$-symmetric gapped phases with modular $\mc B$, and we can apply our classification results. Below, we study massive RG flows triggered by relevant deformations of the non-unitary minimal models.\newline

\paragraph{$M(8,15)+\phi_{1,2}$.} The relevant deformation preserves rank seven MFC $\mc B$ with simple objects $\{\mc L_{1,1},\mc L_{3,1},\dots,\mc L_{13,1}\}$. They form $\mc B\simeq psu(2)_{13}$ with identifications
\[ X\cong\mc L_{13,1},\quad Y\cong\mc L_{3,1},\quad Z\cong\mc L_{11,1},\quad U\cong\mc L_{5,1},\quad V\cong\mc L_{9,1},\quad W\cong\mc L_{7,1}. \]
They have the first (non-unitary) quantum dimensions. Their conformal dimensions
\[ (h_{13,1},h_{3,1},h_{11,1},h_{5,1},h_{9,1},h_{7,1})=(\frac{82}5,\frac1{15},11,\frac65,\frac{20}3,\frac{17}5) \]
match our first conformal dimensions mod 1. Having specified the symmetry, we immediately learn from our results that the massive RG flow described by $\mcal B_A$ has $\text{GSD}=7$ and $\mcal B$ symmetry is spontaneously broken.

\paragraph{$M(8,17)+\phi_{5,1}$.} The relevant deformation preserves rank seven MFC $\mc B$ with simple objects $\{\mc L_{1,1},\mc L_{1,2},\dots,\mc L_{1,7}\}$. They form $\mc B\simeq su(2)_6$ with identifications
\[ X\cong\mc L_{1,7},\quad Y\cong\mc L_{1,2},\quad Z\cong\mc L_{1,6},\quad U\cong\mc L_{1,5},\quad V\cong\mc L_{1,3},\quad W\cong\mc L_{1,4}. \]
They have the third (non-unitary) quantum dimensions. Their conformal dimensions
\[ (h_{1,7},h_{1,2},h_{1,6},h_{1,5},h_{1,3},h_{1,4})=(\frac{45}2,\frac{35}{32},\frac{515}{32},\frac{43}4,\frac{13}4,\frac{207}{32}) \]
match our first conformal dimensions mod 1. With this knowledge on symmetry, our classification result immediately implies that the massive RG flow described by $\mc B_A$ has $\text{GSD}=7$ and $\mc B$ symmetry is spontaneously broken.

\paragraph{$M(9,17)+\phi_{1,2}$.} The relevant deformation preserves rank eight MFC $\mc B$ with simple objects $\{\mc L_{1,1},\mc L_{3,1},\dots,\mc L_{15,1}\}$. They form $\mc B\simeq psu(2)_{15}$ with identifications
\[ X\cong\mc L_{15,1},\quad Y\cong\mc L_{3,1},\quad Z\cong\mc L_{13,1},\quad T\cong\mc L_{5,1},\quad U\cong\mc L_{11,1},\quad V\cong\mc L_{7,1},\quad W\cong\mc L_{9,1}. \]
They have the first (non-unitary) quantum dimensions. Their conformal dimensions
\[ (h_{15,1},h_{3,1},h_{13,1},h_{5,1},h_{11,1},h_{7,1},h_{9,1})=(\frac{385}{17},\frac1{17},\frac{276}{17},\frac{20}{17},\frac{185}{17},\frac{57}{17},\frac{112}{17}) \]
match our second conformal dimensions mod 1. Given the symmetry, we immediately learn the massive RG flow described by $\mc B_A$ has $\text{GSD}=8$ and $\mc B$ symmetry is spontaneously broken.

\paragraph{$M(9,19)+\phi_{5,1}$.} The relevant deformation preserves rank eight MFC $\mc B$ with simple objects $\{\mc L_{1,1},\mc L_{1,2},\dots,\mc L_{1,8}\}$. They form $\mc B\simeq su(2)_7$ with identifications
\[ X\cong\mc L_{1,8},\quad Y\cong\mc L_{1,7},\quad Z\cong\mc L_{1,2},\quad T\cong\mc L_{1,3},\quad U\cong\mc L_{1,6},\quad V\cong\mc L_{1,5},\quad W\cong\mc L_{1,4}. \]
They have the fifth (non-unitary) quantum dimensions. Their conformal dimensions
\[ (h_{1,8},h_{1,7},h_{1,2},h_{1,3},h_{1,6},h_{1,5},h_{1,4})=(\frac{119}4,\frac{67}3,\frac{13}{12},\frac{29}9,\frac{575}{36},\frac{32}3,\frac{77}{12}) \]
match our third conformal dimensions mod 1. Having specified the symmetry, we find the massive RG flow described by $\mc B_A$ has $\text{GSD}=8$ and $\mc B$ symmetry is spontaneously broken.

\paragraph{$M(10,19)+\phi_{1,2}$.} The relevant deformation preserves rank nine MFC $\mc B$ with simple objects $\{\mc L_{1,1},\mc L_{3,1},\dots,\mc L_{17,1}\}$. They form $\mc B\simeq psu(2)_{17}$ with identifications
\[ X\cong\mc L_{17,1},\quad Y\cong\mc L_{3,1},\quad Z\cong\mc L_{15,1},\quad S\cong\mc L_{5,1},\quad T\cong\mc L_{13,1},\quad U\cong\mc L_{7,1},\quad V\cong\mc L_{11,1},\quad W\cong\mc L_{9,1}. \]
They have the first (non-unitary) quantum dimensions. Their conformal dimensions
\[ (h_{17,1},h_{3,1},h_{15,1},h_{5,1},h_{13,1},h_{7,1},h_{11,1},h_{9,1})=(\frac{568}{19},\frac1{19},\frac{427}{19},\frac{22}{19},\frac{306}{19},\frac{63}{19},\frac{205}{19},\frac{124}{19}) \]
match our second conformal dimensions mod 1. With the symmetry, we immediately learn the massive RG flow described by $\mc B_A$ has $\text{GSD}=9$ and $\mc B$ symmetry is spontaneously broken.

\paragraph{$M(10,21)+\phi_{5,1}$.} The relevant deformation preserves rank nine MFC $\mc B$ with simple objects $\{\mc L_{1,1},\mc L_{1,2},\dots,\mc L_{1,9}\}$. They form $\mc B\simeq su(2)_8$ with identifications
\[ X\cong\mc L_{1,9},\quad Y\cong\mc L_{1,2},\quad Z\cong\mc L_{1,8},\quad S\cong\mc L_{1,7},\quad T\cong\mc L_{1,3},\quad U\cong\mc L_{1,4},\quad V\cong\mc L_{1,6},\quad W\cong\mc L_{1,5}. \]
They have the third (non-unitary) quantum dimensions. Their conformal dimensions
\[ (h_{1,9},h_{1,2},h_{1,8},h_{1,7},h_{1,3},h_{1,4},h_{1,6},h_{1,5})=(38,\frac{43}{40},\frac{1183}{40},\frac{111}5,\frac{16}5,\frac{51}8,\frac{127}8,\frac{53}5) \]
match our first conformal dimensions mod 1. With the symmetry, we find the massive RG flow described by $\mc B_A$ has $\text{GSD}\in\{6,9\}$ and $\mc B$ symmetry is spontaneously broken.

\section*{Acknowledgment}
We thank Gert Vercleyen for correspondences on typos in AnyonWiki. We also thank Terry Gannon for teaching us the classification results in $\mc C(B_4,2)$, which made us realize a mistake in the previous method.

\appendix
\setcounter{section}{0}
\renewcommand{\thesection}{\Alph{section}}
\setcounter{equation}{0}
\renewcommand{\theequation}{\Alph{section}.\arabic{equation}}

\end{document}